\documentclass{amsart}
\usepackage[utf8]{inputenc}
\usepackage{amsmath} 
\usepackage{amssymb,amsthm} 
\usepackage{dsfont}
\usepackage{graphicx}
\usepackage{enumitem}
\usepackage{tasks}
\usepackage{float}
\usepackage[lofdepth,lotdepth]{subfig}
\usepackage{color}
\usepackage{tikz} 
\usetikzlibrary{matrix,arrows} 
\usepackage{braids} 
\usepackage{enumitem} 

\newcommand{\C}{\mathcal{C}}
\newcommand{\D}{\mathcal{D}}

\newcommand{\T}{\mathcal{T}}
\newcommand{\sym}{\mathcal{S}}
\newcommand{\un}{\mathds{1}}
\newcommand{\kk}{\Bbbk}
\newcommand{\Vect}{\mathcal{V}ect_{\kk}}

\newcommand{\vect}{\mathcal{V}ect_{\kk}^{f}}
\newcommand{\ung}{\underline{g}}
\newcommand{\unh}{\underline{h}}
\newcommand{\unk}{\underline{k}}

\newcommand{\unx}{\underline{X}}
\newcommand{\uny}{\underline{Y}}
\newcommand{\unX}{\underline{X}}
\newcommand{\unY}{\underline{Y}}
\newcommand{\unA}{\underline{A}}
\newcommand{\unB}{\underline{B}}

\newcommand{\unD}{\underline{D}}
\newcommand{\unZ}{\underline{Z}}

\newcommand{\V}{\mr{V}}
\newcommand{\W}{\mr{W}}
\newcommand\w[1]{\W_{\C,\alpha}(#1)}
\newcommand\halo[1]{\overset{\circ}{#1}}
\newcommand\til[1]{\widetilde{#1}}
\newcommand{\Cob}{\mathcal{C}ob_{3}}
\newcommand{\Cobn}{\mathcal{C}ob_{n}}
\newcommand{\Cobp}{\mathcal{C}ob^{p}_{3}}
\newcommand{\Cobc}{\mathcal{C}ob^{p,cd}_{3}}
\newcommand{\Sp}{\mathbb{S}}
\newcommand{\DD}{\mathbb{D}}
\newcommand{\id}{\mathrm{id}}
\newcommand{\mr}[1]{\mathrm{#1}}
\newcommand{\dimm}{\mathrm{dim}}
\newcommand{\dimq}{\mathrm{dim}_{q}}
\newcommand{\ww}{\W_{\C}(\phantom{o};\alpha)}

\newcommand\pa[1]{p_{\alpha,#1}}
\newcommand\qa[1]{q_{\alpha,#1}}
\newcommand\pia[1]{\Pi_{\alpha,#1}}
\newcommand\tqft[1]{\V_{\C}(#1;\alpha)}
\newcommand\inv[1]{|#1|_{\C,\alpha}}
\newcommand{\iso}{\overset{\sim}{\longrightarrow}}
\newcommand\exitbox[1]{
\begin{tikzpicture}[anchor=base,baseline]
\draw (0.5,0.2) node[]{#1};
\draw (0.5,-0.2) node[]{\tiny{$\ldots$}};
\draw (0,0) rectangle (1,0.6);
\draw[->] (0.1,0)--(0.1,-0.2);
\draw[->] (0.2,-0.4)--(0.2,-0.2);
\draw[-] (0.1,-0.2)--(0.1,-0.4);
\draw[-] (0.2,-0.2)--(0.2,0);
\draw[->] (0.8,0)--(0.8,-0.2);
\draw[->] (0.9,-0.4)--(0.9,-0.2);
\draw[-] (0.8,-0.2)--(0.8,-0.4);
\draw[-] (0.9,-0.2)--(0.9,0);
\end{tikzpicture}
}
\newcommand\handlebox{
\begin{tikzpicture}[scale=0.6]
\draw (0.89,0.55) node[]{\tiny{$\ldots$}};
\draw (0,0) rectangle (1.8,0.4);
\draw[->] (0.7,0.4) arc (0:180:0.3) ;
\draw[->] (1.7,0.4) arc (0:180:0.3) ;
\end{tikzpicture}
}
\newcommand\enterbox[1]{
\begin{tikzpicture}[baseline={([yshift=-18pt]current bounding box.north)}]
\draw (0.5,0.2) node[]{#1};
\draw (0.5,0.8) node[]{\tiny{$\ldots$}};
\draw (0,0) rectangle (1,0.6);
\draw[->] (0.2,0.6)--(0.2,0.8);
\draw[->] (0.1,1)--(0.1,0.8);
\draw[-] (0.2,0.8)--(0.2,1);
\draw[-] (0.1,0.8)--(0.1,0.6);
\draw[->] (0.9,0.6)--(0.9,0.8);
\draw[->] (0.8,1)--(0.8,0.8);
\draw[-] (0.9,0.8)--(0.9,1);
\draw[-] (0.8,0.8)--(0.8,0.6);
\end{tikzpicture}
}

\newcommand{\languette}{
\begin{tikzpicture}[baseline]
\draw[->] (0,0.25) .. controls (0.15,-0.25) and (0.3,0.25) .. (0.3,0.25);
\end{tikzpicture}}
\newcommand{\psii}{
\begin{tikzpicture}
[scale=0.3,description/.style={fill=white,inner sep=2pt},baseline=(current  bounding  box.center)]
\draw (1,0) arc (0:-180:1) ;
\draw (0,0)--(0,-0.85);
\draw (0,-1.15)--(0,-1.5);
\draw (0,-1.5) node[below] {\footnotesize $C$};
\draw (1,-1) node[right] {\footnotesize $Y$};
\end{tikzpicture}}
\newcommand{\psiii}{
\begin{tikzpicture}
[scale=0.3,description/.style={fill=white,inner sep=2pt},baseline=(current  bounding  box.center)]
\draw (1,0) arc (0:-80:1) ;
\draw (-1,0) arc (-180:-100:1) ;
\draw (0,0)--(0,-1.5);
\draw (0,-1.5) node[below] {\footnotesize $C$};
\draw (1,-1) node[right] {\footnotesize $Y$};
\end{tikzpicture}}

\newtheoremstyle{th}{8pt}{8pt}{\itshape}{}{\bfseries}{ ---}{5pt}{\thmname{#1}\thmnumber{ #2} \thmnote{\bfseries (#3).}}
\theoremstyle{th}
\newtheorem{theorem}{Theorem}[section]
\newtheorem{proposition}[theorem]{Proposition}
\newtheorem{definition}[theorem]{Definition}
\newtheorem{lemma}[theorem]{Lemma}
\newtheorem{example}[theorem]{Example}
\newenvironment{remarks}{\addtocounter{theorem}{1}\noindent\textit{\\\bfseries{Remarks \thetheorem} --}\begin{enumerate}}{\end{enumerate}}

\title{Internal Reshetikhin-Turaev TQFT}
\author{mickael lallouche}
\address{Université Lyon 2, Laboratoire ERIC, 5 av. Pierre Mendès-France, 69676 Bron, FRANCE}

\begin{document}
\maketitle



\begin{abstract}
A 3-dimensional topological quantum field theory (TQFT) is a symmetric monoidal functor from the category of 3-cobordisms to the category of vector spaces. Such TQFTs provide in particular numerical invariants of closed 3-manifolds such as the Reshetikhin-Turaev invariants and representations of the mapping class group of closed surfaces. In 1994,  using a modular category, Turaev explains how to construct a TQFT. 
In this article, we describe a generalization of this construction starting from a ribbon category $\C$ with coend. We present a cobordism by a special kind of tangle and we associate to the latter a morphism defined between tensorial products of the coend as described by Lyubashenko in 1994. Composing with an \emph{admissible} color and using extension of Kirby calculus on 3-cobordisms, this morphism gives rise to an \emph{internal} TQFT which takes values in the symmetric monoidal subcategory of transparent objects of $\C$. When the category $\C$ is modular, this subcategory is equivalent to the category of vector spaces. When the category $\C$ is premodular and normalizable with invertible dimension, our TQFT is a lift of Turaev's one associated to the modularization of $\C$.
\end{abstract}

\section{Introduction and overview}

During last twenty years, deep connections were highlighted between low-dimen\-sional topology and braided categories leading to the construction of so-called \emph{quantum invariants}. Some of them are controlled by a \emph{Topological Quantum Field Theory} (for short TQFT), concept introduced in 1988 by Witten \cite{Wit}, before being formalized by Atiyah \cite{Ati}. A $n$-dimensional TQFT is classically defined as a symmetric monoidal functor defined from the category of $n$-dimensional cobordisms to the category of $\kk$-linear vector spaces. 1-dimensional TQFTs are classified by finite dimensional vector spaces. 2-dimensional TQFTs are classified by commutative Frobenius algebras (see \cite{Kock}). According to our knowledge,  such a classification hasn't been found yet for $3$-dimensional TQFTs. However, several constructions exist and they gives rise in particular to 3-manifolds scalar invariants and representations of mapping class groups of closed surfaces.

In 1991, Reshetikhin and Turaev \cite{RT} give the first rigorous construction of such scalar invariants. Their construction is based on surgery along ribbon links and on the use of links invariants. After that, Blanchet, Habegger, Masbaum et Vogel \cite{BHMV} extend this scalar invariant in a 3-dimensional TQFT using Kauffman bracket, before Turaev \cite{Turaev1994} formalizes the key ingredient of \emph{modular categories} to construct families of TQFTs. Such categories are semisimple ribbon categories, with finite number of isomorphism classes of simple objects, such that a certain \emph{S}-matrix, modelling some link crossings colored by $\C$, is invertible.
In \cite{Brug}, Bruguières shows that the relax condition of \emph{modularizability} for \emph{premodular categories}, where $S$-matrix is not necessarily invertible, still allows constructions of TQFTs.
 
Alternatively, in 1994, Lyubashenko \cite{Lyu} generalizes the scalar invariant of Reshetikhin and Turaev using \emph{quantum groups} and \emph{Hopf algebras} for which representations categories are not necessarily semisimple. The main ingredient of Lyubashenko's work is a ribbon category \emph{coend}, that is a special object of the ribbon category satisfying a universal property and admitting a Hopf algebra structure. A ribbon category with coend $C$ equipped with an \emph{integral} $\un\rightarrow C$, where $\un$ is the monoidal unit, gives rise to a 3-manifold invariant. For example,  using a modular category $\C$, and denoting by $\Lambda_{\C}$ the finite set of simple objects isomorphism classes, the coend of $\C$ is given by
\[C=\sum_{X\in \Lambda_{\C}}X^{*}\otimes X,\]
whereas integral of the coend is given by a morphism $\alpha\colon \un\rightarrow C$  corresponding to the following fusion algebra element of $\C$
\[\Omega=\frac{1}{\mr{dim}(\C)}\sum_{X \in \Lambda_{\C}}\mr{dim}_{q}(X)X.\] 

A first construction of TQFTs starting from a ribbon category with coend is done by Kerler and Lyubashenko in the early 2000's \cite{KL}. For this, they had to modify the source category of the TQFT: boundaries are connected and the monoidal product is the connected sum.
In 2002, Virelizier \cite{Vir} proves that a ribbon category with coend equipped with a \emph{Kirby element}, a more general morphism than an integral, leads to a 3-manifold invariant. 

In this article, our goal is to give an, as general as possible, non-semisimple method to construct 3-dimensional TQFTs, starting from a ribbon category with coend, keeping a classical source category for cobordisms but modifying the target category of vector spaces. By analogy with a homological theory (see \cite{Ati}), changing the target category could be interpreted as a coefficient changing. Note that De Renzi, Gainutdinov, Geer, Patureau-Mirand and Runkel have given a method to construct 3-dimensional TQFTs using non semi-simple modular categories using the so-called \emph{modified traces} to renormalize Lyubashenko invariant in \cite{RGGPR}. For this, they need to use a source category of \emph{admissible} cobordisms which is not rigid in general. 
In our work, we will show that a ribbon category $\C$ with coend $C$ equipped with an \emph{admissible element} $\alpha\colon \un \rightarrow C$ provides a 3-dimensional TQFT which takes values in a symmetric monoidal subcategory $\T$ of $\C$, the category of \emph{transparent} objects. Moreover, we show that our TQFT is explicitly computable using only structural morphisms of the coend $C$, extending framework on \emph{Hopf diagrams} of Bruguière and Virelizier \cite{BrugVirHopf}. As an application of our work, we study cases when $\C$ is \emph{modular} and when $\C$ is \emph{modularizable}. In the first case, we recover Reshetikhin-Turaev TQFT composing our \emph{internal} TQFT $V_{\C,\alpha}$ with the functor $\mathrm{Hom}_{\C}(\un,-)$. In the second case, our internal TQFT is a \emph{transparent lift} of the Reshetikin-Turaev TQFT based on the \emph{modularization} of $\C$. Note that in this case, non-trivial transparent objects of $\C$ represent exactly the obstruction for $\C$ to be modular (see \cite{Brug}).

The paper is organized as follows. In Section 2, we review basics on ribbon categories, coend and TQFTs with anomaly. In Section 3, we define \emph{cobordism} tangles, representing 3-cobordisms, and \emph{opentangles}, allowing us to use universal morphisms coming from coend. We recall equivalence moves on these objects and extension of Kirby calculus on 3-cobordism. In Section 4, we explicit all results to construct our \emph{internal} TQFT. Starting from a ribbon category $\C$ with coend $C$, we define an \emph{admissible element} in Definition~\ref{maindefinition} as a morphism $\alpha\colon\un \rightarrow C$ satisfying five \emph{admissible} conditions allowing topological properties of handles slidings, normalizations and compositions. Using such an admissible element, we are able to construct our internal TQFT $\mathrm{V}_{\C,\alpha}$ as claimed in the first main Theorem~\ref{maintheorem} of the article. Second main Theorem~\ref{maintheorem2} explains we can compute all the internal TQFT using only structural morphisms of the coend $C$. In Section 5, we give two applications of our work expliciting cases when $\C$ is modular and $\C$ is modularizable. In the first case, composing the internal TQFT $\mathrm{V}_{\C,\alpha}$ with functor $\mathrm{Hom}_{\C}(\un,-)$, we recover Reshetikhin-Turaev TQFT. In the second case, we obtain a lift of the Reshetikhin-Turaev TQFT on the modularized category of $\C$ taking values in transparent object of $\C$.

\subsection{Acknowledgements} I am deeply grateful to my PhD supervisors Alain Bruguières and Alexis Virelizier for introducing me to the beautiful field of TQFTs. Precious advices and numerous helps from Alain have contributed to increase the quality of this text.


\section{Categorical preliminaries}
In this section, we review main basic definitions and useful results concerning ribbon categories, (pre)modular categories, split idempotents, coends, universal morphism and TQFTs with anomaly. We denote by $\kk$ a field when nothing else is mentionned.

\subsection{Ribbon categories} Let $\C$ be a \emph{monoidal category} which means $\C$ is equipped with a functor $\otimes \colon \C\times \C \rightarrow \C$, called the monoidal (or tensor) product and a unit object $\un$ with natural identifications which guarantee associativity of the product and unit property. In the sequel, we suppose that every monoidal category is a \emph{strict} monoidal category which means identifications are identities (see \cite{MacL}, Theorem XI.3.1). For example, the category $\Vect$ of $\kk$-vector spaces equipped with the tensor product and the unit $\kk$ is a monoidal category.

A \emph{rigid category} is a monoidal category equipped with \emph{left} and \emph{right duality}. A \emph{left duality} in $\C$ associates to every object $X$ an object $X^{*}$, called \emph{left dual}, and two morphisms $\mathrm{ev}_{X}\colon X^{*}\otimes X\rightarrow \un$ and $\mathrm{coev}_{X}\colon \un\rightarrow X\otimes X^{*}$ satisfying for every object $X$:
\begin{align}(\mathrm{id}_{X}\otimes \mathrm{ev}_{X})(\mathrm{coev}_{X}\otimes \mathrm{id}_{X})=\mathrm{id}_{X} \\(\mathrm{ev}_{X}\otimes \mathrm{id}_{X^{*}})(\mathrm{id}_{X^{*}}\otimes \mathrm{coev}_{X})=\mathrm{id}_{X^{*}} \end{align}
We have similar definition for a \emph{right duality} in $\C$ using a \emph{right dual} ${}^{*}\!X$ and morphisms $\til{\mathrm{ev}}_X\colon X\otimes {}^{*}\!X\rightarrow \un $ and $\til{\mathrm{coev}}_X\colon \un  \rightarrow {}^{*}\!X\otimes X$.

A \emph{pivotal category} is a rigid category endowed with a monoidal natural isomorphism between left and right duality. Without loss of generality, we assume that this isomorphism is the identity. In particular, $X^{*}={}^{*}\!X$.

A \emph{braided category} is a monoidal category endowed with a {\it braiding}, that is, a natural system of isomorphisms $\{\tau_{X,Y}\colon X\otimes Y \rightarrow Y\otimes X\}_{X,Y\in \C}$ satisfying for every objects $X,Y,Z$:
\begin{align}
\tau_{X\otimes Y,Z}=(\tau_{X,Z}\otimes \mr{id}_{Y})(\mr{id}_{X}\otimes \tau_{Y,Z}) \\ \tau_{X,Y\otimes Z}=(\mr{id}_{Y}\otimes \tau_{X,Z})(\tau_{X,Y}\otimes \mr{id}_{Z}) 
\end{align}

A particular case of braided category is a \emph{symmetric monoidal category}: this is a monoidal category endowed with a \emph{symmetric braiding}, that means for all objects $X,Y$ of $\C$, the braiding $\tau$  is such that 
\[\tau_{Y,X}\tau_{X,Y}=\mathrm{id}_{X\otimes Y}\]
The target category in our construction of TQFT leads to such example of category as the {\it full subcategory of transparent objects} of a braided category, where $X\in \C$ is a {\it transparent object} if for every object $Y$ of $\C$,
\begin{align}
\tau_{Y,X}\tau_{X,Y}=\mathrm{id}_{X\otimes Y}.
\end{align}

A \emph{balanced category} is a braided category endowed with a {\it twist}, that is, a natural system of isomorphisms $\{\theta_X\colon X\rightarrow X\}_{X\in \C}$ satisfying for every objects $X,Y$:
\begin{align}
\theta_{X\otimes Y}=\tau_{Y,X}\tau_{X,Y}(\theta_{X}\otimes \theta_{Y})
\end{align}
Note that a pivotal balanced category has a left balanced structure given for every object $X$ by $\theta_{X}^{l}=(\mathrm{id}_{X}\otimes \til{\mathrm{ev}}_{X})(\tau_{X,X}\otimes \mathrm{id}_{X^{*}})(\mathrm{id}_{X}\otimes \mathrm{coev}_{X})$ and a right balanced structure given by $\theta_{X}^{r}=(\mathrm{ev}_{X}\otimes \mathrm{id}_{X})(\mathrm{id}_{X}\otimes \tau_{X,X})(\til{\mathrm{coev}}_{X}\otimes \mathrm{id}_{X})$.

A \emph{ribbon category} is a braided pivotal category such that for all objects X, $\theta_{X}^{l}=\theta_{X}^{r}$. This axiom is equivalent to ask for all objects $X$:
\begin{align}
(\theta_{X}^{l}\otimes \mathrm{id_{X}})\mathrm{coev}_{X}=(\mathrm{id}_{X}\otimes\theta_{X^{*}}^{l} )\mathrm{coev}_{X}\\
(\theta_{X}^{r}\otimes \mathrm{id_{X}})\mathrm{coev}_{X}=(\mathrm{id}_{X}\otimes\theta_{X^{*}}^{r} )\mathrm{coev}_{X}
\end{align}

\begin{figure}[H]
\begin{center}
\captionsetup[subfigure]{margin=1pt}
\subfloat[$f\colon X \rightarrow Y$]{
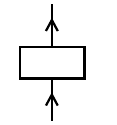}
\subfloat[$g\circ f\colon X\rightarrow Z$]{
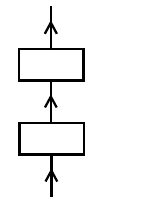}
\subfloat[$f\otimes g \colon X\otimes U\rightarrow Y\otimes V$]{
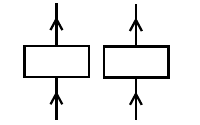}
\subfloat[$\mathrm{id}_{X}$]{
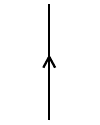}
\subfloat[$\mathrm{id}_{X^{*}}$]{
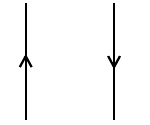}
\end{center}
\begin{center}
\subfloat[Left and right duality morphisms.]{
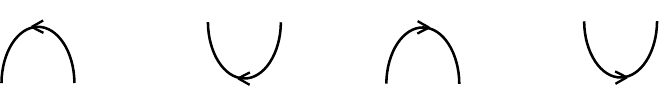}
\end{center}
\begin{center}
\subfloat[Braiding and twist.]{
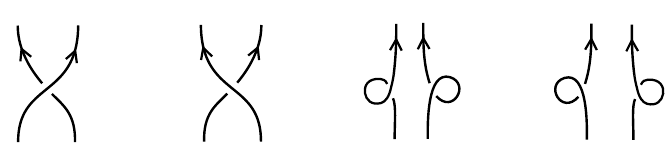}
\end{center}
\caption{Diagrammatic representation of morphisms in a ribbon category.}
\label{ribboncategory}
\end{figure}

\noindent The category of oriented ribbon tangles satisfies a universal property (see \cite{Shu94}), which means in particular that any oriented ribbon tangle $T$ colored by objects of $\C$ defines a morphism in~$\C$. 
Using Penrose graphical calculus in a ribbon category $\C$, we represent morphisms of $\C$ by plane diagrams as in Figure~\ref{ribboncategory}.

\subsection{Premodular and modular categories}

Let $\kk$ be a commutative ring. A \emph{$\kk$-linear category} $\C$ is a category where Hom-sets are $\kk$-modules, the composition of morphisms is $\kk$-bilinear,
and any finite family of objects has a direct sum. In particular, such a category $\C$ has
a zero object.
An object $X$ of  $\C$ is \emph{scalar} if the map $\kk \rightarrow \mathrm{End}_{\C}(X)$,
$\lambda \mapsto \lambda \mr{id}_X$ is bijective.

Let $\C$ and $\D$ be two $\kk$-additive categories. A \emph{$\kk$-linear functor} from $\C$ to $\D$ is a functor which defines a $\kk$-linear map between the $\kk$-modules $\mathrm{Hom}_{\C}(X,Y)$ and $\mathrm{Hom}_{\D}(F(X),F(Y))$ for every objects $X,Y \in \C$.

A \emph{$\kk$-linear monoidal category}\index{\(\kk\)-linear monoidal category} is a monoidal category which is $\kk$-linear in such a way that the monoidal product is $\kk$-bilinear.

 A \emph{fusion category over $\kk$} is a $\kk$-linear rigid category $\C$  such that:
\begin{itemize}
\item each object of $\C$ is a finite direct sum of scalar objects;
\item for any non-isomorphic scalar objects $i$ and $j$ of $\C$, $\mathrm{Hom}_{\C}(i, j) = 0$;
\item the isomorphism classes of scalar objects of $\C$ form a finite set;
\item the unit object $\mathds{1}$ is scalar.
\end{itemize}



A \emph{representative set of scalar objects}\index{representative set of scalar objects} of $\C$ is a set $\Lambda$ of scalar objects such that $\mathds{1}\in \Lambda$ and every scalar object of $\C$ is isomorphic to exactly one element of $\Lambda$.
Note that if $\kk$ is a field, a fusion category over $\kk$ is abelian and semisimple.
Recall that an abelian category is \emph{semisimple} if its objects are direct sums of simple
objects.

A \emph{premodular category} $\C$ over $\kk$ is a ribbon and fusion category over~$\kk$. Pick a representative set $\Lambda$ of scalar objects of $\C$. For $i, j \in \Lambda$, set

\begin{minipage}{0.9\textwidth}
    \raisebox{10mm}{$S_{i,j}= (\mathrm{ev}_{i}\otimes \mathrm{\widetilde{ev}}_{j})(\mathrm{id}_{i^{*}}\otimes \tau_{j,i}\tau_{i,j}\otimes \mathrm{id}_{j^{*}})(\mathrm{\widetilde{coev}}_{i}\otimes \mathrm{coev}_{j}) =$}
\raisebox{7mm}{
\begingroup%
  \makeatletter%
  \providecommand\color[2][]{%
    \errmessage{(Inkscape) Color is used for the text in Inkscape, but the package 'color.sty' is not loaded}%
    \renewcommand\color[2][]{}%
  }%
  \providecommand\transparent[1]{%
    \errmessage{(Inkscape) Transparency is used (non-zero) for the text in Inkscape, but the package 'transparent.sty' is not loaded}%
    \renewcommand\transparent[1]{}%
  }%
  \providecommand\rotatebox[2]{#2}%
  \newcommand*\fsize{\dimexpr\f@size pt\relax}%
  \newcommand*\lineheight[1]{\fontsize{\fsize}{#1\fsize}\selectfont}%
  \ifx\svgwidth\undefined%
    \setlength{\unitlength}{56.69291339bp}%
    \ifx\svgscale\undefined%
      \relax%
    \else%
      \setlength{\unitlength}{\unitlength * \real{\svgscale}}%
    \fi%
  \else%
    \setlength{\unitlength}{\svgwidth}%
  \fi%
  \global\let\svgwidth\undefined%
  \global\let\svgscale\undefined%
  \makeatother%
  \begin{picture}(1,0.5)%
    \lineheight{1}%
    \setlength\tabcolsep{0pt}%
    \put(-0.01143404,0.2635378){\color[rgb]{0,0,0}\makebox(0,0)[lt]{\begin{minipage}{0.24323186\unitlength}\raggedright \small{$i$}\end{minipage}}}%
    \put(0.89683513,0.25657545){\color[rgb]{0,0,0}\makebox(0,0)[lt]{\begin{minipage}{0.29321132\unitlength}\raggedright \small{$j$}\end{minipage}}}%
    \put(0,0){\includegraphics[width=\unitlength,page=1]{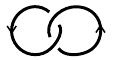}}%
  \end{picture}%
\endgroup%
}
\raisebox{10mm}{$\in \kk$}
\end{minipage}
\vspace{-0.5cm}
\textcolor{white}{OOOOOO} where $\tau$ is the braiding if $\C$. The matrix $S=[S_{i,j}]_{i,j \in I}$ is called the $\emph{S-matrix}$ of $\C$.

 A \emph{modular category} over $\kk$ is a premodular category over $\kk$ for which the $S$-matrix is invertible. This definition is equivalent to ask that the Hopf pairing $\omega$ defined in (\ref{hopfpairing}) is non-degenerate.

There is a way to embed a premodular category into a modular category, assuming some conditions. Let $\C$ and $\D$ be two categories. For $X$ and $Y$ two objects of $\C$, the object $X$ is a \emph{retract} of the object $Y$ if there exist two morphisms $i\colon X \rightarrow Y $ and $p \colon Y \rightarrow X$ such that $p\circ i=\mathrm{id}_{X}$. A functor $F \colon \C \rightarrow \D$ is \emph{dominant}\index{dominant} if for every object $Z\in \D$, there exists an object $X\in \C$ such that $Z$ is a retract of $F(X)$ in $\D$.

A \emph{modularization} of a premodular category $\C$ over $\kk$ is a dominant strong monoidal ribbon $\kk$-linear functor $F\colon \C\rightarrow \widetilde{\C}$, where $\widetilde{\C}$ is a modular category.  A premodular category $\C$ is \emph{modularizable} if it admits a modularization. Note that this modularization is unique up to equivalence in characteristic zero. See \cite{Brug} for more details on modularization.

\subsection{Category with split idempotents}\label{splitidempotent}
An \emph{idempotent}\index{idempotent} of a category $\C$ is an endomorphism $\Pi$ of an object of $\C$ such that $\Pi^{2}=\Pi$.
A \emph{split decomposition}\index{split decomposition} of an idempotent $\Pi$ of an object $X$ of $\C$ is a triple $(A,p,q)$
where $A$ is an object of $\C$ and $p\colon X\rightarrow A$ and $q\colon A\rightarrow X$ are two morphisms of $\C$ satisfying:
\[pq=\mathrm{id}_{A} \quad \text{ and }\quad qp=\Pi.\]
We say that $\C$ is a  \emph{category with split idempotents} if  any idempotent of $\C$ admits a split decomposition. We will always assume that in a category with split idempotents, for each idempotent $\Pi$, a splitting 
$(\mr{Im}(\Pi), p_\Pi,q_\Pi)$ has been chosen. Below, we give some useful results on idempotents.

\begin{proposition}\label{lemmasplit}
Let $\C$ be category with split idempotents. 
Let $X$, $X'$ be two objects of $\C$ and $\Pi$, $\Pi'$ be idempotents of $X$ and $X'$ respectively. 
\begin{enumerate}[label=(\roman*)]
\item If $f\colon X\rightarrow Y$ is a morphism such that $\Pi'f=f\Pi$, then there exists a unique morphism
$g\colon \mr{Im}(\Pi)\rightarrow \mr{Im}(\Pi')$ such that the following diagram \ref{commuteidempotent} commutes: 
\begin{align}\label{commuteidempotent}
\begin{tikzpicture}[description/.style={fill=white,inner sep=2pt},baseline=(current  bounding  box.center)]
\matrix (m) [ampersand replacement=\&,matrix of math nodes, row sep=3em,column sep=2.5em, text height=1.5ex, text depth=0.25ex]
{ X\& \& \mr{Im}(\Pi) \& \& X\\
 Y\& \& \mr{Im}(\Pi') \& \& Y\\ };
\path[->,font=\scriptsize]
(m-1-1) edge node[auto=left] {$p_{\Pi}$} (m-1-3)
		edge node[auto=right] {$f$} (m-2-1)
(m-2-1) edge node[auto=right] {$p_{\Pi'}$} (m-2-3)
(m-1-3) edge node[auto=left] {$g$} (m-2-3)
(m-1-3) edge node[auto=left] {$q_{\Pi}$} (m-1-5)
(m-2-3) edge node[auto=right] {$q_{\Pi'}$} (m-2-5)
(m-1-5) edge node[auto=left] {$f$} (m-2-5);
\end{tikzpicture}
\end{align}
The morphism $g$ is given by $p_{\Pi'}fq_{\Pi}$. We call $g$ the \emph{restriction of $f$ to the images of the idempotents}\index{restriction to the image of idempotents} $\Pi$ and $\Pi'$.
\item This construction induces a bijection
\[\{f\in \mr{Hom}_{\C}(X,X')\;|\; \Pi'f \Pi=f\}\to \mr{Hom}_{\C}\big(\mr{Im}(\Pi),\mr{Im}(\Pi')\big)\] by sending $f$ to its restriction.
\item If $X''$ is a third object and $\Pi''$ an idempotent of $X''$, and if $f\colon X \to X'$ and $f'\colon X' \to X''$ satisfy $\Pi'f = f\Pi$ and $\Pi''f' = f'\Pi'$ then the restriction of $f'f$ on images of idempotents is the compositum of the restrictions of $f$ and $f'$. 
\end{enumerate}
\end{proposition}

In this way, a morphism between images of idempotents $\mr{Im}(\Pi)$ and $\mr{Im}(\Pi')$ can be represented by a morphism $f$ between $X$ and $X'$ satisfying $\Pi'f \Pi=f$. We will be able to use this representation during the construction of our TQFT, when we will project our 3-cobordisms invariant on transparent objects (see definition equality~(\ref{hallowedmorphism})) and when we will unitalize our braided monoidal functor (see Section~\ref{unitalization}).

\subsection{Coend of a ribbon category}
Let $\C$ be a ribbon category. Consider the functor $F\colon \C^{op}\times \C \rightarrow \C$ defined by
\begin{align}\label{dualityfunctor}
F(X,Y)=X^{*}\otimes Y\;
\text{ and } \;
F(f,g)=f^{*}\otimes g
\end{align}
for all objects $X,Y \in \C$ and all morphisms $f,g \in \C$.
The \emph{coend of a ribbon category}\index{coend!of a ribbon category}, when it exists, is the coend of the functor $F$, defined in (\ref{dualityfunctor}). For example, the category $\mathcal{R}ep_{H}$ of left $H$-modules of a finite-dimensional Hopf algebra $H$ over a field $\kk$ possesses a coend $(C,\iota)$ where $C=H^{*}=\mathrm{Hom}_{\kk}(H,\kk)$ and is endowed with the coadjoint left $H$-action given by 
\[(h\otimes f)\in H\otimes H^{*}\mapsto f(S(h_{(1)})\cdot\_\cdot h_{(2)})\in H^{*} \]
and, for a left $H$-module $M$, the dinatural map $\iota_{M}\colon M^{*}\otimes M\rightarrow H^{*}$ is given by
\[(l\otimes m)\in M^{*}\otimes M \mapsto l(\_\cdot m)\in H^{*}.\]

Recall a fundamental result on the ribbon category coend algebraic structure below.

\begin{theorem}\label{thcoendhopf}
Let $\C$ be a ribbon category with a coend $(C,\iota)$. Then $C$ is a Hopf algebra in the category $\C$.
\end{theorem}
See \cite{Lyubhopf} for a proof of Theorem~\ref{thcoendhopf}.

In order to explicit structural morphisms of $C$, recall fundamental consequences of the universal property of the coend $C$ and of the Fubini theorem (see \cite{MacL}, Section IX.4 to IX.7) in the case of a ribbon category:

\begin{proposition}\label{coendpropexistence}
Let $\C$ be a ribbon category with coend $(C,\iota)$, $A$ be an object of $\C$ and 
\[d=\{d_{X_1,\ldots,X_{n}}\colon X_1^{*}\otimes X_1\otimes \ldots\otimes X_n^{*}\otimes X_n \rightarrow A \}_{X_1,\ldots,X_n \in \C}\] be a system of morphisms of $\C$ which is dinatural in every $X_i$ for $1\leq i \leq n$. \newline Then there exists a unique morphism $\phi\colon C^{\otimes n}\rightarrow A$ such that
\[d_{X_1,\ldots,X_n}=\phi\circ(\iota_{X_1}\otimes \ldots\otimes \iota_{X_n})\]
\end{proposition}

\begin{figure}[!h]
\centering
\resizebox{0.99\linewidth}{!}{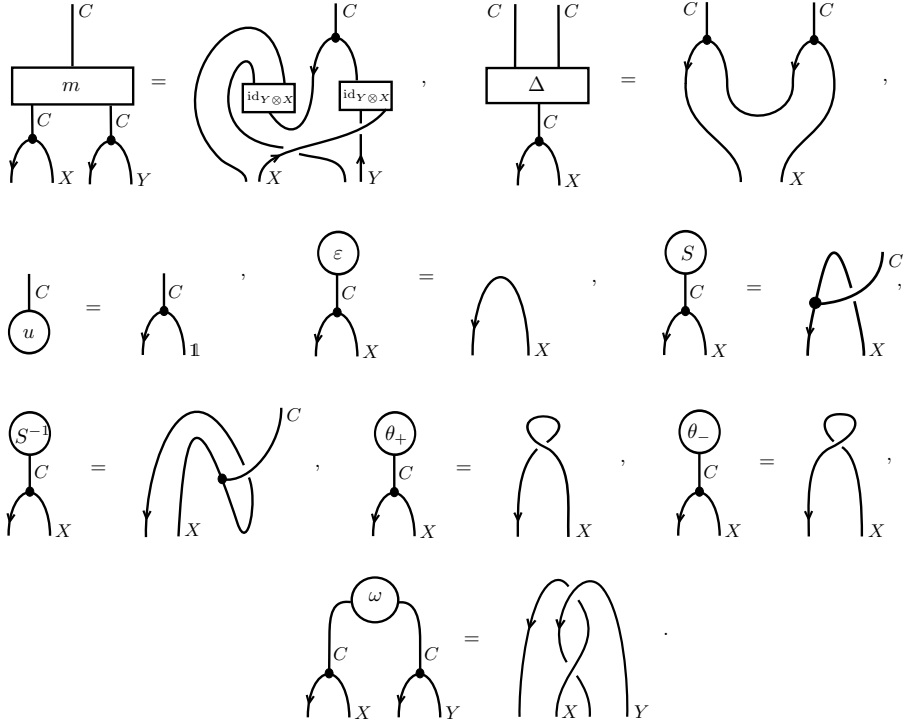}
\caption{Structural morphisms of $C$.}
\label{structural2}
\end{figure}

All structural morphisms of $C$ are defined using Proposition~\ref{coendpropexistence}. Let us define the product $m$, the coproduct $\Delta$, the unit $u$, the counit $\varepsilon$ and the antipode $S$ of $C$ by:
\begin{align}
\iota_{Y\otimes X}(\zeta_{X,Y}\otimes \mathrm{id}_{Y\otimes X})(\mathrm{id}_{X^{*}}\otimes \tau_{X,Y^{*}\otimes Y}) =m(\iota_X \otimes \iota_Y)\colon X^*\otimes X \otimes Y^*\otimes Y \rightarrow C\\
\iota_{\mathds{1}}=u \colon \mathds{1}=\mathds{1}^{*}\otimes \mathds{1}\rightarrow \mathds{1}\\
(\iota_{X} \otimes \iota_{X})(\mathrm{id}_{X^{*}}\otimes \mathrm{coev}_{X}\otimes \mathrm{id}_{X}) = \Delta\iota_X \colon X^{*}\otimes X \rightarrow C \otimes C \\
\mathrm{ev}_{X}=\varepsilon\iota_X \colon X^{*}\otimes X \rightarrow \mathds{1}\\
(\mathrm{ev}_{X}\otimes \iota_{X^{*}})(\mathrm{id}_{X^{*}}\otimes \tau_{X^{**},X}\otimes \mathrm{id}_{X*})(\mathrm{coev}_{X^{*}}\otimes \tau_{X^{*},X})=S\iota_X\colon X^{*}\otimes X\rightarrow C
\end{align}
where equalities are satisfied for every objects $X,Y$ of $\C$ and ${\zeta_{X,Y}\colon X^{*}\otimes Y^{*} \xrightarrow{\sim} (Y \otimes X)^{*}}$ is the isomorphism defined by \[\zeta_{X,Y}=(\mathrm{ev}_{X}(\mathrm{id}_{X^*}\otimes \mathrm{ev}_{Y}\otimes \mathrm{id}_X)\otimes \mathrm{id}_{(Y\otimes X)^*})(\mathrm{id}_{X^*\otimes Y^*}\otimes \mathrm{coev}_{Y\otimes X})\] The antipode $S$ is invertible, with inverse defined via: 
\begin{align}
(ev_{X}\otimes \iota_{X^{*}})(\mathrm{id}_{X^{*}}\otimes \tau_{X^{**},X}^{-1}\otimes \mathrm{id}_{X^*})(\mathrm{coev}_{X^{*}}\otimes \tau_{X^{*},X}^{-1})=&S^{-1}\iota_X\colon X^{*}\otimes X\rightarrow C
\end{align}
and it can be shown that $S^{2}=\theta_{C}$ (see \cite{Lyubhopf}). The coend $C$ is equipped with three additionnal structural morphims, $\theta_{+}\colon C \rightarrow \mathds{1}$\index{$\theta_{+}$}, $\theta_{-}\colon C \rightarrow \mathds{1}$\index{$\theta_{-}$} and $\omega \colon C\otimes C\rightarrow \mathds{1}$ defined using still the universal property of the coend:
\begin{align}
\mathrm{ev}_{X}(\mathrm{id}_{X^*}\otimes \theta_{C})=&\theta_{+}\iota_X\colon X^{*}\otimes X\rightarrow \mathds{1}\\
\mathrm{ev}_{X}(\mathrm{id}_{X^*}\otimes \theta_{C}^{-1})=&\theta_{-}\iota_X\colon X^{*}\otimes X\rightarrow \mathds{1}\\
(\mathrm{ev}_X\otimes \mathrm{ev}_Y)(\mathrm{id}_{X^*}\otimes \tau_{Y^*,X}\tau_{X,Y^*}\otimes \mathrm{id}_Y)=&\omega(\iota_X\otimes \iota_Y)\colon X^*\otimes X \otimes Y^*\otimes Y \rightarrow \mathds{1} \label{hopfpairing}
\end{align}
The morphism $\omega \colon C\otimes C\rightarrow \mathds{1}$ is a {\it Hopf pairing}. Recall that a pairing $\omega$ is said to be non-degenerate if there exists of a morphism $\Omega \colon \mathds{1}\rightarrow C\otimes C$ such that ${(\omega\otimes \mathrm{id}_{C})(\mathrm{id}_C\otimes \Omega)=\mathrm{id}_{C}=(\mathrm{id}_{C}\otimes \omega)(\Omega\otimes \mathrm{id}_C)}$. This is equivalent to say
that $(\omega\otimes \mathrm{id}_{C^*})(\mathrm{id}_C\otimes \mathrm{coev}_C)\colon C \rightarrow C^{*}$ and $(\mathrm{id}_{C^*}\otimes \omega)(\widetilde{\mathrm{coev}}_C \otimes \mathrm{id}_C)\colon C \rightarrow C^{*}$ are isomorphisms. See Section 3A in \cite{BrugVircenter} for more details on Hopf pairing.
The \emph{universal dinatural transformation} of the coend $C$ on an object $X$, denoted by $\iota_{X}\colon X^{*}\otimes X\rightarrow \mathds{1}$\index{$\iota_X$}, is depicted as in Figure~\ref{dinatural}.
Using the graphical calculus illustrated in Figure~\ref{ribboncategory}, we depict all equalities defining the structural morphisms of $C$ in Figure~\ref{structural2}. Note that there is a natural version of all these equalities, using the \emph{universal coaction} of $C$ on objects $X\in \C$ defined by \[\delta_{X}=(\mathrm{id}_{X}\otimes\iota_{X})\circ(\mathrm{coev}_{X}\otimes \mathrm{id}_X)\colon X\rightarrow X\otimes C,\]
and depicted as in Figure~\ref{natural}.

\begin{minipage}{0.5\textwidth}
\label{universalcoaction}
\begin{figure}[H]
\centering
\resizebox{.5\linewidth}{!}{\raisebox{10mm}{$\iota_X=$}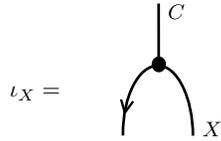}
\caption{Universal dinatural transformation of coend.}
\label{dinatural}
\end{figure}
\end{minipage}
\begin{minipage}{0.5\textwidth}
\begin{figure}[H]
\centering
\resizebox{0.5\linewidth}{!}{\raisebox{10mm}{$\delta_X=$}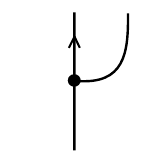}
\caption{Universal natural transformation of coend.}
\label{natural}
\end{figure}
\end{minipage}

\begin{figure}[H]
\centering
\resizebox{.7\linewidth}{!}{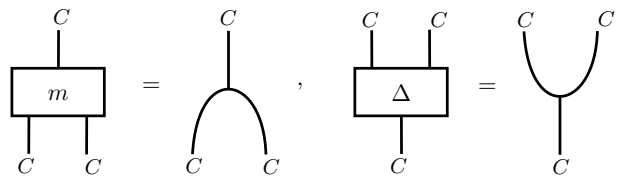}
\caption{The product and the coproduct of $C$.}
\label{pc}
\end{figure}

\phantom{OO}\\
The product $m$ and the coproduct $\Delta$ could be another way depicted as in Figure~\ref{pc}.

\subsection{Universal morphism} Thanks to the universal property of the coend, we define the first algebraic tool we will need to start the construction of our 3-cobordism invariant, and call it {\it the universal morphism}.
Let $\C$ be a ribbon category admitting a coend $(C,\iota)$.
Denote by $\underline{X}=X_1,\ldots,X_n$ an ordered collection of $n$ objects of $\C$ and $\underline{Y}=Y_1,Y_2,\ldots,Y_{2m-1}, Y_{2m}$ an ordered collection of $2m$ objects of $\C$.
\begin{lemma}\label{lemma0}
Consider a system of morphisms of $\C$
\[f=\left
\{f_{\unx,\uny} 
\colon \bigotimes_{i=1}^{n} (X_{i}^{*}\otimes X_{i}) \otimes \bigotimes_{j=1}^{2m} Y_{j} \rightarrow \bigotimes_{j=1}^{2m} Y_{j}
\right\}_{\unx,\uny \in \C}
\]
which is:
\begin{itemize}
    \item dinatural between the functor $(X,Y)\mapsto X^{*} \otimes Y$ and one object $B\in\C$ on each component $X_{i}$ for $1\leq i \leq n$
    \item natural between endofunctors $Id_{\C}$ on each component $Y_{j}$ for $1\leq j \leq 2m$.
\end{itemize} 
Then there exists a unique morphism \[|f| \colon C^{\otimes n+m} \rightarrow C^{\otimes m}\] such that for all $ X_1,\dots,X_n\in\C$ and for all  $Y_1,\ldots,Y_m \in \C$:
\[(\iota_{Y1}\otimes\ldots\otimes\iota_{Y_m})f_{X_1,\ldots,X_n,Y_1^{*},Y_1,\ldots,Y_{m}^{*},Y_{m}} =|f|\circ (\iota_{X_1}\otimes \ldots\otimes \iota_{X_n}\otimes \iota_{Y_{1}}\otimes\ldots\otimes \iota_{Y_m}).
\]
We can depict last equality as:\\
\begin{center}
\resizebox{.95\linewidth}{!}{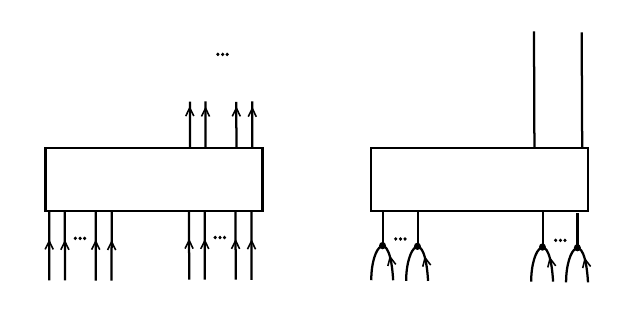}
\end{center}
\end{lemma}
The morphism $|f|$\index{$\lvert f \rvert $} is called the {\it universal morphism}\index{universal morphism} associated to the system of morphisms $f$.
\begin{proof}The system of morphisms
$(\iota_{Y1}\otimes\ldots\otimes\iota_{Y_m})f_{X_1,\ldots,X_{n},Y_1^*,Y_1,\ldots,Y_{m}^*,Y_{m}}$ is dinatural for each component $X_i$ by assumptions. Moreover the system is dinatural in every $Y_j^{*}\otimes Y_j$ because each $Y$ component is supposed to be natural and $\iota$ is dinatural.  The result is then the consequence of Fubini and parameters theorems for coends. For details of the proof, see \cite{MacL}, Section IX.7 and \cite{FS}.
\end{proof}

\subsection{TQFT with anomaly} The aim of this article is to give a construction, as general as possible, of 3-dimensional TQFTs. More precisely, our construction gives raise to a {\it TQFT with anomaly}. In this subsection, we will give the definition of such an object and describe a general algebraic process to lift the \emph{anomaly}. 

A pair of morphisms $(g,f)$ of the category $\C$ is \emph{composable} if the source of $f$ is the target of $g$. A \emph{$2$-cocycle} $\gamma$ for a category $\C$ associates to all pairs $(g,f)$ of composable morphisms of $\C$ a scalar $\gamma_{g,f} \in \kk^ {\times}$\index{$\gamma_{g,f}$} such that for all composable pairs of morphisms $(h,g)$ and $(g,f)$,
 \begin{align}\label{cocycle}
 \gamma_{hg,f}\gamma_{h,g}=\gamma_{h,gf}\gamma_{g,f}.
 \end{align}
For every object $X \in \C$, denote by $\gamma_{X,X}$ the scalar $\gamma_{\mathrm{id}_X,\mathrm{id}_X}$.
\begin{lemma}\label{lemmacocycle}
 Let $\gamma$ be a $2$-cocycle of $\C$ and $f\in \mathrm{Hom}_{\C}(X,Y)$.
Then 
 \begin{align*}\gamma_{f,\mathrm{id}_{X}}=\gamma_{X,X} \text{ and } \gamma_{\mathrm{id}_{Y},f}=\gamma_{Y,Y}.
\end{align*}
 \end{lemma}
\begin{proof}
According to the equality (\ref{cocycle}) defining a $2$-cocycle , we have $\gamma_{f\mathrm{id}_{X},\mathrm{id}_X}\gamma_{f,\mathrm{id}_X}=\gamma_{f,\mathrm{id}_X\mathrm{id}_X}\gamma_{X,X}$, and as $\gamma_{f,\mathrm{id}_X}$ is invertible, $\gamma_{f,\mathrm{id}_{X}}=\gamma_{X,X}$.

 In the same way, $\gamma_{\mathrm{id}_{Y}\mathrm{id}_{Y},f}\gamma_{\mathrm {id}_{Y},\mathrm{id}_{Y}}=\gamma_{\mathrm{id}_{Y},\mathrm{id}_{Y}f}\gamma_{\mathrm{id}_{Y},f}$ and, as $\gamma_{\mathrm{id}_{Y},f}$ is invertible, $\gamma_{Y,Y}=\gamma_{\mathrm{id}_{Y},f}$.
 \end{proof}

\begin{definition}
     A \emph{functor with anomaly} from the category $\C$ to a $\kk$-linear category $\D$ is a pair $(F, \gamma)$ where:
 \begin{itemize}
 \item $F$ associates to each object $X$ of $\C$ an object $F(X)$ of $\D$ and to each morphism $f\in \mathrm{Hom}_{\C}(X,Y)$ a morphism of $D$ in $\mathrm{Hom}_{\D}(F(X),F(Y))$;
 \item $\gamma$ is a $2$-cocycle of $\C$ called the \emph{anomaly}\index{anomaly};
 \end{itemize}
such that for all pairs $(g,f)$ of composable morphisms of $\C$, \[F(g\circ f)=\gamma_{g,f}F(g)\circ F(f).\]
\end{definition}

We notice that we don't know what happens on identities. Next result claims that identities are sent on idempotent morphisms. For every object $X \in \C$, denote by $\Pi_{X}$ the morphism $\gamma_{X,X}F(\mathrm{id}_{X})$ of $\D$.
 \begin{lemma}\label{lemmaidempipo}
 For every object $X \in \C$, the morphism $\Pi_{X}$ is an idempotent of $\D$. 
 \end{lemma}
 \begin{proof}
 We have $\Pi_{X}^{2}=\gamma_{X,X}(\gamma_{X,X}F(\mathrm{id}_X)\circ F(\mathrm{id}_X))=\gamma_{X,X}F(\mathrm{id}_X\circ \mathrm{id}_X)=\gamma_{X,X}F(\mathrm{id}_X)$ so $\Pi_X^2=\Pi_X$.
 \end{proof}
 

 \begin{definition}
A functor with anomaly $(F, \gamma) \colon \C \rightarrow \D$ is \emph{unital} if for all $X \in Ob(\C)$, \[\Pi_{X}=\mathrm{id}_{F(X)}.\]
A unital functor with anomaly is \emph{strict}\index{unital functor with anomaly!strict} if $\gamma_{X,X}=1$.   
 \end{definition}

Even if we can guarantee that identities are sent on identities, there is a natural process to send them on identities up to a scalar. Suppose that the category $\D$ has splitting idempotents and let $(F,\gamma)\colon \C\rightarrow \D$ be a functor with anomaly. Then for all $X \in Ob(\C)$, since $\Pi_X=\gamma_{X,X}F(\mathrm{id}_X)$ is an idempotent (see Lemma~\ref{lemmacocycle}), there exist an object $\mathrm{Im}(\Pi_X) \in \D$ and two morphisms $p_{X} \colon F(X) \rightarrow \mathrm{Im}(\Pi_X)$ and $q_{X} \colon \mathrm{Im}(\Pi_X) \rightarrow F(X)$ such that $p_Xq_X=\mathrm{id}_{\mathrm{Im}(\Pi_X)}$ and $q_{X}p_{X}=\Pi_X$. Then define the \emph{unitalized functor with anomaly}.

\begin{definition}\label{unitalization}
Let $(F,\gamma)$ be a functor with anomaly from $\C$ to $\D$.
The \emph{unitalized functor with anomaly $(\widetilde{F},\gamma)\colon \C \rightarrow \D$ of $(F,\gamma)$} is defined by:
\begin{itemize}
\item for all $X \in \mathrm{Ob}(\C)$, $\widetilde{F}(X)=\mathrm{Im}(\Pi_{X})$;
\item for all $f\in \mathrm{Hom}_{\C}(X,Y)$, $\widetilde{F}(f)=p_{Y}\circ F(f)\circ q_{X}$.
\end{itemize}\label{unital functor}
\end{definition}

The unitalized functor is still a functor with anomaly as proved below.
\begin{lemma}
Let $(F,\gamma)\colon \C \rightarrow \D$ a functor with anomaly. Then the unitalized functor with anomaly $(\widetilde{F},\gamma)$ is a unital functor with anomaly.
\end{lemma}

\begin{proof}
For $X \in \C$, compute $\gamma_{X,X}\widetilde{F}(f)$:
\begin{align*}
\gamma_{X,X}\widetilde{F}(\id_{X}) & =\gamma_{X,X}p_{X}F(\id_{X})q_{X}=p_{X}(q_{X}p_{X})q_{X}=(p_Xq_X)(p_Xq_X)=\id_{\mathrm{Im}(\Pi_{X})}\id_{\mathrm{Im}(\Pi_X)}\\ & =\id_{\widetilde{F}(X)}
\end{align*}
so the unitary condition is satisfied for $\widetilde{F}$.
Next, for $f\colon X \rightarrow Y$ and $g\colon Y \rightarrow Z$, compute $\widetilde{F}(gf)$:
\begin{align*}
\widetilde{F}(gf)=p_{Z}F(gf)q_X & =\gamma_{g,f}F(g)F(f)q_{X}=\gamma_{g,f}F(g)F(\id_{Y}f)q_X\\
&=\gamma_{g,f}p_ZF(g)(\gamma_{\id_{Y},f}F(\id_{Y})F(f))q_{X}
\end{align*}
But, as $\gamma$ is a 2-cocyle and according to the Lemma \ref{lemmacocycle}, $\gamma_{id_{Y},f}=\gamma_{Y,Y}$. Thus:
\begin{align*}
\gamma_{g,f}p_ZF(g)(\gamma_{\id_{Y},f}F(\id_{Y})F(f))q_{X}&=\gamma_{g,f}p_ZF(g)(\gamma_{Y,Y}F(\id_{Y})F(f))q_{X}\\
&=\gamma_{g,f}p_ZF(g)\Pi_{Y}F(f)q_{X}=\gamma_{g,f}p_ZF(g)q_Yp_YF(f)q_{X}\\
&=\gamma_{g,f}(p_ZF(g)q_Y)(p_YF(f)q_{X})=\gamma_{g,f}\widetilde{F}(g)\widetilde{F}(f)
\end{align*}
Then $\widetilde{F}(gf)=\gamma_{g,f}\widetilde{F}(g)\widetilde{F}(f)$.
\end{proof}

To define a general notion of TQFT with anomaly, we need to generalize the notion of natural transformation between functors with anomaly.

\begin{definition}
Let be $(F,\gamma)$ and $(G,\eta)$ two functors with anomaly  from $\C$ to $\D$. \emph{A natural transformation with anomaly} between $F$ and $G$ is a pair ${(\rho \colon (F,\gamma)\rightarrow (G,\eta),\nu)}$ where:
\begin{itemize}
\item for all morphism $f\in Hom_{\C}(X,Y)$, $\nu_{f} \in \kk^{\times}$ and for all pairs ${(g\colon Y\rightarrow Z, f \colon X \rightarrow Y)}$ of composable morphisms, $\nu_{gf}\eta_{g,f}=\nu_{g}\nu_{f}\gamma_{g,f}$;
\item for all $X \in Ob(\C)$, $\rho_{X} \colon F(X) \rightarrow G(X)$ is a morphism of $\D$ such that for every morphism $f\colon X \rightarrow Y$ of $\C$, the following diagram commutes:
\end{itemize}

\begin{center}
\begin{tikzpicture}
\matrix[ampersand replacement=\&,matrix of math nodes,row sep=1cm,column sep=1cm]{
|(C)| F(X) \& |(D)| G(X) \\
|(C')| F(Y) \& |(D')| G(Y) \\
};
\begin{scope}[every node/.style={pos=0.5,font=\scriptsize}] 
\draw[->] (C) -- (D) node[auto=left] {$\rho_{X}$};
\draw[->] (C') -- (D') node[auto=right] {$\rho_{Y}$} ;
\draw[->] (C) -- (C') node[auto=right] {$F(f)$};
\draw[->] (D) -- (D') node[auto=left] {$\nu_{f}G(f)$};
\end{scope}
\end{tikzpicture}
\end{center}
\end{definition}

\begin{remarks}
\item The map $\nu$ is called the \emph{anomaly} of the natural transformation with anomaly $(\rho,\nu)$. 
\item Note that when the anomaly is the constant map equal to $1$, both functors $F$ and $G$ have the same anomaly.
\item When we consider natural transformations between functors with anomaly, we shall mean natural transformations with anomaly.
\end{remarks}

Let's remark that the unizalized functor is universal among a certain class of functors with anomaly. Indeed, recall the definition of unitalized functor $(\widetilde{F},\gamma)$ of a functor with anomaly $(F,\gamma)$ (see Definition~\ref{unital functor}). Note that $p_{X}\colon F(X)\rightarrow \widetilde{F}(X)$ and $q_{X}\colon \widetilde{F}(X)\rightarrow F(X)$ define two natural transformations $p\colon (F,\gamma) \rightarrow (\widetilde{F},\gamma)$ and $ q \colon (\widetilde{F},\gamma) \rightarrow (F,\gamma)$. Then $\Pi_{X}$ defines also a natural transformation ${\Pi \colon (F,\gamma)\rightarrow (F,\gamma)}$. 

The next Lemma~\ref{universalunitalized} explicits the universal property satisfied by the unitalized functor $(\widetilde{F},\gamma)$.

\begin{lemma}\label{universalunitalized}
Let $(F,\gamma)\colon \C \rightarrow \D$ be a functor with anomaly.
Then for any functor with anomaly $(G,\gamma)\colon \C\rightarrow \D$ equipped with two natural transformations $\alpha\colon (F,\gamma)\rightarrow (G,\gamma)$ and $\beta\colon (G,\gamma)\rightarrow (F,\gamma)$ such that $\alpha\beta=\mathrm{id}_{G}$ and $\beta\alpha=\Pi$, there exists a unique natural isomorphism \[\eta=\{\eta_X \colon \widetilde{F}(X)\rightarrow G(X)\}_{X \in \C}\] such that the  following diagram commutes for all objects $X \in \C$:
\begin{center}
\begin{tikzpicture}
\matrix[ampersand replacement=\&,matrix of math nodes,row sep=1cm,column sep=1cm]{
 |(C)| F(X) \& |(D)| \widetilde{F}(X) \& |(E)| F(X)\\
 |(C')|  \& |(D')| G(X) \& |(E')|  \\
};
\begin{scope}[every node/.style={pos=0.5,font=\scriptsize}] 
\draw[->] (C) -- (D) node[auto=left] {$p_{X}$};
\draw[->] (C) -- (D') node[auto=right] {$\alpha_X$};
\draw[->] (D) -- (D') node[auto=left] {$\eta_X$};
\draw[->] (D') -- (E) node[auto=right] {$\beta_X$};
\draw[->] (D) -- (E) node[auto=left] {$q_{X}$};
\end{scope}
\end{tikzpicture}
\end{center}
\end{lemma}
\begin{proof}
Let us define $\eta=\alpha q$. Then $\eta$ is a natural isomorphism that satisfies the conditions (its inverse is $p\beta$).
\end{proof}

We need to extend now the notion of \emph{monoidality} to functors with anomaly. Suppose that the category $\C$ and $\D$ are monoidal. Let $(F,\gamma)$ be a functor with anomaly from $\C$ to $\D$. The anomaly $\gamma$ is \emph{monoidal} if for all morphisms $f\colon X \rightarrow Y$, $g \colon Y\rightarrow Z$, $f'\colon X' \rightarrow Y'$, $g' \colon Y'\rightarrow Z'$ of $\C$, 
\begin{align}
\gamma_{g\otimes g',f\otimes f'}=\gamma_{g,f}\gamma_{g',f'}.
\end{align}

\begin{definition}
A \emph{monoidal functor with anomaly} from the category $\C$ to the category $\D$ is a quadruplet $(F, F_{2},F_{0},\gamma)$ where $(F,\gamma)$ is a functor with monoidal anomaly $\gamma$ from $\C$ to $\D$, $F_{2}$ is a natural transformation between $(F\otimes F,\gamma_{g,f}\gamma_{g',f'})$ and $(F\otimes, \gamma_{g\otimes g',f\otimes f'})$ denoted by  
\[F_{2}=\{F_{2}(X,Y)\colon F(X)\otimes F(Y)\rightarrow F(X\otimes Y)\}_{X,Y\in \C}\]
 and $F_{0}\colon \mathds{1}\rightarrow F(\mathds{1})$ is a morphism of $\D$ such that, for every objects $X$, $Y$ and $Z$ of $\C$, the following diagrams (\ref{diagmonofunctor1}), (\ref{diagmonofunctor2}) and (\ref{diagmonofunctor3}) commute:

\begin{align}\label{diagmonofunctor1}
\begin{tikzpicture}[description/.style={fill=white,inner sep=2pt},baseline=(current  bounding  box.center)]
\matrix (m) [ampersand replacement=\&,matrix of math nodes, row sep=3em,column sep=2.5em, text height=1.5ex, text depth=0.25ex]
{ F(X)\otimes F(Y)\otimes F(Z)\& \& \&F(X\otimes Y)\otimes F(Z) \\
F(X)\otimes F(Y\otimes Z)\& \& \& F(X\otimes Y\otimes Z)\\ };
\path[->,font=\scriptsize]
(m-1-1) edge node[auto=left] {$F_{2}(X,Y)\otimes \mathrm{id}_{F(Z)}$} (m-1-4)
		edge node[auto=right] {$\mathrm{id}_{F(X)}\otimes F_{2}(Y,Z)$} (m-2-1)
(m-2-1) edge node[auto=right] {$F_{2}(X,Y\otimes Z)$} (m-2-4)
(m-1-4) edge node[auto] {$F_{2}(X\otimes Y,Z)$} (m-2-4);
\end{tikzpicture}
\end{align}

\begin{minipage}[t]{0.5\linewidth}
\begin{align}\label{diagmonofunctor2}
\begin{tikzpicture}[description/.style={fill=white,inner sep=2pt},baseline=(current  bounding  box.center)]
\matrix (m) [ampersand replacement=\&,matrix of math nodes, row sep=3em,column sep=2.5em, text height=1.5ex, text depth=0.25ex]
{F(X)\otimes \mathds{1} \& \&F(X)\otimes F(\mathds{1}) \\
\& F(X)\& \\ };
\path[->,font=\scriptsize]
(m-1-1) edge node[auto=left] {$\mathrm{id}_{F(X)}\otimes F_{0}$} (m-1-3)
		edge node[auto=right] {$\mathrm{id}_{F(X)}$} (m-2-2)
(m-1-3) edge node[auto] {$F_{2}(X,\mathds{1})$} (m-2-2);
\end{tikzpicture}
\end{align}
\end{minipage}
\begin{minipage}[t]{0.5\linewidth}
\begin{align}\label{diagmonofunctor3}
\begin{tikzpicture}[description/.style={fill=white,inner sep=2pt},baseline=(current  bounding  box.center)]
\matrix (m) [ampersand replacement=\&,matrix of math nodes, row sep=3em,column sep=2.5em, text height=1.5ex, text depth=0.25ex]
{\mathds{1}\otimes F(X) \& \&  F(\mathds{1})\otimes F(X)\\
\&  F(X) \& \\ };
\path[->,font=\scriptsize]
(m-1-1) edge node[auto=left] {$F_{0}\otimes \mathrm{id}_{F(X)}$} (m-1-3)
		edge node[auto=right] {$\mathrm{id}_{F(X)}$} (m-2-2)
(m-1-3) edge node[auto] {$F_{2}(\mathds{1},X)$} (m-2-2);
\end{tikzpicture}
\end{align}\\
\end{minipage}
\end{definition}

If $F_{2}$ and $F_{0}$ are isomorphisms, the monoidal functor with anomaly $F$ is \emph{strong}. If $F_{2}$ and $F_{0}$ are identities, the monoidal functor $F$ is \emph{strict}.

Let $(F,F_2,F_0,\gamma)$ and $(G,G_2,G_0,\gamma)$ be two monoidal functors with anomaly. A \emph{monoidal natural transformation}\index{monoidal natural transformation between functors with the same anomaly} from $F$ to $G$ is a natural transformation $\rho=\{\rho_{X}\colon F(X)\rightarrow G(X)\}_{X\in\C}$ such that for all objects $X, Y\in \C$, the following diagrams (\ref{diagmononaturalfunctor1}) and (\ref{diagmononaturalfunctor2}) commute:

\begin{minipage}[t]{0.4\linewidth}
\begin{align}\label{diagmononaturalfunctor1}
\begin{tikzpicture}[description/.style={fill=white,inner sep=2pt},baseline=(current  bounding  box.center)]
\matrix (m) [ampersand replacement=\&,matrix of math nodes, row sep=3em,column sep=2.5em, text height=1.5ex, text depth=0.25ex]
{ F(X)\otimes F(Y)\& \&  G(X)\otimes G(Y) \\
F(X\otimes Y)\& \&  G(X\otimes Y)\\ };
\path[->,font=\scriptsize]
(m-1-1) edge node[auto=left] {$\rho_{X}\otimes \rho_{Y}$} (m-1-3)
		edge node[auto=right] {$F_2(X,Y)$} (m-2-1)
(m-2-1) edge node[auto=right] {$\rho_{X\otimes Y}$} (m-2-3)
(m-1-3) edge node[auto] {$G_{2}(X,Y)$} (m-2-3);
\end{tikzpicture}
\end{align}
\end{minipage}
\begin{minipage}[t]{0.5\linewidth}
\begin{align}\label{diagmononaturalfunctor2}
\begin{tikzpicture}[description/.style={fill=white,inner sep=2pt},baseline=(current  bounding  box.center)]
\matrix (m) [ampersand replacement=\&,matrix of math nodes, row sep=3em,column sep=2.5em, text height=1.5ex, text depth=0.25ex]
{ \& \mathds{1} \&  \\
F(\mathds{1})\& \&  G(\mathds{1})\\ };
\path[->,font=\scriptsize]
(m-1-2) edge node[auto=right] {$F_0$} (m-2-1)
		edge node[auto=left] {$G_0$} (m-2-3)
(m-2-1) edge node[auto=right] {$\rho_{\mathds{1}}$} (m-2-3);
\end{tikzpicture}
\end{align}
\end{minipage}

\begin{lemma}
Let $F=(F,F_{2},F_{0},\gamma)\colon \C \rightarrow \D$  be a monoidal functor with anomaly. Then 
\begin{enumerate}[label=(\roman*)]
\item \label{1} The unitalized functor $(\widetilde{F},\gamma)$ admits a unique structure of monoidal functor with anomaly such that the natural transformations $p\colon F\rightarrow \widetilde{F}$ and ${q\colon \widetilde{F}\rightarrow F}$ (see Definition~\ref{unital functor}) are monoidal.
\item \label{2} If $(F,\gamma)$ is strict (respectively strong) then $(\widetilde{F},\gamma)$ is strict (respectively strong).
\end{enumerate}
\end{lemma}
\begin{proof}
Define the monoidal structure of $\widetilde{F}=(\widetilde{F},\widetilde{F_2}, \widetilde{F_0},\gamma)$ by \[\widetilde{F_2}(X,Y)=p_{X\otimes Y}F_2(X,Y)(q_{X}\otimes q_{Y})\colon\widetilde{F}(X)\otimes \widetilde{F}(Y)\rightarrow \widetilde{F}(X\otimes Y)\] and \[\widetilde{F_0}=p_{\un}F_0 \colon \un \rightarrow \widetilde{F}(\un) \]
where $X,Y$ are objects of $\C$, $q_{X}p_{X}=\Pi_X=\gamma_{X,X}F(\mathrm{id}_X)$, $\widetilde{F}(X)=\mathrm{Im}(\Pi_X)$ and $p_{X}q_{X}=\mathrm{id}_{\widetilde{F}(X)}$. Set that $p$ and $q$ are monoidal transformations means that for every objects $X, Y\in \C$, 
\[\widetilde{F_2}(X,Y)(p_X\otimes p_Y)=p_{X\otimes Y}F_2(X,Y) \text{ and } q_{X\otimes Y}\widetilde{F_2}(X,Y)=F_2(X,Y)(q_X\otimes q_Y)\]
so the choice of $\widetilde{F_2}$ is uniquely determined. Moreover, it means that $\widetilde{F_0}=p_{\mathds{1}}F_{0}$ and $q_{\mathds{1}}\widetilde{F_0}=F_0$ so the choice of $\widetilde{F_0}$ is uniquely determined.


First, we check the naturality of $\widetilde{F_2}$. Let $f\colon X \rightarrow X'$ and $g\colon Y \rightarrow Y'$ be two morphisms of $\C$. We compute $\widetilde{F_{2}}(X',Y')(\widetilde{F}(f)\otimes\widetilde{F}(g))$:
\begin{align*}
\widetilde{F_{2}}(X',Y')(\widetilde{F}(f)\otimes\widetilde{F}(g))&=p_{X'\otimes Y'}F_2(X',Y')(q_{X'}\otimes q_{Y'})(p_{X'}F(f)q_{X}\otimes p_{Y'}F(g)q_{Y})\\
&=p_{X'\otimes Y'}F_2(X',Y')(q_{X'}p_{X'}F(f)q_{X}\otimes q_{Y'}p_{Y'}F(g)q_{Y})\\
&\stackrel{(3)}{=}p_{X'\otimes Y'}F_2(X',Y')(\gamma_{X',X'}F(\id_{X'})F(f)q_{X}\otimes \gamma_{Y',Y'}F(\id_{Y'})p_{Y'}F(g)q_{Y})\\
&\stackrel{(4)}{=}p_{X'\otimes Y'}F_2(X',Y')(\gamma_{{\id_{X'},f}}F(\id_{X'})F(f)q_{X}\otimes \gamma_{\id_{Y'},g}F(\id_{Y'})F(g)q_{Y})\\
&=p_{X'\otimes Y'}F_2(X',Y')(F(f)\otimes F(g))(q_X\otimes q_Y)\\
&\stackrel{(6)}{=}p_{X'\otimes Y'}F(f\otimes g)F_2(X,Y)(q_X\otimes q_Y)\\
&=p_{X'\otimes Y'}F((f\otimes g)\id_{X\otimes Y})F_2(X,Y)(q_X\otimes q_Y)\\
&=p_{X'\otimes Y'}(\gamma_{f\otimes g,\id_{X\otimes Y}}F(f\otimes g)F(\id_{X\otimes Y}))F_2(X,Y)(q_X\otimes q_Y)\\
&\stackrel{(9)}{=}p_{X'\otimes Y'}(F(f\otimes g)\Pi_{X\otimes Y})F_2(X,Y)(q_X\otimes q_Y)\\
&\stackrel{(10)}{=}p_{X'\otimes Y'}F(f\otimes g)q_{X\otimes Y}p_{X\otimes Y}F_2(X,Y)(q_X\otimes q_Y)\\
&=\widetilde{F}(f\otimes g)\widetilde{F_2}(X,Y)
\end{align*}

Equalities $(3)$ and $(10)$ follow from the definitions of idempotents $\Pi_{X'}=q_{X'}p_{X'}$, $\Pi_{Y'}=q_{Y'}p_{Y'}$ and $\Pi_{X\otimes Y}=q_{X\otimes Y}p_{X\otimes Y}$. As $\gamma$ is a $2$-cocycle, $\gamma_{X',X'}=\gamma_{\id_X',f}$, $\gamma_{Y',Y'}=\gamma_{\id_Y',g}$ and $\gamma_{X\otimes Y,X\otimes Y}=\gamma_{f\otimes g, \id_{X\otimes Y}}$ which imply equality $(4)$ and $(9)$. Equality $(6)$ comes from the naturality of $F_2\colon F\otimes F\rightarrow F\otimes $.

Next, we have to check the coherence axioms of the definition of a monoidal functor with anomaly. We first compute $\widetilde{F_2}(X,Y\otimes Z)(\mathrm{id}_{\widetilde{F}(X)}\otimes \widetilde{F_{2}}(Y,Z))$ for any objects $X,Y$ of $\C$ and show that the axiom (\ref{diagmonofunctor1}) is satisfied for $\widetilde{F}$:
\begin{align*}
&\widetilde{F_2}(X,Y\otimes Z)(\mathrm{id}_{\widetilde{F}(X)}\otimes \widetilde{F_{2}}(Y,Z))\\
&=p_{X \otimes Y \otimes Z}\widetilde{F}_{2}(X,Y\otimes Z)(q_{X}\otimes q_{Y\otimes Z}))(\mathrm{id}_{\widetilde{F}(X)}\otimes (p_{Y\otimes Z}F_2(Y,Z)(q_Y\otimes q_Z)))\\
&=p_{X \otimes Y \otimes Z}F_{2}(X,Y\otimes Z)(q_{X}\otimes (q_{Y\otimes Z}p_{Y\otimes Z}F_2(Y,Z)(q_Y\otimes q_Z)))\\
&=p_{X \otimes Y \otimes Z}F_{2}(X,Y\otimes Z)(q_{X}\otimes (\Pi_{Y\otimes Z}F_2(Y,Z)(q_Y\otimes q_Z)))\\
&\stackrel{(4)}{=}p_{X \otimes Y \otimes Z}F_{2}(X,Y\otimes Z)(q_{X}\otimes (F_2(Y\otimes Z)(\Pi_Y\otimes \Pi_Z)(q_Y\otimes q_Z)))\\
&=p_{X \otimes Y \otimes Z}F_{2}(X,Y\otimes Z)(\mathrm{id}_{F(X)}\otimes F_2(Y\otimes Z))(q_{X}\otimes(\Pi_Y\otimes \Pi_Z)(q_Y\otimes q_Z))\\
&\stackrel{(6)}{=}p_{X \otimes Y \otimes Z}F_{2}(X\otimes Y,Z)(F_2(X\otimes Y)\otimes \mathrm{id}_{F(Z)} )(q_{X}\otimes(\Pi_Y\otimes \Pi_Z)(q_Y\otimes q_Z))\\
&=p_{X \otimes Y \otimes Z}F_{2}(X\otimes Y,Z)(F_2(X\otimes Y)\otimes \mathrm{id}_{F(Z)} )(q_{X}\otimes\Pi_Yq_Y\otimes \Pi_Zq_Z)\\
&\stackrel{(8)}{=}p_{X \otimes Y \otimes Z}F_{2}(X\otimes Y,Z)(F_2(X\otimes Y)\otimes \mathrm{id}_{F(Z)} )(q_{X}\otimes q_Y\otimes q_Z)\\
&\stackrel{(9)}{=}p_{X \otimes Y \otimes Z}F_{2}(X\otimes Y,Z)(F_2(X\otimes Y)\otimes \mathrm{id}_{F(Z)} )(\Pi_{X}q_{X}\otimes \Pi_Yq_Y\otimes q_Z)\\
&=p_{X \otimes Y \otimes Z}F_{2}(X\otimes Y,Z)(F_2(X\otimes Y)(\Pi_X\otimes \Pi_Y)\otimes q_Z)(q_{X}\otimes q_Y\otimes \mathrm{id}_{\widetilde{F}(Z)})\\
&\stackrel{(11)}{=}p_{X \otimes Y \otimes Z}F_{2}(X\otimes Y,Z)(\Pi_{X\otimes Y}F_2(X\otimes Y)\otimes q_Z)(q_{X}\otimes q_Y\otimes \mathrm{id}_{\widetilde{F}(Z)})\\
&\stackrel{(12)}{=}p_{X \otimes Y \otimes Z}F_{2}(X\otimes Y,Z)(q_{X\otimes Y}p_{X\otimes Y}F_2(X\otimes Y)\otimes q_Z)(q_{X}\otimes q_Y\otimes \mathrm{id}_{\widetilde{F}(Z)})\\
&=(p_{X \otimes Y \otimes Z}F_{2}(X\otimes Y,Z)(q_{X\otimes Y}\otimes q_Z))(p_{X\otimes Y}F_2(X\otimes Y)\otimes \mathrm{id}_{\widetilde{F}(Z)})(q_{X}\otimes q_Y\otimes \mathrm{id}_{\widetilde{F}(Z)})\\
&=(p_{X \otimes Y \otimes Z}F_{2}(X\otimes Y,Z)(q_{X\otimes Y}\otimes q_Z))((p_{X\otimes Y}F_2(X\otimes Y)(q_{X}\otimes q_Y))\otimes \mathrm{id}_{\widetilde{F}(Z)})\\
&=\widetilde{F_2}(X\otimes Y,Z)(\widetilde{F_2}(X,Y)\otimes \mathrm{id}_{\widetilde{F}(Z)})\\
\end{align*}

Equalities $(4)$ and $(11)$ are due to the naturality of $F_2\colon F \otimes F \rightarrow F\otimes$ and the monoidal property of $\gamma$.
Equality $(6)$ follows from the axiom (\ref{diagbraidfunctor}) satisfied by $F$.
Equalities $(8)$, $(9)$, $(11)$ and $(12)$ are due to the definition of $\Pi$, $p$ and $q$.

Now, we show that $\widetilde{F}$ the axiom (\ref{diagmonofunctor2}) is satisfied.
For every object $X\in \C$,
\begin{align*}
\widetilde{F}(X,\un)(\mathrm{id}_{\widetilde{F}(X)}\otimes \widetilde{F_0})
&=p_{X\otimes \un}F_{2}(X,\un)(q_{X}\otimes q_{\un})(\mathrm{id}_{\widetilde{F}(X)}\otimes p_{\un}F_0)\\
&\stackrel{(2)}{=}p_{X}F_{2}(X,\un)(q_{X}\otimes q_{\un})(p_Xq_X\otimes p_{\un}F_0)\\
&=p_{X}F_{2}(X,\un)(\mathrm{\Pi}_{X}\otimes \mathrm{\Pi}_{\un})(q_X\otimes F_0)\\
&\stackrel{(4)}{=}p_{X}\mathrm{\Pi}_{X\otimes \un}F_{2}(X,\un)(q_X\otimes F_0)\\
&=p_{X}\mathrm{\Pi}_{X}F_{2}(X,\un)(\mathrm{id}_{F(X)}\otimes F_0)q_X\\
&\stackrel{(6)}{=}p_{X}\mathrm{\Pi}_{X}\mathrm{id}_{F(X)}q_X\\
&\stackrel{(7)}{=}p_{X}q_{X}p_{X}\mathrm{id}_{F(X)}q_X\\
&\stackrel{(8)}{=}\mathrm{id}_{\widetilde{F}(X)}\\
\end{align*}

Equality $(2)$, $(3)$, $(7)$ and $(8)$ are due to the equality $p_Xq_X=\mathrm{id}_{\mathrm{Im}(F(X))}$ and $\mathrm{\Pi}_X$.
Equality $(4)$ follows from the identity $F_2(X,Y)(\Pi_X\otimes \Pi_Y)=\Pi_{X\otimes Y}F_2(X,Y)$ which follows from the naturality of $F_2\colon F\otimes F \rightarrow F(\otimes) $ and the fact that $\gamma$ is monoidal.
Equality $(6)$ follows from the axiom (\ref{diagmonofunctor2}) which is satisfied by $F$.

In the same way, the axiom (\ref{diagmonofunctor3}) is satisfied by $\widetilde{F}$.
\end{proof}

Now, we have to extend notions of braided and symmetric functors to functors with anomaly. Suppose that the categories $\C$ and $\D$ are braided and denote indifferently by $\tau$ the braiding of the two categories.

\begin{definition}
A \emph{braided functor with anomaly}\index{braided functor with anomaly} is a monoidal functor with anomaly $(F,F_{2},F_{0},\gamma)$ such that for all objects $X$ and $Y$ of $\C$, the following diagram commutes:
\begin{align}\label{diagbraidfunctor}
\begin{tikzpicture}[description/.style={fill=white,inner sep=2pt},baseline=(current bounding  box.center)]
\matrix (m) [ampersand replacement=\&,matrix of math nodes, row sep=3em,column sep=2.5em, text height=1.5ex, text depth=0.25ex]
{ F(X)\otimes F(Y)\& \& \& \& \& F(Y)\otimes F(X) \\
F(X\otimes Y)\& \& \& \& \& F(Y\otimes X)\\ };
\path[->,font=\scriptsize]
(m-1-1) edge node[auto=left] {$(F(\mathrm{id}_{Y})\otimes F(\mathrm{id}_{X}))\tau_{F(X),F(Y)}$} (m-1-6)
		edge node[auto=right] {$F_{2}(X,Y)$} (m-2-1)
(m-2-1) edge node[auto=right] {$F(\tau_{X,Y})$} (m-2-6)
(m-1-6) edge node[auto] {$F_{2}(Y,X)$} (m-2-6);
\end{tikzpicture}
\end{align}
A \emph{symmetric monoidal functor with anomaly}\index{symmetric monoidal functor with anomaly} is a braided functor with anomaly between symmetric categories.
\end{definition}

\begin{remarks}
\item This lax definition of a braided functor allows us to consider braided functors with anomaly which don't send isomorphisms on isomorphisms.
\end{remarks}

\begin{lemma}\label{unilemma}
Let $F=(F, F_{2},F_{0},\gamma)$ be a braided functor with anomaly.
Then the unitalized functor $\widetilde{F}$ is a braided functor with anomaly.
\end{lemma}
\begin{proof}
Recall the monoidal structure of the unitalized monoidal functor with anomaly $\widetilde{F}=(F,F_2,F_0,\gamma)$: $\text{for } X,Y \in \C, \;\widetilde{F_2}(X,Y)=p_{X\otimes Y}F_2(X,Y)q_{X}\otimes q_{Y}$, and, as $p$ and $q$ are monoidal natural transformations,
 $\widetilde{F_{2}}(X,Y)p_{X}\otimes p_{Y}=p_{X\otimes Y}F_{2}(X,Y)$ and
$q_{X\otimes Y}\widetilde{F_{2}}(X,Y)=F_{2}(X,Y)q_{X}\otimes q_{Y}$.
Now, we can compute $\widetilde{F}(\tau_{X,Y})\widetilde{F_{2}}(X,Y)$:
\begin{align*}
\widetilde{F}(\tau_{X,Y})\widetilde{F_{2}}(X,Y)&=
p_{Y\otimes X}F(\tau_{X,Y})q_{X\otimes Y}p_{X\otimes Y}F_{2}(X,Y)q_{X}\otimes q_{Y}\\
&=p_{Y\otimes X}F(\tau_{X,Y})q_{X\otimes Y}\widetilde{F_{2}}(X,Y)(p_{X}\otimes p_{Y})(q_{X}\otimes q_{Y})\\
&=p_{Y\otimes X}F(\tau_{X,Y})F_{2}(X,Y)(q_{X}\otimes q_{Y})(p_{X}\otimes p_{Y})(q_{X}\otimes q_{Y})\\
&=p_{Y\otimes X}F(\tau_{X,Y})(q_{X}\otimes q_{Y})\\
&=
p_{Y\otimes X}F_{2}(Y,X)(F(\id_Y)\otimes F(\id_X))\tau_{F(X),F(Y)}(q_{X}\otimes q_{Y})\\
&\stackrel{(5)}{=}
\gamma_{X,X}^{-1}\gamma_{Y,Y}^{-1}p_{Y\otimes X}F_{2}(Y,X)(q_{Y}p_{Y}\otimes q_{X}p_{X})\tau_{F(X),F(Y)}(q_{X}\otimes q_{Y})\\
&=
\gamma_{X,X}^{-1}\gamma_{Y,Y}^{-1}p_{Y\otimes X}F_{2}(Y,X)(q_{Y}\otimes q_{X})(p_{Y}\otimes p_{X})\tau_{F(X),F(Y)}(q_{X}\otimes q_{Y})\\
&=
\gamma_{X,X}^{-1}\gamma_{Y,Y}^{-1}p_{Y\otimes X}F_{2}(Y,X)(q_{Y}\otimes q_{X})(p_{Y}\otimes p_{X})(q_{Y}\otimes q_{X})\tau_{\widetilde{F}(X),\widetilde{F}(Y)}\\
&=
p_{Y\otimes X}F_{2}(Y,X)(q_{Y}\otimes q_{X})\tau_{\widetilde{F}(X),\widetilde{F}(Y)}\\
&=
p_{Y\otimes X}F_{2}(Y,X)(q_{Y}\otimes q_{X})(\gamma_{Y,Y}^{-1}\id_{Y}\otimes \gamma_{X,X}^{-1}\id_{X})\tau_{\widetilde{F}(X),\widetilde{F}(Y)}\\
&=
\widetilde{F_2}(Y,X)(\widetilde{F}(\id_Y)\otimes \widetilde{F}(\id_X))\tau_{\widetilde{F}(X),\widetilde{F}(Y)}\\
\end{align*}
using that $F(\id_X)=\gamma_{X,X}^{-1}q_Xp_X$, $p_Xq_X=\id_{\widetilde{F}(X)}$, the braiding is natural and that $F$ is a braided functor with anomaly in equality $(5)$.
\end{proof}

\label{defcob}
Before defining TQFTs with anomaly, let us briefly recall the definition of the category of cobordims. Let $n$ be a non-negative integer. A \emph{$n$-manifold} means a topological manifold of dimension $n$ and the empty set $\emptyset$ is a $n$-manifold for any $n$. The \emph{category of $n$-dimensional cobordisms} $\Cobn$\index{\(\Cobn\)} is defined as follows. The objects of $\Cobn$ are closed oriented $(n-1)$-manifolds as objects. A morphism from a $(n-1)$-manifold $\Sigma$ to a $(n-1)$-manifold $\Sigma'$ is represented by a pair $(M,h)$ where $M$ is a compact oriented $n$-manifold and $h$ is a orientation preserving homeomorphism between $\overline{\Sigma}\bigsqcup\Sigma'$ and $\partial M$ where $\overline{\Sigma}$ represents the manifold $\Sigma$ with the opposite orientation. Two such pairs $(M,h \colon \overline{\Sigma}\bigsqcup\Sigma' \rightarrow \partial M )$ and $(N,k \colon \overline{\Sigma}\bigsqcup\Sigma' \rightarrow \partial N)$ represent the same morphism in $\Cobn$ if there exists a preserving-orientation homeomorphism $f\colon M \rightarrow N$ such that $k=fh$. The composition in the category $\Cobn$ is given by the gluing of two $n$-cobordisms: the composition of the two morphisms $(M,h) \colon (\Sigma,\phi_{\Sigma}) \rightarrow (\Sigma',\phi_{\Sigma'})$ and $(N,k) \colon (\Sigma',\phi_{\Sigma'}) \rightarrow (\Sigma'',\phi_{\Sigma''})$ is represented by the morphism $(L,g)$ where $M$ is the gluing of $M$ on $N$ along $\Sigma'$ given by the gluing homeomorphism $kh^{-1} \colon h(\Sigma')\rightarrow k(\Sigma')$ and $g$ is the homeomorphism $h_{\Sigma}\sqcup k_{\Sigma''}$. The identity morphism of the $(n-1)$-manifold $\Sigma$ is represented by the $n$-cobordism $(\Sigma\times [0,1], c \colon \overline{\Sigma}\bigsqcup\Sigma \rightarrow \Sigma\times \{0\}\bigsqcup \Sigma\times \{1\})$ where $c|_{\overline{\Sigma}}(x)=(x,0)$ and $c|_{\Sigma}(x)=(x,1)$. The category $\Cobn$ is a symmetric
monoidal category with tensor product given by disjoint union, unit object is the empty $(n-1)$-manifold and the symmetric braiding $\tau_{\Sigma,\Sigma'}$ between two $(n-1)$-manifolds $\Sigma$ and $\Sigma'$ is represented by the $n$-cobordism $\big((\Sigma\sqcup \Sigma')\times [0,1],d\colon\overline{\Sigma\sqcup \Sigma'}\bigsqcup \Sigma'\sqcup \Sigma\rightarrow (\Sigma\sqcup\Sigma')\times\{0\}\bigsqcup (\Sigma'\sqcup \Sigma)\times\{1\}\big)$ where $d|_{\overline{\Sigma\sqcup \Sigma'}}(x)=(x,0)$ and $d|_{\Sigma'\sqcup \Sigma}(x)=(x,1)$.

Now, we are able to give the main object definition we will construct for $n=3$.

\begin{definition}
A \emph{$n$-dimensional Topological Quantum Field Theory with anomaly} (TQFT with anomaly) is a symmetric monoidal unital functor with anomaly from $\Cobn$ to a symmetric monoidal category $\sym$.
\end{definition}

\begin{remarks}
\item A classical $n$-dimensional TQFT is a particular case of a TQFT with anomaly where $\sym=\Vect$ and anomaly equal to $1_{\kk}$.
\end{remarks}

\subsection{Anomaly lifting} If someone needs to work with a more classical TQFT, there is a general process to lift the anomaly from functors with anomaly.

Let $(F,\gamma)$ be a functor with anomaly from a category $\C$ to a category $\D$ and let $\nu$ be a map that associates to any morphism $f$ of $\C$ an invertible scalar $\nu_{f}\in \kk$. First, define a new functor with anomaly $(F^{\nu},\gamma^{\nu})$ called the \emph{$\nu$-rescaling of $(F,\gamma)$} by:
\[F^{\nu}(X)=X \text{ and } F^{\nu}(f)=\nu_{f}F(f),\]
for any objects $X$ of $\C$ and any morphism $f$ of $\C$. The anomaly $\gamma^{\nu}$ is then given by $\frac{\nu_{gf}}{\nu_f\nu_g}\nu_{g,f}$ on a pair of composable morphisms $(g,f)$ of $\C$.

There is a canonical way to remove the anomaly of $(F,\gamma)$ by modifying less as possible the category $\C$.
Define by $\overline{\C}$ the category:
\begin{itemize} 
\item whose objects are the same than those of $\C$;
\item a morphism of $\mathrm{Hom}_{\overline{\C}}(X,Y)$ is a pair $(f,x)$ where $f\in \mathrm{Hom}_{\C}(X,Y)$ and $x\in \kk^{\times}$;
\item the composition between $(f,x)\colon X\rightarrow Y$ and $(g,y)\colon Y\rightarrow Z$ is given by $(g\circ_{\C}f,\gamma_{g,f}^{-1}xy)$.
\end{itemize} 
Now, define the functor (without anomaly) $\overline{F}\colon \overline{\C}\rightarrow \D$ by:
\[\overline{F}(X)=F(X) \text{ and } \overline{F}(f,x)=xF(f)\]
where $X$ is any object of $\overline{\C}$ and $(f,x)$ is any morphism of $\overline{\C}$.
If $(f,x)$ is a morphism of $\overline{\C}$, set $\nu_{(f,x)}=x$. 
Denote by $U\colon \overline{\C}\rightarrow \C$ the forgetful functor between $\overline{\C}$ and $\C$.
Then we have the following commutative diagram of functors with anomaly:
\begin{align}\label{diagliftanomaly}
\begin{tikzpicture}[description/.style={fill=white,inner sep=2pt},baseline=(current bounding  box.center)]
\matrix (m) [ampersand replacement=\&,matrix of math nodes, row sep=3em,column sep=2.5em, text height=1.5ex, text depth=0.25ex]
{ \overline{\C}\& \&   \\
\C\& \&  \D  \\ };
\path[->,font=\scriptsize]
(m-1-1) edge node[auto=left] {$\overline{F}^{\nu}$} (m-2-3)
(m-1-1)	edge node[auto=right] {$U$} (m-2-1)
(m-2-1) edge node[auto=right] {$(F,\gamma)$} (m-2-3);
\end{tikzpicture}
\end{align}

Note that the triplet $(\overline{\C},U,\overline{F}^{\omega})$ has the following property:
for all triplets $(\widetilde{\C},E,G^{\alpha})$ such that $E\colon \widetilde{\C}\rightarrow \C$ is an equivalence of category, $G^{\alpha}$ is a functor with anomaly which is the $\alpha$-rescaling of a functor $G\colon \widetilde{\C}\rightarrow \D$ and $FE=G^{\alpha}$, there exists an equivalence of category $H$ such that the following diagram (\ref{diagunivanomaly}) commutes

\begin{align}\label{diagunivanomaly}
\begin{tikzpicture}[description/.style={fill=white,inner sep=2pt},baseline=(current bounding  box.center)]
\matrix (m) [ampersand replacement=\&,matrix of math nodes, row sep=3em,column sep=2.5em, text height=2ex, text depth=0.25ex]
{ \& \overline{\C} \&  \\
 \& \widetilde{\C}\&   \\
\C \& \&  \D  \\ };
\path[->,font=\scriptsize]
(m-1-2)	edge[bend right] node[auto=right] {$U$} (m-3-1)
(m-1-2)	edge[bend left] node[auto=left] {$\overline{F}^{\nu}$} (m-3-3)
(m-1-2)	edge node[auto=right] {$H$} (m-2-2)
(m-2-2)	edge node[auto=right] {$E$} (m-3-1)
(m-2-2)	edge node[auto=left] {$G^{\alpha}$} (m-3-3)
(m-3-1) edge node[auto=right] {$(F,\gamma)$} (m-3-3);
\end{tikzpicture}
\end{align}

\section{Topological preliminaries}
In this section, we define all topological objects we need to define our TQFT. We briefly recall the definition of Turaev of
a ribbon graph (see \cite{Turaev1994}) before specifying the main combinatorial object of this thesis, that is, \emph{ribbon cobordism tangles}. These tangles give us a way to present 3-cobordisms by surgery and we recall generalized Kirby moves to explicit when two different ribbon cobordism tangles present two homeomorphic cobordisms.

\subsection{Ribbon graphs}
An \emph{arc} is an embedding of the square $[0,1]\times[0,1]$ in $\mathbb{R}^{3}$. The image of $[0,1]\times \{0\}$ and of $[0,1]\times \{1\}$ are called \emph{bases} of the arc whereas the image of $[0,1]\times \{\frac{1}{2}\}$ is called the \emph{core} of the arc.
A \emph{coupon} is an arc with a distinguished base called the \emph{bottom base} (the other base is the \emph{top base}). 
A \emph{closed component} is an embedding of the the surface $\mathbb{S}^{1}\times [0,1]$ in $\mathbb{R}^{3}$. The image of $\mathbb{S}^{1}\times \{\frac{1}{2}\}$ is the \emph{core} of the closed component.

A \emph{ribbon graph} with $k$ bottom endpoints and $l$ top endpoints is an oriented surface $G$ embedded in $\mathbb{R}^{2}\times [0,1]$ which is a finite disjoint union of arcs, coupons, and closed components such that
\begin{itemize}
\item the set $G\cap \mathbb{R}^{2}\times \{0\}$ (respectively $G\cap \mathbb{R}^{2}\times \{1\}$)  is the union of the $k$ (respectively $l$) disjoint segments $[(1,1,0),(1,2,0)],\ldots,[(1,2k-1,0),(1,2k,0)]$ (respectively $[(1,1,1),(1,2,1)],\ldots,[(1,2l-1,1),(1,2l,1)]$)  which belong to some arcs of $G$ and the orientation of $G$ near these segments is given by the normal vector $(1,0,0)$;
\item all other bases of arcs lie on bases of coupons;
\item the core of arcs and closed components are oriented.
\end{itemize}
For more details on ribbon graphs, see \cite{Turaev1994}, Chapter I.2.1.

A \emph{ribbon tangle} is a ribbon graph with $k$ bottom endpoints and $l$ top endpoints without coupons.

\begin{figure}[H]
\centering
\resizebox{.8\linewidth}{!}{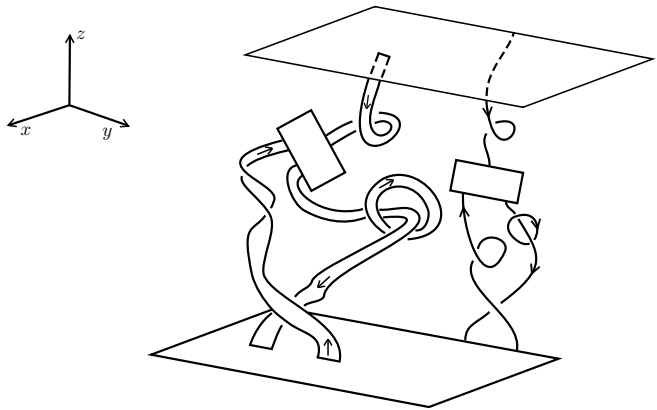}
\caption{A ribbon graph and one of its diagrams in the plane $\{0\}\times \mathbb{R}\times \mathbb{R}$.}
\label{ribbonprojection}
\end{figure}

A \emph{diagram} of a ribbon graph $G$ is a projection of the coupons and the core of arcs and closed components of $G$ in the plane $\{0\}\times \mathbb{R}\times \mathbb{R}$ such that the crossing have only double points and the orientation of a coupon is the orientation of $\{0\}\times\mathbb{R}\times \mathbb{R}$; we distinguish the overcrossing and the undercrossing in such a case. Except on Figure~\ref{examplegraph}, the bottom base of a coupon is parallel to the line $\{0\}\times\mathbb{R}\times \{0\}$ and will be considered as lower than the top base. An example of a projection of a ribbon graph is given in Figure~\ref{ribbonprojection}.

To rebuild the ribbon graph starting from one of its diagrams, just thicken the cores of arcs and closed components of the diagram in the plane $\{0\}\times \mathbb{R}\times \mathbb{R}$. We say that diagrams are represented with convention of \emph{blackboard framing}.

By an \emph{isotopy of ribbon graphs}, we mean an orientation preserving isotopy in  $\mathbb{R}^{2}\times [0,1]$ constant on the boundary segments and preserving the splitting into arcs, coupons and closed components as well as the orientation of cores.
Recall that two diagrams represent the same isotopy class of a ribbon tangle if and only if one can be obtained from the other by deformation (planar isotopies) and a finite sequence of ribbon Reidemeister moves. See \cite{Turaev1994} for more details on isotopies and ribbon Reidemeister moves.

Finally, remember that you can compose two isotopy classes of ribbon graphs by juxtaposing them when the number of top endpoints of the first coincides with the number of bottom endpoints of the second. Moreover, you can obtain a new ribbon graph by putting two ribbon graphs side by side. These operations turn the set of isotopy class of ribbon graphs into a monoidal category (see \cite{Turaev1994}).

Now, we start to define all objects we need to represent a 3-cobordism with entrance boundary of multigenus $g$ and exit boundary of multigenus $h$. In such representations, we denote by $n$ the number of closed components; we will use them as surgery components.

\subsection{Ribbon (\underline{g},n,\underline{h})-graphs}

Let $\ung=(g_{1},\ldots,g_{r})$ and $\unh=(h_{1},\ldots,h_{s})$ be two tuples of non-negative integers, $n$ be a non-negative integer, and denote by ${|\ung|=\sum_{i=1}^{r}g_i}$.
By a \emph{ribbon $(\ung,n,\unh)$-graph}, we shall mean a ribbon graph $G \subset \mathbb{R}^{2}\times [0,1]$ consisting of $n$ closed components, $r$ ordered coupons called \emph{entrance coupons}, $s$ ordered coupons called \emph{exit coupons}, and $|\ung|+|\unh|$ arcs based on coupons. Moreover we assume that for all $1\leq i\leq r$ and all $1\leq j\leq s$ the $i$th entrance coupon has $2g_{i}$ top endpoints and no bottom endpoints and the $j$th exit coupon has $2h_{j}$ bottom endpoints and no top endpoints such that:
\begin{itemize}
\item for $1\leq k \leq g_i$, an arc joins the $(2k-1)$th and the $2k$th top endpoints of the $i$th entrance coupon and its core is oriented from the $2k$th top endpoint to the $(2k-1)$th top endpoint;
\item for $1\leq k \leq h_j$, an arc joins the $(2k-1)$th and the $2k$th bottom endpoints of the $j$th exit coupon and its core is oriented from the $(2k-1)$th bottom endpoint to the $2k$th bottom endpoint.
\end{itemize}
A closed component of a ribbon $(\ung,n,\unh)$-graph is called a \emph{sugery component}. An arc based on an entrance coupon (respectively exit coupon) is an \emph{entrance component}(respectively \emph{exit component}). A connected component of an entrance (resp. exit) coupon constitutes an \emph{entrance boundary component}(resp. \emph{exit boundary component}).
A ribbon $(\ung,n,\unh)$-graph is represented by a diagram with blackboard framing.

\begin{figure}[H]
  \centering
  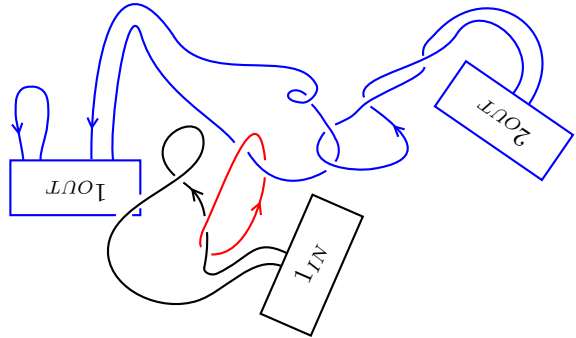
 \caption{A ribbon $((1),1,(2,1))$-graph.}
 \label{examplegraph}
\end{figure}

For such an example of a ribbon $((1),1,(2,1))$-graph, see Figure~\ref{examplegraph} where there is one surgery component (red), one entrance component which is a part of the unique entrance boundary component (black), and three exit components separated on the two exit boundary components (blue). 

We denote by 
$Graph(\underline{g},n,\underline{h})$\index{$Graph(\underline{g},n,\underline{h})$} the set of all isotopy classes of $(\ung,n,\unh)$-graphs, and by \[Graph=\bigsqcup_{(\underline{g},n,\underline{h})}Graph(\underline{g},n,\underline{h}).\]

\subsection{Ribbon cobordism tangles}\label{defcobtangle}

Let $\ung=(g_{1},\ldots,g_{r})$ and $\unh=(h_{1},\ldots,h_{s})$ be two tuples of non-negative integers and $n$ be a non-negative integer.
By a \emph{ribbon $(\ung,n,\unh)$-cobordism tangle}, we shall mean a ribbon tangle $T \subset \mathbb{R}^{2}\times [0,1]$ with $2|\ung|$ bottom endpoints and $2|\unh|$ top endpoints consisting of $n$ closed oriented components called \emph{surgery components}, $|\ung|$ arcs called \emph{entrance components}\ and $|\unh|$ arcs called \emph{exit components} such that:
\begin{itemize}
\item for $1\leq k \leq |\ung|$, the $k$th entrance component joins the $(2k-1)$th and the $2k$th bottom endpoints and its core is oriented from the $2k$th bottom endpoint to the $(2k-1)$th bottom endpoint;
\item for $1\leq k \leq |\unh|$, the $k$th exit component joins the $(2k-1)$th and the $2k$th top endpoints and its core is oriented from the $(2k-1)$th top endpoint to the $2k$th top endpoint.
\end{itemize}
For every $1\leq i\leq r$, the disjoint union of the $k$th entrance component for $\left(\sum\limits_{j=1}^{i-1}g_{j}\right)+1 \leq k \leq \sum\limits_{j=1}^{i}g_{j}$ is called the $i$th \emph{entrance boundary component}\index{cobordism tangle!entrance boundary component} of the cobordism tangle. For every $1\leq i\leq s$, the disjoint union of the $k$th top arcs for $\left(\sum\limits_{j=1}^{i-1}h_j\right)+1 \leq k \leq \left(\sum\limits_{j=1}^{i}h_j\right)$ is called the $i$th \emph{exit boundary component}\index{cobordism tangle!exit boundary component} of the cobordism tangle.

A ribbon cobordism tangle is represented by a diagram with blackboard framing.

\begin{figure}[H]
  \centering
  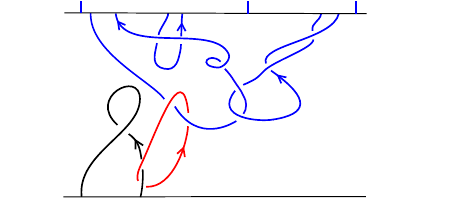
 \caption{A ribbon $((1),1,(2,1))$-cobordism tangle.}
 \label{exampletangle}
\end{figure}

For such an example of a ribbon $((1),1,(2,1))$-cobordism tangle, see Figure~\ref{exampletangle} where there is one surgery component (red), one entrance component which forms the only one entrance boundary component (black) and three exit components separated on two exit boundary components materialized by vertical little segments (blue). 

We denote 
by $Tang^{Cob}(\underline{g},n,\underline{h})$\index{$Tang^{Cob}(\underline{g},n,\underline{h})$} the set of all isotopy classes of $(\ung,n,\unh)$-cobordism tangles, and by \[Tang^{Cob}=\bigsqcup_{(\underline{g},n,\underline{h})}Tang^{Cob}(\underline{g},n,\underline{h}).\]

There is a surjective map $Gr \colon Tang^{Cob}\rightarrow Graph$ defined on an isotopy class $T$ of a $(\ung,n,\unh)$-cobordism tangle where $\ung=(g_1,\ldots,g_r)$ and $\unh=(h_1,\ldots,h_s)$ by:

\begin{align}\label{surjectivemapgr}
Gr (T)= \left(\exitbox{$h_1$}\otimes \ldots\otimes \exitbox{$h_s$} \right)\circ T\circ \left(\enterbox{$g_1$}\otimes \ldots\otimes \enterbox{$g_r$} \right)
\end{align}
where a coupon colored by a positive integer $m$ means that there is $m$ couples of oriented arrows attached on the coupon as specified in the definition (\ref{surjectivemapgr}) of the map $Gr$.

\begin{figure}[!h]
  \centering
  \resizebox{.8\linewidth}{!}{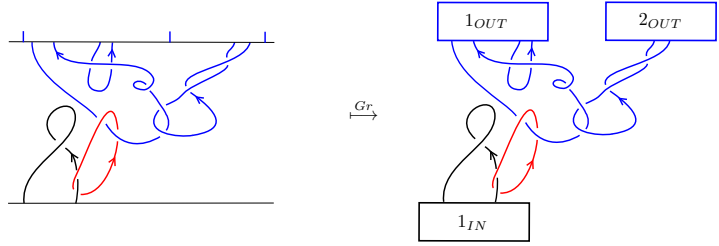}
 \caption{The action of map $Gr$ from $Tang^{Cob}$ to $Graph$ on a $((1),1,(2,1))$-cobordism tangle}
 \label{mapgraph}
\end{figure}

In order to define an isotopy invariant of cobordism tangles, we define a useful ribbon tool called \emph{ribbon opentangles}.

\subsection{Ribbon opentangles}\label{opentangle}

Let $\ung=(g_{1},\ldots,g_{r})$ and $\unh=(h_{1},\ldots,h_{s})$ be two tuples of non-negative integers and $n$ be a non-negative integer.
By a \emph{ribbon $(\ung,n,\unh)$-opentangle}, we shall mean a ribbon tangle $T \subset \mathbb{R}^{2}\times [0,1]$ with $2(|\ung|+n+|\unh|)$ bottom endpoints and $2|\unh|$ top endpoints, consisting of $N=|\ung|+n+2|\unh|$ arcs components without any closed component, such that:
\begin{itemize}
\item for $1\leq k \leq |\ung|+n$, the $k$th arc joins the $(2k-1)$th and the $2k$th bottom endpoints and its core is oriented from the $(2k-1)$th bottom endpoint to the $2k$th bottom endpoint;
\item for $|\ung|+n+1\leq k \leq N$, the $k$th arc joins the $(k+|\ung|+n)$th bottom endpoint with the $(k-|\ung|-n)$th top endpoint and its core is oriented upwards if $k-|\ung|-n$ is odd and donwards if not.
\end{itemize}
The $|\ung|$ first arcs are the \emph{entrance components}, the $n$ following arcs are the \emph{surgery components} and the last $|\unh|$ couples of consecutive arcs are the \emph{exit components}.
For every $1\leq i\leq r$, the disjoint union of the $k$th arcs for $\left(\sum\limits_{j=1}^{i-1}g_{j}\right)+1 \leq k \leq \sum\limits_{j=1}^{i}g_{j}$ is called the $i$th \emph{entrance boundary component} of the opentangle. For every $1\leq i\leq s$, the disjoint union of the $k$th arcs for $\left(|\ung|+n+\sum\limits_{j=1}^{i-1}h_j\right)+1 \leq k \leq |\ung|+n+2\left(\sum\limits_{j=1}^{i}h_j\right)$ is called the $i$th \emph{exit boundary component} of the opentangle.

A ribbon opentangle is represented by a diagram with blackboard framing.

\begin{figure}[!h]
 \centering
 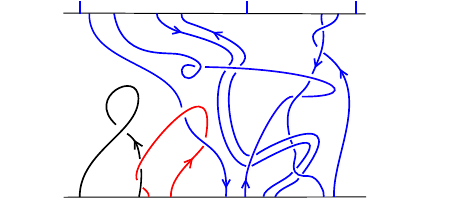
 \caption{A ribbon $((1),1,(2,1))$-opentangle.}
 \label{exampleopentangle}
\end{figure}

For such an example of a ribbon $((1),1,(2,1))$-opentangle, see Figure~\ref{exampleopentangle} where there is one surgery component (red), one entrance component which forms the only entrance boundary component (black)  and three exit components separated on two exit boundary components materialized by vertical little segments (blue).

We denote 
by $Otang(\underline{g},n,\underline{h})$\index{$Otang(\underline{g},n,\underline{h})$} the set of all isotopy classes of $(\ung,n,\unh)$-opentangles, and by \[Otang=\bigsqcup_{(\underline{g},n,\underline{h})}Otang(\underline{g},n,\underline{h}).\]

There is a surjective map $U \colon Otang\rightarrow Tang^{Cob}$ whose restriction on $Otang(\ung,n,\unh)\rightarrow Tang^{Cob}(\ung,n,\unh)$ is also surjective and is given by the closure of surgery components and the bottom closure of exit components of a class of opentangles, that is,
\begin{align}\label{surjectivemap}
U(O)=O\circ (\downarrow \uparrow ^{\otimes|\ung|}\otimes\languette^{\phantom{o}\otimes n+|\unh|})
\end{align}
where $O$ is an isotopy class of a $(\ung,n,\unh)$-opentangle.
\begin{figure}[!h]
  \centering
  \resizebox{.8\linewidth}{!}{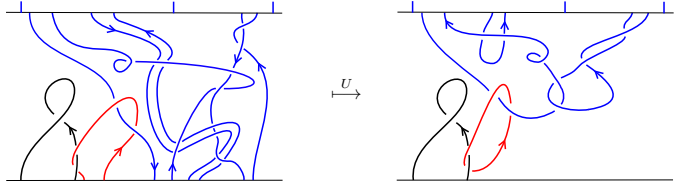}
 \caption{The map $U$ from $Otang$ to $Tang^{Cob}$.}
 \label{mapgraph}
\end{figure}

\begin{figure}[!h]\index{$BA$}
\centering
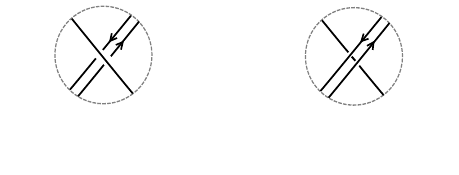
\caption{The move $BA$ on diagrams of opentangles.}
\label{equivalence}
\end{figure}

\begin{figure}[!h]\index{$ESC$}
\centering
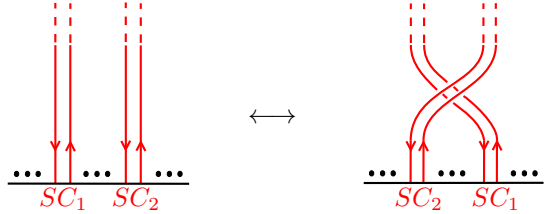
\caption{The move $ESC$ on diagrams of opentangles.}
\label{moveesc}
\end{figure}

\begin{figure}[!h]\index{$ROT$}
\centering
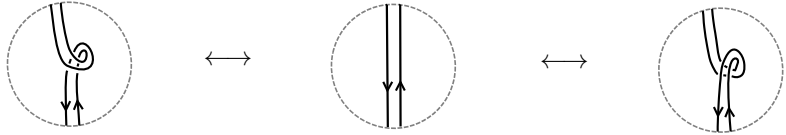
\caption{The move $ROT$ on diagrams of opentangles.}
\label{moverot}
\end{figure}

Since the map $U$ is surjective, it induces a bijection between the set $Tang^{Cob}(\ung,n,\unh)$ and the quotient set $Otang(\ung,n,\unh)/\sim$ where two isotopy classes of opentangles are equivalent if and only if they have the same image under $U$. There is a diagrammatical characterization of this equivalence relation: two isotopy classes of opentangles $O_1$ and $O_2$ are equivalent if and only if a diagram of $O_1$ and a diagram of $O_2$ are related by a finite sequence of planar isotopies, ribbon Reidemeister moves (see \cite{Turaev1994}) and three additionnal moves (see \cite{Lyu}):  moves of type $BA$ ("below-above") defined on Figure \ref{equivalence}, moves of type $ESC$ ("exchange-surgery-components") defined on Figure \ref{moveesc} and moves of type "ROT" ("rotation") defined on Figure \ref{moverot}.

To construct our TQFT, we need to start with a presentation of $3$-cobordisms by cobordism tangles, generalizing surgery of 3-manifolds to 3-cobordisms. Then, we will be able to extend Kirby calculus.

\subsection{Surgery of $3$-cobordisms and presentation by cobordism tangles}
\label{parametrized}
For the definition of the category of 3-cobordims $\Cob$, see Section~\ref{defcob}.
Let $\ung=(g_1,\ldots,g_r)$ be a $r$-tuple of non-negative integers and $g$ be a non-negative integer. If $\Sigma_{g_1},\ldots,\Sigma_{g_r}$ are $r$ connected surfaces of respective genus $g_1,\ldots,g_r$, the \emph{multigenus} of the surface $\Sigma_{\ung}:=\bigsqcup_{i=1}^{r}\Sigma_{g_i}$ is $\ung$. Denote by $S_{g}$\index{\(S_{g}\)} the canonical oriented connected and closed surface of genus $g$ and by $S_{\ung}$ the ordered disjoint union of connected canonical surfaces $\bigsqcup_{i=1}^{r} S_{g_{i}}$. Then we define the \emph{category of $3$-dimensional parametrized cobordisms} $\Cobp$  which contains parametrized surfaces as objects that are pairs $(\Sigma,\phi_{\Sigma} \colon \Sigma \rightarrow S_{\ung})$ with $\Sigma$ a closed oriented surface of multigenus $\ung$ and $\phi_{\Sigma}$ a orientation preserving homeomorphism. Morphisms and composition are defined exactly in the same way that in $\Cob$ and there exists an equivalence of category given by the forgetful functor $\Cobp \rightarrow \Cob$. We denote by $\Cobp(\ung,\unh)$ the set of all parametrized cobordisms from a surface of multigenus $\ung$ to a surface of genus $\unh$ and by $\Cobc(\ung,\unh)$ the subset of $\Cobp(\ung,\unh)$ of connected parametrized cobordisms.

A closed $3$-manifold can be presented by a framed link. Using this result, we describe how to generalize this combinatorial presentation to $3$-cobordims (detailed results and proofs can be found in \cite{Turaev1994}, Chapter I). First, we recall how to {\it present a 3-manifold by a link} and what is {\it the surgery on this link}. Let $L$ be a $n$-components framed link embedded in $\Sp^{3}$ and denote by $V(L)$ a tubular neighborhood of $L$ in $\Sp^{3}$. The boundary of the $3$-manifold $\Sp^{3}\backslash V(L)$ is then homeomorphic to $n$ disjoint canonical tori (or handlebodies) $\DD^{2}\times \Sp^{1}$. We define the $3$-manifold $\Sp^3_{L}$\index{$\Sp^3_{L}$} obtained by surgery of $\Sp^{3}$ along the link $L$ by:
\[\Sp^3_{L}=(\Sp^{3}\backslash V(L)) \bigsqcup_{\phi \colon \partial(\Sp^{3}\backslash V(L)) \rightarrow \sqcup_{i=1}^{n}(\Sp^{1}\times \Sp^{1})_{i}} \sqcup_{i=1}^{n}(\DD^{2}\times \Sp^{1})_{i}\]
where $\phi$ is a homeomorphism between tori that exchanges meridian and parallel and for every $1\leq i\leq n$, the $3$-manifold $(\DD^{2}\times \Sp^{1})_{i}$ is a copy of the canonical torus $\DD^{2}\times \Sp^{1}$.
A result of Lickorish proved in \cite{Lickorish1997} claims that if $M$ be a closed oriented connected $3$-manifold, then there exists a framed link $L$ in $\Sp^{3}$ such that $M$ is homeomorphic to $\Sp^{3}_{L}$.\label{Lickorish}

This last result will allow us to {\it present $3$-cobordisms by ribbon graphs and tangles}. Let $(M, h \colon \overline{\Sigma}\bigsqcup\Sigma' \rightarrow \partial M)$ be a connected $3$-cobordism between the parametrized surfaces $(\Sigma,\phi_{\Sigma} \colon \Sigma \rightarrow S_{\ung})$ and $(\Sigma',\phi_{\Sigma'} \colon \Sigma' \rightarrow S_{\unk})$ where $\ung$ and $\unk$ are respectively a $r$-tuple and a $s$-tuple of nonnegative integers. Denote by $H_{g}$ the $3$-dimensional handlebody of genus $g$ bounded by the canonical closed surface $S_{g}$ and when $\ung$ is a $r$-tuple of nonnegative integers, denote by $H_{\ung}$ the ordered disjoint union of handlebodies $H_{g_{1}},\ldots,H_{g_{r}}$. Then, define by $\widetilde{M}$ the closed $3$-manifold:
\[ \widetilde{M}=H_{\ung}\bigsqcup_{h_{-}\circ\phi_{\Sigma}^{-1}} M \bigsqcup_{\phi_{\Sigma'}\circ h_{+}} H_{\unk} \]
where $h_{-}=h|_{\overline{\Sigma}}$ and $h_{+}=h|_{\Sigma'}$.
Now, according to the surgery theorem of Lickorish, there exist a $n$-component framed link $L$ and a homeomorphism $f$ between $\widetilde{M}$ and $\Sp^{3}_{L}$. Moreover, every canonical handlebody $H_{g}$ of genus $g$ is the tubular neighborhood of an oriented surface which is a ribbon graph $G_{g}$ composed by one coupon and $g$ handles based on the same side of the coupon. We can suppose that, for every $1 \leq i \leq r$ and every $1\leq i \leq s$ , $f(G_{g_{i}})\cap L = \emptyset$ and $f(G_{k_{j}})\cap L = \emptyset$ because the ribbon link $L$ and the image by $f$ of graphs $G_{g_i}$ and $G_{g_k}$ are objects of codimension 1 in $\Sp^{3}$. Consider the ribbon graph in $\Sp^{3}$
\[ G_{\ung,n,\unh}=G_{g_1}\sqcup \ldots \sqcup G_{g_r}\bigsqcup L \bigsqcup G_{k_1}\sqcup \ldots \sqcup G_{k_s}.\]
By isotopy, pull down the graphs $G_{g_1},\ldots,G_{g_r}$ and pull up the graphs $G_{k_1},\ldots,G_{k_s}$ and then cut them: the resulting ribbon $(\ung,n,\unk)$-cobordism tangle is a presentation of the $3$-cobordism $(M, h \colon \overline{\Sigma}\bigsqcup\Sigma' \rightarrow \partial M)$.

\begin{figure}[H]
\centering
\resizebox{.4\linewidth}{!}{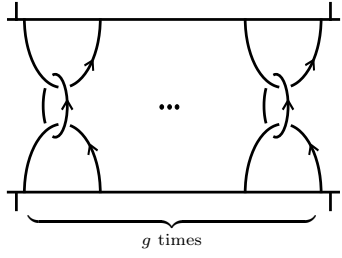}\\
\caption{A cobordism tangle which presents the cylinder of a surface of genus $g$.}
\label{cylinder}
\end{figure}

For example, the cylinder $\id_{\Sigma_g}=(\Sigma_g\times [0,1],\id_{\overline{\Sigma_g}\sqcup \Sigma_g})$ over a surface $\Sigma_g$ of genus $g$ is presented by the ribbon $(g,g,g)$-cobordism tangle given on Figure~\ref{cylinder}.

Henceforth, we know how to present a 3-cobordism by a 2-dimensional object as a cobordism tangle. To know when two cobordism tangles present the same cobordisms, we need to {\it extend Kirby calculus}.
Let $T$ be a $(\ung,n,\unh)$-cobordism tangle in $\Sp^{3}$ with $\ung=(g_1,\ldots,g_r)$ and $\unh=(h_1,\ldots,h_s)$ and denote by $L$ the link defined by the disjoint union of the $n$ surgery components of $T$. Consider the ribbon $(\ung,n,\unh)$-graph $G=Gr(T)$ as defined in map (\ref{surjectivemapgr}) where $r+s$ coupons have been attached to entrance and exit components of $T$. Then do the surgery on $\Sp^{3}$ along $L$ to obtain a closed connected oriented 3-manifold $\Sp^{3}_{L}$ with $r+s$ disjoint embedded ribbon graphs of type \handlebox. Take tubular neighborhoods $N_1,\ldots,N_r$ of the $r$ entrance boundary components of $G$ 
and tubular neighborhoods $N_{r+1},\ldots,N_{r+s}$ of the $s$ exit boundary components of $G$ in $\Sp^{3}_{L}$.
Then we obtain a 3-cobordism \[M_{T}=\Sp^{3}_{L}\backslash (N_1\sqcup\ldots \sqcup N_{r+s})\] from $\partial(N_1\sqcup \ldots \sqcup N_r)$ to $\partial(N_{s+1}\sqcup\ldots\sqcup N_{r+s})$ with a parametrization of these two boundary components by canonical surfaces. 

\begin{figure}[H]
\centering
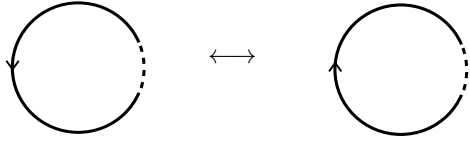
\caption{The move $SO$ on cobordism tangles.}
\label{so}
\end{figure}

\begin{figure}[H]
\centering
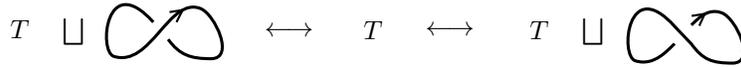
\caption{The move $KI$ on cobordism tangles.}
\label{k1}
\end{figure}

\begin{figure}[H]
\centering
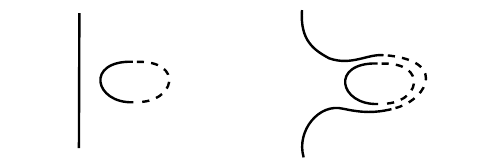
\caption{The move $KII^{g}$ on cobordism tangles.}
\label{k2}
\end{figure}

\begin{figure}[H]
\centering
\begingroup%
  \makeatletter%
  \providecommand\color[2][]{%
    \errmessage{(Inkscape) Color is used for the text in Inkscape, but the package 'color.sty' is not loaded}%
    \renewcommand\color[2][]{}%
  }%
  \providecommand\transparent[1]{%
    \errmessage{(Inkscape) Transparency is used (non-zero) for the text in Inkscape, but the package 'transparent.sty' is not loaded}%
    \renewcommand\transparent[1]{}%
  }%
  \providecommand\rotatebox[2]{#2}%
  \newcommand*\fsize{\dimexpr\f@size pt\relax}%
  \newcommand*\lineheight[1]{\fontsize{\fsize}{#1\fsize}\selectfont}%
  \ifx\svgwidth\undefined%
    \setlength{\unitlength}{368.50393701bp}%
    \ifx\svgscale\undefined%
      \relax%
    \else%
      \setlength{\unitlength}{\unitlength * \real{\svgscale}}%
    \fi%
  \else%
    \setlength{\unitlength}{\svgwidth}%
  \fi%
  \global\let\svgwidth\undefined%
  \global\let\svgscale\undefined%
  \makeatother%
  \begin{picture}(1,0.23076923)%
    \lineheight{1}%
    \setlength\tabcolsep{0pt}%
    \put(0,0){\includegraphics[width=\unitlength,page=1]{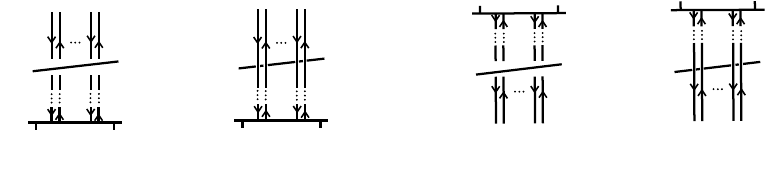}}%
    \put(-0.0387142,0.04838014){\color[rgb]{0,0,0}\makebox(0,0)[lt]{\begin{minipage}{0.64286082\unitlength}\raggedright $\text{All components of the same entrance or exit boundary component can cross any component.}$\end{minipage}}}%
    \put(0.48620335,0.13615092){\color[rgb]{0,0,0}\makebox(0,0)[lt]{\lineheight{1.25}\smash{\begin{tabular}[t]{l}$\text{and}$\end{tabular}}}}%
    \put(0.20733019,0.1403901){\color[rgb]{0,0,0}\makebox(0,0)[lt]{\lineheight{1.25}\smash{\begin{tabular}[t]{l}$\longleftrightarrow$\end{tabular}}}}%
    \put(0.78519342,0.13474994){\color[rgb]{0,0,0}\makebox(0,0)[lt]{\lineheight{1.25}\smash{\begin{tabular}[t]{l}$\longleftrightarrow$\end{tabular}}}}%
  \end{picture}%
\endgroup%

\caption{The move $COUPON$ on cobordism tangles.}
\label{coupon}
\end{figure}

\begin{figure}[H]
\centering
\resizebox{.6\linewidth}{!}{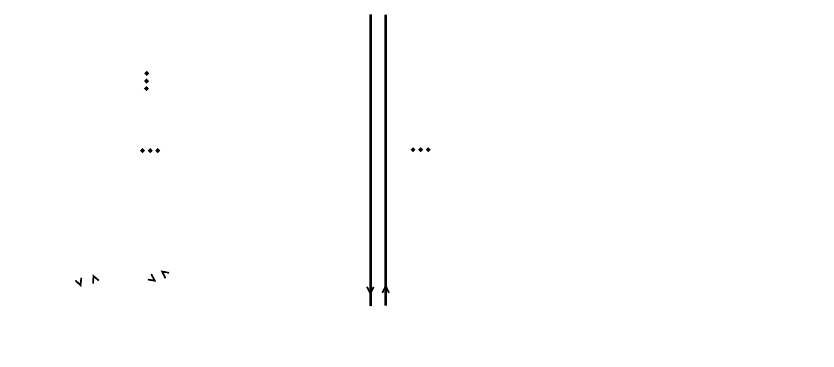}
\begin{center}and\end{center}
\resizebox{.6\linewidth}{!}{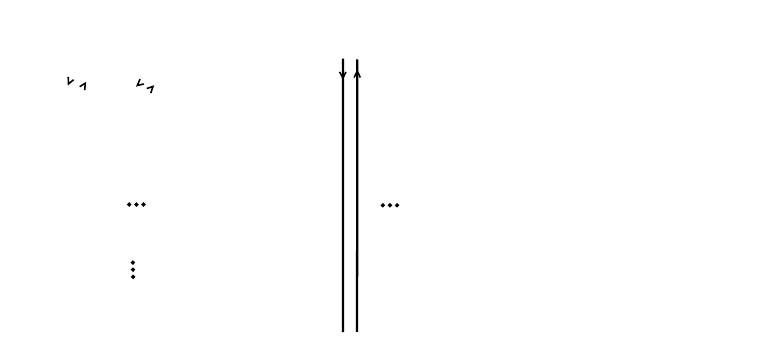}
\begin{center}The simultaneous twist of all components of the same entrance or exit boundary component is considered as nothing happened.\end{center}
\caption{The move $TWIST$ on cobordism tangles.}
\label{twist}
\end{figure}

This construction gives a surjective map $N$ (for "neighborhood") from $\bigsqcup\limits_{n\in \mathbb{N}}^{}Tang^{Cob}(\ung,n,\unh)$ to $\Cobc(\ung,\unh)$ defined on an isotopy class $T$ of a $(\ung,n,\unh)$-cobordism tangle by:
\begin{align}\label{surjectivemapn}
N( T)=M_{T}
\end{align}
The surjectivity of the map $N$ comes from the existence of a presentation for any connected $3$-cobordism by a ribbon $(\ung,n,\unh)$-cobordism tangle.
Since $N$ is surjective, it defines a bijection between the set $\Cobc(\ung,\unh)$ and the quotient set of $\left(\bigsqcup\limits_{n\in \mathbb{N}}Tang^{Cob}(\ung,n,\unh)\right)/\sim$ where two cobordim tangles $T_1$ and $T_2$ are equivalent if and only if $M_{T_1}=M_{T_2}$. In order to give a characterization in terms of cobordism tangles diagrams, we define the moves $SO$ ("surgery orientation"), $KI$ (Kirby I), $KII^{g}$ ("generalized Kirby II"), $COUPON$ and $TWIST$ illustrated respectively in Figures \ref{so}, \ref{k1}, \ref{k2}, \ref{coupon}, and \ref{twist}.


\begin{theorem}\label{theoremkirby}
Two isotopy classes of cobordism tangles $T_1$ and $T_2$ are equivalent if and only if a diagram of $T_1$ and a diagram of $T_2$ differ only by a finite sequence of planar isotopies, ribbon Reidemeister moves, and moves of type $SO$, $KI$, $KII^{g}$, $COUPON$, and $TWIST$. 
\end{theorem}

For details concerning extended Kirby calculus, see \cite{Turaev1994}, Chapter II.3.1, and for a proof of Theorem~\ref{theoremkirby}, see \cite{RT}, Section 7.2.

We need two more topological operations before starting the construction of the TQFT. To define composition of cobordisms and guarantee that the TQFT will be symmetric, we need to add some circles components on cobordism tangles as explain in the next subsection.

\subsection{Two useful operations on cobordism tangles}

\label{operationstar}First we define the {\it encercling composition} denoted by $\star$ between two cobordim tangles. Let $\ung=(g_1,\ldots,g_r)$, $\unh=(h_1,\ldots,h_s)$ and $\unk=(k_1,\ldots,k_t)$ be three tuples of positive integers. Let $M_{S}\in \Cobc(\ung,\unh)$ and $M_{T}\in \Cobc(\unh,\unk)$ two cobordisms represented respectively by a $(\ung,n,\unh)$-cobordism tangle $S$ and a $(\unh,m,\unk)$-cobordism tangle $T$. Then, according to Turaev (see \cite{Turaev1994}), the connected parametrized $3$-cobordism $M_{T}\circ M_{S}\in\Cobc(\ung,\unk)$ is represented by a tangle $T\star S$ obtained by adding $s-1$ trivial knots respectively surrounded the $s-1$ first boundary components of $T$ and juxtaposing this new tangle over $S$ as it is illustrated on Figure \ref{stargluingtangles}.

\label{hallowed} Secondly, we define a special kind of cobordism tangles called {\it hallowed tangles}. Let $T$ be a $(\ung,n,\unh)$-cobordism tangle $T$. We define the \emph{hallowed cobordism tangle} $\halo{T}$ by adding closed components around each connected enter and exit components as depicted in Figure~\ref{figurehallowed}.

\begin{minipage}{0.5\textwidth}
\begin{figure}[H]
\resizebox{.9\linewidth}{!}{\resizebox{!}{0.3\textheight}{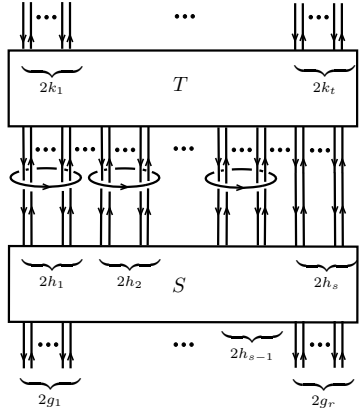}}
\caption{Cobordism tangle $T\star S$.}
\label{stargluingtangles}
\end{figure}
\end{minipage}
\begin{minipage}{0.5\textwidth}
\begin{figure}[H]
\resizebox{.9\linewidth}{!}{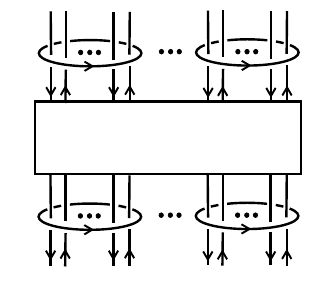}
\caption{Hallowed cobordism tangle $\halo{T}$.}
\label{figurehallowed}
\end{figure}
\end{minipage}


\section{Internal TQFT}
\subsection{Organization and main results} Consider $\C$ is a ribbon category with a coend $C$ and split idempotents. In this section, we start giving the definition of an \emph{admissible element}, a special kind of morphism $\un\rightarrow C$, allowing us to formulate the first main Theorem~\ref{maintheorem} of this article which gives general sufficient conditions to build a three dimensional TQFT with anomaly. In this section we explain all the construction of our \emph{internal} TQFT and give a proof of Theorem~\ref{maintheorem} and Theorem~~\ref{maintheorem2} which precise that internal TQFT is computable just using structural morphisms of the coend. The goal of our work is to extend 3-manifolds quantum invariants described by Virelizier in \cite{Vir} in 3-cobordisms invariants satisfying axioms of a TQFT. To do this, we need to add two conditions (see (Ad4) and (Ad5) below) to the definition of the classical \emph{Kirby element}. That leads to the next definition of an \emph{admissible color}:

\begin{definition}[Admissible element]\label{maindefinition}
Let $\alpha \in \mathrm{Hom}_{\C}(\un,C)$ where $C$ is the coend of the ribbon category $\C$. 
The morphism $\alpha$ is an \emph{admissible element} if it satisfies the following conditions, where $m$, $\Delta$, $\varepsilon$, $S$, $\omega$, and $\theta_{\pm}$  denote respectively the multiplication, the comultiplication, the counit, the antipode, the pairing and the linear forms coming from twists of the coend $C$:\\

\begin{tasks}[label=(\text{Ad}\arabic*)](3)
\task \;\;\;\;\;\;$\varepsilon\alpha=\mathrm{id}_{\mathds{1}}$;\label{ad1}
\task \;\;\;\;\;\; $S\alpha=\alpha$;\label{ad2}
\task \;\;\;\;\;\;$\theta_{\pm}\alpha\in \mathrm{End}_{\C}(\un)^{\times}$;\label{ad3}\\
\task \;\;\;\;\;\;$\forall n\in \mathbb{N},$\\\raisebox{-13mm}{
\begingroup%
  \makeatletter%
  \providecommand\color[2][]{%
    \errmessage{(Inkscape) Color is used for the text in Inkscape, but the package 'color.sty' is not loaded}%
    \renewcommand\color[2][]{}%
  }%
  \providecommand\transparent[1]{%
    \errmessage{(Inkscape) Transparency is used (non-zero) for the text in Inkscape, but the package 'transparent.sty' is not loaded}%
    \renewcommand\transparent[1]{}%
  }%
  \providecommand\rotatebox[2]{#2}%
  \ifx\svgwidth\undefined%
    \setlength{\unitlength}{44bp}%
    \ifx\svgscale\undefined%
      \relax%
    \else%
      \setlength{\unitlength}{\unitlength * \real{\svgscale}}%
    \fi%
  \else%
    \setlength{\unitlength}{\svgwidth}%
  \fi%
  \global\let\svgwidth\undefined%
  \global\let\svgscale\undefined%
  \makeatother%
  \begin{picture}(1,2.18181818)%
    \put(0,0){\includegraphics[width=\unitlength]{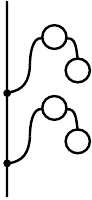}}%
    \put(0.17204626,0.21298827){\color[rgb]{0,0,0}\makebox(0,0)[lt]{\begin{minipage}{0.38151823\unitlength}\raggedright $C^{\otimes n}$\end{minipage}}}%
    \put(0.52439923,0.95317934){\color[rgb]{0,0,0}\makebox(0,0)[lb]{\smash{$\omega$}}}%
    \put(0.77561277,0.67519353){\color[rgb]{0,0,0}\makebox(0,0)[lt]{\begin{minipage}{0.15435022\unitlength}\raggedright $\alpha$\end{minipage}}}%
    \put(0.52445246,1.72663923){\color[rgb]{0,0,0}\makebox(0,0)[lb]{\smash{$\omega$}}}%
    \put(0.775666,1.44865444){\color[rgb]{0,0,0}\makebox(0,0)[lt]{\begin{minipage}{0.15435022\unitlength}\raggedright $\alpha$\end{minipage}}}%
  \end{picture}%
\endgroup%
} = \raisebox{-13mm}{
\begingroup%
  \makeatletter%
  \providecommand\color[2][]{%
    \errmessage{(Inkscape) Color is used for the text in Inkscape, but the package 'color.sty' is not loaded}%
    \renewcommand\color[2][]{}%
  }%
  \providecommand\transparent[1]{%
    \errmessage{(Inkscape) Transparency is used (non-zero) for the text in Inkscape, but the package 'transparent.sty' is not loaded}%
    \renewcommand\transparent[1]{}%
  }%
  \providecommand\rotatebox[2]{#2}%
  \ifx\svgwidth\undefined%
    \setlength{\unitlength}{44bp}%
    \ifx\svgscale\undefined%
      \relax%
    \else%
      \setlength{\unitlength}{\unitlength * \real{\svgscale}}%
    \fi%
  \else%
    \setlength{\unitlength}{\svgwidth}%
  \fi%
  \global\let\svgwidth\undefined%
  \global\let\svgscale\undefined%
  \makeatother%
  \begin{picture}(1,2.18181818)%
    \put(0,0){\includegraphics[width=\unitlength]{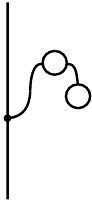}}%
    \put(0.18909928,0.22365554){\color[rgb]{0,0,0}\makebox(0,0)[lt]{\begin{minipage}{0.38151823\unitlength}\raggedright $C^{\otimes n}$\end{minipage}}}%
    \put(0.52789474,1.44536242){\color[rgb]{0,0,0}\makebox(0,0)[lb]{\smash{$\omega$}}}%
    \put(0.77910831,1.1673759){\color[rgb]{0,0,0}\makebox(0,0)[lt]{\begin{minipage}{0.15435022\unitlength}\raggedright $\alpha$\end{minipage}}}%
  \end{picture}%
\endgroup%
}\label{ad4};
\task \;\;\;\;\;\;$\forall n\in \mathbb{N},$\\ \raisebox{-12mm}{
\begingroup%
  \makeatletter%
  \providecommand\color[2][]{%
    \errmessage{(Inkscape) Color is used for the text in Inkscape, but the package 'color.sty' is not loaded}%
    \renewcommand\color[2][]{}%
  }%
  \providecommand\transparent[1]{%
    \errmessage{(Inkscape) Transparency is used (non-zero) for the text in Inkscape, but the package 'transparent.sty' is not loaded}%
    \renewcommand\transparent[1]{}%
  }%
  \providecommand\rotatebox[2]{#2}%
  \ifx\svgwidth\undefined%
    \setlength{\unitlength}{80bp}%
    \ifx\svgscale\undefined%
      \relax%
    \else%
      \setlength{\unitlength}{\unitlength * \real{\svgscale}}%
    \fi%
  \else%
    \setlength{\unitlength}{\svgwidth}%
  \fi%
  \global\let\svgwidth\undefined%
  \global\let\svgscale\undefined%
  \makeatother%
  \begin{picture}(1,1.3)%
    \put(0,0){\includegraphics[width=\unitlength]{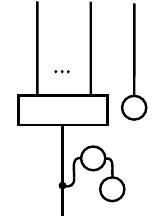}}%
    \put(0.24582802,0.59982379){\color[rgb]{0,0,0}\makebox(0,0)[lb]{\smash{$\mathrm{id}_{C^{\otimes n}}$}}}%
    \put(0.07192465,0.10560186){\color[rgb]{0,0,0}\makebox(0,0)[lt]{\begin{minipage}{0.1306367\unitlength}\raggedright $C^{\otimes n}$\end{minipage}}}%
    \put(0.52132706,0.32150713){\color[rgb]{0,0,0}\makebox(0,0)[lb]{\smash{$\omega$}}}%
    \put(0.63877106,0.13147853){\color[rgb]{0,0,0}\makebox(0,0)[lb]{\smash{$\alpha$}}}%
    \put(0.77012303,0.6207609){\color[rgb]{0,0,0}\makebox(0,0)[lb]{\smash{$\alpha$}}}%
  \end{picture}%
\endgroup%
}=\raisebox{-12mm}{
\begingroup%
  \makeatletter%
  \providecommand\color[2][]{%
    \errmessage{(Inkscape) Color is used for the text in Inkscape, but the package 'color.sty' is not loaded}%
    \renewcommand\color[2][]{}%
  }%
  \providecommand\transparent[1]{%
    \errmessage{(Inkscape) Transparency is used (non-zero) for the text in Inkscape, but the package 'transparent.sty' is not loaded}%
    \renewcommand\transparent[1]{}%
  }%
  \providecommand\rotatebox[2]{#2}%
  \ifx\svgwidth\undefined%
    \setlength{\unitlength}{80bp}%
    \ifx\svgscale\undefined%
      \relax%
    \else%
      \setlength{\unitlength}{\unitlength * \real{\svgscale}}%
    \fi%
  \else%
    \setlength{\unitlength}{\svgwidth}%
  \fi%
  \global\let\svgwidth\undefined%
  \global\let\svgscale\undefined%
  \makeatother%
  \begin{picture}(1,1.3)%
    \put(0,0){\includegraphics[width=\unitlength]{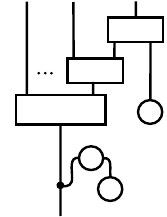}}%
    \put(0.51756469,0.83511152){\color[rgb]{0,0,0}\makebox(0,0)[lb]{\smash{$\Delta$}}}%
    \put(0.76463648,1.09522654){\color[rgb]{0,0,0}\makebox(0,0)[lb]{\smash{$m$}}}%
    \put(0.23298355,0.60120045){\color[rgb]{0,0,0}\makebox(0,0)[lb]{\smash{$\mathrm{id}_{C^{\otimes n}}$}}}%
    \put(0.05908014,0.10697906){\color[rgb]{0,0,0}\makebox(0,0)[lt]{\begin{minipage}{0.1306367\unitlength}\raggedright $C^{\otimes n}$\end{minipage}}}%
    \put(0.50848257,0.32288419){\color[rgb]{0,0,0}\makebox(0,0)[lb]{\smash{$\omega$}}}%
    \put(0.62592646,0.13285484){\color[rgb]{0,0,0}\makebox(0,0)[lb]{\smash{$\alpha$}}}%
    \put(0.86576817,0.59525704){\color[rgb]{0,0,0}\makebox(0,0)[lb]{\smash{$\alpha$}}}%
  \end{picture}%
\endgroup%
}.\label{ad5}
\end{tasks}
\end{definition}

Using an admissible element$\alpha$, we are able to compute our internal three dimensional TQFT as described in the next theorem, which is the first main result of this article. 
Recall that $\Cobp$ denotes the category of parametrized 3-cobordims (see Section~\ref{parametrized}).

\begin{theorem}\label{maintheorem}
Let $\alpha \in  \mathrm{Hom}_{\C}(\mathds{1},C)$ be an admissible element.
Then for every connected parametrized surface $\Sigma_{g}$ of genus $g$ and for every connected $3$-cobordism $M_{T}$ of $\Cobp(\ung,\unh)$ represented by a $(\ung,n,\unh)$-cobordism tangle $T$, the assignment:
\begin{align}
\mathrm{W}_{\C}(\Sigma_{g};\alpha)&=(C^{\otimes g})_{\alpha}\\ 
\mathrm{W}_{\C}(M_{T};\alpha)&=\nu_{\alpha}(T)\;|\halo{T}|_{\C,\alpha}
\end{align}
defines a braided monoidal functor with anomaly $(\mathrm{W_{\C,\alpha}},\gamma)$ between the category $\Cobp$ and $\C$.

Moreover the associated TQFT with anomaly $(\mathrm{V_{\C,\alpha}},\gamma)$ takes values in the full subcategory of transparent objects $\T \subset \C$.
\end{theorem}

\begin{remarks}
\item The object $(C^{\otimes n})_{\alpha}$ is the image of $C^{\otimes n}$ by idempotent defined in equality (\ref{hallowedmorphism}). The tangle invariant $|\phantom{O}|_{\C,\alpha}$ is defined in Lemma~\ref{lemmaisotopyinv}.
\item Note that $\halo{T}$ means we add some circle components on the cobordism tangle $T$ to specify connected (enter and exit) components of the represented cobordism (see Section~\ref{hallowed}). 
\item Note that the anomaly $\gamma$ is given by:
\[\gamma_{M_{T},M_{T'}} = (\theta_{+}\alpha)^{b_{+}(T)+b_{+}(T')-b_{+}(T\star T')}(\theta_{-}\alpha)^{b_{-}(T)+b_{-}(T')-b_{-}(T\star T')}
\]
where $T\star T'$ encodes the cobordism $M'\circ M$ (see Section~\ref{operationstar}), and $b_+(T)$ and $b_{-}(T)$ are respectively numbers of positive and negatives eigenvalues of the linking matrix of $T$.
\item In the case of premodular categories, the Kirby element from \cite{Vir} corresponding to the Kirby color in the Turaev TQFT construction is an admissible element. We will prove it in Section~\ref{modular}.
\end{remarks}


According our next second main result, the internal TQFT is completely computable, starting from a tangle presentation of a 3-corbordims, using only structural morphisms of the coend $C$.

\begin{theorem}\label{maintheorem2}
Let $\alpha\colon \un \rightarrow C$ be an admissible element.
The TQFT $\mathrm{V}_{\C,\alpha}$ can be expressed entirely in terms of $\circ$, $\otimes$, $+$, $\alpha$ and the structural morphisms
\[m,\; \Delta,\; \varepsilon,\; u,\; S,\; S^{-1},\; \theta_{+},\; \theta_{-},\; \omega,\; \tau_{C,C},\; \tau_{C,C}^{-1},\; \mathrm{id}_{C}\]
of the coend $C$.
\end{theorem}

\subsection{Construction of internal TQFT} We start now the steps to construct our internal TQFT and will prove the two main Theorems~\ref{maintheorem} and \ref{maintheorem2}.
First, we define an isotopy invariant of opentangle:
\begin{lemma}\label{propinvopen}
The universal morphism associated to an opentangle defines a map $|\phantom{O}|_{\C} \colon Otang \rightarrow \mathrm{Hom}_{\C}$ such that every isotopy class of $(\ung,n,\unh)$-opentangle $O$ is mapped to the morphism of $\C$:
\begin{align}
|O|_{\C}\colon C^{\otimes |\ung|}\otimes C^{\otimes n}\otimes C^{\otimes |\unh|}\rightarrow C^{|\unh|}.
\end{align}
\end{lemma}

Let's illustrate the process on the example of the 1-genus surface cylinder $\Sigma_1\times [0,1]$ before starting the proof. Recall that a tangle presentation of cylinder $\Sigma_1\times [0,1]$ is given in Figure~\ref{cylinder}.

\begin{example} Construction of the opentangle invariant based on a tangle presentation of $\Sigma_1\times [0,1]$:

\begingroup%
  \makeatletter%
  \providecommand\color[2][]{%
    \errmessage{(Inkscape) Color is used for the text in Inkscape, but the package 'color.sty' is not loaded}%
    \renewcommand\color[2][]{}%
  }%
  \providecommand\transparent[1]{%
    \errmessage{(Inkscape) Transparency is used (non-zero) for the text in Inkscape, but the package 'transparent.sty' is not loaded}%
    \renewcommand\transparent[1]{}%
  }%
  \providecommand\rotatebox[2]{#2}%
  \newcommand*\fsize{\dimexpr\f@size pt\relax}%
  \newcommand*\lineheight[1]{\fontsize{\fsize}{#1\fsize}\selectfont}%
  \ifx\svgwidth\undefined%
    \setlength{\unitlength}{59.52755906bp}%
    \ifx\svgscale\undefined%
      \relax%
    \else%
      \setlength{\unitlength}{\unitlength * \real{\svgscale}}%
    \fi%
  \else%
    \setlength{\unitlength}{\svgwidth}%
  \fi%
  \global\let\svgwidth\undefined%
  \global\let\svgscale\undefined%
  \makeatother%
  \begin{picture}(1,1.33333333)%
    \lineheight{1}%
    \setlength\tabcolsep{0pt}%
    \put(0,0){\includegraphics[width=\unitlength,page=1]{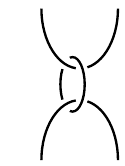}}%
    \put(-0.09510629,0.70935924){\color[rgb]{0,0,0}\makebox(0,0)[lt]{\begin{minipage}{0.93056371\unitlength}\raggedright $T=$\end{minipage}}}%
    \put(0,0){\includegraphics[width=\unitlength,page=2]{cylindertangle2.pdf}}%
  \end{picture}%
\endgroup%
\quad
\raisebox{13mm}{$\longrightarrow$}\quad
\begingroup%
  \makeatletter%
  \providecommand\color[2][]{%
    \errmessage{(Inkscape) Color is used for the text in Inkscape, but the package 'color.sty' is not loaded}%
    \renewcommand\color[2][]{}%
  }%
  \providecommand\transparent[1]{%
    \errmessage{(Inkscape) Transparency is used (non-zero) for the text in Inkscape, but the package 'transparent.sty' is not loaded}%
    \renewcommand\transparent[1]{}%
  }%
  \providecommand\rotatebox[2]{#2}%
  \ifx\svgwidth\undefined%
    \setlength{\unitlength}{80bp}%
    \ifx\svgscale\undefined%
      \relax%
    \else%
      \setlength{\unitlength}{\unitlength * \real{\svgscale}}%
    \fi%
  \else%
    \setlength{\unitlength}{\svgwidth}%
  \fi%
  \global\let\svgwidth\undefined%
  \global\let\svgscale\undefined%
  \makeatother%
  \begin{picture}(1,1)%
    \put(0,0){\includegraphics[width=\unitlength]{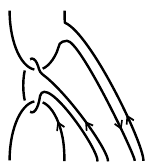}}%
    \put(0.40122499,0.08784635){\color[rgb]{0,0,0}\makebox(0,0)[lt]{\begin{minipage}{0.19604363\unitlength}\raggedright {\footnotesize $X$}\end{minipage}}}%
    \put(0.66159199,0.08784435){\color[rgb]{0,0,0}\makebox(0,0)[lt]{\begin{minipage}{0.19604363\unitlength}\raggedright {\footnotesize $Y$}\end{minipage}}}%
    \put(0.88163381,0.08784435){\color[rgb]{0,0,0}\makebox(0,0)[lt]{\begin{minipage}{0.19604363\unitlength}\raggedright {\footnotesize $Z$}\end{minipage}}}%
    \put(0.09339668,0.93651475){\color[rgb]{0,0,0}\makebox(0,0)[lt]{\begin{minipage}{0.19604363\unitlength}\raggedright {\footnotesize $Z$}\end{minipage}}}%
    \put(-0.04970302,-0.03358647){\color[rgb]{0,0,0}\makebox(0,0)[lt]{\begin{minipage}{1.60312821\unitlength}\raggedright $\text{Morphism $O_{X,Y,Z}$ of $\C$}$\end{minipage}}}%
  \end{picture}%
\endgroup%
\quad
\raisebox{13mm}{$\longrightarrow$}\quad
\begingroup%
  \makeatletter%
  \providecommand\color[2][]{%
    \errmessage{(Inkscape) Color is used for the text in Inkscape, but the package 'color.sty' is not loaded}%
    \renewcommand\color[2][]{}%
  }%
  \providecommand\transparent[1]{%
    \errmessage{(Inkscape) Transparency is used (non-zero) for the text in Inkscape, but the package 'transparent.sty' is not loaded}%
    \renewcommand\transparent[1]{}%
  }%
  \providecommand\rotatebox[2]{#2}%
  \newcommand*\fsize{\dimexpr\f@size pt\relax}%
  \newcommand*\lineheight[1]{\fontsize{\fsize}{#1\fsize}\selectfont}%
  \ifx\svgwidth\undefined%
    \setlength{\unitlength}{63.75bp}%
    \ifx\svgscale\undefined%
      \relax%
    \else%
      \setlength{\unitlength}{\unitlength * \real{\svgscale}}%
    \fi%
  \else%
    \setlength{\unitlength}{\svgwidth}%
  \fi%
  \global\let\svgwidth\undefined%
  \global\let\svgscale\undefined%
  \makeatother%
  \begin{picture}(1,1.76470588)%
    \lineheight{1}%
    \setlength\tabcolsep{0pt}%
    \put(0,0){\includegraphics[width=\unitlength,page=1]{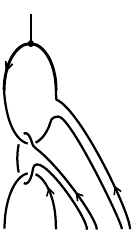}}%
    \put(0.44251184,0.08237947){\color[rgb]{0,0,0}\makebox(0,0)[lt]{\begin{minipage}{0.23063956\unitlength}\raggedright {\footnotesize $X$}\end{minipage}}}%
    \put(0.74882596,0.08237712){\color[rgb]{0,0,0}\makebox(0,0)[lt]{\begin{minipage}{0.23063956\unitlength}\raggedright {\footnotesize $Y$}\end{minipage}}}%
    \put(1.00769869,0.08237712){\color[rgb]{0,0,0}\makebox(0,0)[lt]{\begin{minipage}{0.23063956\unitlength}\raggedright {\footnotesize $Z$}\end{minipage}}}%
    \put(0.28623663,1.61603336){\color[rgb]{0,0,0}\makebox(0,0)[lt]{\begin{minipage}{0.23063956\unitlength}\raggedright {\footnotesize $C$}\end{minipage}}}%
    \put(-0.07258931,-0.01854738){\color[rgb]{0,0,0}\makebox(0,0)[lt]{\begin{minipage}{1.6158275\unitlength}\raggedright $\text{Composition by $\iota_Z$}$\end{minipage}}}%
  \end{picture}%
\endgroup%

\begingroup%
  \makeatletter%
  \providecommand\color[2][]{%
    \errmessage{(Inkscape) Color is used for the text in Inkscape, but the package 'color.sty' is not loaded}%
    \renewcommand\color[2][]{}%
  }%
  \providecommand\transparent[1]{%
    \errmessage{(Inkscape) Transparency is used (non-zero) for the text in Inkscape, but the package 'transparent.sty' is not loaded}%
    \renewcommand\transparent[1]{}%
  }%
  \providecommand\rotatebox[2]{#2}%
  \ifx\svgwidth\undefined%
    \setlength{\unitlength}{80bp}%
    \ifx\svgscale\undefined%
      \relax%
    \else%
      \setlength{\unitlength}{\unitlength * \real{\svgscale}}%
    \fi%
  \else%
    \setlength{\unitlength}{\svgwidth}%
  \fi%
  \global\let\svgwidth\undefined%
  \global\let\svgscale\undefined%
  \makeatother%
  \begin{picture}(1,1)%
    \put(0,0){\includegraphics[width=\unitlength,page=1]{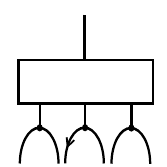}}%
    \put(0.38237688,0.57498363){\color[rgb]{0,0,0}\makebox(0,0)[lt]{\begin{minipage}{0.43583263\unitlength}\raggedright $|O|_{\C}$\end{minipage}}}%
    \put(0,0){\includegraphics[width=\unitlength,page=2]{construction6.pdf}}%
    \put(0.21642184,0.08439455){\color[rgb]{0,0,0}\makebox(0,0)[lt]{\begin{minipage}{0.19604362\unitlength}\raggedright {\footnotesize $X$}\end{minipage}}}%
    \put(0.48738252,0.08287655){\color[rgb]{0,0,0}\makebox(0,0)[lt]{\begin{minipage}{0.19604362\unitlength}\raggedright {\footnotesize $Y$}\end{minipage}}}%
    \put(0.76053304,0.08135955){\color[rgb]{0,0,0}\makebox(0,0)[lt]{\begin{minipage}{0.19604362\unitlength}\raggedright {\footnotesize $Z$}\end{minipage}}}%
    \put(0.36628517,0.90379174){\color[rgb]{0,0,0}\makebox(0,0)[lt]{\begin{minipage}{0.18586327\unitlength}\raggedright \end{minipage}}}%
    \put(0.09811091,0.34379144){\color[rgb]{0,0,0}\makebox(0,0)[lt]{\begin{minipage}{0.18586327\unitlength}\raggedright \end{minipage}}}%
    \put(0.37059002,0.34379244){\color[rgb]{0,0,0}\makebox(0,0)[lt]{\begin{minipage}{0.18586327\unitlength}\raggedright \end{minipage}}}%
    \put(0.64848262,0.34740144){\color[rgb]{0,0,0}\makebox(0,0)[lt]{\begin{minipage}{0.18586327\unitlength}\raggedright \end{minipage}}}%
    \put(-0.05970911,-0.01668267){\color[rgb]{0,0,0}\makebox(0,0)[lt]{\begin{minipage}{1.62247787\unitlength}\raggedright \end{minipage}}}%
    \put(-0.03302936,0.54639723){\color[rgb]{0,0,0}\makebox(0,0)[lt]{\begin{minipage}{0.43583263\unitlength}\raggedright $=$\end{minipage}}}%
    \put(0.35491933,0.91862723){\color[rgb]{0,0,0}\makebox(0,0)[lt]{\begin{minipage}{0.43583263\unitlength}\raggedright $C$\end{minipage}}}%
    \put(0.0743998,0.34415713){\color[rgb]{0,0,0}\makebox(0,0)[lt]{\begin{minipage}{0.43583263\unitlength}\raggedright $C$\end{minipage}}}%
    \put(0.338904,0.34415723){\color[rgb]{0,0,0}\makebox(0,0)[lt]{\begin{minipage}{0.43583263\unitlength}\raggedright $C$\end{minipage}}}%
    \put(0.62820546,0.34415723){\color[rgb]{0,0,0}\makebox(0,0)[lt]{\begin{minipage}{0.43583263\unitlength}\raggedright $C$\end{minipage}}}%
  \end{picture}%
\endgroup%
\quad\vspace{0.5cm}
\end{example}

\begin{proof}
Choosing a $(\ung,n,\unh)$-opentangle O, the construction of the universal morphism $|O|_{\C}$ is well-defined by Lemma~\ref{lemma0}. Morevover, the construction does not depend on the choice of an element in the isotopy class of O because of Shum’s result (see \cite{Shu94}): two isotopic opentangles define the same morphism in the ribbon category $\C$.
\end{proof}

Then, we can define an isotopy invariant of cobordisms tangles:
\begin{lemma}\label{lemmaisotopyinv} 
Let $\alpha \in \mathrm{Hom}_{\C}(\mathds{1},C)$.
There is a map ${|\phantom{O}|_{\C,\alpha}\colon Tang^{Cob}(\ung,n,\unh) \rightarrow \mathrm{Hom}_{\C}}$ defined on every isotopy class of $(\ung,n,\unh)$-cobordism tangle $T$ by
\begin{align}
|T|_{\C,\alpha}=|O|_{\C}\circ(\mathrm{id}_{C^{\otimes |\ung|}}\otimes \alpha^{\otimes n+|\unh|})
\end{align}
where $O$ is any preimage of $T$ by the map $U\colon Otang \rightarrow Tang^{Cob}$ defined in formula  (\ref{surjectivemap}).
\end{lemma}

\begin{example} Construction of the cobordism tangle invariant based on a tangle presentation $T$ of $\Sigma_1\times [0,1]$ :
\begin{center}
\begingroup%
  \makeatletter%
  \providecommand\color[2][]{%
    \errmessage{(Inkscape) Color is used for the text in Inkscape, but the package 'color.sty' is not loaded}%
    \renewcommand\color[2][]{}%
  }%
  \providecommand\transparent[1]{%
    \errmessage{(Inkscape) Transparency is used (non-zero) for the text in Inkscape, but the package 'transparent.sty' is not loaded}%
    \renewcommand\transparent[1]{}%
  }%
  \providecommand\rotatebox[2]{#2}%
  \newcommand*\fsize{\dimexpr\f@size pt\relax}%
  \newcommand*\lineheight[1]{\fontsize{\fsize}{#1\fsize}\selectfont}%
  \ifx\svgwidth\undefined%
    \setlength{\unitlength}{75bp}%
    \ifx\svgscale\undefined%
      \relax%
    \else%
      \setlength{\unitlength}{\unitlength * \real{\svgscale}}%
    \fi%
  \else%
    \setlength{\unitlength}{\svgwidth}%
  \fi%
  \global\let\svgwidth\undefined%
  \global\let\svgscale\undefined%
  \makeatother%
  \begin{picture}(1,1)%
    \lineheight{1}%
    \setlength\tabcolsep{0pt}%
    \put(0,0){\includegraphics[width=\unitlength,page=1]{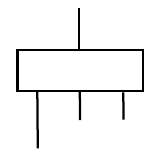}}%
    \put(0.4119965,0.61420333){\color[rgb]{0,0,0}\makebox(0,0)[lt]{\begin{minipage}{0.59194691\unitlength}\raggedright $|O|_{\C}$\end{minipage}}}%
    \put(0,0){\includegraphics[width=\unitlength,page=2]{construction777.pdf}}%
    \put(0.47974024,0.16761728){\color[rgb]{0,0,0}\makebox(0,0)[lt]{\begin{minipage}{0.16512439\unitlength}\raggedright \small{$\alpha$}\end{minipage}}}%
    \put(0.76017024,0.16567628){\color[rgb]{0,0,0}\makebox(0,0)[lt]{\begin{minipage}{0.16512439\unitlength}\raggedright \small{$\alpha$}\end{minipage}}}%
    \put(0.09713,0.05543484){\color[rgb]{0,0,0}\makebox(0,0)[lt]{\lineheight{0}\smash{\begin{tabular}[t]{l}\small{$C$}\end{tabular}}}}%
    \put(0.37045,0.87100952){\color[rgb]{0,0,0}\makebox(0,0)[lt]{\lineheight{0}\smash{\begin{tabular}[t]{l}\small{$C$}\end{tabular}}}}%
    \put(-0.77473121,0.60409163){\color[rgb]{0,0,0}\makebox(0,0)[lt]{\begin{minipage}{0.90000055\unitlength}\raggedright $|T|_{\C,\alpha}=$\end{minipage}}}%
  \end{picture}%
\endgroup%

\end{center}
\end{example}

 \begin{proof} 
Consider an isotopy class $T$ of a $(\ung,n,\unh)$-cobordism tangle and two $(\ung,n,\unh)$-opentangles $T^{o}_1$ and $T^{o}_2$ such that $U(T^{o}_1)=T =U(T^{o}_2)$ where the map $U$ is defined in (\ref{surjectivemap}). Then, as it is explained in Section \ref{opentangle}, there exists a finite sequence of planar isotopies and ribbon Reidemeister moves, moves $BA$ (see Figure~\ref{equivalence}), moves $ESC$ (see Figure~\ref{moveesc}), and moves $ROT$ (see Figure~\ref{moverot}) between diagrams of $T^{o}_1$ and $T^{o}_2$. 
We have to show that 
\begin{align}\label{identity1}
|T^{o}_1|_{\C}(\id_{C^{\otimes |\ung|}}\otimes \alpha^{\otimes n+|\unh|})=|T^{o}_2|_{\C}(\id_{C^{\otimes |\ung|}}\otimes \alpha^{\otimes n+|\unh|}).
\end{align}
If $T^{o}_1$ and $T^{o}_2$ differ by planar isotopies and ribbon Reidemeister moves, the equality~(\ref{identity1}) is obvious because $|\phantom{O}|_{\C}$ is an isotopy invariant of opentangles (see Lemma~\ref{propinvopen}). Suppose now that $T^{o}_1$ and $T^{o}_2$, considered as diagrams, differ only from one move $BA$ as it is illustrated in Figure~\ref{differentcrossings}.

\begin{figure}[H]
\centering
\resizebox{!}{0.2\textheight}{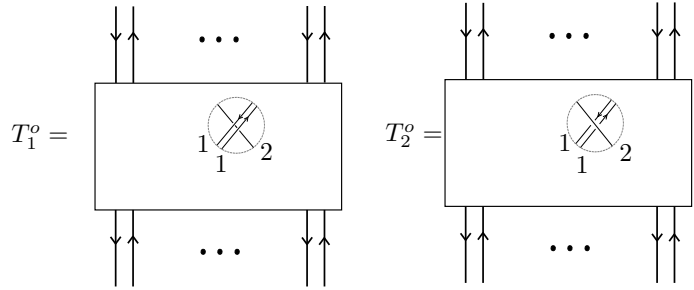}
\caption{One different crossing between $T_{1}^{o}$ and $T_{2}^{o}$.}
\label{differentcrossings}
\end{figure}

By isotopy, we modify the diagrams of opentangles $T^{o}_1$ and $T^{o}_2$ such that the difference of crossings is at the bottom of the diagram as shown in Figure~\ref{differentbottom}.

\begin{figure}[H]
\centering
\resizebox{!}{0.15\textheight}{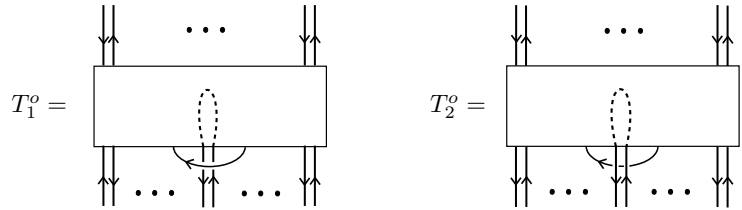}
\caption{One different bottom crossing between $T_{1}^{o}$ and $T_{2}^{o}$.}
\label{differentbottom}
\end{figure}

Now, let $\unX=X_1,\ldots,X_{|\ung|+n+|\unh|}$ be objects of $\C$ and colore $T^{o}_1$ and $T^{o}_2$ by those objects to obtain morphisms of $\C$, $T^{o}_{1,\unX}$ and $T^{o}_{2,\unX}$. Suppose that the crossing which is different between $T^{o}_1$ and $T^{o}_2$ affects the $i$th component such that $|\ung|+1\leq i \leq |\ung|+n+|\unh|$ (a surgery component or an exit component) and denote by $Y$ the color of the other component in the crossing and by\\
$f \colon X_1^{*}\otimes X_1\otimes\ldots\otimes Y^{*}\otimes X_{i}^{*}\otimes X_{i}\otimes Y\otimes \ldots \otimes X^{*}_{|\ung|+n+|\unh|}\otimes X_{|\ung|+n+|\unh|}\rightarrow$\\$X^{*}_{|\ung|+n+1}\otimes X_{|\ung|+n+1}\otimes \ldots\otimes X^{*}_{|\ung|+n+|\unh|}\otimes X_{|\ung|+n+|\unh|}$
the morphism of $\C$ defined in Figure~\ref{morphimst12}.

\begin{figure}[H]
\centering
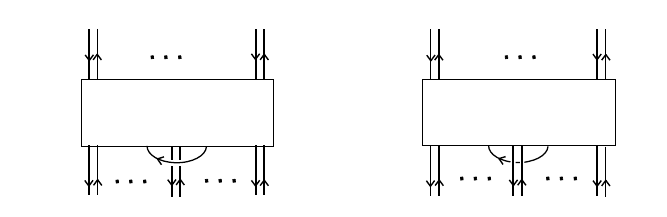
\caption{The morphims $T^{o}_{1,\unX}$ and $T^{o}_{2,\unX}$ of $\C$ associated to tangles $T^{o}_1$ and $T^{o}_2$.}
\label{morphimst12}
\end{figure} 

Note that there exists $j\in \mathbb{N}$ such that $Y=X_j$.

Suppose that $j\neq i$ and without loss of generality that $j=|\ung|+n+|\unh|$.
Using dinaturality of $(\iota_{X_{\ung+n+1}}\otimes\ldots\otimes \iota_{X_{|\ung|+n+|\unh|}})\circ f$ in $X_{1},\ldots,X_{|\ung| +n+ |\unh|-1}$ and using Fubini theorem (see~\cite{MacL} and see Lemma~\ref{lemma0}) with parameters $Y^{*}$ and $Y$, by universal property of the coend $C$, there exists a unique morphism \[\psi_{f} \colon C^{\otimes i-1}\otimes Y^{*}\otimes C\otimes Y\otimes C^{\otimes |\ung|+n+|\unh|-i-1}\otimes Y^{*}\otimes Y\rightarrow  C^{\otimes |\unh|}\] such that, for all objects $X_1, \ldots, X_{|g|+n+|h|-1}$ of $\C$, 
\[
(\iota_{X_{|\ung|+n+1}}\otimes\ldots\otimes \iota_{X_{|\ung|+n+|\unh|}})\circ f=\psi_{f}\circ (\iota_{X_1}\otimes \ldots\otimes \id_{Y^{*}}\otimes \iota_{X_i}\otimes \id_{Y}\otimes\iota_{X_{i+1}}\otimes \ldots\otimes \iota_{X_{|\ung|+n+|\unh|-1}}\otimes \id_{Y^{*}\otimes Y})  
\]
Morevover, we know by Lemma~\ref{propinvopen} that 
\[
(\iota_{X_{|\ung|+n+1}}\otimes\ldots\otimes \iota_{X_{|\ung|+n+|\unh|}})\circ T^{o}_{1,\unX}=|T^{o}_{1}|_{\C}\circ (\iota_{X_1}\otimes \ldots\otimes \iota_{X_{|\ung|+n+|\unh|}}) 
\] and
\[
(\iota_{X_{|\ung|+n+1}}\otimes\ldots\otimes \iota_{X_{|\ung|+n+|\unh|}})\circ T^{o}_{2,\unX}=|T^{o}_{2}|_{\C}\circ (\iota_{X_1}\otimes \ldots\otimes \iota_{X_{|\ung|+n+|\unh|}}).
\]
Thus we have the two identities showed on Figure~\ref{universal}.

\begin{figure}[H]
\resizebox{!}{0.18\textheight}{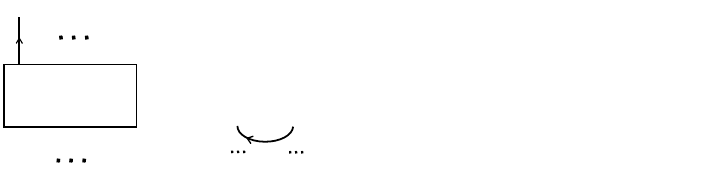}
\caption{Two factorizations of $(\iota_{X_{|\ung|+n+1}}\otimes\ldots\otimes \iota_{X_{|\ung|+n+|\unh|}})T_{1,\unX}^{o}$ (resp. $T_{2,\unX}^{o}$)}
\label{universal}
\end{figure} 

Remark that in one hand, we have for all objects $Y$ of $\C$ 
\begin{align}\label{eq1}
|T^{o}_1|_{\C}\circ (\id_{C^{\otimes |\ung|+n+|\unh|-1}}\otimes \iota_{Y})=\psi_{f}\circ (\id_{C^{\otimes i-1}}\otimes \psii\otimes \id_{C^{\otimes |\ung|+n+|\unh|-i-1}}\otimes \id_{Y^{*}\otimes Y})
\end{align}
and in the other hand, we have 
\begin{align}\label{eq2}
|T^{o}_2|_{\C}\circ (\id_{C^{\otimes |\ung|+n+|\unh|-1}}\otimes \iota_{Y})=\psi_{f}\circ (\id_{C^{\otimes i-1}}\otimes \psiii\otimes \id_{C^{\otimes |\ung|+n+|\unh|-i-1}}\otimes \id_{Y^{*}\otimes Y})
\end{align}
Indeed, as the right member of equation (\ref{eq1}) is dinatural in $Y$ with parameters $C$, there exists a unique morphism $\phi \colon C^{\otimes |\ung|+n+|\unh|}\rightarrow C^{\otimes |\unh|}$ such that, for all $Y\in \C$, 
\[\psi_{f}\circ (\id_{C^{\otimes i-1}}\otimes \psii\otimes \id_{C^{\otimes |\ung|+n+|\unh|-i-1}}\otimes \id_{Y^{*}\otimes Y})=\phi\circ(\id_{C^{\otimes |\ung|+n+|\unh|-1}}\otimes \iota_{Y})\]
and composing with morphism $\iota_{X_1}\otimes \ldots\otimes \iota_{X_{|\ung|+n+|\unh|-1}}$, 
\[\psi_{f}\circ (\iota_{X_1}\otimes \ldots\otimes \psii\circ \iota_{X_i}\otimes \ldots \otimes \iota_{X_{|\ung|+n+|\unh|-1}}\otimes \id_{Y^{*}\otimes Y})=\phi\circ(\iota_{X_1}\otimes \ldots \otimes \iota_{X_{|\ung|+n+|\unh|-1}}\otimes \iota_{Y})\]
so, using the first equality of Figure \ref{universal},
\[|T_1^{o}|_{\C}\circ(\iota_{X_1}\otimes \ldots\otimes \iota_{X_{|\ung|+n+|\unh|}})=\phi\circ(\iota_{X_1}\otimes \ldots\otimes \iota_{X_{|\ung|+n+|\unh|}})\]
and by unicity of this factorization
\[|T_1^{o}|_{\C}=\phi.\]

Now, we compute the invariant $|T|_{\C,\underline{\alpha}}$ using $T^{o}_1$ and $T^{o}_2$. Then, we get equalities showed on Figure~\ref{invariant}.


\begin{figure}[H]
\centering
\resizebox{0.8\linewidth}{!}{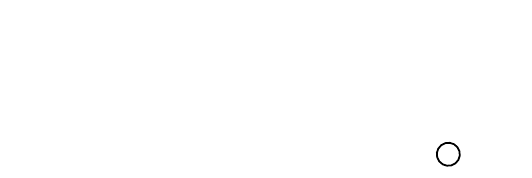}\vspace{0.5cm}
\resizebox{0.8\linewidth}{!}{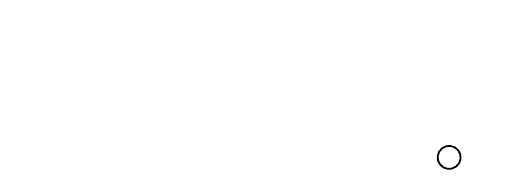}
\caption{Toward the computation of the invariant $|T|_{\C,\underline{\alpha}}$ using $T_{1}^{o}$ and $T_{2}^{o}$.}
\label{invariant}
\end{figure} 



But the braiding of $C$ is natural and $c_{Y,\mathds{1}}=\id_{Y}$ as shown in Figure~\ref{naturalproof}.

\begin{figure}[H]
\centering
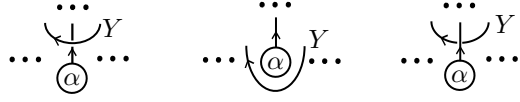
\caption{Naturality of the braiding and trivial braiding on $\mathds{1}$.}
\label{naturalproof}
\end{figure} 

Thus, as the second hand side of the equality in Figure~\ref{invariant} doesn't depend on the crossing between the braid colored by $Y$ and the $\alpha$ colored component, we have identities of Figure \ref{equalinv}:

\begin{figure}[H]
\centering
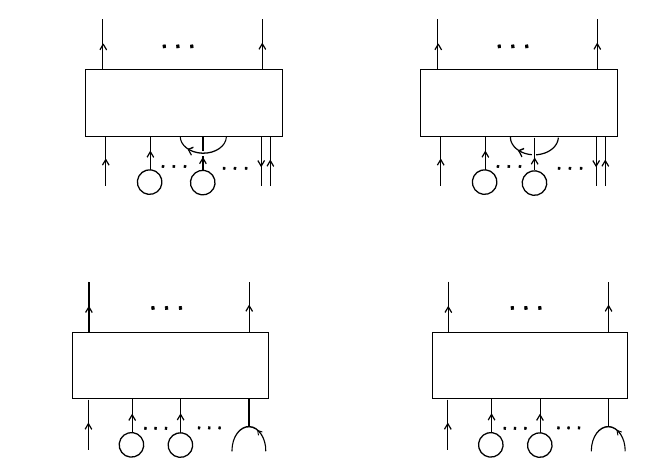
\caption{Equal invariants.}
\label{equalinv}
\end{figure} 
and using factorization property of the coend $C$, 
\[|T^{o}_1|_{\C}(\id_{C^{\otimes |\ung|}}\otimes \alpha^{\otimes n+|\unh|-1}\otimes \id_{C})=|T^{o}_2|_{\C}(\id_{C^{\otimes |\ung|}}\otimes \alpha^{\otimes n+|\unh|-1}\otimes \id_{C})\]
so composing with the missing $\alpha$, we have proved the result:
\[|T_1|_{\C,\alpha}=|T_2|_{\C,\alpha}.\]

Suppose now that $Y=X_i$ (see Figure~\ref{morphimst12}) that means a surgery component is auto-intersecting. In this case, without loss of generality, after sliding the concerning crossing to the bottom of the concerned tangle, it is enough to show that 
\begin{center}
$|
\begingroup%
  \makeatletter%
  \providecommand\color[2][]{%
    \errmessage{(Inkscape) Color is used for the text in Inkscape, but the package 'color.sty' is not loaded}%
    \renewcommand\color[2][]{}%
  }%
  \providecommand\transparent[1]{%
    \errmessage{(Inkscape) Transparency is used (non-zero) for the text in Inkscape, but the package 'transparent.sty' is not loaded}%
    \renewcommand\transparent[1]{}%
  }%
  \providecommand\rotatebox[2]{#2}%
  \newcommand*\fsize{\dimexpr\f@size pt\relax}%
  \newcommand*\lineheight[1]{\fontsize{\fsize}{#1\fsize}\selectfont}%
  \ifx\svgwidth\undefined%
    \setlength{\unitlength}{22.67716535bp}%
    \ifx\svgscale\undefined%
      \relax%
    \else%
      \setlength{\unitlength}{\unitlength * \real{\svgscale}}%
    \fi%
  \else%
    \setlength{\unitlength}{\svgwidth}%
  \fi%
  \global\let\svgwidth\undefined%
  \global\let\svgscale\undefined%
  \makeatother%
  \begin{picture}(1,0.75)%
    \lineheight{1}%
    \setlength\tabcolsep{0pt}%
    \put(0,0){\includegraphics[width=\unitlength,page=1]{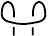}}%
  \end{picture}%
\endgroup%
|_{\C,\alpha}=|
\begingroup%
  \makeatletter%
  \providecommand\color[2][]{%
    \errmessage{(Inkscape) Color is used for the text in Inkscape, but the package 'color.sty' is not loaded}%
    \renewcommand\color[2][]{}%
  }%
  \providecommand\transparent[1]{%
    \errmessage{(Inkscape) Transparency is used (non-zero) for the text in Inkscape, but the package 'transparent.sty' is not loaded}%
    \renewcommand\transparent[1]{}%
  }%
  \providecommand\rotatebox[2]{#2}%
  \newcommand*\fsize{\dimexpr\f@size pt\relax}%
  \newcommand*\lineheight[1]{\fontsize{\fsize}{#1\fsize}\selectfont}%
  \ifx\svgwidth\undefined%
    \setlength{\unitlength}{22.67716535bp}%
    \ifx\svgscale\undefined%
      \relax%
    \else%
      \setlength{\unitlength}{\unitlength * \real{\svgscale}}%
    \fi%
  \else%
    \setlength{\unitlength}{\svgwidth}%
  \fi%
  \global\let\svgwidth\undefined%
  \global\let\svgscale\undefined%
  \makeatother%
  \begin{picture}(1,0.75)%
    \lineheight{1}%
    \setlength\tabcolsep{0pt}%
    \put(0,0){\includegraphics[width=\unitlength,page=1]{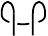}}%
  \end{picture}%
\endgroup%
|_{\C,\alpha}$
\end{center}

After applying the factorization property of the coend, we just have to prove that
\begin{center}
   \raisebox{-5mm}{\resizebox{0.25\linewidth}{!}{
\begingroup%
  \makeatletter%
  \providecommand\color[2][]{%
    \errmessage{(Inkscape) Color is used for the text in Inkscape, but the package 'color.sty' is not loaded}%
    \renewcommand\color[2][]{}%
  }%
  \providecommand\transparent[1]{%
    \errmessage{(Inkscape) Transparency is used (non-zero) for the text in Inkscape, but the package 'transparent.sty' is not loaded}%
    \renewcommand\transparent[1]{}%
  }%
  \providecommand\rotatebox[2]{#2}%
  \newcommand*\fsize{\dimexpr\f@size pt\relax}%
  \newcommand*\lineheight[1]{\fontsize{\fsize}{#1\fsize}\selectfont}%
  \ifx\svgwidth\undefined%
    \setlength{\unitlength}{85.03937008bp}%
    \ifx\svgscale\undefined%
      \relax%
    \else%
      \setlength{\unitlength}{\unitlength * \real{\svgscale}}%
    \fi%
  \else%
    \setlength{\unitlength}{\svgwidth}%
  \fi%
  \global\let\svgwidth\undefined%
  \global\let\svgscale\undefined%
  \makeatother%
  \begin{picture}(1,0.66666667)%
    \lineheight{1}%
    \setlength\tabcolsep{0pt}%
    \put(0,0){\includegraphics[width=\unitlength,page=1]{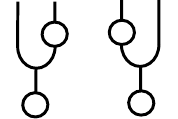}}%
    \put(0.26645968,0.4604015){\color[rgb]{0,0,0}\makebox(0,0)[lt]{\lineheight{1.25}\smash{\begin{tabular}[t]{l}\tiny{$\theta_C$}\end{tabular}}}}%
    \put(0.46212172,0.28714962){\color[rgb]{0,0,0}\makebox(0,0)[lt]{\lineheight{1.25}\smash{\begin{tabular}[t]{l}$=$\end{tabular}}}}%
    \put(0.6444841,0.47145541){\color[rgb]{0,0,0}\makebox(0,0)[lt]{\lineheight{1.25}\smash{\begin{tabular}[t]{l}\tiny{$\theta_C$}\end{tabular}}}}%
    \put(0.17143133,0.06185527){\color[rgb]{0,0,0}\makebox(0,0)[lt]{\lineheight{1.25}\smash{\begin{tabular}[t]{l}\tiny{$\alpha$}\end{tabular}}}}%
    \put(0.77001647,0.06921574){\color[rgb]{0,0,0}\makebox(0,0)[lt]{\lineheight{1.25}\smash{\begin{tabular}[t]{l}\tiny{$\alpha$}\end{tabular}}}}%
  \end{picture}%
\endgroup%
}}
\end{center}
Indeed, we get:
\begin{center}
   \resizebox{1.1\linewidth}{!}{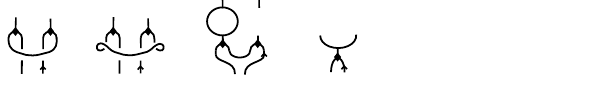}
\end{center}

Equality (2) is due to the definition of the coproduct and antipode for coend $C$ (see \cite{BrugVir}). Composing by $S\otimes S$ and colored by $\alpha$ the two morphisms obtained, and recalling that $S^2=\theta_C$, we get the wanted result using naturality and self-duality of $\theta_C$.

Secondly, suppose that diagrams of opentangles $T^{o}_1$ and $T^{o}_2$ differ only by one move $ESC$ (see Figure~\ref{moveesc}) as it is illustrated in Figure~\ref{proofmoveesc1}:
\begin{figure}[H]
\centering
\resizebox{1\linewidth}{!}{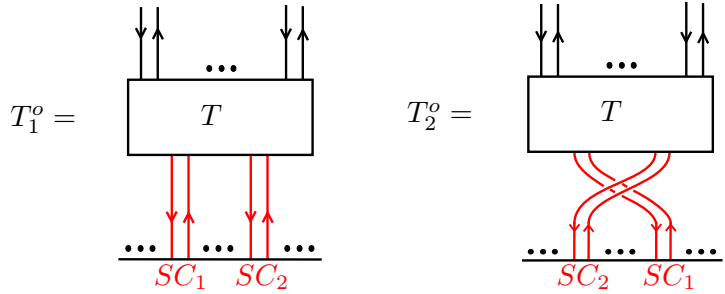}
\caption{One permutation in the ordered surgery components of an opentangle.}
\label{proofmoveesc1}
\end{figure} 
where $T$ is a $(\ung,n,\unh)$-opentangle. Applying Lemma~\ref{lemma0} to the $(\ung,n,\unh)$-opentangle $T$, there exists a morphism $|T|_{\C}\colon C^{\otimes |\ung|+n+|\unh|}\rightarrow C^{|\unh|}$ such that we have identities of Figure~\ref{proofmoveesc2}.
\begin{figure}[H]
\centering
\resizebox{1\linewidth}{!}{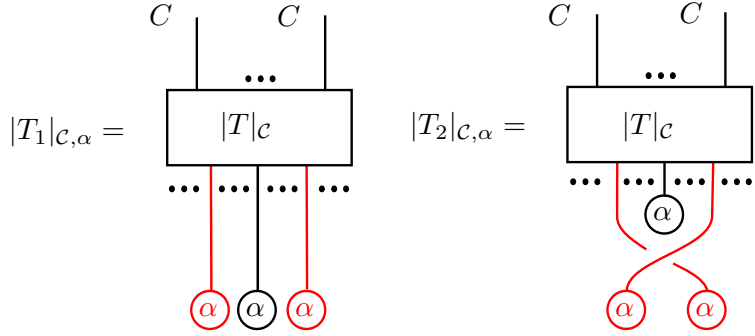}
\caption{Computation of the invariant $|\;\;|_{\C,\alpha}$ with $T^{o}_1$ and $T^{o}_2$.}
\label{proofmoveesc2}
\end{figure} 
As the braiding $\tau$ is natural and $\tau_{\mathds{1},\mathds{1}}=\id_{\mathds{1}}=\id_{\mathds{1}\otimes \mathds{1}}$, we have equalities of Figure~\ref{proofmoveesc3}.
\begin{figure}[H]
\centering
\resizebox{1\linewidth}{!}{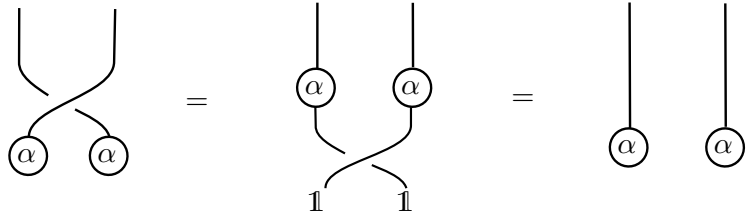}
\caption{Natural braiding and transparent object $\mathds{1}$ }
\label{proofmoveesc3}
\end{figure} 
And so $|T_1|_{\C,\alpha}=|T_2|_{\C,\alpha}$.

Lastly, suppose that diagrams of opentangles $T^{o}_1$ and $T^{o}_2$ differ only by one move $ROT$ (see Figure~\ref{moverot}). The invariance by this last move comes from naturality of $\theta$. Indeed we have then  $\theta_{\pm}\alpha=\alpha\theta_{\un}=\alpha$ since $\theta_{\un}=\id_{\un}$.
\end{proof}

Note that the isotopy invariant of cobordism tangles $|\phantom{O}|_{\C,\alpha}$ is multiplicative for the disjoint union as claimed in the following Lemma.
\begin{lemma}\label{lemmamultiplicative}
Let $\alpha \in \mr{Hom}_{\C}(\mathds{1},C)$ and $T_1$,$T_2$ be two cobordism tangles. Then \[|T_1\sqcup T_2|_{\C,\alpha}=|T_1|_{\C,\alpha}\otimes |T_2|_{\C,\alpha}.\]
\end{lemma}
\begin{proof}

Let us give the idea of the proof on simple cobordism tangles. Suppose that $T_1$ is a $((1),2,(1))$-cobordism tangle, $T_2$ is a $((1),1,(1,1))$-cobordism tangle. Denote by $T_1^{o}$ a $((1),2,(1))$-opentangle and by $T^{o}_2$ a $((1),1,(1,1))$-opentangle such that $U(T^{o}_1)=T_1$ and $U(T^{o}_2)=T_2$ (see \ref{surjectivemap}) that means we have equalities of Figure~\ref{proofconnectedsum1}.
\begin{figure}[!h]
\centering
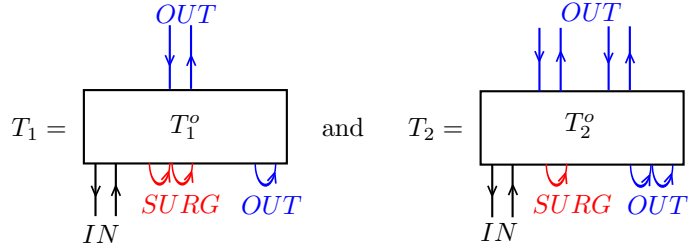
\caption{Opentangles $T_1^{o}$ and $T_2^{o}$ such that $T_1=U(T^{o}_{1})$ and $T_2=U(T^{o}_{2})$.}
\label{proofconnectedsum1}
\end{figure}

Choose the following $((1,1),3,(1,2))$-opentangle $O$ drawn on Figure~\ref{proofconnectedsum2}.
\begin{figure}[!h]
\centering
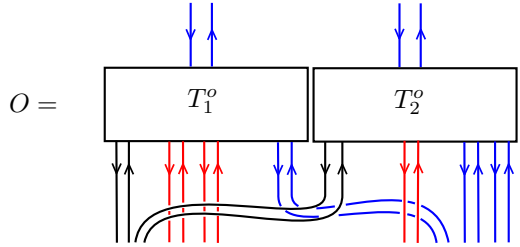
\caption{Opentangle $O$ such that $T_1\sqcup T_2=U(O)$.}
\label{proofconnectedsum2}
\end{figure}

Remark that the universal morphism $|O|_{\C}$ is equal to the morphism \[(|T^{o}_1|_{\C}\otimes |T^{o}_2|_{\C})\left(\begin{tikzpicture}[scale=0.3,description/.style={fill=white,inner sep=2pt},baseline=(current  bounding  box.center)]
\braid[number of strands=8] a_4^{-1} a_3^{-1} a_2^{-1} a_5^{-1} ;
\draw (1,-5.5) node[] {$C$};
\draw (2,-5.5) node[] {$C$};
\draw (3,-5.5) node[] {$C$};
\draw (4,-5.5) node[] {$C$};
\draw (5,-5.5) node[] {$C$};
\draw (6,-5.5) node[] {$C$};
\draw (7,-5.5) node[] {$C$};
\draw (8,-5.5) node[] {$C$};
\end{tikzpicture}\right)\]
and then
\[|T_1\sqcup T_2|_{\C,\alpha}=(|T^{o}_1|_{\C}\otimes|T^{o}_2|_{\C})\left(\begin{tikzpicture}[scale=0.3,description/.style={fill=white,inner sep=2pt},baseline=(current  bounding  box.center)]
\braid[number of strands=8] a_4^{-1} a_3^{-1} a_2^{-1} a_5^{-1} ;
\draw (1,-5.5) node[] {$C$};
\draw (2,-5.5) node[] {$C$};
\draw (3,-5.5) node[] {$C$};
\draw (4,-5.5) node[] {$C$};
\draw (5,-5.5) node[] {$C$};
\draw (6,-5.5) node[] {$C$};
\draw (7,-5.5) node[] {$C$};
\draw (8,-5.5) node[] {$C$};
\end{tikzpicture}\right)(\id_{C^{\otimes 2}}\otimes \alpha^{\otimes 6})=|T^{o}|\]
As the braiding $\tau$ is natural, $C\otimes \mathds{1}=\mathds{1}\otimes C=C$, $\tau_{\mathds{1},C}=\id_{C}$ and $\tau_{\mathds{1},\mathds{1}}=\id_{\mathds{1}}$, we have
\[(|T^{o}_1|_{\C}\otimes|T^{o}_2|_{\C})\left(\begin{tikzpicture}[scale=0.3,description/.style={fill=white,inner sep=2pt},baseline=(current  bounding  box.center)]
\braid[number of strands=8] a_4^{-1} a_3^{-1} a_2^{-1} a_5^{-1} ;
\draw (1,-5.5) node[] {$C$};
\draw (2,-5.5) node[] {$C$};
\draw (3,-5.5) node[] {$C$};
\draw (4,-5.5) node[] {$C$};
\draw (5,-5.5) node[] {$C$};
\draw (6,-5.5) node[] {$C$};
\draw (7,-5.5) node[] {$C$};
\draw (8,-5.5) node[] {$C$};
\end{tikzpicture}\right)(\id_{C^{\otimes 2}}\otimes \alpha^{\otimes 6})=|T^{o}_1|_{\C}(\id_{C}\otimes \alpha^{\otimes 3})\otimes |T^{o}_2|_{\C}(\id_{C}\otimes \alpha^{\otimes 3})\]
and we get the expected result.
\end{proof}

We need now to extend the isotopy invariant in a homeomorphism 3-cobordims invariant. In Section~\ref{hallowed}, we have defined hallowed tangles $\halo{T}$. This topological operation is equivalent to compose and precompose invariant $|T|_{\C,\alpha}$ with a particular morphism encoded by the pairing of the coend $\omega$. Let's us explicit this morphism. Let $X,Y$ be any objects of $\C$ and consider the morphism $\Pi_{X,Y}$ defined in Figure~\ref{usefulmorphism}:
\begin{figure}[H]
\centering
\resizebox{0.3\linewidth}{!}{\raisebox{15mm}{$\Pi_{X,Y}=$}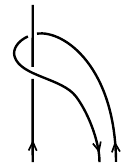}
\caption{The $\Pi_{X,Y}$ morphism}
\label{usefulmorphism}
\end{figure}
\noindent Let us express this morphism $\Pi_{X,Y}$ using structural morphisms of the coend $C$ of $\C$:
\begin{center}
\resizebox{0.2\linewidth}{!}{
\begingroup%
  \makeatletter%
  \providecommand\color[2][]{%
    \errmessage{(Inkscape) Color is used for the text in Inkscape, but the package 'color.sty' is not loaded}%
    \renewcommand\color[2][]{}%
  }%
  \providecommand\transparent[1]{%
    \errmessage{(Inkscape) Transparency is used (non-zero) for the text in Inkscape, but the package 'transparent.sty' is not loaded}%
    \renewcommand\transparent[1]{}%
  }%
  \providecommand\rotatebox[2]{#2}%
  \ifx\svgwidth\undefined%
    \setlength{\unitlength}{64bp}%
    \ifx\svgscale\undefined%
      \relax%
    \else%
      \setlength{\unitlength}{\unitlength * \real{\svgscale}}%
    \fi%
  \else%
    \setlength{\unitlength}{\svgwidth}%
  \fi%
  \global\let\svgwidth\undefined%
  \global\let\svgscale\undefined%
  \makeatother%
  \begin{picture}(1,1.25)%
    \put(0,0){\includegraphics[width=\unitlength]{usefulmorphism.pdf}}%
    \put(-0.01973751,0.15713871){\color[rgb]{0,0,0}\makebox(0,0)[lt]{\begin{minipage}{0.42182498\unitlength}\raggedright $X$\end{minipage}}}%
    \put(0.95828518,0.15118309){\color[rgb]{0,0,0}\makebox(0,0)[lt]{\begin{minipage}{0.42182498\unitlength}\raggedright $Y$\end{minipage}}}%
    \put(0.00677807,1.23110171){\color[rgb]{0,0,0}\makebox(0,0)[lt]{\begin{minipage}{0.42182498\unitlength}\raggedright $X$\end{minipage}}}%
  \end{picture}%
\endgroup%
}
\raisebox{9mm}{$=$}
\resizebox{0.2\linewidth}{!}{
\begingroup%
  \makeatletter%
  \providecommand\color[2][]{%
    \errmessage{(Inkscape) Color is used for the text in Inkscape, but the package 'color.sty' is not loaded}%
    \renewcommand\color[2][]{}%
  }%
  \providecommand\transparent[1]{%
    \errmessage{(Inkscape) Transparency is used (non-zero) for the text in Inkscape, but the package 'transparent.sty' is not loaded}%
    \renewcommand\transparent[1]{}%
  }%
  \providecommand\rotatebox[2]{#2}%
  \ifx\svgwidth\undefined%
    \setlength{\unitlength}{64bp}%
    \ifx\svgscale\undefined%
      \relax%
    \else%
      \setlength{\unitlength}{\unitlength * \real{\svgscale}}%
    \fi%
  \else%
    \setlength{\unitlength}{\svgwidth}%
  \fi%
  \global\let\svgwidth\undefined%
  \global\let\svgscale\undefined%
  \makeatother%
  \begin{picture}(1,1.25)%
    \put(0,0){\includegraphics[width=\unitlength]{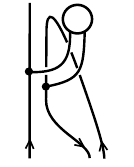}}%
    \put(0.01330559,0.18364521){\color[rgb]{0,0,0}\makebox(0,0)[lt]{\begin{minipage}{0.42182498\unitlength}\raggedright $X$\end{minipage}}}%
    \put(0.81629933,0.16654434){\color[rgb]{0,0,0}\makebox(0,0)[lt]{\begin{minipage}{0.42182498\unitlength}\raggedright $Y$\end{minipage}}}%
    \put(0.02080357,1.23580621){\color[rgb]{0,0,0}\makebox(0,0)[lt]{\begin{minipage}{0.42182498\unitlength}\raggedright $X$\end{minipage}}}%
    \put(0.54048296,1.13874938){\color[rgb]{0,0,0}\makebox(0,0)[lt]{\begin{minipage}{0.17580648\unitlength}\raggedright $\omega$\end{minipage}}}%
  \end{picture}%
\endgroup%
}
\raisebox{9mm}{$=$}
\resizebox{0.2\linewidth}{!}{
\begingroup%
  \makeatletter%
  \providecommand\color[2][]{%
    \errmessage{(Inkscape) Color is used for the text in Inkscape, but the package 'color.sty' is not loaded}%
    \renewcommand\color[2][]{}%
  }%
  \providecommand\transparent[1]{%
    \errmessage{(Inkscape) Transparency is used (non-zero) for the text in Inkscape, but the package 'transparent.sty' is not loaded}%
    \renewcommand\transparent[1]{}%
  }%
  \providecommand\rotatebox[2]{#2}%
  \ifx\svgwidth\undefined%
    \setlength{\unitlength}{64bp}%
    \ifx\svgscale\undefined%
      \relax%
    \else%
      \setlength{\unitlength}{\unitlength * \real{\svgscale}}%
    \fi%
  \else%
    \setlength{\unitlength}{\svgwidth}%
  \fi%
  \global\let\svgwidth\undefined%
  \global\let\svgscale\undefined%
  \makeatother%
  \begin{picture}(1,1.25)%
    \put(0,0){\includegraphics[width=\unitlength]{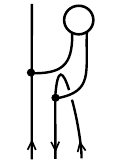}}%
    \put(0.02893698,0.17567384){\color[rgb]{0,0,0}\makebox(0,0)[lt]{\begin{minipage}{0.42182498\unitlength}\raggedright $X$\end{minipage}}}%
    \put(0.65489727,0.17856146){\color[rgb]{0,0,0}\makebox(0,0)[lt]{\begin{minipage}{0.42182498\unitlength}\raggedright $Y$\end{minipage}}}%
    \put(0.03643495,1.22783534){\color[rgb]{0,0,0}\makebox(0,0)[lt]{\begin{minipage}{0.42182498\unitlength}\raggedright $X$\end{minipage}}}%
    \put(0.55025061,1.13198388){\color[rgb]{0,0,0}\makebox(0,0)[lt]{\begin{minipage}{0.17580648\unitlength}\raggedright $\omega$\end{minipage}}}%
  \end{picture}%
\endgroup%
}
\raisebox{9mm}{$=$}
\resizebox{0.2\linewidth}{!}{
\begingroup%
  \makeatletter%
  \providecommand\color[2][]{%
    \errmessage{(Inkscape) Color is used for the text in Inkscape, but the package 'color.sty' is not loaded}%
    \renewcommand\color[2][]{}%
  }%
  \providecommand\transparent[1]{%
    \errmessage{(Inkscape) Transparency is used (non-zero) for the text in Inkscape, but the package 'transparent.sty' is not loaded}%
    \renewcommand\transparent[1]{}%
  }%
  \providecommand\rotatebox[2]{#2}%
  \ifx\svgwidth\undefined%
    \setlength{\unitlength}{64bp}%
    \ifx\svgscale\undefined%
      \relax%
    \else%
      \setlength{\unitlength}{\unitlength * \real{\svgscale}}%
    \fi%
  \else%
    \setlength{\unitlength}{\svgwidth}%
  \fi%
  \global\let\svgwidth\undefined%
  \global\let\svgscale\undefined%
  \makeatother%
  \begin{picture}(1,1.25)%
    \put(0,0){\includegraphics[width=\unitlength]{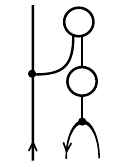}}%
    \put(0.03695818,0.16657484){\color[rgb]{0,0,0}\makebox(0,0)[lt]{\begin{minipage}{0.42182498\unitlength}\raggedright $X$\end{minipage}}}%
    \put(0.7710588,0.17505596){\color[rgb]{0,0,0}\makebox(0,0)[lt]{\begin{minipage}{0.42182498\unitlength}\raggedright $Y$\end{minipage}}}%
    \put(0.04445615,1.21873634){\color[rgb]{0,0,0}\makebox(0,0)[lt]{\begin{minipage}{0.42182498\unitlength}\raggedright $X$\end{minipage}}}%
    \put(0.57027701,0.70170739){\color[rgb]{0,0,0}\makebox(0,0)[lt]{\begin{minipage}{0.27345618\unitlength}\raggedright $S$\end{minipage}}}%
    \put(0.55066138,1.1173885){\color[rgb]{0,0,0}\makebox(0,0)[lt]{\begin{minipage}{0.17580648\unitlength}\raggedright $\omega$\end{minipage}}}%
  \end{picture}%
\endgroup%
}
\end{center}

Suppose that $\alpha$ satisfies \ref{ad2}. Then, for $n$ a non negative integer, we define the \emph{hallowed morphism}, denoted by $\Pi_{\alpha,n}$, we will use to compose our invariant $|T|_{\C,\alpha}$:
\begin{align}\label{hallowedmorphism}
\Pi_{\alpha,n}=[\mr{id}_{C}\otimes\omega(\mr{id}_{C}\otimes \alpha)]\delta_{C^{\otimes n}}
=
\raisebox{-7mm}{
\begingroup%
  \makeatletter%
  \providecommand\color[2][]{%
    \errmessage{(Inkscape) Color is used for the text in Inkscape, but the package 'color.sty' is not loaded}%
    \renewcommand\color[2][]{}%
  }%
  \providecommand\transparent[1]{%
    \errmessage{(Inkscape) Transparency is used (non-zero) for the text in Inkscape, but the package 'transparent.sty' is not loaded}%
    \renewcommand\transparent[1]{}%
  }%
  \providecommand\rotatebox[2]{#2}%
  \newcommand*\fsize{\dimexpr\f@size pt\relax}%
  \newcommand*\lineheight[1]{\fontsize{\fsize}{#1\fsize}\selectfont}%
  \ifx\svgwidth\undefined%
    \setlength{\unitlength}{56.69291339bp}%
    \ifx\svgscale\undefined%
      \relax%
    \else%
      \setlength{\unitlength}{\unitlength * \real{\svgscale}}%
    \fi%
  \else%
    \setlength{\unitlength}{\svgwidth}%
  \fi%
  \global\let\svgwidth\undefined%
  \global\let\svgscale\undefined%
  \makeatother%
  \begin{picture}(1,1)%
    \lineheight{1}%
    \setlength\tabcolsep{0pt}%
    \put(0,0){\includegraphics[width=\unitlength,page=1]{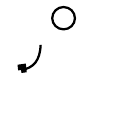}}%
    \put(0.48818057,0.808867){\color[rgb]{0,0,0}\makebox(0,0)[lt]{\lineheight{0}\smash{\begin{tabular}[t]{l}$\omega$\end{tabular}}}}%
    \put(0,0){\includegraphics[width=\unitlength,page=2]{projection235.pdf}}%
    \put(0.66816854,0.53919594){\color[rgb]{0,0,0}\makebox(0,0)[lt]{\lineheight{1.25}\smash{\begin{tabular}[t]{l}$\alpha$\end{tabular}}}}%
    \put(0,0){\includegraphics[width=\unitlength,page=3]{projection235.pdf}}%
    \put(0.22005785,0.04220458){\color[rgb]{0,0,0}\makebox(0,0)[lt]{\lineheight{1.25}\smash{\begin{tabular}[t]{l}$C^{\otimes n}$\end{tabular}}}}%
  \end{picture}%
\endgroup%
}
\end{align}

In the Reshetikhin-Turaev TQFT, this morphism is the \emph{transparent projector} which sends every object on its transparent part. To construct the internal TQFT, we will need to ask this morphism to be a projector. For now, trying to construct a 3-cobordisms invariant, we only need the hallowed morphism to make a link between invariant $|T|_{\C,\alpha}$ and the invariant $|\halo{T}|_{\C,\alpha}$.
\begin{lemma}\label{lemmahalonothalo}
Let $\alpha\colon \un \rightarrow C$ and $T$ be a $(\ung,n,\unh)$-cobordism tangle. Then 
\begin{align}
\Pi_{\alpha,\unh}|T|_{\C,\alpha}\Pi_{\alpha,\ung}=|\halo{T}|_{\C,\alpha}
\end{align}
\end{lemma}
\begin{proof}
Just observe the topological operation that consists to add encircling entrance and exit components is exactly algebraically encoded by composing with hallowed morphism (\ref{hallowedmorphism}).
\end{proof}

Before defining a Kirby II move invariant, we need a last result concerning \ref{ad5}.
\begin{lemma}\label{lemmanaturalkirby}\textcolor{white}{o}
\begin{flushleft}
\raisebox{-10mm}{If \ref{ad5} holds,}
\raisebox{-10mm}{then for any $n\in \mathbb{N}$,}
\raisebox{-25mm}{
\begingroup%
  \makeatletter%
  \providecommand\color[2][]{%
    \errmessage{(Inkscape) Color is used for the text in Inkscape, but the package 'color.sty' is not loaded}%
    \renewcommand\color[2][]{}%
  }%
  \providecommand\transparent[1]{%
    \errmessage{(Inkscape) Transparency is used (non-zero) for the text in Inkscape, but the package 'transparent.sty' is not loaded}%
    \renewcommand\transparent[1]{}%
  }%
  \providecommand\rotatebox[2]{#2}%
  \ifx\svgwidth\undefined%
    \setlength{\unitlength}{160bp}%
    \ifx\svgscale\undefined%
      \relax%
    \else%
      \setlength{\unitlength}{\unitlength * \real{\svgscale}}%
    \fi%
  \else%
    \setlength{\unitlength}{\svgwidth}%
  \fi%
  \global\let\svgwidth\undefined%
  \global\let\svgscale\undefined%
  \makeatother%
  \begin{picture}(1,0.5)%
    \put(0,0){\includegraphics[width=\unitlength]{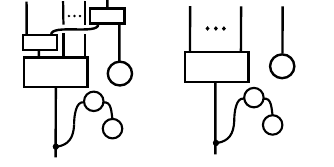}}%
    \put(0.09558464,0.30368775){\color[rgb]{0,0,0}\makebox(0,0)[lt]{\begin{minipage}{0.13708045\unitlength}\raggedright $\mathrm{id}_{C^{\otimes n}}$\end{minipage}}}%
    \put(0.26792941,0.18677044){\color[rgb]{0,0,0}\makebox(0,0)[lb]{\smash{\tiny{$\omega$}}}}%
    \put(0.34286215,0.26155201){\color[rgb]{0,0,0}\makebox(0,0)[lb]{\smash{$\alpha$}}}%
    \put(0.31782106,0.09905628){\color[rgb]{0,0,0}\makebox(0,0)[lb]{\smash{$\alpha$}}}%
    \put(0.02587217,0.05823856){\color[rgb]{0,0,0}\makebox(0,0)[lt]{\begin{minipage}{0.06531835\unitlength}\raggedright $C^{\otimes n}$\end{minipage}}}%
    \put(0.10233746,0.35822136){\color[rgb]{0,0,0}\makebox(0,0)[lb]{\smash{\tiny{$\Delta$}}}}%
    \put(0.30253616,0.44221373){\color[rgb]{0,0,0}\makebox(0,0)[lb]{\smash{\tiny{$m$}}}}%
    \put(0.58092074,0.27900679){\color[rgb]{0,0,0}\makebox(0,0)[lb]{\smash{$\mathrm{id}_{C^{\otimes n}}$}}}%
    \put(0.74832149,0.19689688){\color[rgb]{0,0,0}\makebox(0,0)[lb]{\smash{\tiny{$\omega$}}}}%
    \put(0.82674665,0.28793587){\color[rgb]{0,0,0}\makebox(0,0)[lb]{\smash{$\alpha$}}}%
    \put(0.79868832,0.11412236){\color[rgb]{0,0,0}\makebox(0,0)[lb]{\smash{$\alpha$}}}%
    \put(0.50279458,0.07354756){\color[rgb]{0,0,0}\makebox(0,0)[lt]{\begin{minipage}{0.06531835\unitlength}\raggedright $C^{\otimes n}$\end{minipage}}}%
    \put(0.43599374,0.31911206){\color[rgb]{0,0,0}\makebox(0,0)[lt]{\begin{minipage}{0.09094481\unitlength}\raggedright $=$\end{minipage}}}%
  \end{picture}%
\endgroup%
}
\raisebox{-10mm}{.}
\end{flushleft}
\end{lemma}

\begin{proof} 
Set $\Pi_{n}=[\mr{id}_{C^{\otimes n}}\otimes\omega(\mr{id}_{C}\otimes \alpha)]\delta_{C^{\otimes n}}$. We have:
\begin{align*}
&(\mr{id}_{C^{\otimes n}}\otimes m(\mr{id}_{C}\otimes \alpha))(\mr{id}_{C^{\otimes n-1}}\otimes \tau_{C,C})(\mr{id}_{C}\otimes \tau_{C,C^{\otimes n-2}}\otimes \mr{id}_{C})\Pi_n \\
&=(\tau_{C,C^{\otimes n-2}}\otimes \mr{id_{C^{\otimes 2}}})(\mr{id}_{C}\otimes \tau_{C^{\otimes n-2},C}\otimes m(\mr{id}_{C}\otimes \alpha)))(\mr{id}_{C^{\otimes n-1}}\otimes \Delta)(\mr{id}_{C}\otimes \tau_{C^{\otimes n-2},C})(\tau_{C,C^{\otimes n-2}}\otimes \mr{id}_{C})\Pi_n\\
&\stackrel{(2)}{=}(\tau_{C,C^{\otimes n-2}}\otimes \mr{id_{C^{\otimes 2}}})(\mr{id}_{C}\otimes \tau_{C^{\otimes n-2},C}\otimes m(\mr{id}_{C}\otimes \alpha)))(\mr{id}_{C^{\otimes n-1}}\otimes \Delta)\Pi_n (\mr{id}_{C}\otimes \tau_{C^{\otimes n-2},C})(\tau_{C,C^{\otimes n-2}}\otimes \mr{id}_{C})\\
&\stackrel{(3)}{=}(\mr{\tau_{C,C^{n-2}}}\otimes \mr{id}_{C}\otimes \alpha)(\mr{id_{C}}\otimes \tau_{C^{\otimes n-2,C}})\Pi_n (\mr{id}_{C}\otimes \tau_{C^{\otimes n-2},C})(\tau_{C,C^{\otimes n-2}}\otimes \mr{id}_{C})\\
&\stackrel{(4)}{=}(\mr{\tau_{C,C^{n-2}}}\otimes \mr{id}_{C}\otimes \alpha)(\mr{id_{C}}\otimes \tau_{C^{\otimes n-2,C}})(\mr{id}_{C}\otimes \tau_{C^{\otimes n-2},C})(\tau_{C,C^{\otimes n-2}}\otimes \mr{id}_{C})\Pi_n \\
&=\Pi_n\otimes \alpha\\
\end{align*}
Equalities $(2)$ and $(4)$ are due to the naturality of $\Pi_n$ between identity functors whereas equality $(3)$ comes from the "If" part of the Lemma.
\end{proof}

In order to construct a 3-cobordisms invariant, we first develop a generalized Kirby II move $KII^{g}$ invariant (see Figure~\ref{k2}).

\begin{lemma}\label{lemmainvariantkirby2}
Let $\alpha\colon \un \rightarrow C$ that satisfies \ref{ad1}, \ref{ad2}, \ref{ad5} and $T$ be a $(\ung,n,\unh)$-cobordism tangle. Then the morphism
\[|\halo{T}|_{\C,\alpha}\colon C^{\otimes \ung}\rightarrow C^{\otimes \unh}\]
is invariant by the generalized Kirby move $KII^{g}$ (see Figure~\ref{k2}) on cobordism tangle $T$.
\end{lemma}
\begin{proof}
Let $T_1$ and $T_{2}$ be two $(\ung,n,\unh)$-cobordism tangles that differ by one move $KII^{g}$ . Let us compute $|\halo{T_1}|_{\C,\alpha}$ and $|\halo{T_2}|_{\C,\alpha}$.

\noindent \textbullet\; First, suppose that an entrance component slides over a surgery component and without loss of generality, we suppose that the entrance component and the surgery component are the first ones as illustrated on the following picture:
\begin{center}
\begingroup%
  \makeatletter%
  \providecommand\color[2][]{%
    \errmessage{(Inkscape) Color is used for the text in Inkscape, but the package 'color.sty' is not loaded}%
    \renewcommand\color[2][]{}%
  }%
  \providecommand\transparent[1]{%
    \errmessage{(Inkscape) Transparency is used (non-zero) for the text in Inkscape, but the package 'transparent.sty' is not loaded}%
    \renewcommand\transparent[1]{}%
  }%
  \providecommand\rotatebox[2]{#2}%
  \newcommand*\fsize{\dimexpr\f@size pt\relax}%
  \newcommand*\lineheight[1]{\fontsize{\fsize}{#1\fsize}\selectfont}%
  \ifx\svgwidth\undefined%
    \setlength{\unitlength}{307.5bp}%
    \ifx\svgscale\undefined%
      \relax%
    \else%
      \setlength{\unitlength}{\unitlength * \real{\svgscale}}%
    \fi%
  \else%
    \setlength{\unitlength}{\svgwidth}%
  \fi%
  \global\let\svgwidth\undefined%
  \global\let\svgscale\undefined%
  \makeatother%
  \begin{picture}(1,0.26829268)%
    \lineheight{1}%
    \setlength\tabcolsep{0pt}%
    \put(0,0){\includegraphics[width=\unitlength,page=1]{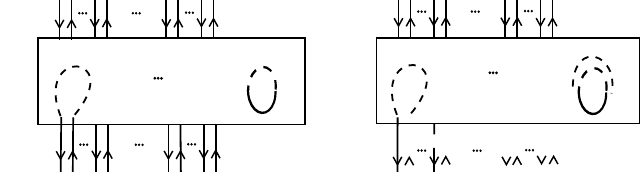}}%
    \put(0.51129599,0.17101844){\color[rgb]{0,0,0}\makebox(0,0)[lt]{\begin{minipage}{0.13082143\unitlength}\raggedright $T_2=$\end{minipage}}}%
    \put(-0.01829204,0.16337342){\color[rgb]{0,0,0}\makebox(0,0)[lt]{\begin{minipage}{0.13082143\unitlength}\raggedright $T_1=$\end{minipage}}}%
    \put(0,0){\includegraphics[width=\unitlength,page=2]{proofkirby12.pdf}}%
  \end{picture}%
\endgroup%

\end{center} 
Note that we can always suppose that the "sliding" part of the entrance component is located on the bottom of the picture (if not, you could "transport" by isotopy the little piece of the entrance component of $T_2$ that had slided over the surgery component to the bottom of the tangle $T_2$). 
We choose two opentangles $O_1$ and $O_2$ coming respectively from $T_{1}$ and $T_{2}$ (that means $U(O_1)=T_1$ and $U(O_2)=T_2$) as shown just below:
\begin{center}
\resizebox{!}{0.2\textheight}{
\begingroup%
  \makeatletter%
  \providecommand\color[2][]{%
    \errmessage{(Inkscape) Color is used for the text in Inkscape, but the package 'color.sty' is not loaded}%
    \renewcommand\color[2][]{}%
  }%
  \providecommand\transparent[1]{%
    \errmessage{(Inkscape) Transparency is used (non-zero) for the text in Inkscape, but the package 'transparent.sty' is not loaded}%
    \renewcommand\transparent[1]{}%
  }%
  \providecommand\rotatebox[2]{#2}%
  \newcommand*\fsize{\dimexpr\f@size pt\relax}%
  \newcommand*\lineheight[1]{\fontsize{\fsize}{#1\fsize}\selectfont}%
  \ifx\svgwidth\undefined%
    \setlength{\unitlength}{425.19685039bp}%
    \ifx\svgscale\undefined%
      \relax%
    \else%
      \setlength{\unitlength}{\unitlength * \real{\svgscale}}%
    \fi%
  \else%
    \setlength{\unitlength}{\svgwidth}%
  \fi%
  \global\let\svgwidth\undefined%
  \global\let\svgscale\undefined%
  \makeatother%
  \begin{picture}(1,0.26666667)%
    \lineheight{1}%
    \setlength\tabcolsep{0pt}%
    \put(0,0){\includegraphics[width=\unitlength,page=1]{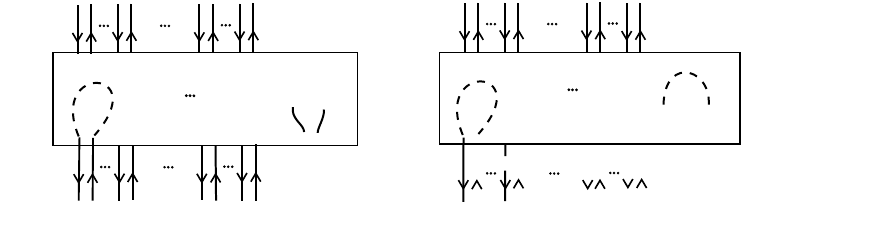}}%
    \put(0.43294185,0.1779724){\color[rgb]{0,0,0}\makebox(0,0)[lt]{\begin{minipage}{0.10803061\unitlength}\raggedright $O_2=$\end{minipage}}}%
    \put(-0.00438495,0.17203105){\color[rgb]{0,0,0}\makebox(0,0)[lt]{\begin{minipage}{0.10803061\unitlength}\raggedright $O_1=$\end{minipage}}}%
    \put(0,0){\includegraphics[width=\unitlength,page=2]{proofkirby2222.pdf}}%
    \put(0.26616961,0.02887037){\color[rgb]{0,0,0}\rotatebox{-34.0118342}{\makebox(0,0)[lt]{\lineheight{1.25}\smash{\begin{tabular}[t]{l}$|\ung|$\end{tabular}}}}}%
    \put(0.69533075,0.02658193){\color[rgb]{0,0,0}\rotatebox{-29.47395042}{\makebox(0,0)[lt]{\lineheight{1.25}\smash{\begin{tabular}[t]{l}$|\ung|$\end{tabular}}}}}%
    \put(0.08022658,0.02181243){\color[rgb]{0,0,0}\rotatebox{-5.13710752}{\makebox(0,0)[lt]{\lineheight{1.25}\smash{\begin{tabular}[t]{l}$1$\end{tabular}}}}}%
    \put(0.51528265,0.01937188){\color[rgb]{0,0,0}\rotatebox{-1.94226635}{\makebox(0,0)[lt]{\lineheight{1.25}\smash{\begin{tabular}[t]{l}$1$\end{tabular}}}}}%
    \put(0.63155938,0.03331022){\color[rgb]{0,0,0}\rotatebox{-28.55723644}{\makebox(0,0)[lt]{\lineheight{1.25}\smash{\begin{tabular}[t]{l}\small{$|\ung|-g_r+1$}\end{tabular}}}}}%
    \put(0.19622118,0.04483511){\color[rgb]{0,0,0}\rotatebox{-35.27591831}{\makebox(0,0)[lt]{\lineheight{1.25}\smash{\begin{tabular}[t]{l}\small{$|\ung|-g_r+1$}\end{tabular}}}}}%
    \put(0.12126327,0.02029409){\color[rgb]{0,0,0}\rotatebox{-3.52623405}{\makebox(0,0)[lt]{\lineheight{1.25}\smash{\begin{tabular}[t]{l}$g_1$\end{tabular}}}}}%
    \put(0.55818236,0.01911952){\color[rgb]{0,0,0}\rotatebox{-0.96030764}{\makebox(0,0)[lt]{\lineheight{1.25}\smash{\begin{tabular}[t]{l}$g_1$\end{tabular}}}}}%
    \put(0.75104061,0.02731332){\color[rgb]{0,0,0}\rotatebox{-27.61687378}{\makebox(0,0)[lt]{\lineheight{1.25}\smash{\begin{tabular}[t]{l}$|\ung|+1$\end{tabular}}}}}%
    \put(0.31119478,0.03668787){\color[rgb]{0,0,0}\rotatebox{-32.16166693}{\makebox(0,0)[lt]{\lineheight{1.25}\smash{\begin{tabular}[t]{l}$|\ung|+1$\end{tabular}}}}}%
  \end{picture}%
\endgroup%
}
\end{center}
 After coloring $O_1$ and $O_2$ by objects of $\C$ (we only particularize the color $X$ corresponding to the first entrance component of $O_1$ and $O_2$ and the color $Y$ on the first surgery component of $O_1$ and $O_2$) and composing this morphism of $\C$ with the universal dinatural action $\iota$ of the coend tensored as many times as the number of exit components $|\unh|$, we obtain a dinatural transformation $d_{X,\ldots,Y,\ldots}$ which is dinatural in every entrance pairs. Thus, 
 
\begin{center}
\resizebox{!}{0.3\textheight}{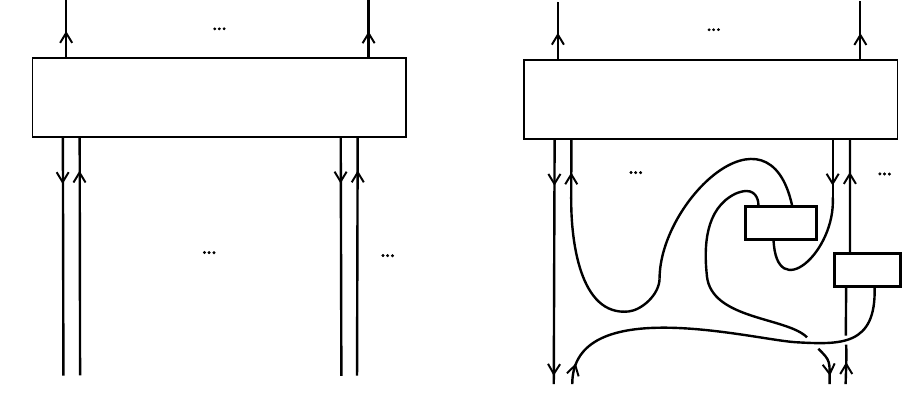}
\end{center} 

By universal property of the coend, there exists a unique morphism $\phi\colon C^{\otimes |\ung|+n+|\unh|}\rightarrow C^{|\unh|}$ such that $\forall X\in \C, \ldots, \forall Y\in \C,\ldots$,
\[
d_{X,\ldots,Y,\ldots}=\phi(\iota_X \otimes \ldots \otimes \iota_Y \otimes \ldots)
\]
and then, we can factorize the two last diagrams using the morphism $\phi$ and the universal action $\iota$:
\begin{center}
\resizebox{!}{0.35\textheight}{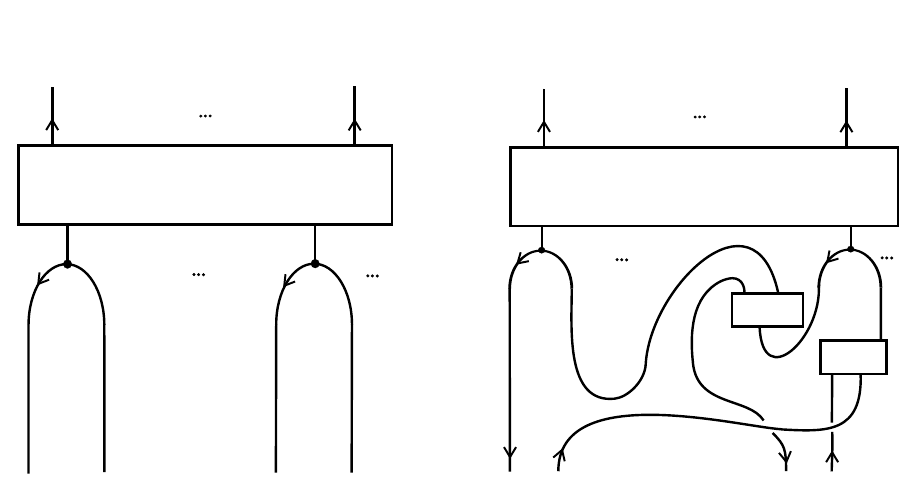}
\end{center}
Using the definition of $\Delta$, the morphism $\iota^{\otimes |\unh|}\circ O_{2;X,\ldots,Y,\ldots}$ is equal to 
\begin{center}
\resizebox{!}{0.25\textheight}{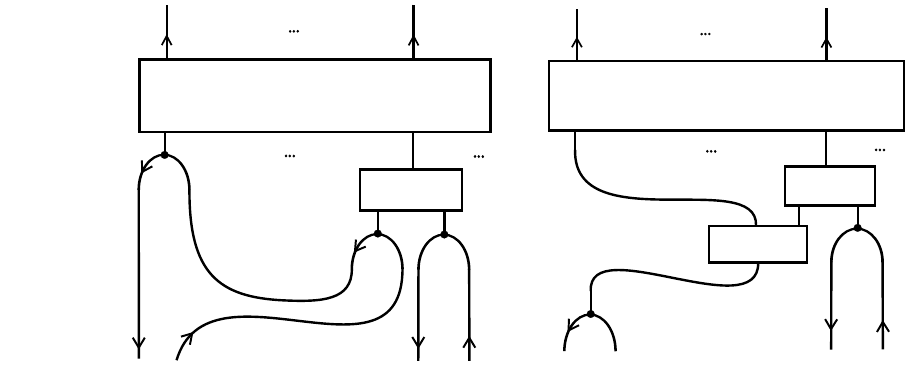}
\end{center}
Consequently, the morphism $|T_2|_{\C,\alpha}$ is given by:
\begin{center}
\resizebox{!}{0.25\textheight}{
\begingroup%
  \makeatletter%
  \providecommand\color[2][]{%
    \errmessage{(Inkscape) Color is used for the text in Inkscape, but the package 'color.sty' is not loaded}%
    \renewcommand\color[2][]{}%
  }%
  \providecommand\transparent[1]{%
    \errmessage{(Inkscape) Transparency is used (non-zero) for the text in Inkscape, but the package 'transparent.sty' is not loaded}%
    \renewcommand\transparent[1]{}%
  }%
  \providecommand\rotatebox[2]{#2}%
  \ifx\svgwidth\undefined%
    \setlength{\unitlength}{256bp}%
    \ifx\svgscale\undefined%
      \relax%
    \else%
      \setlength{\unitlength}{\unitlength * \real{\svgscale}}%
    \fi%
  \else%
    \setlength{\unitlength}{\svgwidth}%
  \fi%
  \global\let\svgwidth\undefined%
  \global\let\svgscale\undefined%
  \makeatother%
  \begin{picture}(1,0.75)%
    \put(0,0){\includegraphics[width=\unitlength]{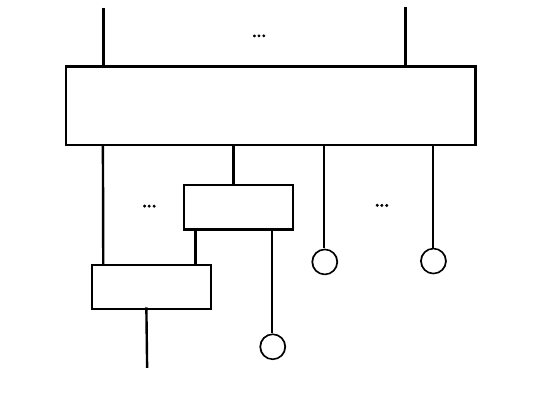}}%
    \put(0.47298964,0.57239576){\color[rgb]{0,0,0}\makebox(0,0)[lt]{\begin{minipage}{0.45164785\unitlength}\raggedright \Large{$\phi$}\end{minipage}}}%
    \put(0.21220398,0.7362417){\color[rgb]{0,0,0}\makebox(0,0)[lt]{\begin{minipage}{0.08772365\unitlength}\raggedright $C$\end{minipage}}}%
    \put(0.78532063,0.73931136){\color[rgb]{0,0,0}\makebox(0,0)[lt]{\begin{minipage}{0.08772365\unitlength}\raggedright $C$\end{minipage}}}%
    \put(0.45796789,0.03128916){\color[rgb]{0,0,0}\makebox(0,0)[lt]{\begin{minipage}{0.09932385\unitlength}\raggedright $|\ung|+1$\end{minipage}}}%
    \put(0.20675225,0.4559451){\color[rgb]{0,0,0}\makebox(0,0)[lt]{\begin{minipage}{0.05396422\unitlength}\raggedright $C$\end{minipage}}}%
    \put(0.46511392,0.45875613){\color[rgb]{0,0,0}\makebox(0,0)[lt]{\begin{minipage}{0.05396422\unitlength}\raggedright $C$\end{minipage}}}%
    \put(0.40193271,0.30128191){\color[rgb]{0,0,0}\makebox(0,0)[lt]{\begin{minipage}{0.05396422\unitlength}\raggedright $C$\end{minipage}}}%
    \put(0.52502942,0.30333032){\color[rgb]{0,0,0}\makebox(0,0)[lt]{\begin{minipage}{0.05396422\unitlength}\raggedright $C$\end{minipage}}}%
    \put(0.42866731,0.3714988){\color[rgb]{0,0,0}\makebox(0,0)[lt]{\begin{minipage}{0.10134591\unitlength}\raggedright $m$\end{minipage}}}%
    \put(0.2647444,0.22580821){\color[rgb]{0,0,0}\makebox(0,0)[lt]{\begin{minipage}{0.12402891\unitlength}\raggedright $\Delta$\end{minipage}}}%
    \put(0.257193,0.03974314){\color[rgb]{0,0,0}\makebox(0,0)[lt]{\begin{minipage}{0.14588983\unitlength}\raggedright $1$\end{minipage}}}%
    \put(0.5015014,0.10501134){\color[rgb]{0,0,0}\makebox(0,0)[lt]{\begin{minipage}{0.08204947\unitlength}\raggedright \small{$\alpha$}\end{minipage}}}%
    \put(0.62781786,0.46362323){\color[rgb]{0,0,0}\makebox(0,0)[lt]{\begin{minipage}{0.05396422\unitlength}\raggedright $C$\end{minipage}}}%
    \put(0.59907003,0.26409969){\color[rgb]{0,0,0}\makebox(0,0)[lt]{\begin{minipage}{0.08204947\unitlength}\raggedright \small{$\alpha$}\end{minipage}}}%
    \put(0.83159564,0.46524054){\color[rgb]{0,0,0}\makebox(0,0)[lt]{\begin{minipage}{0.05396422\unitlength}\raggedright $C$\end{minipage}}}%
    \put(0.80284781,0.265717){\color[rgb]{0,0,0}\makebox(0,0)[lt]{\begin{minipage}{0.08204947\unitlength}\raggedright \small{$\alpha$}\end{minipage}}}%
    \put(-0.15899573,0.55421371){\color[rgb]{0,0,0}\makebox(0,0)[lt]{\begin{minipage}{0.29414213\unitlength}\raggedright $|T_2|_{\C,\alpha}=$\end{minipage}}}%
  \end{picture}%
\endgroup%
}
\end{center}
so, according to Lemma~\ref{lemmahalonothalo}, the morphism $|\halo{T}_{\C,\alpha}|$ is given by:
\begin{center}
\resizebox{!}{0.42\textheight}{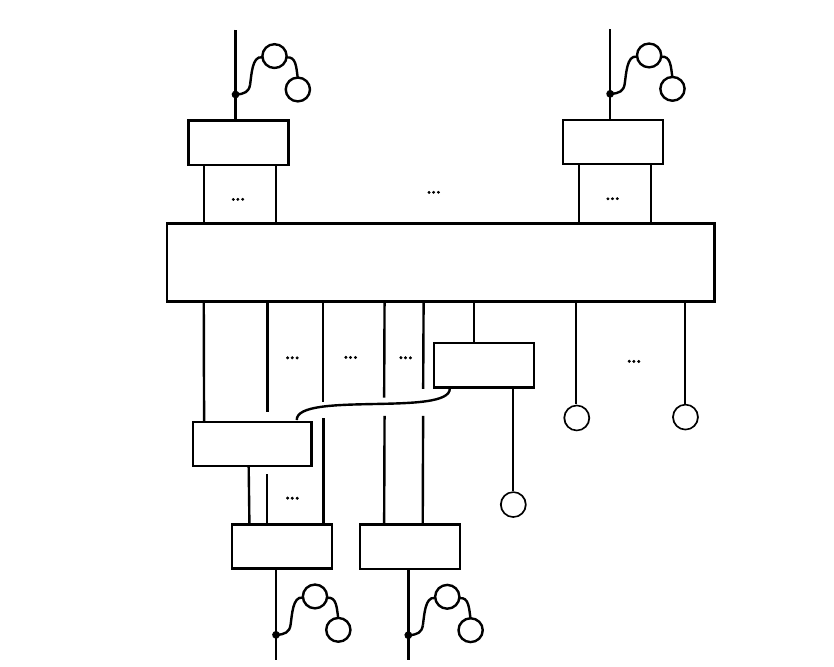}
\end{center}
Since $\alpha$ satisfies \ref{ad5} and using the result of Lemma~\ref{lemmanaturalkirby}:
\vspace{-1.5cm}
\begin{center}
\resizebox{!}{0.55\textheight}{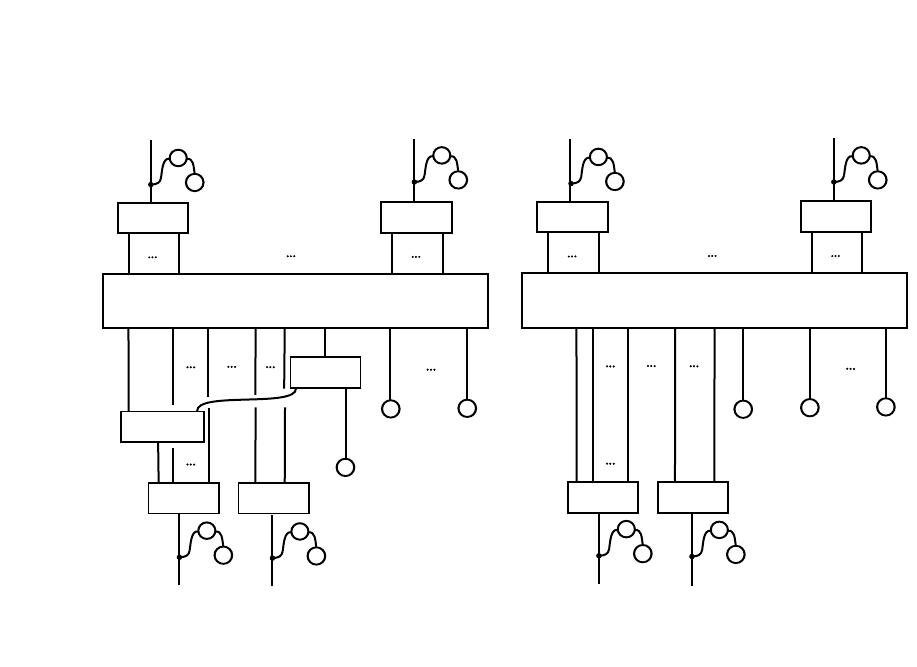}
\end{center}
\vspace{-1cm}

Then $|\halo{T}_1|_{\C,\alpha}=|\halo{T}_2|_{\C,\alpha}$.\\

\noindent\textbullet\;  Secondly, suppose that a surgery component slides over another (or itself) surgery component.
The topological proof is essentially the same than previously. As there is two surgery components, to assure that we have an invariant by the classical move $KII$, this time it is enough to satisfy the axiom
\begin{align}\label{eqkirbysurg}
(\id_{C}\otimes m)(\Delta\otimes \mr{id}_{C})(\alpha\otimes \alpha)=\alpha\otimes \alpha
\end{align}

The latter is a consequence of axiom \ref{ad1} and \ref{ad5}. Indeed, as the we suppose \ref{ad5}, we have 
\[(\id_{C}\otimes m)(\Delta\otimes \mr{id}_{C})(\id_C\otimes \alpha)(\id_C\otimes \omega(\id_C\otimes \alpha))\delta_{C}=[(\id_C\otimes \omega(\id_C\otimes \alpha))\delta_{C}]\otimes \alpha.\]

Composing the last equality by $\alpha$, we obtain

\[(\varepsilon\alpha)(\id_{C}\otimes m)(\Delta\otimes \mr{id}_{C})(\alpha\otimes \alpha)=(\varepsilon\alpha)(\alpha\otimes \alpha)\]

because the transformation $\delta=\{\delta_{X}\colon X\rightarrow X\otimes C\}$ is natural, $u=\delta_{\un}$, ${\omega(u\otimes \id_C)=\varepsilon}$ (properties of the Hopf pairing $\omega$) and $\varepsilon\alpha$ is invertible by \ref{ad1}. Consequently, the identity (\ref{eqkirbysurg}) is true. See \cite{Vir} for more details on this case.

\noindent\textbullet\; Thirdly, suppose that an exit component slides over a surgery component and without loss of generality that the first exit component slides over the first surgery component. Moreover, we can assume that the "sliding part" of the exit component is located at the top of the exit component as it is drawn of the following picture:
\begin{center}
\resizebox{!}{0.15\textheight}{
\begingroup%
  \makeatletter%
  \providecommand\color[2][]{%
    \errmessage{(Inkscape) Color is used for the text in Inkscape, but the package 'color.sty' is not loaded}%
    \renewcommand\color[2][]{}%
  }%
  \providecommand\transparent[1]{%
    \errmessage{(Inkscape) Transparency is used (non-zero) for the text in Inkscape, but the package 'transparent.sty' is not loaded}%
    \renewcommand\transparent[1]{}%
  }%
  \providecommand\rotatebox[2]{#2}%
  \ifx\svgwidth\undefined%
    \setlength{\unitlength}{420bp}%
    \ifx\svgscale\undefined%
      \relax%
    \else%
      \setlength{\unitlength}{\unitlength * \real{\svgscale}}%
    \fi%
  \else%
    \setlength{\unitlength}{\svgwidth}%
  \fi%
  \global\let\svgwidth\undefined%
  \global\let\svgscale\undefined%
  \makeatother%
  \begin{picture}(1,0.28571429)%
    \put(0,0){\includegraphics[width=\unitlength]{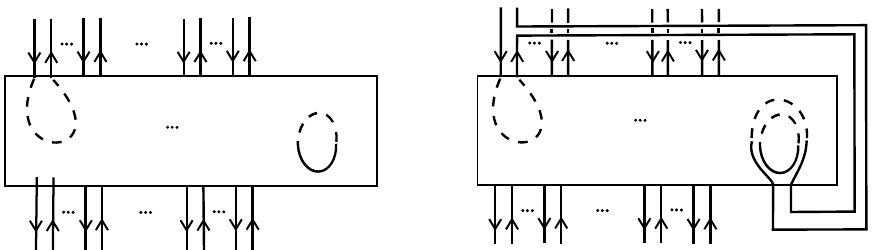}}%
    \put(0.46736314,0.16358964){\color[rgb]{0,0,0}\makebox(0,0)[lt]{\begin{minipage}{0.13383753\unitlength}\raggedright $T_2=$\end{minipage}}}%
    \put(-0.07443455,0.15642545){\color[rgb]{0,0,0}\makebox(0,0)[lt]{\begin{minipage}{0.13383753\unitlength}\raggedright $T_1=$\end{minipage}}}%
  \end{picture}%
\endgroup%
}
\end{center}

We choose two opentangles $P_1$ and $P_2$ associated respectively to $T_1$ and $T_2$ as shown on the picture:
\begin{center}
\resizebox{!}{0.15\textheight}{
\begingroup%
  \makeatletter%
  \providecommand\color[2][]{%
    \errmessage{(Inkscape) Color is used for the text in Inkscape, but the package 'color.sty' is not loaded}%
    \renewcommand\color[2][]{}%
  }%
  \providecommand\transparent[1]{%
    \errmessage{(Inkscape) Transparency is used (non-zero) for the text in Inkscape, but the package 'transparent.sty' is not loaded}%
    \renewcommand\transparent[1]{}%
  }%
  \providecommand\rotatebox[2]{#2}%
  \ifx\svgwidth\undefined%
    \setlength{\unitlength}{440bp}%
    \ifx\svgscale\undefined%
      \relax%
    \else%
      \setlength{\unitlength}{\unitlength * \real{\svgscale}}%
    \fi%
  \else%
    \setlength{\unitlength}{\svgwidth}%
  \fi%
  \global\let\svgwidth\undefined%
  \global\let\svgscale\undefined%
  \makeatother%
  \begin{picture}(1,0.29090909)%
    \put(0,0){\includegraphics[width=\unitlength]{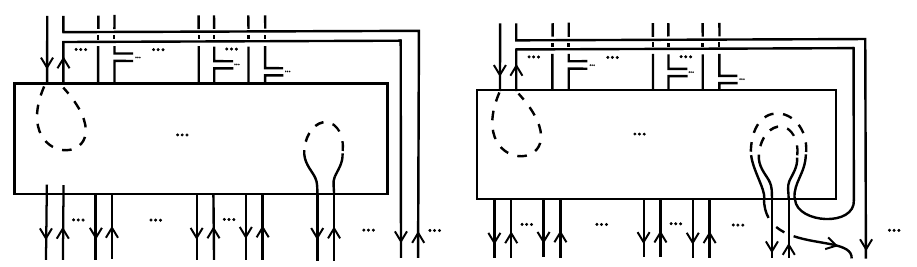}}%
    \put(0.46772918,0.15534755){\color[rgb]{0,0,0}\makebox(0,0)[lt]{\begin{minipage}{0.127754\unitlength}\raggedright $P_2=$\end{minipage}}}%
    \put(-0.04225354,0.14437708){\color[rgb]{0,0,0}\makebox(0,0)[lt]{\begin{minipage}{0.127754\unitlength}\raggedright $P_1=$\end{minipage}}}%
  \end{picture}%
\endgroup%
}
\end{center}

The sequel of the reasoning is the same as in the first case. As $\iota$ is dinatural, note that multiplication $m$ of the coend could be defined graphically by the two following forms:
\begin{center}
\resizebox{!}{0.12\textheight}{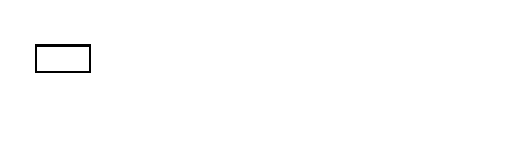}
\end{center}
We use the second form in this case and to conclude, we observe that we only need to satisfy the axiom:
\begin{align}\label{eqkirbysurg2}
(m\otimes \id_{C})(\mr{id}_{C}\otimes \Delta)(\alpha\otimes \alpha)=\alpha\otimes \alpha
\end{align}

Note that in the case where \ref{ad2} is satisfied, the two identities (\ref{eqkirbysurg}) and (\ref{eqkirbysurg2}) are equivalent as it is proved in \cite{Vir}. As we have proved in the second case that (\ref{eqkirbysurg}) is true, equality (\ref{eqkirbysurg2}) is satisfied and the third and last case are then proved.
\end{proof}

Denote by $b_{+}(T)$ (respectively $b_{-}(T)$)\index{$b_{+}(T)$}\index{$b_{-}(T)$} the number of positive (respectively negative) eigenvalues of the linking matrix of the link composed by all the surgery components of $T$ and set
\begin{align}\label{normalization}
\nu_{\alpha}(T)=(\theta_{+}\alpha)^{-b_{+}(T)}(\theta_{-}\alpha)^{-b_{-}(T)}\in \mathrm{End}_{\C}(\un).\end{align}
Obviously, $\nu_{\alpha}$ is multiplicative for disjoint union of tangles:
\[\nu_{\alpha}(T\sqcup T')=\nu_{\alpha}(T)\nu_{\alpha}(T').\] 
Then we are able to define a 3-cobordim invariant.

\begin{lemma}\label{lemmainv}
Let $\alpha\in \mathrm{Hom}_{\C}(\un,C)$ satisfying \ref{ad1}, \ref{ad2}, \ref{ad3}, \ref{ad5} and let $M_{T}$ be a connected cobordism represented by a cobordism tangle $T$. 
Then \[\W_{\C,\alpha}(M_T)=\nu_{\alpha}(T)|\halo{T}|_{\C,\alpha}\] is a topological $3$-cobordism invariant.
\end{lemma}

\begin{example} Construction of the cobordism homeomorphism invariant based on a tangle presentation T of $\Sigma_1\times [0,1]$:
\begin{center}
\resizebox{!}{0.20\textheight}{\raisebox{10mm}{
\begingroup%
  \makeatletter%
  \providecommand\color[2][]{%
    \errmessage{(Inkscape) Color is used for the text in Inkscape, but the package 'color.sty' is not loaded}%
    \renewcommand\color[2][]{}%
  }%
  \providecommand\transparent[1]{%
    \errmessage{(Inkscape) Transparency is used (non-zero) for the text in Inkscape, but the package 'transparent.sty' is not loaded}%
    \renewcommand\transparent[1]{}%
  }%
  \providecommand\rotatebox[2]{#2}%
  \newcommand*\fsize{\dimexpr\f@size pt\relax}%
  \newcommand*\lineheight[1]{\fontsize{\fsize}{#1\fsize}\selectfont}%
  \ifx\svgwidth\undefined%
    \setlength{\unitlength}{202.5bp}%
    \ifx\svgscale\undefined%
      \relax%
    \else%
      \setlength{\unitlength}{\unitlength * \real{\svgscale}}%
    \fi%
  \else%
    \setlength{\unitlength}{\svgwidth}%
  \fi%
  \global\let\svgwidth\undefined%
  \global\let\svgscale\undefined%
  \makeatother%
  \begin{picture}(1,0.68518519)%
    \lineheight{1}%
    \setlength\tabcolsep{0pt}%
    \put(0,0){\includegraphics[width=\unitlength,page=1]{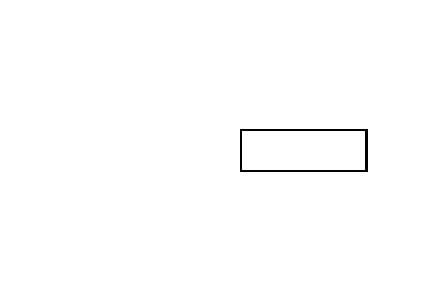}}%
    \put(0.68482568,0.35200171){\color[rgb]{0,0,0}\makebox(0,0)[lt]{\begin{minipage}{0.2192396\unitlength}\raggedright $|T|_{\C,\alpha}$\end{minipage}}}%
    \put(0.72472079,0.01436262){\color[rgb]{0,0,0}\makebox(0,0)[lt]{\lineheight{0}\smash{\begin{tabular}[t]{l}\small{$C$}\end{tabular}}}}%
    \put(-0.05364174,0.35138697){\color[rgb]{0,0,0}\makebox(0,0)[lt]{\begin{minipage}{1.21491301\unitlength}\raggedright $W_{\C,\alpha}(\Sigma_1\times[0,1]_T)= \nu_{\alpha}(T)$\end{minipage}}}%
    \put(0,0){\includegraphics[width=\unitlength,page=2]{construction778.pdf}}%
    \put(0.84806524,0.2270267){\color[rgb]{0,0,0}\makebox(0,0)[lt]{\lineheight{0}\smash{\begin{tabular}[t]{l}$\omega$\end{tabular}}}}%
    \put(0,0){\includegraphics[width=\unitlength,page=3]{construction778.pdf}}%
    \put(0.91438032,0.11866649){\color[rgb]{0,0,0}\makebox(0,0)[lt]{\lineheight{0}\smash{\begin{tabular}[t]{l}$\alpha$\end{tabular}}}}%
    \put(0,0){\includegraphics[width=\unitlength,page=4]{construction778.pdf}}%
    \put(0.73310823,0.65581438){\color[rgb]{0,0,0}\makebox(0,0)[lt]{\lineheight{0}\smash{\begin{tabular}[t]{l}\small{$C$}\end{tabular}}}}%
    \put(0,0){\includegraphics[width=\unitlength,page=5]{construction778.pdf}}%
    \put(0.84966637,0.59558581){\color[rgb]{0,0,0}\makebox(0,0)[lt]{\lineheight{0}\smash{\begin{tabular}[t]{l}$\omega$\end{tabular}}}}%
    \put(0,0){\includegraphics[width=\unitlength,page=6]{construction778.pdf}}%
    \put(0.91598146,0.48722495){\color[rgb]{0,0,0}\makebox(0,0)[lt]{\lineheight{0}\smash{\begin{tabular}[t]{l}$\alpha$\end{tabular}}}}%
    \put(0,0){\includegraphics[width=\unitlength,page=7]{construction778.pdf}}%
  \end{picture}%
\endgroup%
}}
\end{center}
where $\nu_{\alpha}(T)=1$ (see crossings on Figure~\ref{cylinder}).
\end{example}

\begin{proof}
Let $M$ be a connected $3$-cobordism of $\Cobp(\ung,\unh)$ and let $\left[T_{1}\right]$ and $\left[T_{2}\right]$ be two cobordism tangles of $\bigsqcup\limits_{n\in \mathbb{N}}^{}Tang^{Cob}(\ung,n,\unh)$ which represent the cobordism $M$. That is  $N(\left[T_{1}\right])=M=N(\left[T_{2}\right])$ where the map $N$ is defined in equality~(\ref{surjectivemapn}). Our goal is to prove that $\W_{\C}(M_{T_1};\alpha)=\W_{\C}(M_{T_2};\alpha)$.\\

As indicated in Subsection~\ref{parametrized}, $\left[T_1\right]$ and $\left[T_2\right]$ are equivalent if and only if a diagram of $\left[T_1\right]$ and a diagram of $\left[T_2\right]$ are related by a planar isotopy with ribbon Reidemeister moves and a finite sequence of moves of type $SO$, $KI$, $KII^{g}$, $COUPON$, and $TWIST$. In the sequel, we denote indifferently by $T_1$ and $T_2$ the tangles and their diagrams. We can suppose without loss of generality that $T_1$ and $T_2$ are only isotopic or only differ by only one of the four moves $SO$, $KI$, $KII^{g}$, $COUPON$, and $TWIST$.\\

\noindent \textbullet \;\;  If $T_1$ and $T_2$ are isotopic in sphere $\Sp^{3}$, then $\nu_{\alpha}(T_1)=\nu_{\alpha}(T_2)$ cause it is well-known that the linking number is an isotopy invariant and so it is for the linking matrix. Moreover, if $T_1$ and $T_2$ are isotopic then $\halo{T}_1$ and $\halo{T}_2$ are obviously isotopic by construction. We have already shown that $|\phantom{o}|_{\C,\alpha}$ is an isotopy invariant (see Lemma \ref{lemmaisotopyinv}) and as a consequence, $\W_{\C}(M_{T_1};\alpha)=\W_{\C}(M_{T_2};\alpha)$.\\

\noindent \textbullet \;\; If one surgery component of $T_1$ and $T_{2}$ differs by its orientation:\\
In order to simplify notations on the proof, we assume that the surgery component is isolated from others components of the tangle. The general case is a direct rewriting of this particular case. Suppose that $T_1=T\sqcup L$ and $T_2=T\sqcup \overline{L}$ where $T$ is a tangle of cobordism and $L$ and $\overline{L}$ are respectively the surgery components of $T_1$ and $T_2$ that differ by their orientation. Cause $L$ and $\overline{L}$ are links, $\halo{T}_1=\halo{T}\sqcup L$ and $\halo{T}_2=\halo{T}\sqcup \overline{L}$.
As $|\halo{T}_1|_{\C,\alpha}=|\halo{T}\sqcup L|_{\C,\alpha}=|\halo{T}|_{\C,\alpha}\otimes |L|_{\C,\alpha}$ and $|\halo{T}_2|_{\C,\alpha}=|\halo{T}\sqcup \overline{L}|_{\C,\alpha}=|\halo{T}|_{\C,\alpha}\otimes |\overline{L}|_{\C,\alpha}$ by Lemma \ref{lemmamultiplicative}, it is sufficient to prove that $|L|_{\C,\alpha}=|\overline{L}|_{\C,\alpha}$.

First, compute $|L|_{\C,\alpha}=\Bigg|\raisebox{-6mm}{
\begingroup%
  \makeatletter%
  \providecommand\color[2][]{%
    \errmessage{(Inkscape) Color is used for the text in Inkscape, but the package 'color.sty' is not loaded}%
    \renewcommand\color[2][]{}%
  }%
  \providecommand\transparent[1]{%
    \errmessage{(Inkscape) Transparency is used (non-zero) for the text in Inkscape, but the package 'transparent.sty' is not loaded}%
    \renewcommand\transparent[1]{}%
  }%
  \providecommand\rotatebox[2]{#2}%
  \ifx\svgwidth\undefined%
    \setlength{\unitlength}{40bp}%
    \ifx\svgscale\undefined%
      \relax%
    \else%
      \setlength{\unitlength}{\unitlength * \real{\svgscale}}%
    \fi%
  \else%
    \setlength{\unitlength}{\svgwidth}%
  \fi%
  \global\let\svgwidth\undefined%
  \global\let\svgscale\undefined%
  \makeatother%
  \begin{picture}(1,1)%
    \put(0,0){\includegraphics[width=\unitlength]{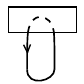}}%
  \end{picture}%
\endgroup%
}\Bigg|_{\C,\alpha}$.
By the factorization property of the coend $C$ of $\C$, there exists a unique morphism $\phi \colon C \rightarrow \mathds{1}$ such that for all objects $X$ of $\C$, 
\begin{center}
\resizebox{!}{0.08\textheight}{
\begingroup%
  \makeatletter%
  \providecommand\color[2][]{%
    \errmessage{(Inkscape) Color is used for the text in Inkscape, but the package 'color.sty' is not loaded}%
    \renewcommand\color[2][]{}%
  }%
  \providecommand\transparent[1]{%
    \errmessage{(Inkscape) Transparency is used (non-zero) for the text in Inkscape, but the package 'transparent.sty' is not loaded}%
    \renewcommand\transparent[1]{}%
  }%
  \providecommand\rotatebox[2]{#2}%
  \ifx\svgwidth\undefined%
    \setlength{\unitlength}{216bp}%
    \ifx\svgscale\undefined%
      \relax%
    \else%
      \setlength{\unitlength}{\unitlength * \real{\svgscale}}%
    \fi%
  \else%
    \setlength{\unitlength}{\svgwidth}%
  \fi%
  \global\let\svgwidth\undefined%
  \global\let\svgscale\undefined%
  \makeatother%
  \begin{picture}(1,0.37037037)%
    \put(0,0){\includegraphics[width=\unitlength]{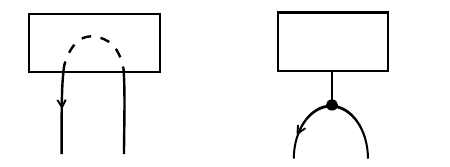}}%
    \put(0.71983035,0.30101477){\color[rgb]{0,0,0}\makebox(0,0)[lt]{\begin{minipage}{0.25931042\unitlength}\raggedright \Large{$\phi$}\end{minipage}}}%
    \put(0.47398377,0.24443384){\color[rgb]{0,0,0}\makebox(0,0)[lt]{\begin{minipage}{0.0659557\unitlength}\raggedright $=$\end{minipage}}}%
    \put(0.29308029,0.05383934){\color[rgb]{0,0,0}\makebox(0,0)[lt]{\begin{minipage}{0.08046049\unitlength}\raggedright $X$\end{minipage}}}%
    \put(0.84492569,0.04803508){\color[rgb]{0,0,0}\makebox(0,0)[lt]{\begin{minipage}{0.08046049\unitlength}\raggedright $X$\end{minipage}}}%
  \end{picture}%
\endgroup%
}\label{lemmainvproofantipode1}
\end{center}
and $|L|_{\C,\alpha}=\phi\alpha$.

Then, compute $|\overline{L}|_{\C,\alpha}=\Bigg|\raisebox{-6mm}{
\begingroup%
  \makeatletter%
  \providecommand\color[2][]{%
    \errmessage{(Inkscape) Color is used for the text in Inkscape, but the package 'color.sty' is not loaded}%
    \renewcommand\color[2][]{}%
  }%
  \providecommand\transparent[1]{%
    \errmessage{(Inkscape) Transparency is used (non-zero) for the text in Inkscape, but the package 'transparent.sty' is not loaded}%
    \renewcommand\transparent[1]{}%
  }%
  \providecommand\rotatebox[2]{#2}%
  \ifx\svgwidth\undefined%
    \setlength{\unitlength}{40bp}%
    \ifx\svgscale\undefined%
      \relax%
    \else%
      \setlength{\unitlength}{\unitlength * \real{\svgscale}}%
    \fi%
  \else%
    \setlength{\unitlength}{\svgwidth}%
  \fi%
  \global\let\svgwidth\undefined%
  \global\let\svgscale\undefined%
  \makeatother%
  \begin{picture}(1,1)%
    \put(0,0){\includegraphics[width=\unitlength]{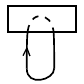}}%
  \end{picture}%
\endgroup%
}\Bigg|_{\C,\alpha}$. Since 
\raisebox{-6mm}{
\begingroup%
  \makeatletter%
  \providecommand\color[2][]{%
    \errmessage{(Inkscape) Color is used for the text in Inkscape, but the package 'color.sty' is not loaded}%
    \renewcommand\color[2][]{}%
  }%
  \providecommand\transparent[1]{%
    \errmessage{(Inkscape) Transparency is used (non-zero) for the text in Inkscape, but the package 'transparent.sty' is not loaded}%
    \renewcommand\transparent[1]{}%
  }%
  \providecommand\rotatebox[2]{#2}%
  \ifx\svgwidth\undefined%
    \setlength{\unitlength}{40bp}%
    \ifx\svgscale\undefined%
      \relax%
    \else%
      \setlength{\unitlength}{\unitlength * \real{\svgscale}}%
    \fi%
  \else%
    \setlength{\unitlength}{\svgwidth}%
  \fi%
  \global\let\svgwidth\undefined%
  \global\let\svgscale\undefined%
  \makeatother%
  \begin{picture}(1,1)%
    \put(0,0){\includegraphics[width=\unitlength]{lemmainvproofantipode3.pdf}}%
  \end{picture}%
\endgroup%
}=\raisebox{-6mm}{
\begingroup%
  \makeatletter%
  \providecommand\color[2][]{%
    \errmessage{(Inkscape) Color is used for the text in Inkscape, but the package 'color.sty' is not loaded}%
    \renewcommand\color[2][]{}%
  }%
  \providecommand\transparent[1]{%
    \errmessage{(Inkscape) Transparency is used (non-zero) for the text in Inkscape, but the package 'transparent.sty' is not loaded}%
    \renewcommand\transparent[1]{}%
  }%
  \providecommand\rotatebox[2]{#2}%
  \ifx\svgwidth\undefined%
    \setlength{\unitlength}{40bp}%
    \ifx\svgscale\undefined%
      \relax%
    \else%
      \setlength{\unitlength}{\unitlength * \real{\svgscale}}%
    \fi%
  \else%
    \setlength{\unitlength}{\svgwidth}%
  \fi%
  \global\let\svgwidth\undefined%
  \global\let\svgscale\undefined%
  \makeatother%
  \begin{picture}(1,1)%
    \put(0,0){\includegraphics[width=\unitlength]{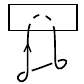}}%
  \end{picture}%
\endgroup%
}, compute $|\overline{L}|_{\C,\alpha}$ using this second diagram of $\overline{L}$. For all objects $X$ of $\C$,
\begin{center}
\resizebox{!}{0.10\textheight}{
\begingroup%
  \makeatletter%
  \providecommand\color[2][]{%
    \errmessage{(Inkscape) Color is used for the text in Inkscape, but the package 'color.sty' is not loaded}%
    \renewcommand\color[2][]{}%
  }%
  \providecommand\transparent[1]{%
    \errmessage{(Inkscape) Transparency is used (non-zero) for the text in Inkscape, but the package 'transparent.sty' is not loaded}%
    \renewcommand\transparent[1]{}%
  }%
  \providecommand\rotatebox[2]{#2}%
  \ifx\svgwidth\undefined%
    \setlength{\unitlength}{320bp}%
    \ifx\svgscale\undefined%
      \relax%
    \else%
      \setlength{\unitlength}{\unitlength * \real{\svgscale}}%
    \fi%
  \else%
    \setlength{\unitlength}{\svgwidth}%
  \fi%
  \global\let\svgwidth\undefined%
  \global\let\svgscale\undefined%
  \makeatother%
  \begin{picture}(1,0.3)%
    \put(0,0){\includegraphics[width=\unitlength]{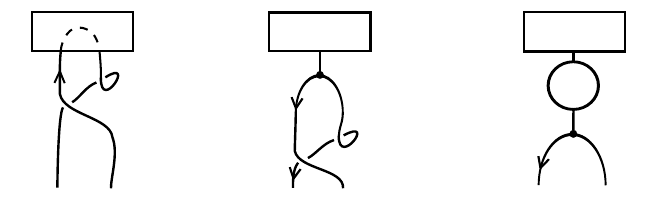}}%
    \put(0.28016405,0.16008309){\color[rgb]{0,0,0}\makebox(0,0)[lt]{\begin{minipage}{0.03689916\unitlength}\raggedright $=$\end{minipage}}}%
    \put(0.17656785,0.03145015){\color[rgb]{0,0,0}\makebox(0,0)[lt]{\begin{minipage}{0.0344391\unitlength}\raggedright $X$\end{minipage}}}%
    \put(0.67984579,0.16897512){\color[rgb]{0,0,0}\makebox(0,0)[lt]{\begin{minipage}{0.03689916\unitlength}\raggedright $=$\end{minipage}}}%
    \put(0.52301901,0.0312125){\color[rgb]{0,0,0}\makebox(0,0)[lt]{\begin{minipage}{0.0344391\unitlength}\raggedright $X$\end{minipage}}}%
    \put(0.91605493,0.0361838){\color[rgb]{0,0,0}\makebox(0,0)[lt]{\begin{minipage}{0.0344391\unitlength}\raggedright $X$\end{minipage}}}%
    \put(0.46638827,0.27195698){\color[rgb]{0,0,0}\makebox(0,0)[lt]{\begin{minipage}{0.13345611\unitlength}\raggedright \Large{$\phi$}\end{minipage}}}%
    \put(0.84982452,0.27169649){\color[rgb]{0,0,0}\makebox(0,0)[lt]{\begin{minipage}{0.13345611\unitlength}\raggedright \Large{$\phi$}\end{minipage}}}%
    \put(0.52018681,0.16040505){\color[rgb]{0,0,0}\makebox(0,0)[lt]{\begin{minipage}{0.05689169\unitlength}\raggedright $X^{*}$\end{minipage}}}%
    \put(0.24516036,0.13316215){\color[rgb]{0,0,0}\makebox(0,0)[lt]{\begin{minipage}{0.1416911\unitlength}\raggedright $X=X^{**}$\end{minipage}}}%
    \put(0.60098062,0.14760209){\color[rgb]{0,0,0}\makebox(0,0)[lt]{\begin{minipage}{0.21330658\unitlength}\raggedright $\text{by definition of }$\end{minipage}}}%
    \put(0.83904312,0.18048){\color[rgb]{0,0,0}\makebox(0,0)[lt]{\begin{minipage}{0.0796377\unitlength}\raggedright \small{$S^{-1}$}\end{minipage}}}%
    \put(0.67943541,0.11843498){\color[rgb]{0,0,0}\makebox(0,0)[lt]{\begin{minipage}{0.03557525\unitlength}\raggedright $S^{-1}$\end{minipage}}}%
  \end{picture}%
\endgroup%
}\label{lemmainvproofantipode5.pdf_tex}
\end{center}

and by the factorization property of the coend, we have $|\overline{L}|_{\C,\alpha}=\phi S^{-1}\alpha$.
Thus, as $S\alpha=\alpha$ (the morphism $\alpha$ satisfies \ref{ad2}), we have $S^{-1}\alpha=\alpha$ too so $|\overline{L}|_{\C,\alpha}=\phi\alpha$ and we conclude that $|L|_{\C,\alpha}=|\overline{L}|_{\C,\alpha}$. The number of positive (respectively negative) eigenvalues is invariant by changing the orientation of a link component: indeed, changing the orientation of one component goes back to change the sign of exactly the $k$th-line and the $k$th-column for a certain integer $k$ of the linking matrix and this new matrix is similar to the first one. Then we get that $\W_{\C}(M_{T_1};\alpha)=\W_{\C}(M_{T_2};\alpha)$.\\

\noindent \textbullet \indent If $T_1=T_{2}$ $\sqcup$ \raisebox{-1.5mm}{
\begingroup%
  \makeatletter%
  \providecommand\color[2][]{%
    \errmessage{(Inkscape) Color is used for the text in Inkscape, but the package 'color.sty' is not loaded}%
    \renewcommand\color[2][]{}%
  }%
  \providecommand\transparent[1]{%
    \errmessage{(Inkscape) Transparency is used (non-zero) for the text in Inkscape, but the package 'transparent.sty' is not loaded}%
    \renewcommand\transparent[1]{}%
  }%
  \providecommand\rotatebox[2]{#2}%
  \ifx\svgwidth\undefined%
    \setlength{\unitlength}{16bp}%
    \ifx\svgscale\undefined%
      \relax%
    \else%
      \setlength{\unitlength}{\unitlength * \real{\svgscale}}%
    \fi%
  \else%
    \setlength{\unitlength}{\svgwidth}%
  \fi%
  \global\let\svgwidth\undefined%
  \global\let\svgscale\undefined%
  \makeatother%
  \begin{picture}(1,1)%
    \put(0,0){\includegraphics[width=\unitlength]{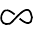}}%
  \end{picture}%
\endgroup%
},
we have $|\halo{T}_1|_{\C,\alpha}=|\halo{T}_2|_{\C,\alpha}\otimes |\raisebox{-1.5mm}{
\begingroup%
  \makeatletter%
  \providecommand\color[2][]{%
    \errmessage{(Inkscape) Color is used for the text in Inkscape, but the package 'color.sty' is not loaded}%
    \renewcommand\color[2][]{}%
  }%
  \providecommand\transparent[1]{%
    \errmessage{(Inkscape) Transparency is used (non-zero) for the text in Inkscape, but the package 'transparent.sty' is not loaded}%
    \renewcommand\transparent[1]{}%
  }%
  \providecommand\rotatebox[2]{#2}%
  \ifx\svgwidth\undefined%
    \setlength{\unitlength}{16bp}%
    \ifx\svgscale\undefined%
      \relax%
    \else%
      \setlength{\unitlength}{\unitlength * \real{\svgscale}}%
    \fi%
  \else%
    \setlength{\unitlength}{\svgwidth}%
  \fi%
  \global\let\svgwidth\undefined%
  \global\let\svgscale\undefined%
  \makeatother%
  \begin{picture}(1,1)%
    \put(0,0){\includegraphics[width=\unitlength]{thetaplus.pdf}}%
  \end{picture}%
\endgroup%
}|_{\C,\alpha}$ by Lemma \ref{lemmamultiplicative}. Compute $|\raisebox{-1.5mm}{
\begingroup%
  \makeatletter%
  \providecommand\color[2][]{%
    \errmessage{(Inkscape) Color is used for the text in Inkscape, but the package 'color.sty' is not loaded}%
    \renewcommand\color[2][]{}%
  }%
  \providecommand\transparent[1]{%
    \errmessage{(Inkscape) Transparency is used (non-zero) for the text in Inkscape, but the package 'transparent.sty' is not loaded}%
    \renewcommand\transparent[1]{}%
  }%
  \providecommand\rotatebox[2]{#2}%
  \ifx\svgwidth\undefined%
    \setlength{\unitlength}{16bp}%
    \ifx\svgscale\undefined%
      \relax%
    \else%
      \setlength{\unitlength}{\unitlength * \real{\svgscale}}%
    \fi%
  \else%
    \setlength{\unitlength}{\svgwidth}%
  \fi%
  \global\let\svgwidth\undefined%
  \global\let\svgscale\undefined%
  \makeatother%
  \begin{picture}(1,1)%
    \put(0,0){\includegraphics[width=\unitlength]{thetaplus.pdf}}%
  \end{picture}%
\endgroup%
}|_{\C,\alpha}$: by definition of $\theta_{+}\colon C \rightarrow \mathds{1}$ (see Figure~\ref{structural2}), for all objects $X\in \C$,
\begin{center}
\resizebox{!}{0.06\textheight}{
\begingroup%
  \makeatletter%
  \providecommand\color[2][]{%
    \errmessage{(Inkscape) Color is used for the text in Inkscape, but the package 'color.sty' is not loaded}%
    \renewcommand\color[2][]{}%
  }%
  \providecommand\transparent[1]{%
    \errmessage{(Inkscape) Transparency is used (non-zero) for the text in Inkscape, but the package 'transparent.sty' is not loaded}%
    \renewcommand\transparent[1]{}%
  }%
  \providecommand\rotatebox[2]{#2}%
  \ifx\svgwidth\undefined%
    \setlength{\unitlength}{160bp}%
    \ifx\svgscale\undefined%
      \relax%
    \else%
      \setlength{\unitlength}{\unitlength * \real{\svgscale}}%
    \fi%
  \else%
    \setlength{\unitlength}{\svgwidth}%
  \fi%
  \global\let\svgwidth\undefined%
  \global\let\svgscale\undefined%
  \makeatother%
  \begin{picture}(1,0.45)%
    \put(0,0){\includegraphics[width=\unitlength]{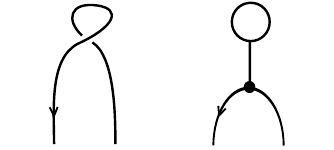}}%
    \put(0.72134049,0.37145159){\color[rgb]{0,0,0}\makebox(0,0)[lb]{\smash{$\theta_{+}$}}}%
    \put(0.4823867,0.26562268){\color[rgb]{0,0,0}\makebox(0,0)[lt]{\begin{minipage}{0.16601929\unitlength}\raggedright $=$\end{minipage}}}%
    \put(0.36899185,0.0585508){\color[rgb]{0,0,0}\makebox(0,0)[lt]{\begin{minipage}{0.10699021\unitlength}\raggedright $X$\end{minipage}}}%
    \put(0.87310834,0.0548626){\color[rgb]{0,0,0}\makebox(0,0)[lt]{\begin{minipage}{0.10699021\unitlength}\raggedright $X$\end{minipage}}}%
  \end{picture}%
\endgroup%
}\label{lemmainvproofnormal1.pdf_tex}
\end{center} 

then $|\raisebox{-1.5mm}{
\begingroup%
  \makeatletter%
  \providecommand\color[2][]{%
    \errmessage{(Inkscape) Color is used for the text in Inkscape, but the package 'color.sty' is not loaded}%
    \renewcommand\color[2][]{}%
  }%
  \providecommand\transparent[1]{%
    \errmessage{(Inkscape) Transparency is used (non-zero) for the text in Inkscape, but the package 'transparent.sty' is not loaded}%
    \renewcommand\transparent[1]{}%
  }%
  \providecommand\rotatebox[2]{#2}%
  \ifx\svgwidth\undefined%
    \setlength{\unitlength}{16bp}%
    \ifx\svgscale\undefined%
      \relax%
    \else%
      \setlength{\unitlength}{\unitlength * \real{\svgscale}}%
    \fi%
  \else%
    \setlength{\unitlength}{\svgwidth}%
  \fi%
  \global\let\svgwidth\undefined%
  \global\let\svgscale\undefined%
  \makeatother%
  \begin{picture}(1,1)%
    \put(0,0){\includegraphics[width=\unitlength]{thetaplus.pdf}}%
  \end{picture}%
\endgroup%
}|_{\C,\alpha}=\theta_{+}\alpha$. Since $\nu_{\alpha}(T_1)=\nu_{\alpha}(T_2\sqcup\raisebox{-1.5mm}{
\begingroup%
  \makeatletter%
  \providecommand\color[2][]{%
    \errmessage{(Inkscape) Color is used for the text in Inkscape, but the package 'color.sty' is not loaded}%
    \renewcommand\color[2][]{}%
  }%
  \providecommand\transparent[1]{%
    \errmessage{(Inkscape) Transparency is used (non-zero) for the text in Inkscape, but the package 'transparent.sty' is not loaded}%
    \renewcommand\transparent[1]{}%
  }%
  \providecommand\rotatebox[2]{#2}%
  \ifx\svgwidth\undefined%
    \setlength{\unitlength}{16bp}%
    \ifx\svgscale\undefined%
      \relax%
    \else%
      \setlength{\unitlength}{\unitlength * \real{\svgscale}}%
    \fi%
  \else%
    \setlength{\unitlength}{\svgwidth}%
  \fi%
  \global\let\svgwidth\undefined%
  \global\let\svgscale\undefined%
  \makeatother%
  \begin{picture}(1,1)%
    \put(0,0){\includegraphics[width=\unitlength]{thetaplus.pdf}}%
  \end{picture}%
\endgroup%
})=\nu_{\alpha}(T_2)(\theta_{+}\alpha)^{-1}$, 

\[\nu_{\alpha}(T_1)|\halo{T}_1|_{\C,\alpha}=\nu_{\alpha}(T_2)(\theta_{+}\alpha)^{-1}\left(|\halo{T}_2|_{\C,\alpha}\otimes\theta_+\alpha \right)=\nu_{\alpha}(T_2)|\halo{T}_2|_{\C,\alpha}\]

that is $\W_{\C}(M_{T_1};\alpha)=\W_{\C}(M_{T_2};\alpha)$.
The case where $T_1=T_{2}$ $\sqcup$ \raisebox{-1.5mm}{
\begingroup%
  \makeatletter%
  \providecommand\color[2][]{%
    \errmessage{(Inkscape) Color is used for the text in Inkscape, but the package 'color.sty' is not loaded}%
    \renewcommand\color[2][]{}%
  }%
  \providecommand\transparent[1]{%
    \errmessage{(Inkscape) Transparency is used (non-zero) for the text in Inkscape, but the package 'transparent.sty' is not loaded}%
    \renewcommand\transparent[1]{}%
  }%
  \providecommand\rotatebox[2]{#2}%
  \ifx\svgwidth\undefined%
    \setlength{\unitlength}{16bp}%
    \ifx\svgscale\undefined%
      \relax%
    \else%
      \setlength{\unitlength}{\unitlength * \real{\svgscale}}%
    \fi%
  \else%
    \setlength{\unitlength}{\svgwidth}%
  \fi%
  \global\let\svgwidth\undefined%
  \global\let\svgscale\undefined%
  \makeatother%
  \begin{picture}(1,1)%
    \put(0,0){\includegraphics[width=\unitlength]{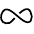}}%
  \end{picture}%
\endgroup%
} is similar.\\ 

\noindent\textbullet \indent If $T_1$ and $T_{2}$ differ by one move $KII^{g}$, since $\alpha$ satisfies conditions \ref{ad1} et \ref{ad5}, the result is the consequence of Lemma \ref{lemmainvariantkirby2} and the fact that $b_{\pm}(T)$ is invariant by classical Kirby move II (handle sliding).\\

\noindent \textbullet \indent If $T_1$ and $T_2$ differ by one move $COUPON$, without loss of generality, suppose that the move $COUPON$ concerns the first boundary component of $T_1$ and $T_2$ as shown on the left part on the following diagram where $T_1$ and $T_2$ are $(\ung,n,\unh)$-cobordism tangles and $\ung=(g_1,\ldots,g_r)$. Then $\halo{T}_1$ and $\halo{T}_2$ are given in the right part of the following diagram.

\begin{center}
    \resizebox{!}{0.12\textheight}{
\begingroup%
  \makeatletter%
  \providecommand\color[2][]{%
    \errmessage{(Inkscape) Color is used for the text in Inkscape, but the package 'color.sty' is not loaded}%
    \renewcommand\color[2][]{}%
  }%
  \providecommand\transparent[1]{%
    \errmessage{(Inkscape) Transparency is used (non-zero) for the text in Inkscape, but the package 'transparent.sty' is not loaded}%
    \renewcommand\transparent[1]{}%
  }%
  \providecommand\rotatebox[2]{#2}%
  \newcommand*\fsize{\dimexpr\f@size pt\relax}%
  \newcommand*\lineheight[1]{\fontsize{\fsize}{#1\fsize}\selectfont}%
  \ifx\svgwidth\undefined%
    \setlength{\unitlength}{651.96850394bp}%
    \ifx\svgscale\undefined%
      \relax%
    \else%
      \setlength{\unitlength}{\unitlength * \real{\svgscale}}%
    \fi%
  \else%
    \setlength{\unitlength}{\svgwidth}%
  \fi%
  \global\let\svgwidth\undefined%
  \global\let\svgscale\undefined%
  \makeatother%
  \begin{picture}(1,0.17391304)%
    \lineheight{1}%
    \setlength\tabcolsep{0pt}%
    \put(0,0){\includegraphics[width=\unitlength,page=1]{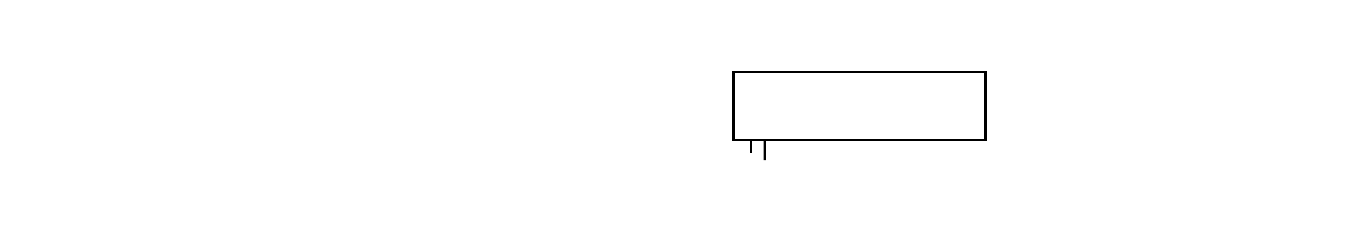}}%
    \put(0.49004863,0.10394502){\color[rgb]{0,0,0}\makebox(0,0)[lt]{\begin{minipage}{0.13484639\unitlength}\raggedright $\halo{T}_1=$\end{minipage}}}%
    \put(0,0){\includegraphics[width=\unitlength,page=2]{test23.pdf}}%
    \put(0.57019053,0.03903261){\color[rgb]{0,0,0}\makebox(0,0)[lt]{\begin{minipage}{0.02292813\unitlength}\raggedright $L$\end{minipage}}}%
    \put(0.61579509,0.06723578){\color[rgb]{0,0,0}\makebox(0,0)[lt]{\begin{minipage}{0.02292813\unitlength}\raggedright $A$\end{minipage}}}%
    \put(0,0){\includegraphics[width=\unitlength,page=3]{test23.pdf}}%
    \put(0.73439081,0.10667293){\color[rgb]{0,0,0}\makebox(0,0)[lt]{\begin{minipage}{0.13484639\unitlength}\raggedright $\halo{T}_2=$\end{minipage}}}%
    \put(0,0){\includegraphics[width=\unitlength,page=4]{test23.pdf}}%
    \put(0.8145327,0.0417605){\color[rgb]{0,0,0}\makebox(0,0)[lt]{\begin{minipage}{0.02292813\unitlength}\raggedright $L$\end{minipage}}}%
    \put(0.86013722,0.06996367){\color[rgb]{0,0,0}\makebox(0,0)[lt]{\begin{minipage}{0.02292813\unitlength}\raggedright $A$\end{minipage}}}%
    \put(0,0){\includegraphics[width=\unitlength,page=5]{test23.pdf}}%
    \put(0.44714334,0.08305589){\color[rgb]{0,0,0}\makebox(0,0)[lt]{\lineheight{1.25}\smash{\begin{tabular}[t]{l}$\text{and}$\end{tabular}}}}%
    \put(0,0){\includegraphics[width=\unitlength,page=6]{test23.pdf}}%
    \put(-0.03993243,0.098959){\color[rgb]{0,0,0}\makebox(0,0)[lt]{\begin{minipage}{0.08195193\unitlength}\raggedright $T_1=$\end{minipage}}}%
    \put(0,0){\includegraphics[width=\unitlength,page=7]{test23.pdf}}%
    \put(0.20366726,0.09506589){\color[rgb]{0,0,0}\makebox(0,0)[lt]{\begin{minipage}{0.08195193\unitlength}\raggedright $T_2=$\end{minipage}}}%
    \put(0,0){\includegraphics[width=\unitlength,page=8]{test23.pdf}}%
  \end{picture}%
\endgroup%
}  
\end{center}

Just remark that you can transform $\halo{T}_1$ into $\halo{T}_2$ by a move $KII^{g}$ considering that the arc $A$ slides on the link $L$ as indicated on the diagram just above. We have already shown that $ |\phantom{o}|_{\C,\alpha}$ is invariant by move $KII^{g}$ when $\alpha$ verifies \ref{ad1} and \ref{ad5} (see Lemma~\ref{lemmainvariantkirby2}) so $|\halo{T}_1|_{\C,\alpha}=|\halo{T}_2|_{\C,\alpha}$ and since the move $COUPON$ doesn't affect the normalization coefficient $\nu_{\alpha}(T)$, we conclude that $\mr{W}_{\C}(M_{T_1};\alpha)=\mr{W}_{\C}(M_{T_2};\alpha)$.\\

\noindent \textbullet \indent If $T_1$ and $T_2$ differ by one move $TWIST$, then, thanks to Fenn and Rourke moves (see \cite{FR}), note that you can eliminate a twist of a boundary component adding an encircling closed twisted component. So add such a surgery component and remark that it could slides on halo of $\halo{T}_1$ or $\halo{T}_2$. Since $\theta_{\pm}\alpha$ is invertible (see \ref{ad3}), we get easily the result.
\end{proof}

\begin{remarks}
\item We extend the invariant $\W_{\C}(\phantom{O};\alpha)$ to non-connected cobordisms. If $M$ is a $3$-cobordism, denote by $M^{\#}$ the $3$-cobordism obtained as the connected sum of connected components of $M$. Then, for the non-connected cobordism $M$, we just set $W_{\C,\alpha}(M):=W_{\C,\alpha}(M^{\#})$.
\end{remarks}

To define a TQFT, we still have to understand what is going on for the composition of cobordisms. Recall that the cobordism tangle which encodes the compositum of two cobordisms is not the compositum of the tangles: we have to add several closed components (see Section~\ref{operationstar}). This leads to the definition of the hallowed morphism $\Pi_{\alpha,n}$ in equality (\ref{hallowedmorphism}). To guarantee the composition is well defined, we need to ask that this morphism is an idempotent in the category $\C$, which is given by the admissibility condition \ref{ad4}. As $\C$ has splitting idempotents then $\Pi_{\alpha,n}$ has a split decomposition,  that is, there are an object $(C^{\otimes n})_{\alpha} \in \C$ and two morphisms ${\pa{n}\colon C^{\otimes n}\rightarrow (C^{\otimes n})_{\alpha}}$, $\qa{n}\colon (C^{\otimes n})_{\alpha}\rightarrow C^{\otimes n}$ such that 
\[\pa{n}\qa{n}=\mr{id}_{(C^{\otimes n})_{\alpha}} \text{ and } \qa{n}\pa{n}=\Pi_{\alpha,n}.\]
For any tuple $\ung=(g_1,\ldots,g_r)$ of integers, set
\[\pia{\ung}=\pia{g_1}\otimes\ldots\otimes \pia{g_r},\; \pa{\ung}=p_{\alpha,g_1}\otimes\ldots\otimes p_{\alpha,g_r} \text{ and } \qa{\ung}=q_{\alpha,g_1}\otimes\ldots\otimes q_{\alpha,g_r}.\] 
Note that in the Reshestikhin-Turaev TQFT, $(C^{\otimes n})_{\alpha}$ plays the role of the transparent part of object $C^{\otimes n}$. 
To prove the main theorem, we will need a result of characterization on transparency of the image of an idempotent. We will use it to show that the internal TQFT takes values in the subcategory of $\C$ of transparent objects.  
\begin{lemma}\label{transparency}
Let $\Pi\colon X\rightarrow X$ be an idempotent of a ribbon category $\C$ with coend.
Then $\mr{Im}(\Pi)$ is transparent if and only if 
\[\Pi(\id_{X}\otimes \omega(\id_{C}\otimes S))\delta_{X}=\Pi\otimes \varepsilon.\]
\end{lemma}
\begin{proof}
Let $(\mr{Im}(\Pi),p,q)$ a decomposition of the idempotent $\Pi$ that means $pq=\id_{\mr{Im}(\Pi)}$ and $qp=\Pi$.
$\mr{Im}(\Pi)$ is transparent if and only if
\begin{align*}
\forall Y\in \C, &(\id_{\mr{Im}(\Pi)}\otimes \mr{coev}_Y)(\tau^{-1}_{\mr{Im}(\Pi),Y}\tau^{-1}_{Y,\mr{Im}(\Pi)}\otimes \id_Y) =\id_{\mr{Im}(\Pi)}\otimes \mr{coev}_Y \\
\Longleftrightarrow \quad & \forall Y\in \C, (\id_{\mr{Im}(\Pi)}\otimes \omega(\id_C\otimes S))\delta_{\mr{Im}(\Pi)}(\id_{\mr{Im}(\Pi)}\otimes \iota_Y) =\id_{\mr{Im(\Pi)}}\otimes \mr{coev}_Y \\
\stackrel{(2)}{\Longleftrightarrow} \quad & (\id_{\mr{Im}(\Pi)}\otimes \omega(\id_C\otimes S))\delta_{\mr{Im}(\Pi)} =\id_{\mr{Im(\Pi)}}\otimes \varepsilon \\
\stackrel{(3)}{\Longleftrightarrow} \quad &   q(\id_{\mr{Im}(\Pi)}\otimes \omega(\id_C\otimes S))\delta_{\mr{Im}(\Pi)}p =\Pi\otimes \varepsilon\\
\stackrel{(4)}{\Longleftrightarrow} \quad & qp(\id_{\mr{Im}(\Pi)}\otimes \omega(\id_C\otimes S))\delta_{X} =\Pi\otimes \varepsilon\\
\Longleftrightarrow \quad&\Pi(\id_{X}\otimes \omega(\id_{C}\otimes S))\delta_{X}=\Pi\otimes \varepsilon
\end{align*}

Equivalence $(2)$ is due to the universal property of the coend, equivalence $(3)$ comes from the identities $pq=\id_{\mr{Im}(\Pi)}$ et $qp=\Pi$ and equivalence $(4)$ is the consequence of naturality of $\delta$.
\end{proof}

\subsection{Proof of the two main theorems}
Finally, using the 3-cobordisms invariant $\W_{\C,\alpha}$, we are able to construct the \emph{internal} TQFT given in Theorem~\ref{maintheorem}. Before proving this Theorem, just give an example of the computation of the internal TQFT on the cylinder of $\Sigma_1$.



\begin{example} Computations of the internal TQFT on the one genus surface cylinder $\Sigma_1\times [0,1]$, using unizalisation (see Definition ~\ref{unitalization}) of the braided monoidal functor $W_{\C,\alpha}$:
\[V_{\C,\alpha}(\Sigma_1)=\mathrm{Im}(\Pi_{\Sigma_{1}}) \text{ and } V_{\C,\alpha}(\Sigma_1\times [0,1])=\mathrm{id}_{\mathrm{Im}(\Pi_{\Sigma_{1}})}\] 

where 

\begin{center}
\raisebox{20mm}{$\Pi_{\Sigma_1}=$}\resizebox{!}{0.21\textheight}{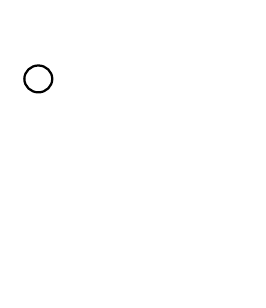}
\end{center}
\end{example}

\begin{proof} We prove Theorem~\ref{maintheorem}.\\ 
Recall that for any $3$-cobordism, we set $\gamma_{M,N}:=\gamma_{M^{\#},N^{\#}}$ where $M^{\#}$ is the connected sum of connected components of $M$.
First, let us check that $(\W_{\C},\gamma)$  is a functor with anomaly.\\\\
\textbullet\; Let us see what's going on objects of $\Cobp$.
Let $(g_1,\ldots,g_r)$ be a $r$-tuple of integers and denote by $\Sigma_{\ung}$ a surface of multigenus $\ung$ that means there exists, for $1\leq i \leq r$, a connected surface $\Sigma_{g_i}$ of genus $g_i$ such that
\[\Sigma_{\ung}=\Sigma_{g_1}\sqcup\ldots\sqcup\Sigma_{g_r}.\]
Note that we will forget parametrizations of surfaces in this proof. For a surface of multigenus $\ung=(g_1,\ldots,g_r)$, set:
\[\W_{\C,\alpha}(\Sigma_{\ung}):=(C^{\otimes g_1})_{\alpha}\otimes \ldots\otimes (C^{\otimes  g_r})_{\alpha}.\]
Thus we have a canonical identification $\W_2(\Sigma_{\ung},\Sigma_{\unh})\colon \W_{\C}(\Sigma_{\ung})\otimes \W_{\C}(\Sigma_{\unh})\iso \W_{\C}(\Sigma_{\unh}\otimes \Sigma_{\unh})$. Moreover, we set $\W_{\C,\alpha}(\emptyset;\alpha)=\un$.


\noindent \textbullet\;
Suppose that $M_{T}\in \Cobp(\ung,\unh)$ and $M_{T'}\in \Cobp(\unh,\unk)$ are two connected $3$-cobordims represented respectively by a $(\ung,n,\unh)$-cobordism tangle $T$ and by a $(\unh,m,\unk)$-cobordism tangle $T'$.
We want to compare $\W_{\C,\alpha}(M_{T'}\circ M_{T})$ and $\W_{\C,\alpha}(M_{T'})\circ \W_{\C,\alpha}(M_{T})$.

According to Turaev (see \cite{Turaev1994}), the cobordim tangle $T'\star T$ (defined in Section~\ref{operationstar}) encodes the compositum of $3$-cobordims $M_{T'}\circ M_{T}$ that means $M_{T'}\circ M_{T}$ is homeomorphic (as $3$-cobordims) to $M_{T'\star T}$. Recall that the tangle $T'\star T$ is defined as illustrated below. We want to show that $|\halo{T'\star T}|_{\C,\alpha}=|\halo{T'}|_{\C,\alpha}\circ |\halo{T}|_{\C,\alpha}$. To be more symmetric, we define a new operation $\square$ on cobordism tangles $T'$ and $T$ as specified below.
\begin{center}
     \resizebox{!}{0.28\textheight}{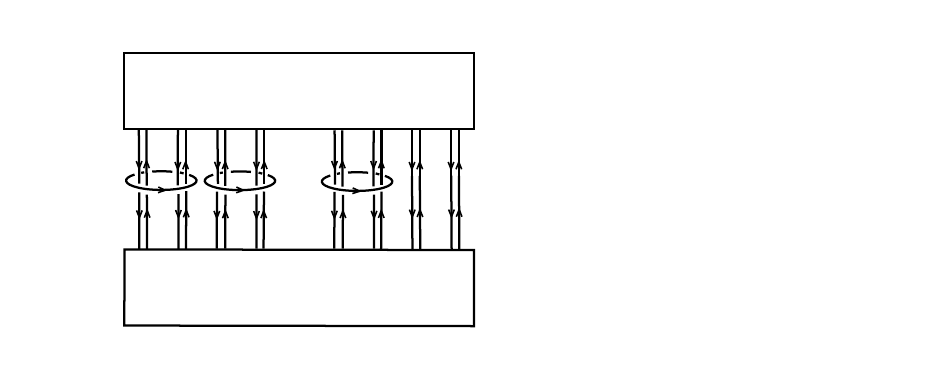}
\end{center}


Note that $|\halo{T'\square T}|_{\C,\alpha}=|\halo{T'\star T}|_{\C,\alpha}$. Indeed, observe that
\begin{center}
\begin{align*}
\raisebox{45mm}{$|\halo{T'\square T}|_{\C,\alpha}=$}&\scalebox{0.75}{
\begingroup%
  \makeatletter%
  \providecommand\color[2][]{%
    \errmessage{(Inkscape) Color is used for the text in Inkscape, but the package 'color.sty' is not loaded}%
    \renewcommand\color[2][]{}%
  }%
  \providecommand\transparent[1]{%
    \errmessage{(Inkscape) Transparency is used (non-zero) for the text in Inkscape, but the package 'transparent.sty' is not loaded}%
    \renewcommand\transparent[1]{}%
  }%
  \providecommand\rotatebox[2]{#2}%
  \ifx\svgwidth\undefined%
    \setlength{\unitlength}{220bp}%
    \ifx\svgscale\undefined%
      \relax%
    \else%
      \setlength{\unitlength}{\unitlength * \real{\svgscale}}%
    \fi%
  \else%
    \setlength{\unitlength}{\svgwidth}%
  \fi%
  \global\let\svgwidth\undefined%
  \global\let\svgscale\undefined%
  \makeatother%
  \begin{picture}(1,1.27272727)%
    \put(0,0){\includegraphics[width=\unitlength]{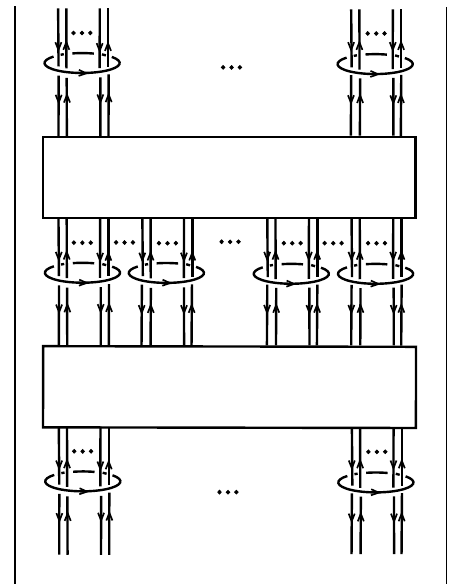}}%
    \put(0.47740203,0.91328986){\color[rgb]{0,0,0}\makebox(0,0)[lt]{\begin{minipage}{0.07638696\unitlength}\raggedright $T'$\end{minipage}}}%
    \put(0.47503629,0.44254782){\color[rgb]{0,0,0}\makebox(0,0)[lt]{\begin{minipage}{0.07638696\unitlength}\raggedright $T$\end{minipage}}}%
    \put(0.98419306,0.0293634){\color[rgb]{0,0,0}\makebox(0,0)[lt]{\begin{minipage}{0.22084041\unitlength}\raggedright $\C,\alpha$\end{minipage}}}%
  \end{picture}%
\endgroup%
}
\raisebox{45mm}{$\stackrel{(2)}{=}$}
\scalebox{0.75}{
\begingroup%
  \makeatletter%
  \providecommand\color[2][]{%
    \errmessage{(Inkscape) Color is used for the text in Inkscape, but the package 'color.sty' is not loaded}%
    \renewcommand\color[2][]{}%
  }%
  \providecommand\transparent[1]{%
    \errmessage{(Inkscape) Transparency is used (non-zero) for the text in Inkscape, but the package 'transparent.sty' is not loaded}%
    \renewcommand\transparent[1]{}%
  }%
  \providecommand\rotatebox[2]{#2}%
  \ifx\svgwidth\undefined%
    \setlength{\unitlength}{220bp}%
    \ifx\svgscale\undefined%
      \relax%
    \else%
      \setlength{\unitlength}{\unitlength * \real{\svgscale}}%
    \fi%
  \else%
    \setlength{\unitlength}{\svgwidth}%
  \fi%
  \global\let\svgwidth\undefined%
  \global\let\svgscale\undefined%
  \makeatother%
  \begin{picture}(1,1.27272727)%
    \put(0,0){\includegraphics[width=\unitlength]{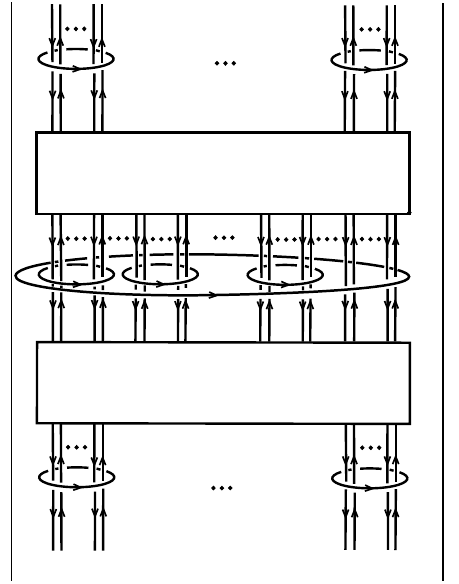}}%
    \put(0.46411013,0.92222582){\color[rgb]{0,0,0}\makebox(0,0)[lt]{\begin{minipage}{0.07638696\unitlength}\raggedright $T'$\end{minipage}}}%
    \put(0.46174439,0.45148379){\color[rgb]{0,0,0}\makebox(0,0)[lt]{\begin{minipage}{0.07638696\unitlength}\raggedright $T$\end{minipage}}}%
    \put(0.97591356,0.03828452){\color[rgb]{0,0,0}\makebox(0,0)[lt]{\begin{minipage}{0.22084041\unitlength}\raggedright $\C,\alpha$\end{minipage}}}%
  \end{picture}%
\endgroup%
}\\
\raisebox{45mm}{$\stackrel{(3)}{=}$}&\scalebox{0.75}{
\begingroup%
  \makeatletter%
  \providecommand\color[2][]{%
    \errmessage{(Inkscape) Color is used for the text in Inkscape, but the package 'color.sty' is not loaded}%
    \renewcommand\color[2][]{}%
  }%
  \providecommand\transparent[1]{%
    \errmessage{(Inkscape) Transparency is used (non-zero) for the text in Inkscape, but the package 'transparent.sty' is not loaded}%
    \renewcommand\transparent[1]{}%
  }%
  \providecommand\rotatebox[2]{#2}%
  \ifx\svgwidth\undefined%
    \setlength{\unitlength}{220bp}%
    \ifx\svgscale\undefined%
      \relax%
    \else%
      \setlength{\unitlength}{\unitlength * \real{\svgscale}}%
    \fi%
  \else%
    \setlength{\unitlength}{\svgwidth}%
  \fi%
  \global\let\svgwidth\undefined%
  \global\let\svgscale\undefined%
  \makeatother%
  \begin{picture}(1,1.27272727)%
    \put(0,0){\includegraphics[width=\unitlength]{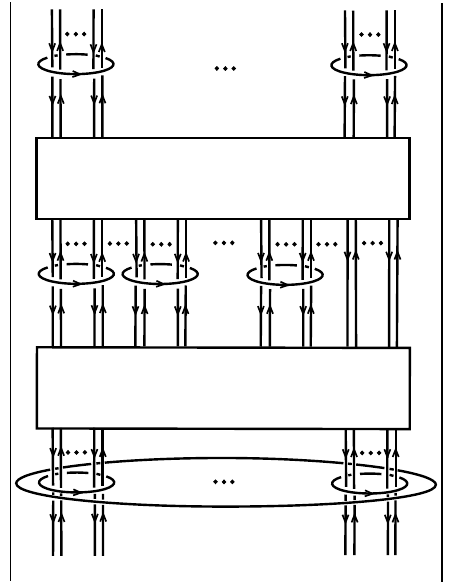}}%
    \put(0.46388093,0.9108868){\color[rgb]{0,0,0}\makebox(0,0)[lt]{\begin{minipage}{0.07638696\unitlength}\raggedright $T'$\end{minipage}}}%
    \put(0.46151519,0.44014477){\color[rgb]{0,0,0}\makebox(0,0)[lt]{\begin{minipage}{0.07638696\unitlength}\raggedright $T$\end{minipage}}}%
    \put(0.97370554,0.03764329){\color[rgb]{0,0,0}\makebox(0,0)[lt]{\begin{minipage}{0.22084041\unitlength}\raggedright $\C,\alpha$\end{minipage}}}%
  \end{picture}%
\endgroup%
}
\raisebox{45mm}{$\stackrel{(4)}{=}$}
\scalebox{0.75}{
\begingroup%
  \makeatletter%
  \providecommand\color[2][]{%
    \errmessage{(Inkscape) Color is used for the text in Inkscape, but the package 'color.sty' is not loaded}%
    \renewcommand\color[2][]{}%
  }%
  \providecommand\transparent[1]{%
    \errmessage{(Inkscape) Transparency is used (non-zero) for the text in Inkscape, but the package 'transparent.sty' is not loaded}%
    \renewcommand\transparent[1]{}%
  }%
  \providecommand\rotatebox[2]{#2}%
  \ifx\svgwidth\undefined%
    \setlength{\unitlength}{260bp}%
    \ifx\svgscale\undefined%
      \relax%
    \else%
      \setlength{\unitlength}{\unitlength * \real{\svgscale}}%
    \fi%
  \else%
    \setlength{\unitlength}{\svgwidth}%
  \fi%
  \global\let\svgwidth\undefined%
  \global\let\svgscale\undefined%
  \makeatother%
  \begin{picture}(1,1.07692308)%
    \put(0,0){\includegraphics[width=\unitlength]{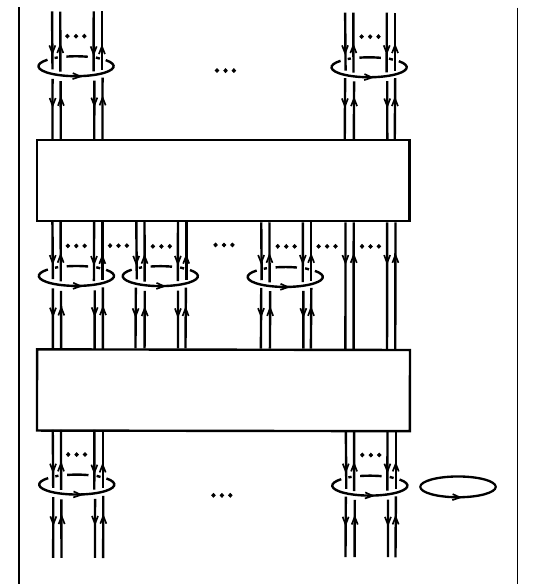}}%
    \put(0.39277997,0.76690277){\color[rgb]{0,0,0}\makebox(0,0)[lt]{\begin{minipage}{0.06463512\unitlength}\raggedright $T'$\end{minipage}}}%
    \put(0.39077819,0.36858259){\color[rgb]{0,0,0}\makebox(0,0)[lt]{\begin{minipage}{0.06463512\unitlength}\raggedright $T$\end{minipage}}}%
    \put(0.96370581,0.02332921){\color[rgb]{0,0,0}\makebox(0,0)[lt]{\begin{minipage}{0.18686496\unitlength}\raggedright $\C,\alpha$\end{minipage}}}%
  \end{picture}%
\endgroup%
}\\
&\raisebox{45mm}{$=|\halo{T'\star T}|_{\C,\alpha}\otimes |\raisebox{-2.5mm}{
\begingroup%
  \makeatletter%
  \providecommand\color[2][]{%
    \errmessage{(Inkscape) Color is used for the text in Inkscape, but the package 'color.sty' is not loaded}%
    \renewcommand\color[2][]{}%
  }%
  \providecommand\transparent[1]{%
    \errmessage{(Inkscape) Transparency is used (non-zero) for the text in Inkscape, but the package 'transparent.sty' is not loaded}%
    \renewcommand\transparent[1]{}%
  }%
  \providecommand\rotatebox[2]{#2}%
  \ifx\svgwidth\undefined%
    \setlength{\unitlength}{20bp}%
    \ifx\svgscale\undefined%
      \relax%
    \else%
      \setlength{\unitlength}{\unitlength * \real{\svgscale}}%
    \fi%
  \else%
    \setlength{\unitlength}{\svgwidth}%
  \fi%
  \global\let\svgwidth\undefined%
  \global\let\svgscale\undefined%
  \makeatother%
  \begin{picture}(1,1)%
    \put(0,0){\includegraphics[width=\unitlength]{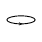}}%
  \end{picture}%
\endgroup%
}|_{\C,\alpha}=|\halo{T'\star T}|_{\C,\alpha}\otimes \varepsilon\alpha\stackrel{(7)}{=}|\halo{T'\star T}|_{\C,\alpha}$}
\end{align*}
\end{center}
Since $\alpha$ is admissible, equalities $(2)$ and $(4)$ are due to Lemma~\ref{lemmainvariantkirby2} which guarantees that $|\phantom{o}|_{\C,\alpha}$ is invariant by generalized Kirby move $KII^{g}$ whereas the third equality is due to Lemma~\ref{lemmaisotopyinv} which assures that $|\phantom{o}|_{\C,\alpha}$ is an isotopy invariant. Finally, equality $(7)$ is the consequence of the first admissibility condition \ref{ad1}.
Then we compute $|T'\star T|_{\C,\alpha}$ using the $(\ung,n+m+|\unh|+s,\unk)$-cobordism tangle $T'\square T$. We choose an opentangle $Q$ such that $U(Q)=T'\square T$. \\Let $X_1,\ldots, X_{|\ung|}, Y_1,\ldots,Y_{n}, Z_1,\ldots, Z_{m},A_1,\ldots, A_{|\unh|},B_1,\ldots, B_{s}, D_1,\ldots, D_{|\unk|}$ be any objects of $\C$ and colore the opentangle $Q$ thanks to these objects:
\vspace{-0.5cm}
\begin{center}
\resizebox{!}{0.3\textheight}{\raisebox{50mm}{$Q_{\unX,\unY,\unZ,\unA,\unB,\unD}=$}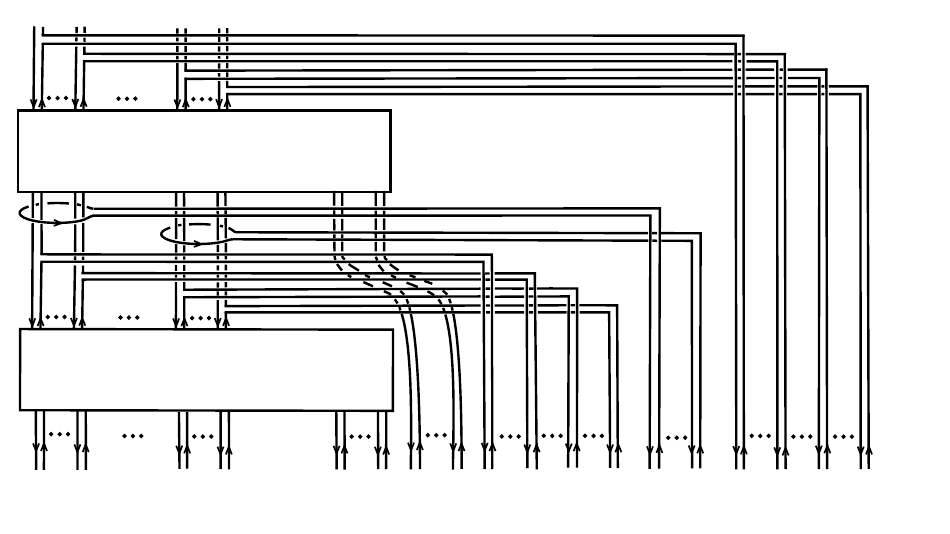}
\end{center}
We are going to factorize this morphism "by part". First, consider one opentangle $O'$ associated to the cobordism tangle $T'$ defined by:
\begin{center}
\resizebox{!}{0.15\textheight}{\raisebox{30mm}{$O'=$}
\begingroup%
  \makeatletter%
  \providecommand\color[2][]{%
    \errmessage{(Inkscape) Color is used for the text in Inkscape, but the package 'color.sty' is not loaded}%
    \renewcommand\color[2][]{}%
  }%
  \providecommand\transparent[1]{%
    \errmessage{(Inkscape) Transparency is used (non-zero) for the text in Inkscape, but the package 'transparent.sty' is not loaded}%
    \renewcommand\transparent[1]{}%
  }%
  \providecommand\rotatebox[2]{#2}%
  \newcommand*\fsize{\dimexpr\f@size pt\relax}%
  \newcommand*\lineheight[1]{\fontsize{\fsize}{#1\fsize}\selectfont}%
  \ifx\svgwidth\undefined%
    \setlength{\unitlength}{453.54330709bp}%
    \ifx\svgscale\undefined%
      \relax%
    \else%
      \setlength{\unitlength}{\unitlength * \real{\svgscale}}%
    \fi%
  \else%
    \setlength{\unitlength}{\svgwidth}%
  \fi%
  \global\let\svgwidth\undefined%
  \global\let\svgscale\undefined%
  \makeatother%
  \begin{picture}(1,0.3125)%
    \lineheight{1}%
    \setlength\tabcolsep{0pt}%
    \put(0,0){\includegraphics[width=\unitlength,page=1]{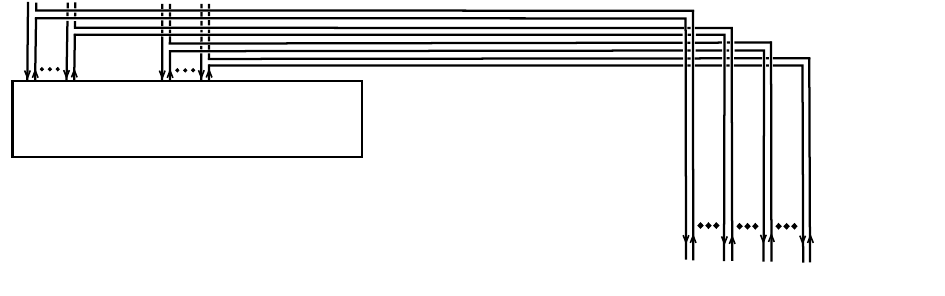}}%
    \put(0.72342272,0.03552635){\color[rgb]{0,0,0}\makebox(0,0)[lt]{\begin{minipage}{0.04278675\unitlength}\raggedright \tiny{$D_{1}$}\end{minipage}}}%
    \put(0.76051728,0.03554166){\color[rgb]{0,0,0}\makebox(0,0)[lt]{\begin{minipage}{0.04278675\unitlength}\raggedright \tiny{$D_{k_1}$}\end{minipage}}}%
    \put(0.8432872,0.03302339){\color[rgb]{0,0,0}\makebox(0,0)[lt]{\begin{minipage}{0.03993109\unitlength}\raggedright \tiny{$D_{|\unk|}$}\end{minipage}}}%
    \put(0,0){\includegraphics[width=\unitlength,page=2]{proof44.pdf}}%
    \put(0.19967287,0.17566077){\color[rgb]{0,0,0}\makebox(0,0)[lt]{\lineheight{1.25}\smash{\begin{tabular}[t]{l}$\text{T'}$\end{tabular}}}}%
  \end{picture}%
\endgroup%
}
\end{center}
\vspace{-0.3cm}
Then applying Lemma~\ref{propinvopen} to the opentangle $O'$, the universal morphism 
$|O'|_{\C}\colon C^{|\unh|+m+|\unk|}\rightarrow C^{\otimes |\unk|}$
is such that for all $\unX,\unY,\unZ,\unA,\unB,\unD$ tuples of objects of $\C$, the morphism $\iota_{D_1\otimes \ldots\otimes \iota_{D_|\unk|}}Q_{\unX,\unY,\unZ,\unA,\unB,\unD}=$ 
 \vspace{-0.5cm}
\begin{center}
\resizebox{!}{0.26\textheight}{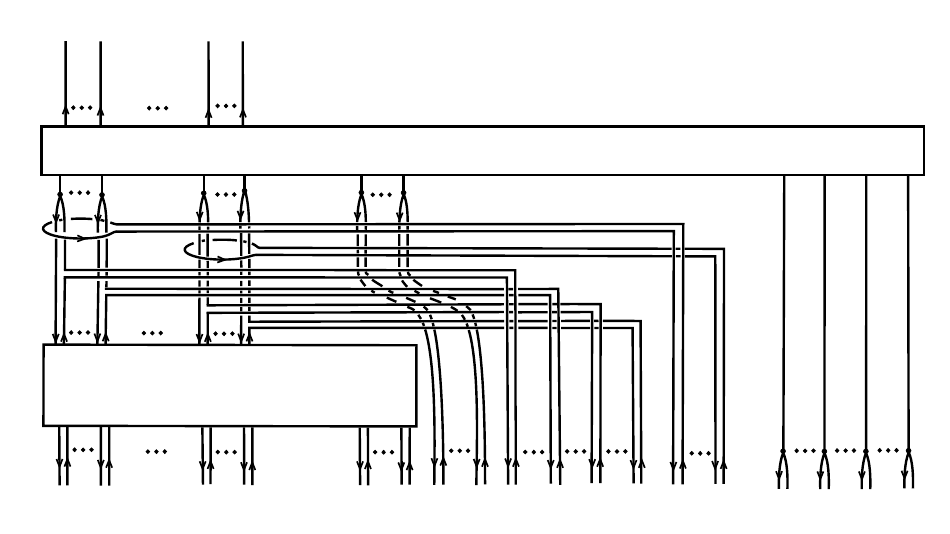}
\end{center}

Now, let $O$ be one opentangle associated to the cobordism tangle $T$ defined by:
\begin{center}
\resizebox{!}{0.15\textheight}{\raisebox{20mm}{$O=$}
\begingroup%
  \makeatletter%
  \providecommand\color[2][]{%
    \errmessage{(Inkscape) Color is used for the text in Inkscape, but the package 'color.sty' is not loaded}%
    \renewcommand\color[2][]{}%
  }%
  \providecommand\transparent[1]{%
    \errmessage{(Inkscape) Transparency is used (non-zero) for the text in Inkscape, but the package 'transparent.sty' is not loaded}%
    \renewcommand\transparent[1]{}%
  }%
  \providecommand\rotatebox[2]{#2}%
  \ifx\svgwidth\undefined%
    \setlength{\unitlength}{320bp}%
    \ifx\svgscale\undefined%
      \relax%
    \else%
      \setlength{\unitlength}{\unitlength * \real{\svgscale}}%
    \fi%
  \else%
    \setlength{\unitlength}{\svgwidth}%
  \fi%
  \global\let\svgwidth\undefined%
  \global\let\svgscale\undefined%
  \makeatother%
  \begin{picture}(1,0.45)%
    \put(0,0){\includegraphics[width=\unitlength]{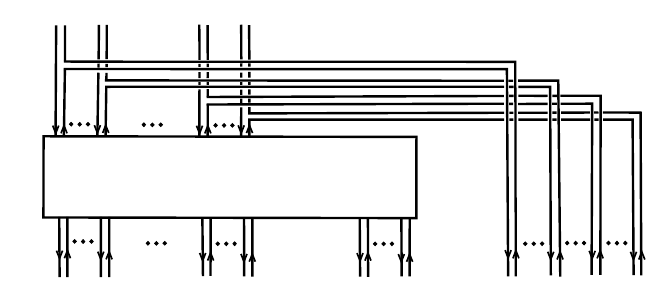}}%
    \put(0.27149054,0.19946678){\color[rgb]{0,0,0}\makebox(0,0)[lt]{\begin{minipage}{0.30703844\unitlength}\raggedright $T$\end{minipage}}}%
  \end{picture}%
\endgroup%
}
\end{center}
\vspace{-0.5cm}
Applying Lemma~\ref{propinvopen} on the opentangle $O$, the universal morphism $|O|_{\C}\colon C^{|\ung|+n+|\unh|}\rightarrow C^{\otimes |\unh|}$ is such that for all $\unX,\unY,\unZ,\unA,\unB,\unD$ tuples of objects of $\C$, the morphism $\iota_{D_1\otimes \ldots\otimes \iota_{D_|\unk|}}Q_{\unX,\unY,\unZ,\unA,\unB,\unD}=$
\vspace{-0.5cm}
\begin{center}
\resizebox{!}{0.25\textheight}{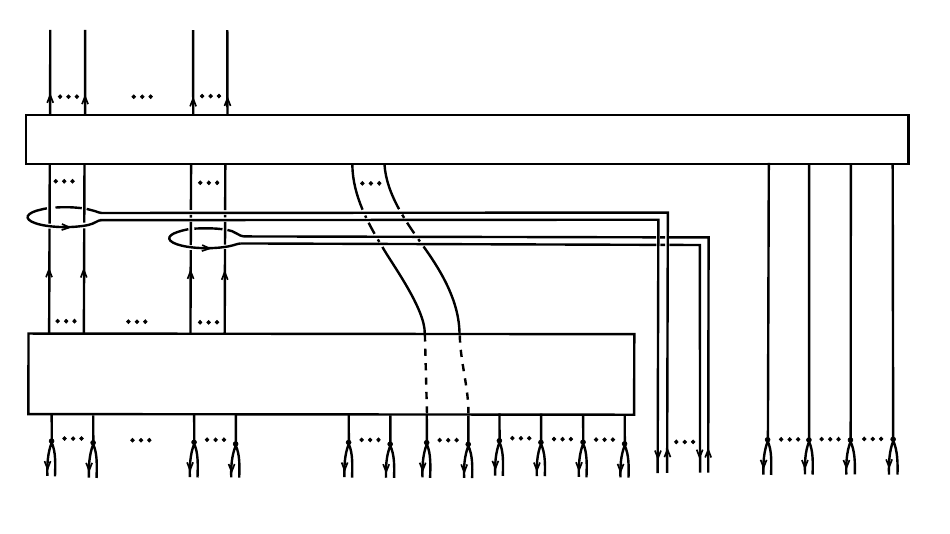}
\end{center}
Note that the mid part of the diagram corresponds to hallowed morphisms (\ref{hallowedmorphism})
and for all $\unX,\unY,\unZ,\unA,\unB,\unD$,$\iota_{D_1\otimes \ldots\otimes \iota_{D_|\unk|}}Q_{\unX,\unY,\unZ,\unA,\unB,\unD}=$
\begin{center}
\resizebox{!}{0.35\textheight}{\raisebox{50mm}{$O'=$}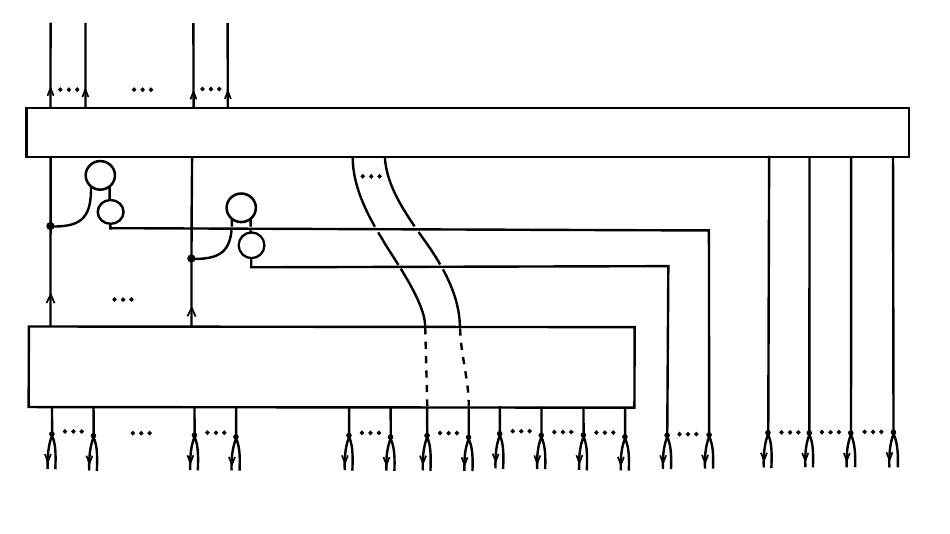}
\end{center}
Thus, we have revealed the universal morphism $|Q|_{\C}$.
As a consequence, and using that $S\alpha=\alpha$, the morphism $|T'\square T|_{\C,\alpha}$ is given by:
\begin{center}
\resizebox{!}{0.38\textheight}{\raisebox{50mm}{$|T'\square T|_{\C,\alpha}=$}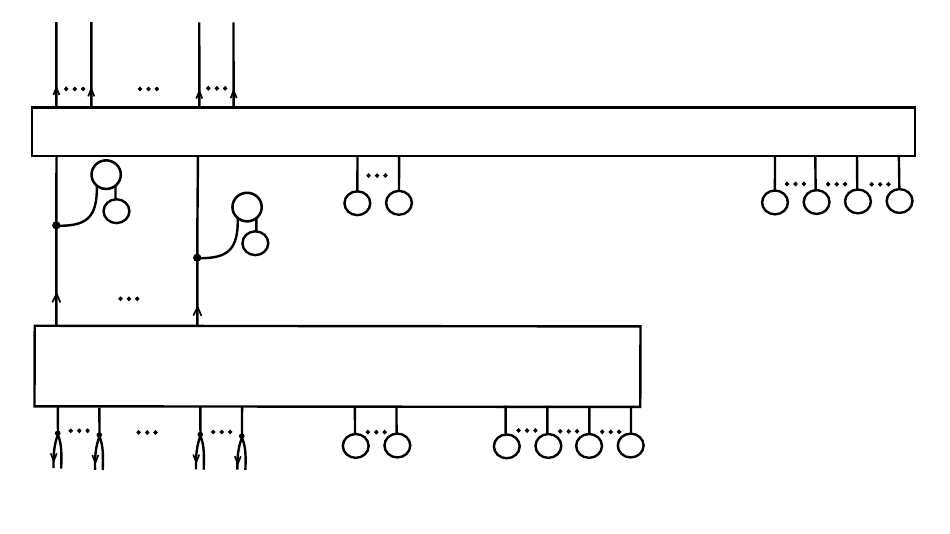}
\end{center}
\begin{align*}
=\inv{T'}\pia{\unh}|T|_{\C,\alpha}
\end{align*}

So $\Pi_{\alpha,\unk}\inv{T'\square T}\Pi_{\alpha,\ung}=\Pi_{\alpha,\unk}\inv{T'}\pia{\unh}\inv{T}\pia{\ung}$
but we have 
\[\pia{\unk}\inv{T'\square T}\pia{\ung}=\pia{\unk}\inv{T'\star T}\pia{\ung}.\]

Indeed, according to Lemma~\ref{lemmahalonothalo}, 
\[\pia{\unk}\inv{T'\square T}\pia{\ung}=|\halo{T'\square T}|_{\C,\alpha} \text{ and } \pia{\unk}\inv{T'\star T}\pia{\ung}=|\halo{T'\star T}|_{\C,\alpha}\]
and as $\halo{T'\square T}=\halo{T'\star T}$, we conclude.

\noindent Thus, 
\[\pia{\unk}\inv{T'\star T}\pia{\ung}=\pia{\unk}\inv{T'}\pia{\unh}\inv{T}\pia{\ung}\]
and as $\pia{\unh}$ is an idempotent,
\[\pia{\unk}\inv{T'\star T}\pia{\ung}=(\pia{\unk}\inv{T'}\pia{\unh})(\pia{\unh}\inv{T}\pia{\ung}).\]
As a consequence, we have shown that:
\[|\halo{T'\star T}|_{\C,\alpha}=|\halo{T'}|_{\C,\alpha}\circ |\halo{T}|_{\C,\alpha}\]

Moreover, $\gamma_{M_{T'},M_T}=\frac{\nu_{\alpha}(T'\star T)}{\nu_{\alpha}(T')\nu_{\alpha}(T)}$ so
\[\W_{\C,\alpha}(M_{T'}\circ M_{T})=\gamma_{M_{T'},M_{T}}\W_{\C,\alpha}(M_T)\W_{\C,\alpha}(M_{T'}).\]

Now, we have to prove this formula for non-connected $3$-cobordisms using the connected case. Let $M$ and $N$ be two composable $3$-cobordisms. Note that the cobordism $(M\circ N)^{\#}$ can differ from the cobordim $M^{\#}\circ N^{\#}$ only by adding a finite number of handles of type $\mathbb{S}^{2}\times [0,1]$. If $M_{T}$ is a cobordism represented by the cobordism tangle $T$, the cobordism $M_{T}\# (\mathbb{S}^{2}\times \mathbb{S}^{1})$ is represented by the cobordism tangle $T\sqcup\raisebox{-1mm}{
\begingroup%
  \makeatletter%
  \providecommand\color[2][]{%
    \errmessage{(Inkscape) Color is used for the text in Inkscape, but the package 'color.sty' is not loaded}%
    \renewcommand\color[2][]{}%
  }%
  \providecommand\transparent[1]{%
    \errmessage{(Inkscape) Transparency is used (non-zero) for the text in Inkscape, but the package 'transparent.sty' is not loaded}%
    \renewcommand\transparent[1]{}%
  }%
  \providecommand\rotatebox[2]{#2}%
  \ifx\svgwidth\undefined%
    \setlength{\unitlength}{12bp}%
    \ifx\svgscale\undefined%
      \relax%
    \else%
      \setlength{\unitlength}{\unitlength * \real{\svgscale}}%
    \fi%
  \else%
    \setlength{\unitlength}{\svgwidth}%
  \fi%
  \global\let\svgwidth\undefined%
  \global\let\svgscale\undefined%
  \makeatother%
  \begin{picture}(1,1)%
    \put(0,0){\includegraphics[width=\unitlength]{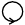}}%
  \end{picture}%
\endgroup%
}$ and $|T\sqcup\raisebox{-1mm}{
\begingroup%
  \makeatletter%
  \providecommand\color[2][]{%
    \errmessage{(Inkscape) Color is used for the text in Inkscape, but the package 'color.sty' is not loaded}%
    \renewcommand\color[2][]{}%
  }%
  \providecommand\transparent[1]{%
    \errmessage{(Inkscape) Transparency is used (non-zero) for the text in Inkscape, but the package 'transparent.sty' is not loaded}%
    \renewcommand\transparent[1]{}%
  }%
  \providecommand\rotatebox[2]{#2}%
  \ifx\svgwidth\undefined%
    \setlength{\unitlength}{12bp}%
    \ifx\svgscale\undefined%
      \relax%
    \else%
      \setlength{\unitlength}{\unitlength * \real{\svgscale}}%
    \fi%
  \else%
    \setlength{\unitlength}{\svgwidth}%
  \fi%
  \global\let\svgwidth\undefined%
  \global\let\svgscale\undefined%
  \makeatother%
  \begin{picture}(1,1)%
    \put(0,0){\includegraphics[width=\unitlength]{circle2.pdf}}%
  \end{picture}%
\endgroup%
}|_{\C,\alpha}=|T|_{ \C,\alpha}\otimes |\raisebox{-1mm}{
\begingroup%
  \makeatletter%
  \providecommand\color[2][]{%
    \errmessage{(Inkscape) Color is used for the text in Inkscape, but the package 'color.sty' is not loaded}%
    \renewcommand\color[2][]{}%
  }%
  \providecommand\transparent[1]{%
    \errmessage{(Inkscape) Transparency is used (non-zero) for the text in Inkscape, but the package 'transparent.sty' is not loaded}%
    \renewcommand\transparent[1]{}%
  }%
  \providecommand\rotatebox[2]{#2}%
  \ifx\svgwidth\undefined%
    \setlength{\unitlength}{12bp}%
    \ifx\svgscale\undefined%
      \relax%
    \else%
      \setlength{\unitlength}{\unitlength * \real{\svgscale}}%
    \fi%
  \else%
    \setlength{\unitlength}{\svgwidth}%
  \fi%
  \global\let\svgwidth\undefined%
  \global\let\svgscale\undefined%
  \makeatother%
  \begin{picture}(1,1)%
    \put(0,0){\includegraphics[width=\unitlength]{circle2.pdf}}%
  \end{picture}%
\endgroup%
}|_{\C,\alpha}=|T|_{\C,\alpha}\otimes \varepsilon\alpha=|T|_{\C,\alpha}$ since $\varepsilon\alpha=1$ so $\W_{\C,\alpha}(M_T\# (\mathbb{S}^{2}\times \mathbb{S}^{1}))=W_{\C,\alpha}(M_T)$. Then
\begin{align*}
\W_{\C,\alpha}(M\circ N)&\stackrel{(1)}{=}\W_{\C,\alpha}((M\circ N)^{\#})
\stackrel{(2)}{=}\W_{\C,\alpha}(M^{\#}\circ N^{\#})\\
& \stackrel{(3)}{=}\gamma_{M^{\#},N^{\#}}\W_{\C,\alpha}(M^{\#})\circ \W_{\C,\alpha}(N^{\#})
\stackrel{(4)}{=}\gamma_{M,N}\W_{\C,\alpha}(M)\circ \W_{\C,\alpha}(N;\alpha)
\end{align*}
Equalities $(1)$ et $(3)$ come from the definition of the invariant $\W_{\C,\alpha}(\phantom{o})$ on any cobordisms, equality $(2)$ is based on the fact that cobordisms $(M\circ N)^{\#}$ and $M^{\#}\circ N^{\#}$ could differ only by adding or suppress handles $\mathbb{S}^{2}\times [0,1]$, operation that is not detected by the invariant $\W$, and equality $(4)$ is true because we have proved it on connected cobordisms.
To conclude that $\W_{\C,\alpha}(\phantom{o})$ is a functor with anomaly, we have to check that $\gamma$ is a $2$-cocycle. We have already explained the difference between the two cobordisms $(M\circ N)^{\#}$ and $M^{\#}\circ N^{\#}$ : they could differ only by handles of type $\mathbb{S}^{2}\times [0,1]$. And as $\nu_{\alpha}(\raisebox{-1mm}{
\begingroup%
  \makeatletter%
  \providecommand\color[2][]{%
    \errmessage{(Inkscape) Color is used for the text in Inkscape, but the package 'color.sty' is not loaded}%
    \renewcommand\color[2][]{}%
  }%
  \providecommand\transparent[1]{%
    \errmessage{(Inkscape) Transparency is used (non-zero) for the text in Inkscape, but the package 'transparent.sty' is not loaded}%
    \renewcommand\transparent[1]{}%
  }%
  \providecommand\rotatebox[2]{#2}%
  \ifx\svgwidth\undefined%
    \setlength{\unitlength}{12bp}%
    \ifx\svgscale\undefined%
      \relax%
    \else%
      \setlength{\unitlength}{\unitlength * \real{\svgscale}}%
    \fi%
  \else%
    \setlength{\unitlength}{\svgwidth}%
  \fi%
  \global\let\svgwidth\undefined%
  \global\let\svgscale\undefined%
  \makeatother%
  \begin{picture}(1,1)%
    \put(0,0){\includegraphics[width=\unitlength]{circle2.pdf}}%
  \end{picture}%
\endgroup%
})=1$, then $\nu_{\alpha}(T\sqcup \raisebox{-1mm}{
\begingroup%
  \makeatletter%
  \providecommand\color[2][]{%
    \errmessage{(Inkscape) Color is used for the text in Inkscape, but the package 'color.sty' is not loaded}%
    \renewcommand\color[2][]{}%
  }%
  \providecommand\transparent[1]{%
    \errmessage{(Inkscape) Transparency is used (non-zero) for the text in Inkscape, but the package 'transparent.sty' is not loaded}%
    \renewcommand\transparent[1]{}%
  }%
  \providecommand\rotatebox[2]{#2}%
  \ifx\svgwidth\undefined%
    \setlength{\unitlength}{12bp}%
    \ifx\svgscale\undefined%
      \relax%
    \else%
      \setlength{\unitlength}{\unitlength * \real{\svgscale}}%
    \fi%
  \else%
    \setlength{\unitlength}{\svgwidth}%
  \fi%
  \global\let\svgwidth\undefined%
  \global\let\svgscale\undefined%
  \makeatother%
  \begin{picture}(1,1)%
    \put(0,0){\includegraphics[width=\unitlength]{circle2.pdf}}%
  \end{picture}%
\endgroup%
})=\nu_{\alpha}(T)\nu_{\alpha}(\raisebox{-1mm}{
\begingroup%
  \makeatletter%
  \providecommand\color[2][]{%
    \errmessage{(Inkscape) Color is used for the text in Inkscape, but the package 'color.sty' is not loaded}%
    \renewcommand\color[2][]{}%
  }%
  \providecommand\transparent[1]{%
    \errmessage{(Inkscape) Transparency is used (non-zero) for the text in Inkscape, but the package 'transparent.sty' is not loaded}%
    \renewcommand\transparent[1]{}%
  }%
  \providecommand\rotatebox[2]{#2}%
  \ifx\svgwidth\undefined%
    \setlength{\unitlength}{12bp}%
    \ifx\svgscale\undefined%
      \relax%
    \else%
      \setlength{\unitlength}{\unitlength * \real{\svgscale}}%
    \fi%
  \else%
    \setlength{\unitlength}{\svgwidth}%
  \fi%
  \global\let\svgwidth\undefined%
  \global\let\svgscale\undefined%
  \makeatother%
  \begin{picture}(1,1)%
    \put(0,0){\includegraphics[width=\unitlength]{circle2.pdf}}%
  \end{picture}%
\endgroup%
})=\nu_{\alpha}(T)$, and it is straightforward to check that $\gamma$ is a $2$-cocycle.\\

\noindent \textbullet\; We show that $\w{\phantom{o}}$ is a strong monoidal functor with anomaly.\\
Let $\Sigma_{\ung}$ and $\Sigma_{\unh}$ be two surfaces of multigenus $\ung$ and $\unh$. We have already seen that we have a canonical identification  $W_{2}(\Sigma_{\ung},\Sigma_{\unh})\colon W_{2}(\Sigma_{\ung})\otimes W_{2}(\Sigma_{\unh})\rightarrow W_{2}(\Sigma_{\ung}\sqcup \Sigma_{\unh})$ and an identity $W_{0}\colon \un \rightarrow W_{\C,\alpha}(\emptyset)$. It remains to be seen if the anomaly of the functor with anomaly $\w{\phantom{o}}$ is monoidal. Let $(M,N)$ and $(M',N')$ be two pairs of composable cobordisms and suppose that there exist $3$-cobordisms $B_{in}^{M}, B_{out}^{M}, B_{in}^{M'},B_{out}^{M'},B_{in}^{N}, B_{out}^{N}, B_{in}^{N'},B_{out}^{N'}$ obtained by composition and juxtaposition of the braiding and its inverse in $\Cobp$ such that 
\[
M=B_{out}^{M}\circ (M_{T_1}\sqcup\ldots\sqcup M_{T_m})\circ B_{in}^{M} \text{ , } 
M'=B_{out}^{M'}\circ (M'_{S_1}\sqcup\ldots\sqcup M_{S_n})\circ B_{in}^{M'} \text{ , }\]
\[N=B_{out}^{N}\circ (N_{R_1}\sqcup\ldots\sqcup N_{R_k})\circ B_{in}^{N} \text{ , } 
N'=B_{out}^{N'}\circ (N'_{O_1}\sqcup\ldots\sqcup N'_{O_l})\circ B_{in}^{N'}.
\]
where $T_1,\ldots,T_m,S_1,\ldots,S_n,R_1,\ldots,R_p,O_1,\ldots,O_l$ are cobordism tangles. Let's denote by $T=T_1\sqcup\ldots \sqcup T_m$, $S=S_1\sqcup\ldots\sqcup S_n$, $R=R_1\sqcup\ldots \sqcup R_k$ et $O_1\sqcup\ldots\sqcup O_l$.
We have 
\begin{align*}
\gamma_{M\sqcup M',N\sqcup N'}&=\gamma_{(M\sqcup M')^{\#},(N\sqcup N')^{\#}}=\dfrac{\nu_{\alpha}((T\sqcup S)\star (R\sqcup O))}{\nu_{\alpha}(T\sqcup S)\nu_{\alpha}(R\sqcup O)}\\
&\stackrel{(3)}{=}\dfrac{\nu_{\alpha}(T\square R)\nu_{\alpha} (S\star O)}{\nu_{\alpha}(T\sqcup S)\nu_{\alpha}(R\sqcup O)} =\dfrac{\nu_{\alpha}(T\star R)\nu_{\alpha} (S\star O)}{\nu_{\alpha}(T\sqcup S)\nu_{\alpha}(R\sqcup O)}=\dfrac{\nu_{\alpha}(T\star R)\nu_{\alpha} (S\star O)}{\nu_{\alpha}(T)\nu_{\alpha}(S)\nu_{\alpha}(R)\nu_{\alpha}( O)}\\
&=\gamma_{M^{\#},N^{\#}}\gamma_{M'^{\#},N'^{\#}}=\gamma_{M,N}\gamma_{M',N'}
\end{align*}
Remember the operation $\square$ on cobordism tangles which is defined above in this proof and remark that if $T,T'$ are two cobordism tangles, the number of positive (respectively negative) eigenvalues $b_{+}$ satisfies $b_+(T\star T')=b_{+}(T\square T')$ (respectively $b_{-}(T\star T')=b_{-}(T\square T'))$. Indeed, the operation $\square$ only adds a circle component which just encircles closed components and so it is not linked with other components. Then we conclude that $\gamma$ is a monoidal anomaly and so $(W_{\C,\alpha},W_2,W_0,\gamma)$ is a monoidal functor with anomaly.

\noindent \textbullet\; We show that $\w{\phantom{o}}$ is a braided functor functor with anomaly.
Let us show this result on connected surfaces $\Sigma_g$ and $\Sigma_h$ respectively of genus $g$ and $h$. The general case is analogous.
Denote by $T_{g,h}$ the following cobordism tangle
\begin{center}
\resizebox{!}{0.15\textheight}{\raisebox{30mm}{$T_{g,h}=$}
\begingroup%
  \makeatletter%
  \providecommand\color[2][]{%
    \errmessage{(Inkscape) Color is used for the text in Inkscape, but the package 'color.sty' is not loaded}%
    \renewcommand\color[2][]{}%
  }%
  \providecommand\transparent[1]{%
    \errmessage{(Inkscape) Transparency is used (non-zero) for the text in Inkscape, but the package 'transparent.sty' is not loaded}%
    \renewcommand\transparent[1]{}%
  }%
  \providecommand\rotatebox[2]{#2}%
  \ifx\svgwidth\undefined%
    \setlength{\unitlength}{280bp}%
    \ifx\svgscale\undefined%
      \relax%
    \else%
      \setlength{\unitlength}{\unitlength * \real{\svgscale}}%
    \fi%
  \else%
    \setlength{\unitlength}{\svgwidth}%
  \fi%
  \global\let\svgwidth\undefined%
  \global\let\svgscale\undefined%
  \makeatother%
  \begin{picture}(1,0.57142857)%
    \put(0,0){\includegraphics[width=\unitlength]{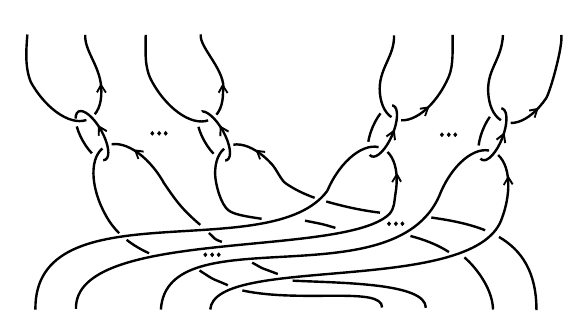}}%
  \end{picture}%
\endgroup%
}
\end{center} 
and remark that $\nu_{\alpha}(T_{g,h})=1$ since its closed components are not linked.
We have
\begin{align*}
\mathrm{W}_{\C,\alpha}\Bigg(\raisebox{-3mm}{
\begingroup%
  \makeatletter%
  \providecommand\color[2][]{%
    \errmessage{(Inkscape) Color is used for the text in Inkscape, but the package 'color.sty' is not loaded}%
    \renewcommand\color[2][]{}%
  }%
  \providecommand\transparent[1]{%
    \errmessage{(Inkscape) Transparency is used (non-zero) for the text in Inkscape, but the package 'transparent.sty' is not loaded}%
    \renewcommand\transparent[1]{}%
  }%
  \providecommand\rotatebox[2]{#2}%
  \ifx\svgwidth\undefined%
    \setlength{\unitlength}{24bp}%
    \ifx\svgscale\undefined%
      \relax%
    \else%
      \setlength{\unitlength}{\unitlength * \real{\svgscale}}%
    \fi%
  \else%
    \setlength{\unitlength}{\svgwidth}%
  \fi%
  \global\let\svgwidth\undefined%
  \global\let\svgscale\undefined%
  \makeatother%
  \begin{picture}(1,1)%
    \put(0,0){\includegraphics[width=\unitlength]{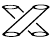}}%
    \put(0.03894843,0.20758066){\color[rgb]{0,0,0}\makebox(0,0)[lt]{\begin{minipage}{0.42250591\unitlength}\raggedright \small{$\Sigma_g$}\end{minipage}}}%
    \put(0.58238388,0.19896399){\color[rgb]{0,0,0}\makebox(0,0)[lt]{\begin{minipage}{0.43980456\unitlength}\raggedright \small{$\Sigma_h$}\end{minipage}}}%
  \end{picture}%
\endgroup%
}\Bigg)&=(\pia{h}\otimes \pia{g})|T_{g,h}|_{\C,\alpha}(\pia{g}\otimes \pia{h})\\
&=(\pia{h}\otimes \pia{g})(\Big|\raisebox{-3mm}{$\underbrace{
\begingroup%
  \makeatletter%
  \providecommand\color[2][]{%
    \errmessage{(Inkscape) Color is used for the text in Inkscape, but the package 'color.sty' is not loaded}%
    \renewcommand\color[2][]{}%
  }%
  \providecommand\transparent[1]{%
    \errmessage{(Inkscape) Transparency is used (non-zero) for the text in Inkscape, but the package 'transparent.sty' is not loaded}%
    \renewcommand\transparent[1]{}%
  }%
  \providecommand\rotatebox[2]{#2}%
  \ifx\svgwidth\undefined%
    \setlength{\unitlength}{24bp}%
    \ifx\svgscale\undefined%
      \relax%
    \else%
      \setlength{\unitlength}{\unitlength * \real{\svgscale}}%
    \fi%
  \else%
    \setlength{\unitlength}{\svgwidth}%
  \fi%
  \global\let\svgwidth\undefined%
  \global\let\svgscale\undefined%
  \makeatother%
  \begin{picture}(1,1)%
    \put(0,0){\includegraphics[width=\unitlength]{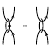}}%
  \end{picture}%
\endgroup%
}_{h\text{ times }}$}\Big|_{\C,\alpha}\otimes \Big|\raisebox{-3mm}{$\underbrace{
\begingroup%
  \makeatletter%
  \providecommand\color[2][]{%
    \errmessage{(Inkscape) Color is used for the text in Inkscape, but the package 'color.sty' is not loaded}%
    \renewcommand\color[2][]{}%
  }%
  \providecommand\transparent[1]{%
    \errmessage{(Inkscape) Transparency is used (non-zero) for the text in Inkscape, but the package 'transparent.sty' is not loaded}%
    \renewcommand\transparent[1]{}%
  }%
  \providecommand\rotatebox[2]{#2}%
  \ifx\svgwidth\undefined%
    \setlength{\unitlength}{24bp}%
    \ifx\svgscale\undefined%
      \relax%
    \else%
      \setlength{\unitlength}{\unitlength * \real{\svgscale}}%
    \fi%
  \else%
    \setlength{\unitlength}{\svgwidth}%
  \fi%
  \global\let\svgwidth\undefined%
  \global\let\svgscale\undefined%
  \makeatother%
  \begin{picture}(1,1)%
    \put(0,0){\includegraphics[width=\unitlength]{smallcylinders.pdf}}%
  \end{picture}%
\endgroup%
}_{g\text{ times }}$}\Big|_{\C,\alpha})\tau_{C^{\otimes g},C^{\otimes h}}(\pia{g}\otimes \pia{h})\\
&=(\pia{h}\Big|\raisebox{-3mm}{$\underbrace{
\begingroup%
  \makeatletter%
  \providecommand\color[2][]{%
    \errmessage{(Inkscape) Color is used for the text in Inkscape, but the package 'color.sty' is not loaded}%
    \renewcommand\color[2][]{}%
  }%
  \providecommand\transparent[1]{%
    \errmessage{(Inkscape) Transparency is used (non-zero) for the text in Inkscape, but the package 'transparent.sty' is not loaded}%
    \renewcommand\transparent[1]{}%
  }%
  \providecommand\rotatebox[2]{#2}%
  \ifx\svgwidth\undefined%
    \setlength{\unitlength}{24bp}%
    \ifx\svgscale\undefined%
      \relax%
    \else%
      \setlength{\unitlength}{\unitlength * \real{\svgscale}}%
    \fi%
  \else%
    \setlength{\unitlength}{\svgwidth}%
  \fi%
  \global\let\svgwidth\undefined%
  \global\let\svgscale\undefined%
  \makeatother%
  \begin{picture}(1,1)%
    \put(0,0){\includegraphics[width=\unitlength]{smallcylinders.pdf}}%
  \end{picture}%
\endgroup%
}_{h\text{ times }}$}\Big|_{\C,\alpha}\pia{g}\otimes \pia{g}\Big|\raisebox{-3mm}{$\underbrace{
\begingroup%
  \makeatletter%
  \providecommand\color[2][]{%
    \errmessage{(Inkscape) Color is used for the text in Inkscape, but the package 'color.sty' is not loaded}%
    \renewcommand\color[2][]{}%
  }%
  \providecommand\transparent[1]{%
    \errmessage{(Inkscape) Transparency is used (non-zero) for the text in Inkscape, but the package 'transparent.sty' is not loaded}%
    \renewcommand\transparent[1]{}%
  }%
  \providecommand\rotatebox[2]{#2}%
  \ifx\svgwidth\undefined%
    \setlength{\unitlength}{24bp}%
    \ifx\svgscale\undefined%
      \relax%
    \else%
      \setlength{\unitlength}{\unitlength * \real{\svgscale}}%
    \fi%
  \else%
    \setlength{\unitlength}{\svgwidth}%
  \fi%
  \global\let\svgwidth\undefined%
  \global\let\svgscale\undefined%
  \makeatother%
  \begin{picture}(1,1)%
    \put(0,0){\includegraphics[width=\unitlength]{smallcylinders.pdf}}%
  \end{picture}%
\endgroup%
}_{g\text{ times }}$}\Big|_{\C,\alpha}\pia{h}) \tau_{(C^{\otimes g})_{\alpha},(C^{\otimes h})_{\alpha}}\\
&=\big(\w{\id_{\Sigma_{h}}}\otimes \w{\id_{\Sigma_{g}}}\big)\tau_{\w{\Sigma_{g}},\w{\Sigma_{h}}}
\end{align*}
and $\w{\phantom{o}}$ is then a braided functor with anomaly.

\noindent \textbullet\; Considering the braided monoidal functor with anomaly $\mathrm{W_{\C,\alpha}}$, remember the associated unitalized functor $\mathrm{V}_{\C,\alpha}$ is still braided monoidal with anomaly (see Lemma~\ref{unilemma}). To show that $\mathrm{V}_{\C,\alpha}$ is a TQFT with anomaly, we just have to prove that objects associated to surfaces are transparent.
Let $\Sigma_{g}$ be the canonical surface of genus $g$. Recall that $\Sigma_{g}\times [0,1]$ is encoded by the $(g,g,g)$-cobordism tangle in Figure~\ref{cylinder}.
Denote by 
\[\Pi_{g}=\omega(\alpha\otimes \alpha)^{-g}\W_{\C}(\Sigma_{g}\times [0,1];\alpha)\]

Remark that $\gamma_{\Sigma_g,\Sigma_g}=\omega(\alpha\otimes\alpha)^{-g}$. Then, as image of an identity by a functor with anomaly up to the anomaly $\gamma_{\Sigma_g,\Sigma_g}$ (see Lemma~\ref{lemmaidempipo}), $\Pi_g$ is an idempotent of $\C$. In order to show that $\mr{Im}(\Pi_g)$ is transparent, we compute\\ ${\omega(\alpha\otimes\alpha)^{g}\Pi_{g}\circ(\id_{C^{\otimes g}}\otimes \omega(\id_{C}\otimes S))\delta_{C^{\otimes g}}}$ which is equal to:\\
\vspace{-0.5cm}
\resizebox{!}{0.18\textheight}{
\begingroup%
  \makeatletter%
  \providecommand\color[2][]{%
    \errmessage{(Inkscape) Color is used for the text in Inkscape, but the package 'color.sty' is not loaded}%
    \renewcommand\color[2][]{}%
  }%
  \providecommand\transparent[1]{%
    \errmessage{(Inkscape) Transparency is used (non-zero) for the text in Inkscape, but the package 'transparent.sty' is not loaded}%
    \renewcommand\transparent[1]{}%
  }%
  \providecommand\rotatebox[2]{#2}%
  \ifx\svgwidth\undefined%
    \setlength{\unitlength}{160bp}%
    \ifx\svgscale\undefined%
      \relax%
    \else%
      \setlength{\unitlength}{\unitlength * \real{\svgscale}}%
    \fi%
  \else%
    \setlength{\unitlength}{\svgwidth}%
  \fi%
  \global\let\svgwidth\undefined%
  \global\let\svgscale\undefined%
  \makeatother%
  \begin{picture}(1,1)%
    \put(0,0){\includegraphics[width=\unitlength,page=1]{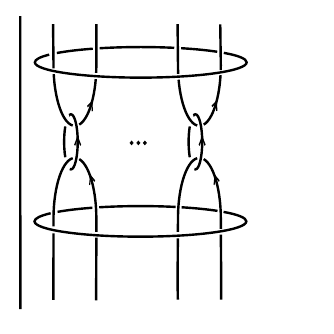}}%
    \put(0.8182777,0.08292474){\color[rgb]{0,0,0}\makebox(0,0)[lt]{\begin{minipage}{0.21030999\unitlength}\raggedright $_{\C,\alpha}$\end{minipage}}}%
    \put(0,0){\includegraphics[width=\unitlength,page=2]{prooftqftunit3.pdf}}%
  \end{picture}%
\endgroup%
}\raisebox{20mm}{$\circ \quad(\id_{C^{\otimes g}}\otimes \omega(\id_{C}\otimes S))\delta_{C^{\otimes g}}$}
\raisebox{20mm}{$=$}\resizebox{!}{0.2\textheight}{
\begingroup%
  \makeatletter%
  \providecommand\color[2][]{%
    \errmessage{(Inkscape) Color is used for the text in Inkscape, but the package 'color.sty' is not loaded}%
    \renewcommand\color[2][]{}%
  }%
  \providecommand\transparent[1]{%
    \errmessage{(Inkscape) Transparency is used (non-zero) for the text in Inkscape, but the package 'transparent.sty' is not loaded}%
    \renewcommand\transparent[1]{}%
  }%
  \providecommand\rotatebox[2]{#2}%
  \ifx\svgwidth\undefined%
    \setlength{\unitlength}{160bp}%
    \ifx\svgscale\undefined%
      \relax%
    \else%
      \setlength{\unitlength}{\unitlength * \real{\svgscale}}%
    \fi%
  \else%
    \setlength{\unitlength}{\svgwidth}%
  \fi%
  \global\let\svgwidth\undefined%
  \global\let\svgscale\undefined%
  \makeatother%
  \begin{picture}(1,1)%
    \put(0,0){\includegraphics[width=\unitlength]{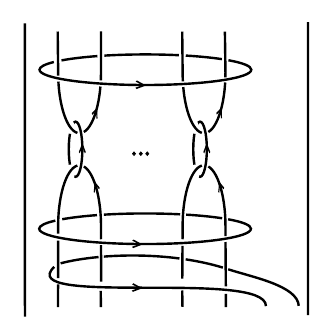}}%
    \put(0.93890959,0.06972814){\color[rgb]{0,0,0}\makebox(0,0)[lt]{\begin{minipage}{0.21031001\unitlength}\raggedright $_{\C,\alpha}$\end{minipage}}}%
  \end{picture}%
\endgroup%
}\\\raisebox{20mm}{$\stackrel{(2)}{=}$}
\resizebox{!}{0.18\textheight}{
\begingroup%
  \makeatletter%
  \providecommand\color[2][]{%
    \errmessage{(Inkscape) Color is used for the text in Inkscape, but the package 'color.sty' is not loaded}%
    \renewcommand\color[2][]{}%
  }%
  \providecommand\transparent[1]{%
    \errmessage{(Inkscape) Transparency is used (non-zero) for the text in Inkscape, but the package 'transparent.sty' is not loaded}%
    \renewcommand\transparent[1]{}%
  }%
  \providecommand\rotatebox[2]{#2}%
  \ifx\svgwidth\undefined%
    \setlength{\unitlength}{160bp}%
    \ifx\svgscale\undefined%
      \relax%
    \else%
      \setlength{\unitlength}{\unitlength * \real{\svgscale}}%
    \fi%
  \else%
    \setlength{\unitlength}{\svgwidth}%
  \fi%
  \global\let\svgwidth\undefined%
  \global\let\svgscale\undefined%
  \makeatother%
  \begin{picture}(1,1)%
    \put(0,0){\includegraphics[width=\unitlength,page=1]{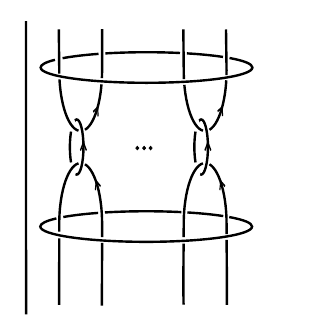}}%
    \put(0.94222294,0.07634804){\color[rgb]{0,0,0}\makebox(0,0)[lt]{\begin{minipage}{0.21030999\unitlength}\raggedright $_{\C,\alpha}$\end{minipage}}}%
    \put(0,0){\includegraphics[width=\unitlength,page=2]{prooftqftunit2.pdf}}%
  \end{picture}%
\endgroup%
}
\raisebox{20mm}{$\stackrel{(4)}{=}$}
\resizebox{!}{0.18\textheight}{}\raisebox{20mm}{$\otimes$}\resizebox{!}{0.18\textheight}{
\begingroup%
  \makeatletter%
  \providecommand\color[2][]{%
    \errmessage{(Inkscape) Color is used for the text in Inkscape, but the package 'color.sty' is not loaded}%
    \renewcommand\color[2][]{}%
  }%
  \providecommand\transparent[1]{%
    \errmessage{(Inkscape) Transparency is used (non-zero) for the text in Inkscape, but the package 'transparent.sty' is not loaded}%
    \renewcommand\transparent[1]{}%
  }%
  \providecommand\rotatebox[2]{#2}%
  \ifx\svgwidth\undefined%
    \setlength{\unitlength}{160bp}%
    \ifx\svgscale\undefined%
      \relax%
    \else%
      \setlength{\unitlength}{\unitlength * \real{\svgscale}}%
    \fi%
  \else%
    \setlength{\unitlength}{\svgwidth}%
  \fi%
  \global\let\svgwidth\undefined%
  \global\let\svgscale\undefined%
  \makeatother%
  \begin{picture}(1,1)%
    \put(0,0){\includegraphics[width=\unitlength]{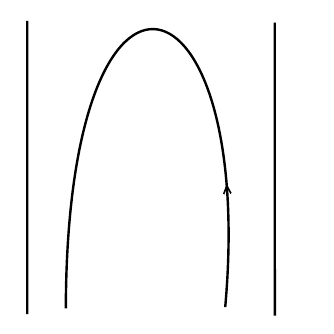}}%
    \put(0.83923638,0.06812534){\color[rgb]{0,0,0}\makebox(0,0)[lt]{\begin{minipage}{0.21031001\unitlength}\raggedright $_{\C,\alpha}$\end{minipage}}}%
  \end{picture}%
\endgroup%
}\\

$=\omega(\alpha\otimes\alpha)^{g}\Pi_{g}\otimes \varepsilon$

Equality $(2)$ uses the fact that $\alpha$ is an admissible element so $|\phantom{O}|_{\C,\alpha}$ is invariant by the generalized Kirby move $KII^{g}$. Equality $(4$) is based on the multiplicativity of $|\phantom{O}|_{\C,\alpha}$ on a cobordism tangle seen as the disjoint union of a $(g,g+2,g)$-cobordism tangle and a $(1,0,0)$-cobordism tangle.
This calculus shows that $\Pi_{g}\circ(\id_{C^{\otimes g}}\otimes \omega(\id_{C}\otimes S))=\Pi_{g}\otimes \varepsilon$. Applying Lemma~\ref{transparency} to $\Pi_g$, we get the result.
\end{proof}

The internal TQFT is now constructed. We will show in Section~\ref{modular} that the construction is a generalization of the Reshetikhin-Turaev one.
Now, we prove the second main Theorem~\ref{maintheorem2} of this paper: the TQFT can be computed using only structural morphisms of the Hopf algebra coend.


\begin{remarks}
    \item To compute the TQFT, note that the product $m$ of the coend $C$ is only used to express the universal coaction on tensorial products of the coend $C$ (see Figure~\ref{deltan}).
\begin{figure}[H]
\begin{center}
\resizebox{!}{0.18\textheight}{\raisebox{20mm}{$\delta_{C^{\otimes n}}=$}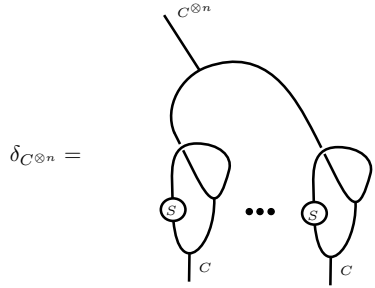}
\end{center}
\caption{The morphism $\delta_{C^{\otimes n}}$.}
\label{deltan}
\end{figure}
\end{remarks}

\begin{proof} We prove Theorem~\ref{maintheorem2}.\\
In \cite{BrugVirHopf}, Bruguières and Virelizier show that one $(|\ung|,n,0)$-cobordism tangle without exit components (called \emph{ribbon handles}) gives a morphism ${|T|_{\C,\alpha}\colon C^{|\ung|}\rightarrow \un}$ which is expressed entirely in terms of structural morphisms (except $m$) of the coend $C$ and $\alpha$. We want to show that it is still the case for any $(\ung,n,\unh)$-cobordism tangle. The product $m$ is needed only to express $\delta_{C^{\otimes n}}$.
Let $T$ be a $(\ung,n,\unh)$-cobordism tangle. We pull down the exit components before pull them up as shown of the following diagram. We denote by $T'$ the $(\ung.\unh,n,0)$-cobordism tangle defined just below.
\begin{center}
\resizebox{!}{0.15\textheight}{
\begingroup%
  \makeatletter%
  \providecommand\color[2][]{%
    \errmessage{(Inkscape) Color is used for the text in Inkscape, but the package 'color.sty' is not loaded}%
    \renewcommand\color[2][]{}%
  }%
  \providecommand\transparent[1]{%
    \errmessage{(Inkscape) Transparency is used (non-zero) for the text in Inkscape, but the package 'transparent.sty' is not loaded}%
    \renewcommand\transparent[1]{}%
  }%
  \providecommand\rotatebox[2]{#2}%
  \newcommand*\fsize{\dimexpr\f@size pt\relax}%
  \newcommand*\lineheight[1]{\fontsize{\fsize}{#1\fsize}\selectfont}%
  \ifx\svgwidth\undefined%
    \setlength{\unitlength}{412.5bp}%
    \ifx\svgscale\undefined%
      \relax%
    \else%
      \setlength{\unitlength}{\unitlength * \real{\svgscale}}%
    \fi%
  \else%
    \setlength{\unitlength}{\svgwidth}%
  \fi%
  \global\let\svgwidth\undefined%
  \global\let\svgscale\undefined%
  \makeatother%
  \begin{picture}(1,0.32727273)%
    \lineheight{1}%
    \setlength\tabcolsep{0pt}%
    \put(0,0){\includegraphics[width=\unitlength,page=1]{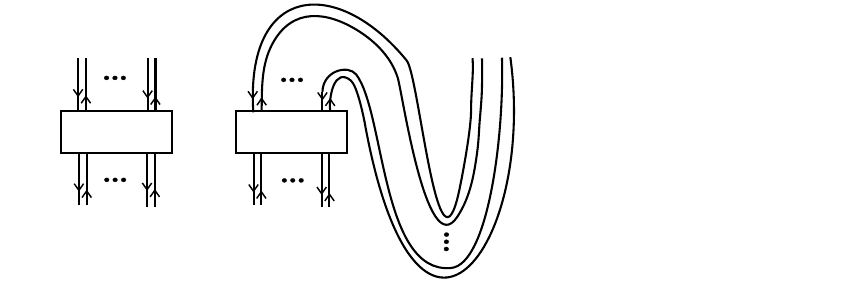}}%
    \put(0.22255912,0.17665779){\color[rgb]{0,0,0}\makebox(0,0)[lt]{\begin{minipage}{0.17644732\unitlength}\raggedright $=$\end{minipage}}}%
    \put(-0.00123371,0.18096039){\color[rgb]{0,0,0}\makebox(0,0)[lt]{\begin{minipage}{0.17644732\unitlength}\raggedright $T=$\end{minipage}}}%
    \put(0.61059203,0.18872571){\color[rgb]{0,0,0}\makebox(0,0)[lt]{\begin{minipage}{0.17644732\unitlength}\raggedright $=$\end{minipage}}}%
    \put(0,0){\includegraphics[width=\unitlength,page=2]{proofhopf11.pdf}}%
    \put(0.79246452,0.19171673){\color[rgb]{0,0,0}\makebox(0,0)[lt]{\begin{minipage}{0.17644732\unitlength}\raggedright $T'$\end{minipage}}}%
    \put(0,0){\includegraphics[width=\unitlength,page=3]{proofhopf11.pdf}}%
  \end{picture}%
\endgroup%
}
\end{center}
where $\ung=(g_1,\ldots,g_r)$, $\unh=(h_1,\ldots,h_s)$ and $\ung.\unh=(g_1,\ldots,g_r,h_1,\ldots,h_s)$.
Then just remark that 
\begin{center}
\resizebox{!}{0.17\textheight}{\raisebox{42mm}{$|T|_{\C,\alpha}=$}
\begingroup%
  \makeatletter%
  \providecommand\color[2][]{%
    \errmessage{(Inkscape) Color is used for the text in Inkscape, but the package 'color.sty' is not loaded}%
    \renewcommand\color[2][]{}%
  }%
  \providecommand\transparent[1]{%
    \errmessage{(Inkscape) Transparency is used (non-zero) for the text in Inkscape, but the package 'transparent.sty' is not loaded}%
    \renewcommand\transparent[1]{}%
  }%
  \providecommand\rotatebox[2]{#2}%
  \ifx\svgwidth\undefined%
    \setlength{\unitlength}{160bp}%
    \ifx\svgscale\undefined%
      \relax%
    \else%
      \setlength{\unitlength}{\unitlength * \real{\svgscale}}%
    \fi%
  \else%
    \setlength{\unitlength}{\svgwidth}%
  \fi%
  \global\let\svgwidth\undefined%
  \global\let\svgscale\undefined%
  \makeatother%
  \begin{picture}(1,1.05)%
    \put(0,0){\includegraphics[width=\unitlength,page=1]{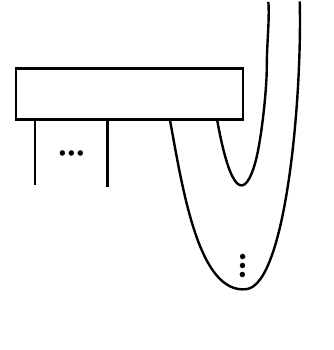}}%
    \put(0.33288133,0.81231713){\color[rgb]{0,0,0}\makebox(0,0)[lt]{\begin{minipage}{0.52924004\unitlength}\raggedright $|T'|_{\C,\alpha}$\end{minipage}}}%
    \put(0,0){\includegraphics[width=\unitlength,page=2]{proofhopf2.pdf}}%
    \put(0.70097177,0.39129513){\color[rgb]{0,0,0}\makebox(0,0)[lt]{\begin{minipage}{0.05796054\unitlength}\raggedright $\alpha$\end{minipage}}}%
    \put(0,0){\includegraphics[width=\unitlength,page=3]{proofhopf2.pdf}}%
    \put(0.70067682,0.08004318){\color[rgb]{0,0,0}\makebox(0,0)[lt]{\begin{minipage}{0.05796054\unitlength}\raggedright $\alpha$\end{minipage}}}%
    \put(0.08030554,0.46644602){\color[rgb]{0,0,0}\makebox(0,0)[lt]{\begin{minipage}{0.20344934\unitlength}\raggedright $C$\end{minipage}}}%
    \put(0.30324363,0.46332747){\color[rgb]{0,0,0}\makebox(0,0)[lt]{\begin{minipage}{0.20344934\unitlength}\raggedright $C$\end{minipage}}}%
    \put(0.72544299,1.01052682){\color[rgb]{0,0,0}\makebox(0,0)[lt]{\begin{minipage}{0.20344934\unitlength}\raggedright $C$\end{minipage}}}%
    \put(0.92264244,1.00983952){\color[rgb]{0,0,0}\makebox(0,0)[lt]{\begin{minipage}{0.20344934\unitlength}\raggedright $C$\end{minipage}}}%
    \put(0.44246189,0.63658342){\color[rgb]{0,0,0}\makebox(0,0)[lt]{\begin{minipage}{0.20344934\unitlength}\raggedright $C$\end{minipage}}}%
    \put(0.57993619,0.63658342){\color[rgb]{0,0,0}\makebox(0,0)[lt]{\begin{minipage}{0.20344934\unitlength}\raggedright $C$\end{minipage}}}%
  \end{picture}%
\endgroup%
}
\end{center}
According to the result of Bruguières and Virelizier, $|T'|_{\C,\alpha}$ can be expressed only thanks to structural morphisms of $C$ (except $m$) and $\alpha$. So is $|T|_{\C,\alpha}$ and then so is the TQFT with anomaly $\mathrm{V}_{\C,\alpha}$ (the normalization coefficient $\nu_{\alpha}$ of the TQFT is expressed using only morphisms $\theta_{+}$, $\theta_{-}$, and $\alpha$; see equality~(\ref{normalization})).
\end{proof}

\begin{remarks}
\item Note that using results on "Hopf diagrams" constructed in \cite{BrugVir}, we are able to compute easily all invariants coming from the internal TQFT, started from a tangle presentation of a cobordism.
\end{remarks}
\section{Applications on modular and premodular cases}\label{modular}

When $\C$ is premodular, using Kirby color, we can build internal TQFTs. In the modular case, we show that our internal TQFT is a transparent lift of Turaev's one before studying the dependancy between our TQFT $\V_{\C}$ and the category $\C$. Using theses results, we are able to compare our TQFT and Turaev's one in the most general case of moduralizable premodular categories: our internal TQFT is still a transparent lift of the Turaev TQFT built on the modularized category of $\C$. We don't know yet if it exists internal TQFTs that are not of this form.

\subsection{Premodular category and internal TQFT} Let $\C$ be a premodular category and denote by $\Lambda_{\C}$ a representative set of simple objects of $\C$. Recall that the category $\C$ has a coend object
\[C=\bigoplus_{\lambda \in \Lambda_{\C}}\lambda^{*}\otimes \lambda.\]

For each $\lambda\in \Lambda_{\C}$, there exist morphisms $p_{\lambda}\colon C\rightarrow\lambda^{*}\otimes \lambda$ and $q_{\lambda}\colon \lambda \rightarrow C$ such that $\mr{id_{\lambda}}=\sum_{\lambda \in \Lambda_{\C}}q_{\lambda}p_{\lambda}$ and $p_{\lambda}q_{\lambda}=\delta_{\lambda,\mu}\mr{id_{\lambda^{*}\otimes \lambda}}$. Any object $X\in \C$ has a decomposition $\bigotimes_{i\in I} \lambda_i$ where $I$ is a finite set and $\lambda_i\in \Lambda_{\C}$. We set :
\[\iota_{X}=\sum_{i\in I}q_{\lambda_i}\circ (Q_i^{*}\otimes P_i)\colon X^{*}\otimes X \rightarrow C\]
where $P_i$ and $Q_i$ are morphisms such that 
\[\mr{id}_X=\sum_{i\in I}Q_iP_i \text{ and } P_iQ_j=\delta_{i,j}\mr{id_{\lambda_i}}\]
Note that $\iota_X$ does not depend on morphisms $P_i$ and $Q_i$ and remark that $\iota_{\lambda}=q_{\lambda}$ for any $\lambda\in \Lambda_{\C}$. It can be easily proved that $(C,\iota)$ is a coend of the functor (\ref{dualityfunctor}). Recall that $C$ has a structure of Hopf algebra (see Theorem~\ref{thcoendhopf}).

Assume that the category is \emph{normalizable}, that is,
\[
\Delta_{\pm}:=\sum_{\lambda\in\Lambda_{\C}}\theta_{\lambda}^{\pm}\mathrm{dim}^{2}_{q} \in \kk^{*}
\]
and suppose that $\C$ has invertible dimension 
\[\mr{dim}(\C):=\sum_{\lambda\in \Lambda_{\C}}\mathrm{dim}^{2}_{q}(\lambda).\]\index{$\mathrm{dim}(\C)$}
The category $\C$ has a Kirby element 
\[\alpha_{K}:=\frac{1}{\mr{dim}(\C)}\sum_{\lambda \in \Lambda_{\C}}\dimq(\lambda)e_{\lambda}\]
where $e_{\lambda}=\iota_{\lambda}\widetilde{\mr{coev}_{\lambda}}$.

\begin{proposition}\label{kirbyadmissible}
The morphism $\alpha_K$ is an admissible element.
\end{proposition}

\begin{proof}
As $\C$ is supposed to be normalizable and has invertible dimension, then $\theta_{\pm}\alpha_K=\Delta_{\pm}$ are invertible. Moreover, as $\epsilon_{\lambda}e_{\lambda}=\mathrm{dim}_{q}(\lambda)$, we have $\varepsilon\alpha_K=1$. Next, according to the Lemma 3.2 (c) from \cite{Vir}, the element $\alpha_K$ is such that $S\alpha_K=\alpha_K$. Thus, admissibility conditions \ref{ad1}, \ref{ad2}, \ref{ad3} are satisfied. 

We have to prove that \ref{ad4} and \ref{ad5} are verified by $\alpha_K$. For this, just remark (see Section 2.5 in \cite{Ker} or Lemma 4.1 in \cite{BrugVircenter} for more general results) that $\alpha_K$ is a 2-sided integral of the coend $C$, that means 
\[m(\alpha_K\otimes\mathrm{id_C})=\alpha_K\epsilon=m(\mathrm{id}_C\otimes \alpha_K)\]
Note that \ref{ad4} is just a special case of Lemma~\ref{lemmaprotran} (1) for $X=C^{\otimes n}$. At last, as $\alpha_K$ is an integral,
\begin{center}
    \resizebox{!}{0.1\textheight}{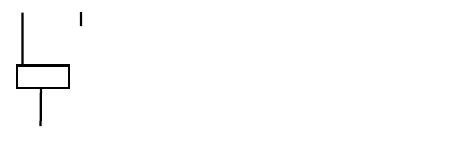}
\end{center}
As a consequence, condition \ref{ad5} is satisfied too.
\end{proof}

We just have proved that $\alpha_K$ is an admissible element and the TQFT with anomaly $V_{\C,\alpha_K}$ is well-defined in the premodular normalizable case.

Denote by $\T$ the subcategory of transparent objects of $\C$. Let $X$ be an object of $\C$ and decomposed as direct sum of $n$ simple objects:
\[X:=\bigoplus_{i=1}^{n}S_i\]
The \emph{transparent part}\index{transparent part} $X_{\T}$\index{\(X_{\T}\)} of $X$ is defined as:
\[X_{\T}:=\bigoplus_{i \in \{k \;|\; S_k \in \T \}}S_i.\]
The next Lemma identifies the natural transformation  
\[\Pi_{\alpha_K,X}=[\id_{X}\otimes \omega(\id_{C}\otimes \alpha_K)]\delta_{X}\]
with the projector on the transparent part of the object $X$. 

\begin{lemma}\label{lemmaprotran}
Let $\C$ be a premodular category with invertible dimension. Then, for all objects $X$ of $\C$:
\begin{enumerate}[label=i)]
\item The morphism $\Pi_{\alpha_K,X}=[\id_{X}\otimes \omega(\id_{C}\otimes \alpha_K)]\delta_{X}$ is an idempotent of $\C$.
\item If $X$ is a transparent object of $\C$, $\Pi_{\alpha_K,X}=\id_X$.
\item If $X$ is such that $\Pi_{\alpha_K,X}=\id_X$, then $X$ is transparent.
\end{enumerate}
\end{lemma}

\begin{proof}
For all objects $X$ of $\C$,
\begin{center}
\begingroup%
  \makeatletter%
  \providecommand\color[2][]{%
    \errmessage{(Inkscape) Color is used for the text in Inkscape, but the package 'color.sty' is not loaded}%
    \renewcommand\color[2][]{}%
  }%
  \providecommand\transparent[1]{%
    \errmessage{(Inkscape) Transparency is used (non-zero) for the text in Inkscape, but the package 'transparent.sty' is not loaded}%
    \renewcommand\transparent[1]{}%
  }%
  \providecommand\rotatebox[2]{#2}%
  \ifx\svgwidth\undefined%
    \setlength{\unitlength}{56bp}%
    \ifx\svgscale\undefined%
      \relax%
    \else%
      \setlength{\unitlength}{\unitlength * \real{\svgscale}}%
    \fi%
  \else%
    \setlength{\unitlength}{\svgwidth}%
  \fi%
  \global\let\svgwidth\undefined%
  \global\let\svgscale\undefined%
  \makeatother%
  \begin{picture}(1,1.42857143)%
    \put(0,0){\includegraphics[width=\unitlength,page=1]{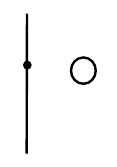}}%
    \put(0.62640174,0.84838759){\color[rgb]{0,0,0}\makebox(0,0)[lt]{\begin{minipage}{0.10782629\unitlength}\raggedright \tiny{$\alpha_K$}\end{minipage}}}%
    \put(0,0){\includegraphics[width=\unitlength,page=2]{trans1.pdf}}%
    \put(0.50243698,1.18118145){\color[rgb]{0,0,0}\makebox(0,0)[lt]{\begin{minipage}{0.10782629\unitlength}\raggedright \tiny{$\omega$}\end{minipage}}}%
    \put(0,0){\includegraphics[width=\unitlength,page=3]{trans1.pdf}}%
    \put(0.61843897,0.27732588){\color[rgb]{0,0,0}\makebox(0,0)[lt]{\begin{minipage}{0.10782629\unitlength}\raggedright \tiny{$\alpha_K$}\end{minipage}}}%
    \put(0,0){\includegraphics[width=\unitlength,page=4]{trans1.pdf}}%
    \put(0.49447421,0.61012017){\color[rgb]{0,0,0}\makebox(0,0)[lt]{\begin{minipage}{0.10782629\unitlength}\raggedright \tiny{$\omega$}\end{minipage}}}%
    \put(0,0){\includegraphics[width=\unitlength,page=5]{trans1.pdf}}%
    \put(0.27741388,0.10528483){\color[rgb]{0,0,0}\makebox(0,0)[lt]{\begin{minipage}{0.39289115\unitlength}\raggedright \tiny{$X$}\end{minipage}}}%
  \end{picture}%
\endgroup%

\raisebox{15mm}{$=$}
\begingroup%
  \makeatletter%
  \providecommand\color[2][]{%
    \errmessage{(Inkscape) Color is used for the text in Inkscape, but the package 'color.sty' is not loaded}%
    \renewcommand\color[2][]{}%
  }%
  \providecommand\transparent[1]{%
    \errmessage{(Inkscape) Transparency is used (non-zero) for the text in Inkscape, but the package 'transparent.sty' is not loaded}%
    \renewcommand\transparent[1]{}%
  }%
  \providecommand\rotatebox[2]{#2}%
  \ifx\svgwidth\undefined%
    \setlength{\unitlength}{56bp}%
    \ifx\svgscale\undefined%
      \relax%
    \else%
      \setlength{\unitlength}{\unitlength * \real{\svgscale}}%
    \fi%
  \else%
    \setlength{\unitlength}{\svgwidth}%
  \fi%
  \global\let\svgwidth\undefined%
  \global\let\svgscale\undefined%
  \makeatother%
  \begin{picture}(1,1.42857143)%
    \put(0,0){\includegraphics[width=\unitlength,page=1]{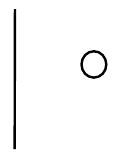}}%
    \put(0.71656417,0.90146588){\color[rgb]{0,0,0}\makebox(0,0)[lt]{\begin{minipage}{0.10782629\unitlength}\raggedright \tiny{$\alpha_K$}\end{minipage}}}%
    \put(0,0){\includegraphics[width=\unitlength,page=2]{trans2.pdf}}%
    \put(0.59259941,1.23426017){\color[rgb]{0,0,0}\makebox(0,0)[lt]{\begin{minipage}{0.10782629\unitlength}\raggedright \tiny{$\omega$}\end{minipage}}}%
    \put(0,0){\includegraphics[width=\unitlength,page=3]{trans2.pdf}}%
    \put(0.68407124,0.36106674){\color[rgb]{0,0,0}\makebox(0,0)[lt]{\begin{minipage}{0.10782629\unitlength}\raggedright \tiny{$\alpha_K$}\end{minipage}}}%
    \put(0,0){\includegraphics[width=\unitlength,page=4]{trans2.pdf}}%
    \put(0.56010648,0.69386102){\color[rgb]{0,0,0}\makebox(0,0)[lt]{\begin{minipage}{0.10782629\unitlength}\raggedright \tiny{$\omega$}\end{minipage}}}%
    \put(0,0){\includegraphics[width=\unitlength,page=5]{trans2.pdf}}%
    \put(0.17440143,0.1430314){\color[rgb]{0,0,0}\makebox(0,0)[lt]{\begin{minipage}{0.39289115\unitlength}\raggedright \tiny{$X$}\end{minipage}}}%
    \put(0,0){\includegraphics[width=\unitlength,page=6]{trans2.pdf}}%
  \end{picture}%
\endgroup%

\raisebox{15mm}{$\stackrel{(2)}{=}$}
\begingroup%
  \makeatletter%
  \providecommand\color[2][]{%
    \errmessage{(Inkscape) Color is used for the text in Inkscape, but the package 'color.sty' is not loaded}%
    \renewcommand\color[2][]{}%
  }%
  \providecommand\transparent[1]{%
    \errmessage{(Inkscape) Transparency is used (non-zero) for the text in Inkscape, but the package 'transparent.sty' is not loaded}%
    \renewcommand\transparent[1]{}%
  }%
  \providecommand\rotatebox[2]{#2}%
  \ifx\svgwidth\undefined%
    \setlength{\unitlength}{56bp}%
    \ifx\svgscale\undefined%
      \relax%
    \else%
      \setlength{\unitlength}{\unitlength * \real{\svgscale}}%
    \fi%
  \else%
    \setlength{\unitlength}{\svgwidth}%
  \fi%
  \global\let\svgwidth\undefined%
  \global\let\svgscale\undefined%
  \makeatother%
  \begin{picture}(1,1.42857143)%
    \put(0,0){\includegraphics[width=\unitlength,page=1]{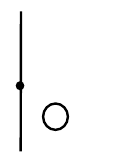}}%
    \put(0.38739575,0.45357774){\color[rgb]{0,0,0}\makebox(0,0)[lt]{\begin{minipage}{0.10782629\unitlength}\raggedright \tiny{$\alpha_K$}\end{minipage}}}%
    \put(0,0){\includegraphics[width=\unitlength,page=2]{trans3.pdf}}%
    \put(0.43981548,1.00475331){\color[rgb]{0,0,0}\makebox(0,0)[lt]{\begin{minipage}{0.10782629\unitlength}\raggedright \tiny{$\omega$}\end{minipage}}}%
    \put(0,0){\includegraphics[width=\unitlength,page=3]{trans3.pdf}}%
    \put(0.72100294,0.44246774){\color[rgb]{0,0,0}\makebox(0,0)[lt]{\begin{minipage}{0.10782629\unitlength}\raggedright \tiny{$\alpha_K$}\end{minipage}}}%
    \put(0.2259914,0.12483954){\color[rgb]{0,0,0}\makebox(0,0)[lt]{\begin{minipage}{0.39289115\unitlength}\raggedright \tiny{$X$}\end{minipage}}}%
    \put(0,0){\includegraphics[width=\unitlength,page=4]{trans3.pdf}}%
  \end{picture}%
\endgroup%

\raisebox{15mm}{$\stackrel{(3)}{=}$}
\begingroup%
  \makeatletter%
  \providecommand\color[2][]{%
    \errmessage{(Inkscape) Color is used for the text in Inkscape, but the package 'color.sty' is not loaded}%
    \renewcommand\color[2][]{}%
  }%
  \providecommand\transparent[1]{%
    \errmessage{(Inkscape) Transparency is used (non-zero) for the text in Inkscape, but the package 'transparent.sty' is not loaded}%
    \renewcommand\transparent[1]{}%
  }%
  \providecommand\rotatebox[2]{#2}%
  \ifx\svgwidth\undefined%
    \setlength{\unitlength}{56bp}%
    \ifx\svgscale\undefined%
      \relax%
    \else%
      \setlength{\unitlength}{\unitlength * \real{\svgscale}}%
    \fi%
  \else%
    \setlength{\unitlength}{\svgwidth}%
  \fi%
  \global\let\svgwidth\undefined%
  \global\let\svgscale\undefined%
  \makeatother%
  \begin{picture}(1,1.42857143)%
    \put(0,0){\includegraphics[width=\unitlength,page=1]{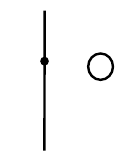}}%
    \put(0.77376692,0.88239002){\color[rgb]{0,0,0}\makebox(0,0)[lt]{\begin{minipage}{0.10782629\unitlength}\raggedright \tiny{$\alpha_K$}\end{minipage}}}%
    \put(0,0){\includegraphics[width=\unitlength,page=2]{trans4.pdf}}%
    \put(0.64980221,1.21518331){\color[rgb]{0,0,0}\makebox(0,0)[lt]{\begin{minipage}{0.10782629\unitlength}\raggedright \tiny{$\omega$}\end{minipage}}}%
    \put(0,0){\includegraphics[width=\unitlength,page=3]{trans4.pdf}}%
    \put(0.4292328,0.13037954){\color[rgb]{0,0,0}\makebox(0,0)[lt]{\begin{minipage}{0.39289115\unitlength}\raggedright \tiny{$X$}\end{minipage}}}%
    \put(0,0){\includegraphics[width=\unitlength,page=4]{trans4.pdf}}%
    \put(0.05868438,0.52621059){\color[rgb]{0,0,0}\makebox(0,0)[lt]{\begin{minipage}{0.10782629\unitlength}\raggedright \tiny{$\alpha_K$}\end{minipage}}}%
    \put(0,0){\includegraphics[width=\unitlength,page=5]{trans4.pdf}}%
    \put(0.11212868,0.94263102){\color[rgb]{0,0,0}\makebox(0,0)[lt]{\begin{minipage}{0.14568266\unitlength}\raggedright \tiny{$\varepsilon$}\end{minipage}}}%
    \put(0,0){\includegraphics[width=\unitlength,page=6]{trans4.pdf}}%
  \end{picture}%
\endgroup%

\raisebox{15mm}{$\stackrel{(4)}{=}$}
\begingroup%
  \makeatletter%
  \providecommand\color[2][]{%
    \errmessage{(Inkscape) Color is used for the text in Inkscape, but the package 'color.sty' is not loaded}%
    \renewcommand\color[2][]{}%
  }%
  \providecommand\transparent[1]{%
    \errmessage{(Inkscape) Transparency is used (non-zero) for the text in Inkscape, but the package 'transparent.sty' is not loaded}%
    \renewcommand\transparent[1]{}%
  }%
  \providecommand\rotatebox[2]{#2}%
  \ifx\svgwidth\undefined%
    \setlength{\unitlength}{56bp}%
    \ifx\svgscale\undefined%
      \relax%
    \else%
      \setlength{\unitlength}{\unitlength * \real{\svgscale}}%
    \fi%
  \else%
    \setlength{\unitlength}{\svgwidth}%
  \fi%
  \global\let\svgwidth\undefined%
  \global\let\svgscale\undefined%
  \makeatother%
  \begin{picture}(1,1.42857143)%
    \put(0,0){\includegraphics[width=\unitlength,page=1]{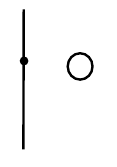}}%
    \put(0.5983683,0.88437288){\color[rgb]{0,0,0}\makebox(0,0)[lt]{\begin{minipage}{0.10782629\unitlength}\raggedright \tiny{$\alpha_K$}\end{minipage}}}%
    \put(0,0){\includegraphics[width=\unitlength,page=2]{trans5.pdf}}%
    \put(0.47440354,1.21716674){\color[rgb]{0,0,0}\makebox(0,0)[lt]{\begin{minipage}{0.10782629\unitlength}\raggedright \tiny{$\omega$}\end{minipage}}}%
    \put(0,0){\includegraphics[width=\unitlength,page=3]{trans5.pdf}}%
    \put(0.24938044,0.14126997){\color[rgb]{0,0,0}\makebox(0,0)[lt]{\begin{minipage}{0.39289115\unitlength}\raggedright \tiny{$X$}\end{minipage}}}%
  \end{picture}%
\endgroup%

\end{center}

Equality $(2)$ uses the axiomatic of the Hopf pairing $\omega$ (see Section 3A in \cite{BrugVircenter}), equality $(3)$ is true because $\alpha_K$ is an integral (see proof of Proposition~\ref{kirbyadmissible}) and equality $(4)$ holds because $\alpha_K$ is an admissible element so the $i-$part of the Lemma is proved.

Now, assume that $X$ is a transparent object. Then for all objects $Y$,
\begin{center}
\resizebox{!}{0.1\textheight}{
\begingroup%
  \makeatletter%
  \providecommand\color[2][]{%
    \errmessage{(Inkscape) Color is used for the text in Inkscape, but the package 'color.sty' is not loaded}%
    \renewcommand\color[2][]{}%
  }%
  \providecommand\transparent[1]{%
    \errmessage{(Inkscape) Transparency is used (non-zero) for the text in Inkscape, but the package 'transparent.sty' is not loaded}%
    \renewcommand\transparent[1]{}%
  }%
  \providecommand\rotatebox[2]{#2}%
  \ifx\svgwidth\undefined%
    \setlength{\unitlength}{56bp}%
    \ifx\svgscale\undefined%
      \relax%
    \else%
      \setlength{\unitlength}{\unitlength * \real{\svgscale}}%
    \fi%
  \else%
    \setlength{\unitlength}{\svgwidth}%
  \fi%
  \global\let\svgwidth\undefined%
  \global\let\svgscale\undefined%
  \makeatother%
  \begin{picture}(1,1.42857143)%
    \put(0,0){\includegraphics[width=\unitlength,page=1]{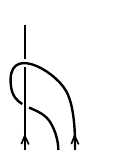}}%
    \put(0.19419065,0.11050182){\color[rgb]{0,0,0}\makebox(0,0)[lt]{\begin{minipage}{0.40525775\unitlength}\raggedright \tiny{$X$}\end{minipage}}}%
    \put(0.62679655,0.10995754){\color[rgb]{0,0,0}\makebox(0,0)[lt]{\begin{minipage}{0.40525775\unitlength}\raggedright \tiny{$Y$}\end{minipage}}}%
  \end{picture}%
\endgroup%
\raisebox{15mm}{$=$}
\begingroup%
  \makeatletter%
  \providecommand\color[2][]{%
    \errmessage{(Inkscape) Color is used for the text in Inkscape, but the package 'color.sty' is not loaded}%
    \renewcommand\color[2][]{}%
  }%
  \providecommand\transparent[1]{%
    \errmessage{(Inkscape) Transparency is used (non-zero) for the text in Inkscape, but the package 'transparent.sty' is not loaded}%
    \renewcommand\transparent[1]{}%
  }%
  \providecommand\rotatebox[2]{#2}%
  \ifx\svgwidth\undefined%
    \setlength{\unitlength}{56bp}%
    \ifx\svgscale\undefined%
      \relax%
    \else%
      \setlength{\unitlength}{\unitlength * \real{\svgscale}}%
    \fi%
  \else%
    \setlength{\unitlength}{\svgwidth}%
  \fi%
  \global\let\svgwidth\undefined%
  \global\let\svgscale\undefined%
  \makeatother%
  \begin{picture}(1,1.42857143)%
    \put(0,0){\includegraphics[width=\unitlength,page=1]{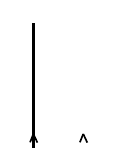}}%
    \put(0.26809264,0.12573325){\color[rgb]{0,0,0}\makebox(0,0)[lt]{\begin{minipage}{0.40525775\unitlength}\raggedright \tiny{$X$}\end{minipage}}}%
    \put(0.70069849,0.12518897){\color[rgb]{0,0,0}\makebox(0,0)[lt]{\begin{minipage}{0.40525775\unitlength}\raggedright \tiny{$Y$}\end{minipage}}}%
    \put(0,0){\includegraphics[width=\unitlength,page=2]{trans7.pdf}}%
  \end{picture}%
\endgroup%
}
\end{center}
so $\Pi_{\alpha_K,X}=\id_X\otimes \varepsilon\alpha_K=\id_X$.

Suppose now $X$ is such that $\Pi_{\alpha_K,X}=\mathrm{id}_X$. For all objects $Y$,
\begin{center}
    \resizebox{!}{0.1\textheight}{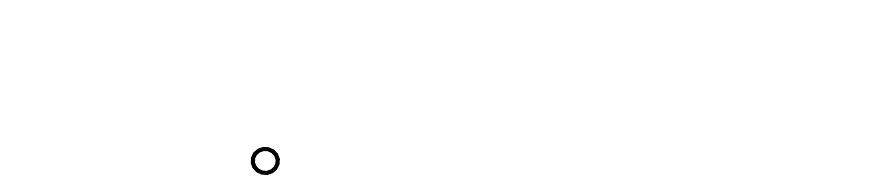}
\end{center}
First and last equalities are due to the statement on the idempotent. Central drawing is a topologival representation of the idempotent: the braid colored by $X$ is encircled by a component colored by $\alpha_K$ and as $\alpha_K$ is admissible, the braid $Y$ can slide along the component colored by $\alpha_K$ so $X$ is transparent.
\end{proof}

\subsection{The modular case : on the Reshetikhin-Turaev TQFTs}
Let $\C$ be a modular category. We have the following comparison result between the TQFT $\mr{RT}_{\C}$\index{$\mr{RT}_{\C}$} of Reshetikhin-Turaev and our TQFT.

\begin{theorem}\label{theoremrt}
Let $\C$ be a modular category and denote by $\alpha_K$ its Kirby element and by $\T$ its subcategory of transparent objects. Then, up to normalization, the following diagram commutes
\begin{align}
\begin{tikzpicture}[description/.style={fill=white,inner sep=2pt},baseline=(current  bounding  box.center)]
\matrix (m) [ampersand replacement=\&,matrix of math nodes, row sep=3em,column sep=2.5em, text height=1.5ex, text depth=0.25ex]
{ \Cob^{p} \& \& \T \\
\& \& \vect \\ };
\path[->,font=\scriptsize]
(m-1-1) edge node[auto=left] {$\mr{V}_{\C,\alpha_{K}}$} (m-1-3)
(m-1-1) edge node[auto=right] {$\mr{RT}_{\C}$} (m-2-3)
(m-1-3) edge node[auto] {$\mr{Hom}_{\C}(\mathds{1},-)$} (m-2-3);
\end{tikzpicture}
\end{align}
where $\vect$ is the category of finite dimensional $\kk$-vector spaces.
\end{theorem}

Before proving this theorem, we need some tools. 
First, recall that the subcategory of transparent objects of a modular category is identified with the category $\vect$ of finite dimensional vector spaces. If $\C$ is a $\kk$-fusion category and $X$ is an object of $\C$, we denote by $<X>$ the smallest monoidal rigid subcategory of $\C$ containing $X$ and stable under direct sums and direct factors. The subcategory $<X>$ is a fusion subcategory of $\C$ such that the simple objects are the direct factors of tensorial products of $X$ and $X^{*}$. 

\begin{lemma}
Let $\C$ be a modular category.
The subcategory $<\un>$ of transparent objects of $\C$ is monoidally equivalent to the category $\vect$.  
\end{lemma}

\begin{proof}
In a modular category, the only simple and transparent object is the monoidal unit $\mathds{1}$. Indeed, if we denote by $S$ a simple and transparent object of $\C$, it satisfies for all simple objects $X$ of $\C$, 
\[\mathrm{tr}_{q}(\tau_{S,X}\tau_{X,S})=\mathrm{dim}_{q}(S)\mathrm{dim}_{q}(X)\]
so the line of the object $S$ in the $S$-matrix is colinear to the line of the object $\mathds{1}$. As the $S$-matrix is invertible, we have \[S=\mathds{1}.\]
Then every transparent object of a modular category is a direct sum of copies of $\mathds{1}$ and the functor defined by $\kk \in \vect\mapsto \mathds{1}\in <\mathds{1}>$ is a $\kk$-linear monoidal equivalence of category.
\end{proof}

Secondly, recall that the cylinder $\Sigma_{g}\times [0,1]$ on a surface of genus $g$ is represented by the $(g,g,g)$-cobordism tangle $T_g$ drawn on Figure~\ref{cylinder}.

\begin{lemma}\label{lemmaunit}
Suppose that $\C$ is a modular category and $\alpha_{K}$ its Kirby element.
Then \[|T_g|_{\C,\alpha_K}=\omega(\alpha_K\otimes \alpha_K)^{g}\mr{id}_{C^{\otimes g}}.\] 
\end{lemma}
\begin{proof}
Let us show that $|T_1|_{\C,\alpha_K}=\omega({\alpha_K\otimes \alpha_K})\id_{C}$; the general result comes from tensorization.
Remark that
\begin{center}
\resizebox{!}{0.1\textheight}{\raisebox{14mm}{$T_1=$}
\begingroup%
  \makeatletter%
  \providecommand\color[2][]{%
    \errmessage{(Inkscape) Color is used for the text in Inkscape, but the package 'color.sty' is not loaded}%
    \renewcommand\color[2][]{}%
  }%
  \providecommand\transparent[1]{%
    \errmessage{(Inkscape) Transparency is used (non-zero) for the text in Inkscape, but the package 'transparent.sty' is not loaded}%
    \renewcommand\transparent[1]{}%
  }%
  \providecommand\rotatebox[2]{#2}%
  \ifx\svgwidth\undefined%
    \setlength{\unitlength}{160bp}%
    \ifx\svgscale\undefined%
      \relax%
    \else%
      \setlength{\unitlength}{\unitlength * \real{\svgscale}}%
    \fi%
  \else%
    \setlength{\unitlength}{\svgwidth}%
  \fi%
  \global\let\svgwidth\undefined%
  \global\let\svgscale\undefined%
  \makeatother%
  \begin{picture}(1,0.55)%
    \put(0,0){\includegraphics[width=\unitlength,page=1]{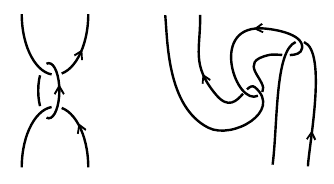}}%
    \put(0.33248259,0.27544457){\color[rgb]{0,0,0}\makebox(0,0)[lt]{\begin{minipage}{0.2465971\unitlength}\raggedright $=$\end{minipage}}}%
  \end{picture}%
\endgroup%
}
\end{center}
We choose the following opentangle $O$ coming from $T_1$:
\begin{center}
\resizebox{!}{0.1\textheight}{\raisebox{20mm}{$O=$}{
\begingroup%
  \makeatletter%
  \providecommand\color[2][]{%
    \errmessage{(Inkscape) Color is used for the text in Inkscape, but the package 'color.sty' is not loaded}%
    \renewcommand\color[2][]{}%
  }%
  \providecommand\transparent[1]{%
    \errmessage{(Inkscape) Transparency is used (non-zero) for the text in Inkscape, but the package 'transparent.sty' is not loaded}%
    \renewcommand\transparent[1]{}%
  }%
  \providecommand\rotatebox[2]{#2}%
  \ifx\svgwidth\undefined%
    \setlength{\unitlength}{160bp}%
    \ifx\svgscale\undefined%
      \relax%
    \else%
      \setlength{\unitlength}{\unitlength * \real{\svgscale}}%
    \fi%
  \else%
    \setlength{\unitlength}{\svgwidth}%
  \fi%
  \global\let\svgwidth\undefined%
  \global\let\svgscale\undefined%
  \makeatother%
  \begin{picture}(1,0.75)%
    \put(0,0){\includegraphics[width=\unitlength,page=1]{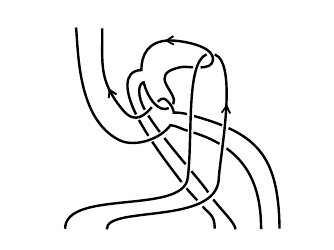}}%
  \end{picture}%
\endgroup%
}}
\end{center}
Then compute $|O|_{\C}$. For $X_1, X_2, X_3$ objects of $\C$,
\begin{center}
\resizebox{!}{0.15\textheight}{\raisebox{21mm}{$\iota_{X_3}O_{X_1,X_2,X_3}=$}{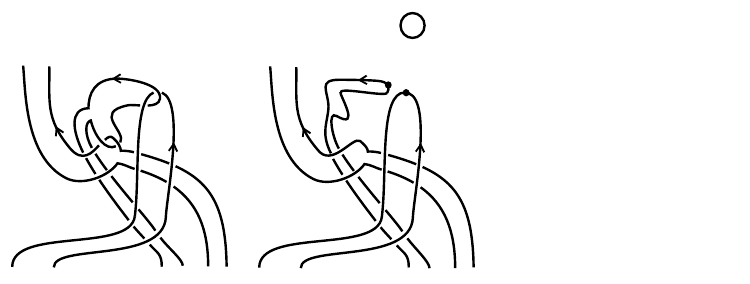}}
\end{center}
so the invariant is
\begin{center}
$|T_1|_{\C,\alpha_K}=|O|_{\C}(\id_{C}\otimes \alpha_K\otimes \alpha_K)=$\resizebox{!}{0.1\textheight}{\raisebox{-20mm}{
\begingroup%
  \makeatletter%
  \providecommand\color[2][]{%
    \errmessage{(Inkscape) Color is used for the text in Inkscape, but the package 'color.sty' is not loaded}%
    \renewcommand\color[2][]{}%
  }%
  \providecommand\transparent[1]{%
    \errmessage{(Inkscape) Transparency is used (non-zero) for the text in Inkscape, but the package 'transparent.sty' is not loaded}%
    \renewcommand\transparent[1]{}%
  }%
  \providecommand\rotatebox[2]{#2}%
  \ifx\svgwidth\undefined%
    \setlength{\unitlength}{160bp}%
    \ifx\svgscale\undefined%
      \relax%
    \else%
      \setlength{\unitlength}{\unitlength * \real{\svgscale}}%
    \fi%
  \else%
    \setlength{\unitlength}{\svgwidth}%
  \fi%
  \global\let\svgwidth\undefined%
  \global\let\svgscale\undefined%
  \makeatother%
  \begin{picture}(1,1)%
    \put(0,0){\includegraphics[width=\unitlength,page=1]{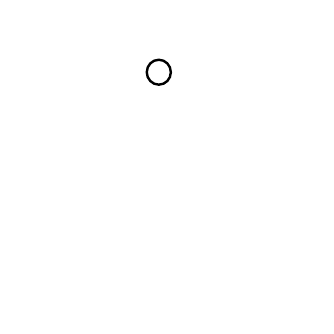}}%
    \put(0.4598393,0.79517764){\color[rgb]{0,0,0}\makebox(0,0)[lt]{\begin{minipage}{0.07380238\unitlength}\raggedright \small{$\omega$}\end{minipage}}}%
    \put(0.70762211,0.48354999){\color[rgb]{0,0,0}\makebox(0,0)[lt]{\begin{minipage}{0.07909745\unitlength}\raggedright \small{$S$}\end{minipage}}}%
    \put(0,0){\includegraphics[width=\unitlength,page=2]{invt1.pdf}}%
    \put(0.46418186,0.60651304){\color[rgb]{0,0,0}\makebox(0,0)[lt]{\begin{minipage}{0.07380238\unitlength}\raggedright \small{$\omega$}\end{minipage}}}%
    \put(0.17410289,0.92669392){\color[rgb]{0,0,0}\makebox(0,0)[lt]{\begin{minipage}{0.20406365\unitlength}\raggedright \small{$C$}\end{minipage}}}%
    \put(0,0){\includegraphics[width=\unitlength,page=3]{invt1.pdf}}%
    \put(0.35051696,0.33727329){\color[rgb]{0,0,0}\makebox(0,0)[lt]{\begin{minipage}{0.10661737\unitlength}\raggedright \tiny{$\alpha_K$}\end{minipage}}}%
    \put(0,0){\includegraphics[width=\unitlength,page=4]{invt1.pdf}}%
    \put(0.76225496,0.07647447){\color[rgb]{0,0,0}\makebox(0,0)[lt]{\begin{minipage}{0.20406365\unitlength}\raggedright \small{$C$}\end{minipage}}}%
    \put(0.34867269,0.12584629){\color[rgb]{0,0,0}\makebox(0,0)[lt]{\begin{minipage}{0.10587358\unitlength}\raggedright \tiny{$\alpha_K$}\end{minipage}}}%
  \end{picture}%
\endgroup%
}}
\end{center}
As $\omega$ is nondegenerate and $\alpha_K$ is an integral of $C$, as shown in \cite{BrugVircenter}, Lemma 3.1, we have the following identity:
\begin{center}
\resizebox{!}{0.2\textheight}{
\begingroup%
  \makeatletter%
  \providecommand\color[2][]{%
    \errmessage{(Inkscape) Color is used for the text in Inkscape, but the package 'color.sty' is not loaded}%
    \renewcommand\color[2][]{}%
  }%
  \providecommand\transparent[1]{%
    \errmessage{(Inkscape) Transparency is used (non-zero) for the text in Inkscape, but the package 'transparent.sty' is not loaded}%
    \renewcommand\transparent[1]{}%
  }%
  \providecommand\rotatebox[2]{#2}%
  \ifx\svgwidth\undefined%
    \setlength{\unitlength}{160bp}%
    \ifx\svgscale\undefined%
      \relax%
    \else%
      \setlength{\unitlength}{\unitlength * \real{\svgscale}}%
    \fi%
  \else%
    \setlength{\unitlength}{\svgwidth}%
  \fi%
  \global\let\svgwidth\undefined%
  \global\let\svgscale\undefined%
  \makeatother%
  \begin{picture}(1,1)%
    \put(0,0){\includegraphics[width=\unitlength,page=1]{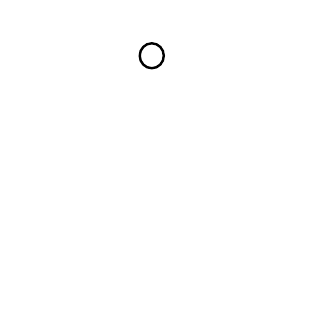}}%
    \put(0.4387145,0.84455279){\color[rgb]{0,0,0}\makebox(0,0)[lt]{\begin{minipage}{0.07380238\unitlength}\raggedright \small{$\omega$}\end{minipage}}}%
    \put(0.08998834,0.83463819){\color[rgb]{0,0,0}\makebox(0,0)[lt]{\begin{minipage}{0.07909745\unitlength}\raggedright \small{$S$}\end{minipage}}}%
    \put(0,0){\includegraphics[width=\unitlength,page=2]{invt2.pdf}}%
    \put(0.4430571,0.65588819){\color[rgb]{0,0,0}\makebox(0,0)[lt]{\begin{minipage}{0.07380238\unitlength}\raggedright \small{$\omega$}\end{minipage}}}%
    \put(0.15297813,0.97606907){\color[rgb]{0,0,0}\makebox(0,0)[lt]{\begin{minipage}{0.20406365\unitlength}\raggedright \small{$C$}\end{minipage}}}%
    \put(0,0){\includegraphics[width=\unitlength,page=3]{invt2.pdf}}%
    \put(0.32939218,0.38664844){\color[rgb]{0,0,0}\makebox(0,0)[lt]{\begin{minipage}{0.10661737\unitlength}\raggedright \tiny{$\alpha_K$}\end{minipage}}}%
    \put(0,0){\includegraphics[width=\unitlength,page=4]{invt2.pdf}}%
    \put(0.74113016,0.12584992){\color[rgb]{0,0,0}\makebox(0,0)[lt]{\begin{minipage}{0.20406365\unitlength}\raggedright \small{$C$}\end{minipage}}}%
    \put(0.32754793,0.17522174){\color[rgb]{0,0,0}\makebox(0,0)[lt]{\begin{minipage}{0.10587358\unitlength}\raggedright \tiny{$\alpha_K$}\end{minipage}}}%
    \put(0,0){\includegraphics[width=\unitlength,page=5]{invt2.pdf}}%
  \end{picture}%
\endgroup%
}\raisebox{20mm}{$=\omega(\alpha_K\otimes \alpha_K)\id_{C}$}
\end{center}
so $|T_1|_{\C,\alpha}=\omega(\alpha_K\otimes \alpha_K)\id_{C}$.
\end{proof}
Thus, in our construction, if $\C$ is modular and $\alpha_K$ is the Kirby element of $\C$:
\[W_{\C,\alpha_K}(\Sigma_{g}\times [0,1])=\id_{(C^{\otimes g}){\alpha_K}}.\]
Then, in this case, the functor with anomaly $\ww$ is already unitalized and $\W_{\C,\alpha_k}(\phantom{o})=\V_{\C,\alpha_k}(\phantom{o})$. Now we are ready to prove Theorem~\ref{theoremrt}.

\begin{proof}
\textbullet \; On a connected surface $\Sigma_{g}$ (we forget the parametrization of the surface) of genus $g$, 
$$\mr{RT}_{\C}(\Sigma_{g})=\mr{Hom}_{\C}\left( \mathds{1}, \bigoplus_{(\lambda_{1},\ldots,\lambda_{g})\in \Lambda_{\C}^{g}}\lambda_1^{*}\otimes \lambda_1\otimes \ldots\otimes \lambda_g^{*}\otimes \lambda_g\right) $$
and 
$$ \V_{\C,\alpha_{K}}(\Sigma_g)=\left(\left( \bigoplus_{\lambda\in \Lambda_{\C}}\lambda^{*}\otimes \lambda\right)^{\otimes g}\right)_{\alpha_K}=\left(\left( \bigoplus_{\lambda\in \Lambda_{\C}}\lambda^{*}\otimes \lambda\right)^{\otimes g}\right)_{\T} $$ according to Lemmas~\ref{lemmaunit} and \ref{lemmaprotran}.
Moreover, for every object $X \in \C$, we have $ \mr{Hom}_{\C}(\mathds{1},X)=\mr{Hom}_{\C}(\mathds{1},X_{\T})$ by Schur lemma. Then we have the result on surfaces.

\noindent \textbullet \;  We show the result now on a connected cobordism $M\colon \Sigma \rightarrow \Sigma'$. For simplicity of notations, assume that $\Sigma$ and $\Sigma'$ have genus $1$. Let $T$ be a $(1,n,1)$-cobordism tangle representing $M$ and denote by $L_1,\ldots, L_n$ the $n$ surgery components of $T$. Colore by $k\in \Lambda_{\C}$ and by $l\in\Lambda_{\C}$ the boundary components of $T$ corresponding respectively to $\Sigma$ and $\Sigma'$ defining a morphism $T_{k}^{l}\colon k^{*}\otimes k\rightarrow l^{*}\otimes l$. Let $c \colon \{L_1,\ldots, L_n\} \rightarrow \Lambda_{\C}$ and denote by $F$ the Shum-Turaev functor from colored ribbon tangles to $\C$ (see~
\cite{Shu94} and \cite{Turaev1994}) and by $e_{\lambda}=\iota_{\lambda}\widetilde{\mr{coev}}_{\lambda}$. Choose an opentangle $O$ such that  $T=U(O)$. Then we have: 

$$ F(T_{k}^{l},c)=p_l\circ |O|_{\C}\circ (\iota_k\otimes e_{c(1)}\otimes \ldots\otimes e_{c(n)}\otimes e_{l})$$
where for all $\lambda\in \Lambda_{\C}$, $p_{\lambda}\colon C \rightarrow  l^{*}\otimes l$ is such that $\id_{C}=\sum_{\lambda \in \Lambda_{\C}}\iota_{\lambda}p_{\lambda}$ and $p_{\lambda}\iota_{\mu}=\delta_{\lambda,\mu}\id_{\lambda^{*}\otimes \lambda}$.
Then
$$ \frac {1}{\mr{dim}(\C)^{n}}\sum_{c} \mr{dim}_{q}(c)F(T_{k}^{l},c)=p_l\circ |O|_{\C}\circ (\iota_k\otimes \alpha_{K}^{\otimes n}\otimes e_{l})$$ 
where $\mr{dim}_{q}(c):=\prod_{i=1}^{n}\mr{dim}_{q}(c(i))$
and so
$$ \frac{\mr{dim}_{q}(l)}{\dimm(\C)^{n}}\sum_{c} \mr{dim}_{q}(c)F(T_{k}^{l},c)=p_l\circ |O|_{\C}\circ (\iota_k\otimes \alpha_{K}^{\otimes n}\otimes \mr{dim}_{q}(l)e_{l}) \stackrel{(2)}{=} p_l\circ |O|_{\C}\circ (\iota_k\otimes \alpha_{K}^{\otimes n}\otimes \alpha_{K})$$

The last equality $(2)$ holds because:

\begin{align*}\dimm(\C)^n(p_l\circ |O|_{\C}\circ (\iota_k\otimes \alpha_{K}^{\otimes n}\otimes \alpha_{RT}))&=
\sum_{\gamma \in \Lambda_{\C}} \mr{dim}_{q}(\gamma) p_l\circ |O|_{\C}\circ (\iota_k\otimes \alpha_{K}^{\otimes n}\otimes e_{\gamma})\\
&=\sum_{\gamma \in \Lambda_{\C}} \mr{dim}_{q}(\gamma) p_l\circ \iota_{\gamma} p_{\gamma}|O|_{\C} (\iota_k\otimes \alpha_{K}^{\otimes n}\otimes e_{\gamma})
\\&=\sum_{\gamma \in \Lambda} \mr{dim}_{q}(\gamma) \delta_{l,\gamma}\id_{l^{*}\otimes l}\circ p_{\gamma}|O|_{\C} (\iota_k\otimes \alpha_{K}^{\otimes n}\otimes e_{\gamma})
\\&=\mr{dim}_{q}(l)p_{l}\circ |O|_{\C}\circ (\iota_k\otimes \alpha_{K}^{\otimes n}\otimes e_{\gamma})
\end{align*}

Thus we have shown that
\[ \frac{\mr{dim}_{q}(l)}{\dimm(\C)^{n}}\sum_{c} \mr{dim}_{q}(c)F(T_{k}^{l},c)= p_l\circ |T|_{\C,\alpha}\circ \iota_{k}.\]

The default of normalization between the two TQFTs is given by
$\D^{-b_{0}(T)-g'-2n}$,
where $\D$ is a square root of $\dimm(\C)$, $b_{0}(T)$ is the number of null eigenvalues of the linking matrix of the surgery components of $T$ and $g'$ is the genus of exit boundary $\Sigma'$ of $M$.

Adding projections and injections on transparent part, we recover exactly \[\mr{RT}_{\C}(M)=\D^{-b_{0}(T)-g'-2n}\mr{Hom}_{\C}(\un,-)\V_{\C,\alpha_K}(\phantom{O}).\]
\end{proof}

\subsection{Functoriality of the construction}

Let $\C$ and $\D$ be two ribbon categories with coend respectively denoted by $(C,\iota)$ and $(D,j)$.
Let $\alpha \colon \mathds{1}\rightarrow C$ and $\beta \colon \mathds{1}\rightarrow D$.
Suppose that $F\colon \C \rightarrow \D$ is a strong monoidal functor which is ribbon such that $(F(C),F(\iota))$ is the coend of the functor $F\colon \C^{op}\otimes \C \rightarrow \D$ defined by 
\begin{align}
F(X\otimes Y)=F(X^{*}\otimes Y)\;
\text{ and }\;
F(f,g)=F(f^{*}\otimes g)
\end{align}

As $F$ is a strong monoidal functor, we have a natural isomorphism $F(X^{*}\otimes Y)\simeq F(X)^{*}\otimes F(Y)$. Then, by the factorization property of the coend $(F(C),F(\iota))$, there exists a unique morphism $\zeta \colon F(C) \rightarrow D$ of $\D$\index{$\zeta$} such that, for every object $X$ of $\C$, the following \ref{imagecoendiscoend} diagram commutes:

\begin{align}\label{imagecoendiscoend}
\begin{tikzpicture}[description/.style={fill=white,inner sep=2pt},baseline=(current  bounding  box.center)]
\matrix (m) [ampersand replacement=\&,matrix of math nodes, row sep=3em,column sep=2.5em, text height=1.5ex, text depth=0.25ex]
{ \& F(X)^{*}\otimes F(X)\& \\
 F(C)\&  \& D \\ };
\path[->,font=\scriptsize]
(m-1-2) edge node[auto=right] {$F(\iota_X)$} (m-2-1)
		edge node[auto=left] {$j_{F(X)}$} (m-2-3)
(m-2-1) edge node[auto=right] {$\zeta $} (m-2-3);
\end{tikzpicture}
\end{align}

If $\alpha \in \mr{Hom}_{\C}(\mathds{1},C)$, denote by $F_{!}\alpha=\zeta F(\alpha)\colon F({\mathds{1}})\simeq \mathds{1}\rightarrow D$. 

\begin{lemma}\label{lemmamainchange}
Let $T$ be a $(\ung,n,\unh)$-cobordism and $\alpha \in \mr{Hom}_{\C}(\un,C)$. Then the following diagram commutes:
\begin{center}
\begin{tikzpicture}[description/.style={fill=white,inner sep=2pt},baseline=(current  bounding  box.center)]
\matrix (m) [ampersand replacement=\&,matrix of math nodes, row sep=3em,column sep=2.5em, text height=1.5ex, text depth=0.25ex]
{ F(C^{\otimes |\ung|})\& \& D^{\otimes |\ung|} \\
 F(C^{\otimes |\unh|})\&  \& D^{\otimes |\unh|} \\ };
\path[->,font=\scriptsize]
(m-1-1) edge node[auto=left] {$\zeta^{\otimes |\ung|}$} (m-1-3)
		edge node[auto=right] {$F(|T|_{\C,\alpha})$} (m-2-1)
(m-2-1) edge node[auto=right] {$\zeta^{\otimes |\unh|}$} (m-2-3)
(m-1-3)	edge node[auto=left] {$|T|_{\D,F_{!}\alpha}$} (m-2-3);
\end{tikzpicture}
\end{center}
\end{lemma}
\begin{proof}
For simplicity of notations, assume that $T$ is a $(g,n,h)$-cobordism tangle where $g$ and $h$ are integers. The case where $T$ is a $(\ung,n,\unh)$-cobordism tangle with multigenus $\ung$ and $\unh$ is similar. 
Let $O$ be a $(\ung,n,\unh)$-opentangle such that $U(O)=T$.
Colore the components of the opentangle $O$ by objects of $\C$: the entrance components are colored by $X_1,\ldots,X_{g}$,  the surgery components are colored by $X_{g+1},\ldots, X_{g+n}$, the exit components are colored by $X_{g+n+1},\ldots,X_{g+n+h}$.

Remark that, since $F$ is ribbon, that:
\begin{align*}
F((\iota_{X_{g+n+1}}\otimes\ldots\otimes \iota_{X_{g+n+h}})\circ O_{X_1,\ldots X_{g+n+h}})
&=F(\iota_{X_{g+n+1}}\otimes\ldots\otimes \iota_{X_{g+n+h}})F(O_{X_1,\ldots X_{g+n+h}})\\
&=F(\iota_{X_{g+n+1}}\otimes \ldots\otimes \iota_{X_{g+n+h}})O_{F(X_1),\ldots,F(X_{g+n+h})}
\end{align*}

And, as $(\iota_{X_{g+n+1}}\otimes\ldots\otimes \iota_{X_{g+n+h}})\circ O_{X_1,\ldots X_{g+n+h}}=
|O|_{\C}(\iota_{X_1}\otimes \ldots\otimes \iota_{X_{g+n+h}})$,

\[F(|O|_{\C}(\iota_{X_1}\otimes \ldots\otimes \iota_{X_{g+n+h}}))=F(\iota_{X_{g+n+1}}\otimes \ldots\otimes \iota_{X_{g+n+h}})O_{F(X_1),\ldots,F(X_{g+n+h})}.\]

so, multiplying by $\zeta^{\otimes h}$,
\begin{align*}
&\zeta^{\otimes h}F(|O|_{\C}(\iota_{X_1}\otimes \ldots\otimes \iota_{X_{g+n+h}}))\\
&=(j_{F(X_{g+n+1})}\otimes \ldots\otimes j_{F(X_{g+n+h})})F(\iota_{X_{g+n+1}}\otimes \ldots\otimes \iota_{X_{g+n+h}})O_{F(X_1),\ldots,F(X_{g+n+h})}\\
&=|O|_{\D}(j_{F(X_1)}\otimes\ldots \otimes j_{F(X_{g+n+h})})\\
&=|O|_{\D}\zeta^{\otimes g+n+h}(F(\iota_{X_1})\otimes \ldots \otimes F(\iota_{X_{g+n+h}})).
\end{align*}

We have \[\zeta^{\otimes h}F(|O|_{\C})(F(\iota_{X_1})\otimes \ldots\otimes F(\iota_{X_{g+n+h}}))=|O|_{\D}\zeta^{\otimes g+n+h}(F(\iota_{X_1})\otimes \ldots \otimes F(\iota_{X_{g+n+h}}))\] so
as $(F(C),F(\iota))$ is a coend,
\[\zeta^{\otimes h}F(|O|_{\C}=|O|_{\D}\zeta^{\otimes g+n+h}\] thus
\[\zeta^{\otimes h}F(|O|_{\C})F(\id_{C^{\otimes g}}\otimes \alpha^{\otimes n+h})=|O|_{\D}\zeta^{\otimes g+n+h}F(\id_{C^{\otimes g}}\otimes \alpha^{\otimes n+h})\] so

\[\zeta^{\otimes h}F(|T|_{\C,\alpha})=F(|T|_{\D,F_{!}\alpha})\zeta^{\otimes g}.\]
\end{proof}

If $\alpha$ is an admissible element of $\C$ and $g$ is a positive integer, then remember the following idempotent of $\C$
\begin{center}
\resizebox{!}{0.15\textheight}{\raisebox{30mm}{$\Pi_{\alpha,g}^{\C}=$}
\begingroup%
  \makeatletter%
  \providecommand\color[2][]{%
    \errmessage{(Inkscape) Color is used for the text in Inkscape, but the package 'color.sty' is not loaded}%
    \renewcommand\color[2][]{}%
  }%
  \providecommand\transparent[1]{%
    \errmessage{(Inkscape) Transparency is used (non-zero) for the text in Inkscape, but the package 'transparent.sty' is not loaded}%
    \renewcommand\transparent[1]{}%
  }%
  \providecommand\rotatebox[2]{#2}%
  \ifx\svgwidth\undefined%
    \setlength{\unitlength}{160bp}%
    \ifx\svgscale\undefined%
      \relax%
    \else%
      \setlength{\unitlength}{\unitlength * \real{\svgscale}}%
    \fi%
  \else%
    \setlength{\unitlength}{\svgwidth}%
  \fi%
  \global\let\svgwidth\undefined%
  \global\let\svgscale\undefined%
  \makeatother%
  \begin{picture}(1,1)%
    \put(0,0){\includegraphics[width=\unitlength]{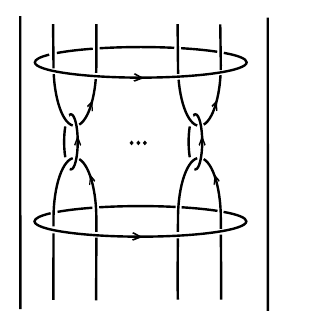}}%
    \put(0.8182777,0.08292474){\color[rgb]{0,0,0}\makebox(0,0)[lt]{\begin{minipage}{0.21031001\unitlength}\raggedright $_{\C,\alpha}$\end{minipage}}}%
    \put(0.16236975,0.06631853){\color[rgb]{0,0,0}\makebox(0,0)[lt]{\begin{minipage}{0.54838828\unitlength}\raggedright $\underbrace{\phantom{ooooooooooooooooo}}_{g}$\end{minipage}}}%
  \end{picture}%
\endgroup%
\raisebox{30mm}{,}}
\end{center}
and  if $\ung=(g_1,\ldots,g_r)$ is a $r$-tuple of positive integers, denote by $\Pi_{\alpha,\ung}^{\C}=\Pi_{\alpha,g_1}^{\C}\otimes \ldots\otimes \Pi_{\alpha,g_r}^{\C}$.

\begin{lemma}\label{lemmapipi}
Let $\alpha\colon \un \rightarrow C$ be an admissible element of $\C$ and suppose that $F_{!}\alpha$ is an admissible element of $\D$. Then the following diagram commutes:
\begin{center}
\begin{tikzpicture}[description/.style={fill=white,inner sep=2pt},baseline=(current  bounding  box.center)]
\matrix (m) [ampersand replacement=\&,matrix of math nodes, row sep=3em,column sep=2.5em, text height=1.5ex, text depth=0.25ex]
{ F(C^{\otimes |\ung|})\& \& D^{\otimes |\ung|} \\
 F(C^{\otimes |\ung|})\&  \& D^{\otimes |\ung|} \\ };
\path[->,font=\scriptsize]
(m-1-1) edge node[auto=left] {$\zeta^{\otimes |\ung|}$} (m-1-3)
		edge node[auto=right] {$F(\Pi_{\alpha,\ung}^{\C})$} (m-2-1)
(m-2-1) edge node[auto=right] {$\zeta^{\otimes |\ung|}$} (m-2-3)
(m-1-3)	edge node[auto=left] {$\Pi_{F_{!}\alpha,\ung}^{\D}$} (m-2-3);
\end{tikzpicture}
\end{center}
\end{lemma}
\begin{proof}
The result is just the consequence of Lemma~\ref{lemmamainchange} applied on the special $(\ung,n,\ung)$-cobordism tangle 
\begin{center}
\resizebox{!}{0.15\textheight}{\raisebox{30mm}{$T_{\ung}=$}
\begingroup%
  \makeatletter%
  \providecommand\color[2][]{%
    \errmessage{(Inkscape) Color is used for the text in Inkscape, but the package 'color.sty' is not loaded}%
    \renewcommand\color[2][]{}%
  }%
  \providecommand\transparent[1]{%
    \errmessage{(Inkscape) Transparency is used (non-zero) for the text in Inkscape, but the package 'transparent.sty' is not loaded}%
    \renewcommand\transparent[1]{}%
  }%
  \providecommand\rotatebox[2]{#2}%
  \ifx\svgwidth\undefined%
    \setlength{\unitlength}{320bp}%
    \ifx\svgscale\undefined%
      \relax%
    \else%
      \setlength{\unitlength}{\unitlength * \real{\svgscale}}%
    \fi%
  \else%
    \setlength{\unitlength}{\svgwidth}%
  \fi%
  \global\let\svgwidth\undefined%
  \global\let\svgscale\undefined%
  \makeatother%
  \begin{picture}(1,0.5)%
    \put(0,0){\includegraphics[width=\unitlength]{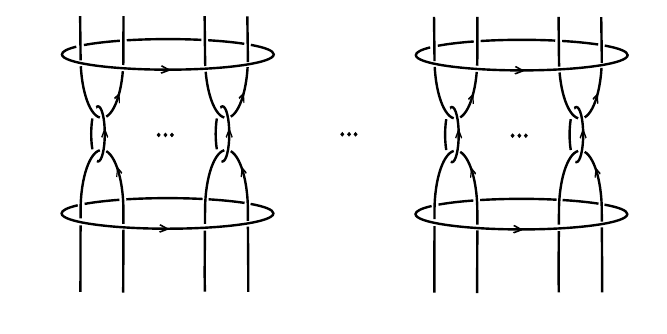}}%
    \put(0.121796,0.04509519){\color[rgb]{0,0,0}\makebox(0,0)[lt]{\begin{minipage}{0.27419414\unitlength}\raggedright $\underbrace{\phantom{ooooooooooooooooo}}_{g_1}$\end{minipage}}}%
    \put(0.65265973,0.04401359){\color[rgb]{0,0,0}\makebox(0,0)[lt]{\begin{minipage}{0.27419414\unitlength}\raggedright $\underbrace{\phantom{ooooooooooooooooo}}_{g_r}$\end{minipage}}}%
  \end{picture}%
\endgroup%
}
\end{center}
where $\ung=(g_1,\ldots,g_r)$.
\end{proof}
Recall that if $\alpha$ is admissible, we have defined a braided functor with anomaly $(\w{\phantom{O}},\gamma)$ and then a TQFT with anomaly $(\tqft{\phantom{O}},\gamma)$. The space associated to the surface of multigenus $\Sigma_{\ung}$ is the image of the idempotent $\W_{\C}(\id_{\Sigma_{\ung}};\alpha)$ denoted by $\Pi_{\alpha}^{\C}$.

\begin{lemma}\label{lemmaisomono}
Let $F=(F,F_2,F_0,\gamma)$ and $G=(G,G_2,G_0,\gamma)$ be two strong monoidal functors with the same anomaly between categories $\C$ and $\D$ where $\C$ is supposed to be rigid.
Then a monoidal natural transformation between $F$ and $G$ with the same anomaly is a natural isomorphism.
\end{lemma}
\begin{proof}
Let
\[\zeta=\{\zeta_{X}\colon F(X)\rightarrow G(X)\}_{X\in \C}\]
be a monoidal natural transformation.
For $X$ an object of $\C$, set
\[\beta_{X}=(G_0G(\widetilde{ev}_{X})G_2(X,X^{*})\otimes \id_{F(X)})(\id_{G(X)}\otimes \zeta_{X^{*}}\otimes \id_{F(X)})(\id_{G(X)}\otimes F_{2}^{-1}(X^{*},X)F(\widetilde{coev}_X)F_0^{-1})\]
and remark that $\zeta_{X}\beta_{X}$ and $\beta_{X}\zeta_{X}$ are identities up to an invertible scalar.
\end{proof}

\begin{lemma}\label{lemmafunctor}
Let $\alpha$ be a morphism of $\mr{Hom}_{\C}(\un,C)$.
If $\alpha$ and $F_{!}\alpha$ are admissible elements then $\zeta$ induces a system of natural isomorphisms
\[\underline{\zeta}=\{\underline{\zeta}\colon F\V_{\C,\alpha}(\Sigma)\rightarrow \V_{\D}(\Sigma,F_{!}\alpha)\}_{\Sigma \in \Cobp}\]
between functors with anomaly $F\V_{\C,\alpha}(\phantom{O})$ and $\V_{\D,F_{!}\alpha}(\phantom{O})$.
\end{lemma}

\begin{proof}
First, note that the two functors $F\V_{\C,\alpha}(\phantom{O})$ and 
$\V_{\D,F_!\alpha}(\phantom{O})$ have the same anomaly. Indeed, the anomaly of $\V_{\C,\alpha}(\phantom{O})$ is given by a product of inverse of $\theta_{\pm}\alpha$ which are morphisms of the form $|T_{\pm}|_{\C,\alpha}$ for some $(0,1,0)$-cobordim tangles $T_{\pm}$ (see Figure \ref{structural2}) so, applying the result of Lemma~\ref{lemmamainchange}, 
$F(\theta_{\pm}^{\C}\alpha)=\theta_{\pm}^{\D}F_!\alpha.$

Let $M_{T}\colon \Sigma_{\ung}\rightarrow \Sigma_{\unh}$ be a cobordism represented by the $(\ung,n,\unh)$-cobordism tangle $T$. As $F\V_{\C,\alpha}(\Sigma_{\ung})$ is a direct factor of $F(C^{\otimes |\ung|})$ and using the commutative diagram of Lemma~\ref{lemmapipi}, we have the following commutative diagram:

\begin{center}
\begin{tikzpicture}[description/.style={fill=white,inner sep=2pt},baseline=(current  bounding  box.center)]
\matrix (m) [ampersand replacement=\&,matrix of math nodes, row sep=3em,column sep=2.5em, text height=1.5ex, text depth=0.25ex]
{F(C^{\otimes |\ung|}) \& \& \& \& D^{\otimes |\ung|} \\ 
 \& F\V_{\C,\alpha}(\Sigma_{\ung}) \& \& V_{\D,F_{!}\alpha}(\Sigma_{\ung}) \&  \\
 \& F\V_{\C,\alpha}(\Sigma_{\unh}) \& \& V_{\D,F_{!}\alpha}(\Sigma_{\unh}) \&  \\
 F(C^{\otimes |\ung|})\& \& \& \& D^{\otimes |\ung|} \\ };
\path[->,font=\scriptsize]
(m-1-1) edge node[auto=left] {$\zeta^{\otimes |\ung|}$} (m-1-5)
		edge node[auto=right] {$F(\w{M_T})$} (m-4-1)
(m-1-5) edge node[auto=left] {$\W_{\D,F_!\alpha}(M_T)$} (m-4-5)
(m-4-1) edge node[auto=right] {$\zeta^{\otimes |\unh|}$} (m-4-5)
(m-2-2) edge node[auto=right] {$F\V_{\C,\alpha}(M_T)$} (m-3-2)
(m-2-2) edge node[auto=right] {} (m-2-4)
(m-2-4) edge node[auto=left]  {$F\V_{\D,F_!\alpha}(M_T)$} (m-3-4)
(m-3-2) edge node[auto=right] {} (m-3-4);
\path[left hook->,font=\scriptsize]
(m-2-2) edge node[auto=left] {} (m-1-1)
(m-3-2) edge node[auto=left] {} (m-4-1);
\path[right hook->,font=\scriptsize]
(m-2-4) edge node[auto=right] {} (m-1-5)
(m-3-4) edge node[auto=right] {} (m-4-5);
\end{tikzpicture}
\end{center}
The induced natural transformation $\underline{\zeta}\colon F\V_{\C,\alpha}\rightarrow V_{\D,F_!\alpha}$ is monoidal by construction. Finally, we can conclude applying Lemma~\ref{lemmaisomono}.
\end{proof}

\begin{remarks}
\item If $\zeta \colon F(C) \rightarrow D$ is an epimorphism then $F_{!}\alpha$ is an admissible element.
\end{remarks}

\subsection{The modularizable case}
\begin{theorem}\label{theoremodu}\textcolor{white}{-}
Let $\C$ be a normalizable premodular category with Kirby element $\alpha_{K}$.
Assume that $\C$ is modularizable, with modularization $F \colon \C \to \widetilde{\C}$. Then $F(\T)$ is a subcategory of the category $\widetilde{T}$  of transparent objects of $\widetilde{\C}$ and there exists a natural isomorphism $\zeta$ such that
\begin{align*}
\begin{tikzpicture}[description/.style={fill=white,inner sep=2pt},baseline=(current  bounding  box.center)]
\matrix (m) [ampersand replacement=\&,matrix of math nodes, row sep=2em,column sep=2.5em, text height=2.2ex, text depth=0.25ex]
{ \Cob \&  \& \T \\
\phantom{0} \& \phantom{0}  \& \widetilde{T}\\
 \phantom{0} \& \phantom{0} \& \vect \\};
\path[->,font=\scriptsize]
(m-1-1) edge node[auto=left] {$\V_{\C,\alpha_{K}}$} (m-1-3)
	    edge node[auto=right] {$\mathrm{RT}_{\widetilde{\C}}$} (m-3-3)
(m-1-3) edge node[auto=left] {$F_{\mid \T}$} (m-2-3)
(m-2-3) edge node[auto=left] {$\mathrm{Hom}_{\widetilde{\C}}(\un,-)$} (m-3-3);
\draw[double,line width=0.8pt] (0.2,0)-- (1.3,1);
\draw (0.5,0.6) node{$\zeta$};
\end{tikzpicture}
\end{align*}
\end{theorem}

\begin{proof}
Apply the result of Lemma~\ref{lemmafunctor} on the modularization functor $F\colon \C\rightarrow \widetilde{\C}$ which is ribbon and so preserves coends. In this case, $F_{!}(\alpha_K)$ is the Kirby element of $\C$ so is admissible. For details, see \cite{Brug}, Section 2.
\end{proof}

%

\bibliographystyle{amsalpha}
\bibliography{bibli}

\newcommand{\etalchar}[1]{$^{#1}$}
\providecommand{\bysame}{\leavevmode\hbox to3em{\hrulefill}\thinspace}
\providecommand{\MR}{\relax\ifhmode\unskip\space\fi MR }
\providecommand{\MRhref}[2]{%
  \href{http://www.ams.org/mathscinet-getitem?mr=#1}{#2}
}
\providecommand{\href}[2]{#2}
\begin{thebibliography}{BCGPM16}

\bibitem[Ati89]{Ati}
M.~Atiyah, \emph{Topological quantum field theories}, Publications
  {M}athématiques de l'{I}{H}{É}{S} \textbf{68} (1989), 175–186.

\bibitem[BCGPM16]{BCGP}
C.~Blanchet, F.~Constantino, N.~Geer, and B.~Patureau-Mirand, \emph{Non
  semi-simple {T}{Q}{F}{T}s, {R}eidemeister torsion and {K}ashaev’s
  invariants}, Advances in Mathematics \textbf{301} (2016), 1–78.

\bibitem[BHMV95]{BHMV}
C.~Blanchet, N~. Habegger, G.~Masbaum, and P.~Vogel, \emph{Topological quantum
  field theories derived from the kauffman bracket}, Topology \textbf{34, No 4}
  (1995), 883–927.

\bibitem[Bru00]{Brug}
A.~Brugui\`eres, \emph{Cat\'egories pr\'emodulaires, modularisations et
  invariants des vari\'et\'es de dimension 3}, Mathematische Annalen
  \textbf{316 (2)} (2000), 215--236.

\bibitem[BV05]{BrugVirHopf}
A.~Brugui\`eres and A.~Virelizier, \emph{Hopf diagrams and quantum invariants},
  Algebraic and Geometric Topology \textbf{5} (2005), 1677--1710.

\bibitem[BV07]{BrugVir}
\bysame, \emph{Hopf monads}, Advances in Mathematics \textbf{215} (2007),
  679--733.

\bibitem[BV13]{BrugVircenter}
\bysame, \emph{On the center of fusion categories}, Pacific Journals of
  Mathematics \textbf{264 , No.1} (2013), 1--30.

\bibitem[DRGG{\etalchar{+}}]{RGGPR}
M.~De~Renzi, A.~Gainutdinov, N.~Geer, B.~Patureau-Mirand, and I.~Runkel,
  \emph{3-dimensional \textsc{TQFTs} from non-semisimple modular categories},
  Selecta Mathematica \textbf{28}, article 42.

\bibitem[FR79]{FR}
R.~Fenn and C.~Rourke, \emph{On {K}irby's calculus of links}, Topology
  \textbf{18} (1979), 1--15.

\bibitem[FS]{FS}
J.~Fuchs and C.~Schweigert, \emph{Coends in conformal field theory},
  ZMP-HH/16-6, Hamburger Beiträge zur Mathematik \textbf{589}, preprint,
  arXiv:1604.01670.

\bibitem[Kas95]{Kassel}
C.~Kassel, \emph{Quantum {G}roups}, Graduate Texts in Mathematics, vol. 155,
  Springer, 1995.

\bibitem[Ker97]{Ker}
T.~Kerler, \emph{Genealogy of nonperturbative quantum-invariants of
  3-manifolds: The surgical family}, Geometry and physics (ed. J. E. Andersen
  et al.; Dekker, New York, 1997) (1997), 503--547.

\bibitem[Kir78]{Kir}
K.~Kirby, \emph{A calculus of framed links in $\mathbb{S}^{3}$}, Inventiones
  mathematicae \textbf{45} (1978), 35--56.

\bibitem[KL01]{KL}
T~Kerler and V.~Lyubashenko, \emph{Non-{S}emisimple {T}opological {Q}uantum
  {F}ield {T}heories for 3-{M}anifolds with {C}orners}, Lecture Notes in
  Mathematics, vol. 1765, Springer, 2001.

\bibitem[Koc03]{Kock}
J.~Kock, \emph{Frobenius algebras and 2{D} topological quantum field theories},
  London Mathematical Society Student Texts, vol.~59, Cambridge University
  Press, 2003.

\bibitem[KRT97]{KRT}
C.~Kassel, M.~Rosso, and V.~Turaev, \emph{Quantum groups and knot invariants},
  Panoramas et Synthèses de la S.M.F., vol.~5, Société Mathématique de
  France, 1997.

\bibitem[Lic97]{Lickorish1997}
W.~B.~R. Lickorish, \emph{An introduction to knot theory}, Springer, 1997.

\bibitem[Lyu95a]{Lyu}
V.~Lyubashenko, \emph{Invariants of 3-manifolds and projective representations
  of mapping class groups via quantum groups at roots of unity}, Communications
  in Mathematical Physics \textbf{172} (1995), 467–516.

\bibitem[Lyu95b]{Lyubhopf}
\bysame, \emph{Tangles and {H}opf algebras in braided categories}, Journal of
  Pure and Applied Algebra \textbf{98} (1995), 245--278.

\bibitem[ML98]{MacL}
S.~Mac~Lane, \emph{Categories for the working mathematician}, 2nd edition,
  1998.

\bibitem[Rad12]{Radford}
D.~E. Radford, \emph{Hopf {A}lgebras}, Series on Knots and Everything, vol.~49,
  World Scientific, 2012.

\bibitem[RT91]{RT}
N.~Reshetikhin and V.~Turaev, \emph{Invariants of 3-manifolds via link
  polynomials and quantum groups}, Inventiones mathematicae \textbf{103}
  (1991), 547--598.

\bibitem[Shu94]{Shu94}
M.~C. Shum, \emph{Tortile {T}ensor {C}ategories}, Journal of Pure and Applied
  Algebra \textbf{93 (1)} (1994), 57--110.

\bibitem[Tur94]{Turaev1994}
V.~Turaev, \emph{Quantum {I}nvariants of {K}nots and 3-{M}anifolds}, Studies in
  Mathematics, vol.~18, W. de Gruyter, 1994.

\bibitem[Vir06]{Vir}
A.~Virelizier, \emph{Kirby elements and quantum invariants}, Proc. London Math.
  Soc. 93 (2006), 474--514.

\bibitem[Wit89]{Wit}
E.~Witten, \emph{Quantum field theory and the {J}ones polynomial},
  Communications in Mathematical Physics \textbf{121, No. 3} (1989), 351–399.

\end{thebibliography}
\nocite{*}
\end{document}